\documentclass[11pt,twoside]{article}
\newcommand{\ifims}[2]{#1} 
\newcommand{\ifAMS}[2]{#1}   
\newcommand{\ifau}[4]{#1}  
\newcommand{\ifbook}[2]{#1}   
\newcommand{\ifLaplace}[2]{#1} 
\newcommand{\ifNL}[2]{#1}  
\newcommand{\ifapp}[2]{#2}  
\newcommand{\ifshort}[2]{#2}  
\newcommand{\iffourG}[2]{#2}  
\newcommand{\ifsqnorm}[2]{#2}  
\newcommand{\ifMLE}[2]{#1}  

\usepackage{tocloft}

\advance\cftsecnumwidth 0.3em\relax
\advance\cftsubsecnumwidth 0.3em\relax

\usepackage{ifthen}
\usepackage[ampersand]{easylist}
\usepackage{amsmath,amssymb,amsthm}
\usepackage[bb=stix]{mathalpha}
\usepackage{mathtools}
\usepackage{natbib}
\usepackage{epsfig,graphicx}
\usepackage{comment}
\usepackage{color}
\usepackage{srcltx}
\usepackage[mathscr]{eucal}
\usepackage[math]{easyeqn}
\usepackage{etoolbox}
\usepackage{hyperref}
\hypersetup{
            colorlinks,
            linkcolor=hookersgreen,
            linktoc=blue,
            citecolor=ultramarine,
            urlcolor=black,
            filecolor=black
            }
\usepackage{nameref}
\usepackage{multirow}

\usepackage{accents}
\usepackage{mathrsfs}
\usepackage{bbm}
\usepackage{algorithm}
\usepackage{algpseudocode}

\ifims{
\textheight=22cm
\textwidth=14.8cm
\topmargin=0pt
\oddsidemargin=1.0cm
\evensidemargin=1.0cm
\linespread{1.3}

}{ 
}

\numberwithin{equation}{section}
\numberwithin{figure}{section}
\newcounter{example}[section]
\numberwithin{example}{section}
\newcounter{remark}[section]
\numberwithin{remark}{section}
\newtheorem{theorem}{Theorem}[section]
\newtheorem{proposition}[theorem]{Proposition}
\newtheorem{lemma}[theorem]{Lemma}

\newtheorem{exmp}[example]{Example}
\newtheorem{rmrk}[remark]{Remark}
\newenvironment{example}{\begin{exmp}\rm}{\end{exmp}}
\newenvironment{remark}{\begin{rmrk}\rm}{\end{rmrk}}

\ifbook{
    \newcommand{\Chapter}[1]{\section{#1}}
    \newcommand{\Section}[1]{\subsection{#1}}
    \newcommand{\Subsection}[1]{\subsubsection{#1}}
    \def\Chname{Section }
    \def\chname{section }
    
  }
  {
    \newcommand{\Chapter}[1]{\chapter{#1}}
    \newcommand{\Section}[1]{\section{#1}}
    \newcommand{\Subsection}[1]{\subsection{#1}}
    \def\Chname{Chapter}
    
  }

\renewcommand{\(}{$\,}
\renewcommand{\)}{\,$}

\def\nquad{\hspace{-1cm}}
\def\eqdef{\stackrel{\operatorname{def}}{=}}

\DeclareMathAlphabet{\mathbbmsl}{U}{bbm}{bx}{sl}

\DeclareMathSymbol{\Alpha}{\mathalpha}{operators}{"41}
\DeclareMathSymbol{\Beta}{\mathalpha}{operators}{"42}
\DeclareMathSymbol{\Epsilon}{\mathalpha}{operators}{"45}
\DeclareMathSymbol{\Zeta}{\mathalpha}{operators}{"5A}
\DeclareMathSymbol{\Eta}{\mathalpha}{operators}{"48}
\DeclareMathSymbol{\Iota}{\mathalpha}{operators}{"49}
\DeclareMathSymbol{\Kappa}{\mathalpha}{operators}{"4B}
\DeclareMathSymbol{\Mu}{\mathalpha}{operators}{"4D}
\DeclareMathSymbol{\Nu}{\mathalpha}{operators}{"4E}
\DeclareMathSymbol{\Omicron}{\mathalpha}{operators}{"4F}
\DeclareMathSymbol{\Rho}{\mathalpha}{operators}{"50}
\DeclareMathSymbol{\Tau}{\mathalpha}{operators}{"54}
\DeclareMathSymbol{\Chi}{\mathalpha}{operators}{"58}
\DeclareMathSymbol{\omicron}{\mathord}{letters}{"6F}

\newcommand{\cc}[1]{\mathscr{#1}}
\newcommand{\bb}[1]{\boldsymbol{#1}}

\DeclareFontFamily{U}{mathx}{\hyphenchar\font45}
\DeclareFontShape{U}{mathx}{m}{n}{
<5><6><7><8><9><10>
<10.95><12><14.4><17.28><20.74><24.88>
mathx10
}{}
\DeclareSymbolFont{mathx}{U}{mathx}{m}{n}
\DeclareFontSubstitution{U}{mathx}{m}{n}
\DeclareMathAccent{\widebar}{0}{mathx}{"73}

\renewcommand{\bar}[1]{\widebar{#1}}
\renewcommand{\hat}[1]{\widehat{#1}}
\renewcommand{\tilde}[1]{\widetilde{#1}}

\makeatletter
\def\mathcenterto#1#2{\mathclap{\phantom{#1}\mathclap{#2}}\phantom{#1}}
\let\old@widetilde\widetilde
\def\widetildeto#1#2{\mathcenterto{#2}{\old@widetilde{\mathcenterto{#1}{#2\,}}}}
\let\old@widehat\widehat
\def\widehatto#1#2{\mathcenterto{#2}{\old@widehat{\mathcenterto{#1}{#2\,}}}}
\makeatother

\newcommand{\thankstitle}[1]{\ifthenelse{\equal{#1}{}}{}{\thanks{#1}}}
\newcommand{\thanksau}[1]{\ifthenelse{\equal{#1}{}}{}{\thanks{#1}}}

\newcommand{\aua}[6]
{\def\authora{#1}
\def\runauthora{#2}
\def\addressa{#3}
\def\emaila{#4}
\def\affiliationa{#5}
\def\thanksa{#6}}

\def\theauthors{
\ifau{ 
  \author{
    \authora
    \thanksau{\thanksa}
    \\[5.pt]
    \addressa \\
    \texttt{ \emaila}
  }
}
{  
  \author{
    \authora
    \thanksau{\thanksa}
    \\[5.pt]
    \addressa \\
    \texttt{ \emaila}
    \and
    \authorb
    \thanksau{\thanksb}
    \\[5.pt]
    \addressb \\
    \texttt{ \emailb}
  }
}
{   
  \author{
    \authora
    \thanksau{\thanksa}
    \\[5.pt]
    \addressa \\
    \texttt{ \emaila}
    \and
    \authorb
    \thanksau{\thanksb}
    \\[5.pt]
    \addressb \\
    \texttt{ \emailb}
    \and
    \authorc
    \thanksau{\thanksc}
    \\[5.pt]
    \addressc \\
    \texttt{ \emailc}
  }
} {   
  \author{
    \authora
    \thanksau{\thanksa}
    \\[5.pt]
    \addressa \\
    \texttt{ \emaila}
    \and
    \authorb
    \thanksau{\thanksb}
    \\[5.pt]
    \addressb \\
    \texttt{ \emailb}
    \and
    \authorc
    \thanksau{\thanksc}
    \\[5.pt]
    \addressc \\
    \texttt{ \emailc}
    \and
    \authord
    \thanksau{\thanksd}
    \\[5.pt]
    \addressd \\
    \texttt{ \emaild}
  }
}
}

\renewcommand{\Gamma}{\varGamma}
\renewcommand{\Pi}{\varPi}
\renewcommand{\Sigma}{\varSigma}
\renewcommand{\Delta}{\varDelta}
\renewcommand{\Lambda}{\varLambda}
\renewcommand{\Psi}{\varPsi}
\renewcommand{\Phi}{\varPhi}
\renewcommand{\Theta}{\varTheta}
\renewcommand{\Omega}{\varOmega}
\renewcommand{\Xi}{\varXi}
\renewcommand{\Upsilon}{\varUpsilon}

\def\argmax{\operatornamewithlimits{argmax}}
\def\argmin{\operatornamewithlimits{argmin}}

\def\av{\bb{a}}
\def\bv{\bb{b}}

\def\hv{\bb{h}}

\def\uv{\bb{u}}

\def\wv{\bb{w}}
\def\xv{\bb{x}}

\def\zv{\bb{z}}

\def\Av{\bb{A}}
\def\Bv{\bb{B}}

\def\Xv{\bb{X}}
\def\Yv{\bb{Y}}
\def\Zv{\bb{Z}}

\def\betav{\bb{\beta}}
\def\chiv{\bb{\chi}}

\def\epsv{\bb{\varepsilon}}
\def\etav{\bb{\eta}}

\def\phiv{\bb{\phi}}

\def\thetav{\bb{\theta}}

\def\xiv{\bb{\xi}}

\def\Psiv{\bb{\Psi}}

\def\sumi{\sum_{i=1}^{n}}

\usepackage{color}

\definecolor{blue(pigment)}{rgb}{0.2, 0.2, 0.6}
\definecolor{ultramarine}{rgb}{0.07, 0.04, 0.56}
\definecolor{darkspringgreen}{rgb}{0.09, 0.45, 0.27}
\definecolor{hookersgreen}{rgb}{0.0, 0.44, 0.0}
\definecolor{hgreen}{rgb}{0.0, 0.44, 0.0}
\definecolor{plum(traditional)}{rgb}{0.56, 0.27, 0.52}
\definecolor{purple(html/css)}{rgb}{0.5, 0.0, 0.5}
\definecolor{magenta(dye)}{rgb}{0.79, 0.08, 0.48}

\newcommand{\scorem}[1]{\score_{\!#1}}

\def\detab{\deta_{0}}

\def\dimX{n_{\xxv}}
\def\dimW{n_{\wwv}}
\def\Mh{\dimp}
\def\SmOpT{\mathbbmsl{A}}
\def\smop{a}
\def\IFTinv{\mathcal{I}}

\def\XXv{\mathbbmsl{X}}
\def\XXvs{\XXv^{*}}
\def\xxv{\mathbbmsl{x}}
\def\wwv{\mathbbmsl{w}}

\def\neff{\mathbbmsl{N}}
\def\hspm{\hspace{1pt}}
\def\AvGP{\Avm_{\!\GP}}
\def\AvGPT{\Avm_{\!\GPT}}

\def\AFN{\mathbbmsl{U}}
\def\Avm{\bb{M}}

\def\DFM{\mathcal{D}}
\def\HPM{H}
\def\DPM{D}
\def\DPMb{\breve{\DPM}}

\def\Avn{\mathcal{M}}

\def\DPNL{D}

\def\DPNLc{\DPNL_{0}}

\def\deta{\tau}
\def\neta{z}
\def\netav{\bb{\neta}}
\def\SmOp{\mathcal{S}}
\def\smop{s}
\def\smopv{\bb{\smop}}
\def\CONSTSO{\CONSTi_{\SmOp}}

\def\F{F}
\def\Fb{\breve{\F}}
\def\Fb{\Phi}
\def\aIF{\alpha}
\def\bIF{\beta}

\def\AFN{\mathbb{Z}}
\def\BFN{\mathbb{B}}
\def\afv{\bb{s}}
\def\afn{\mathbb{z}}
\def\dmax{\kappa}

\def\neta{t}
\def\netav{\bb{\neta}}

\def\bvn{\bb{\mu}}
\def\svn{\bb{\phi}}

\def\AFN{\mathbb{U}}

\def\Tens{\mathcal{T}}

\def\TensU{\Tens}

\def\prmt{\ups}

\def\prmtv{\bb{\prmt}}

\def\prmtvs{\prmtv^{*}}

\def\zvs{\zv^{*}}

\def\targ{x}
\def\targv{\bb{\targ}}

\def\tarp{\theta}
\def\tarpv{\bb{\tarp}}

\def\tarpvs{\tarpv^{*}}

\def\Regrf{M}
\def\Regrfv{\bb{\Regrf}}

\def\Regrfvs{\Regrfv^{*}}

\def\deta{\kappa}

\def\smlc{\varrho}

\def\hmax{\mathsf{c}}

\def\hL{h}

\def\smpa{\smp_{0}}

\def\dagg{\prime}

\def\Etad{\Eta^{\circ}}

\def\amax{\nu}

\def\Matr{\mathfrak{M}}

\def\Eta{\mathcal{H}}

\def\HVB{\mathbbmsl{V}} 

\def\HL{\mathbb{m}}

\def\smp{s}

\def\dltw{\delta}

\def\dltwb{\omega}

\def\dltwu{\tau}
\def\dltwd{\dltw^{\dagg}}
\def\dltwbd{\dltwb^{\dagg}}

\def\dltwbD{\dltwb^{+}}

\def\dltwbD{\dltwb^{+}}

\def\II{\mathcal{I}}

\def\DFN{\DVL}

\def\QF{\mathcal{\QP}}

\def\Projc{\Proj^{\perp}}

\def\R{\mathbbmsl{R}}
\def\E{\mathbbmsl{E}}

\def\P{\mathbbmsl{P}}

\def\kappa{\varkappa}

\def\DNN{\operatorname{DNN}}

\def\blk{\operatorname{block}}

\def\diag{\operatorname{diag}}

\def\ND{\mathcal{N}}

\def\oper{\operatorname{op}}

\def\Var{\operatorname{Var}}

\def\T{\top}
\def\tr{\operatorname{tr}}

\def\nsize{{n}}

\def\sumi{\sum_{i=1}^{\nsize}}

\def\ex{\mathrm{e}}

\def\Id{\mathbbmsl{I}}
\def\Ind{\operatorname{1}\hspace{-4.3pt}\operatorname{I}}

\def\alp{\alpha}

\def\avn{\av}

\def\Avn{\bb{M}}

\def\betavs{\betav^{*}}

\def\bias{\mathsf{b}}

\def\biasD{\bias_{\scalebox{0.666}{$\DPN$}}}
\def\biasDF{\bias_{\scalebox{0.666}{$\DFM$}}}

\def\B{\cc{B}}
\def\BB{\mathbbmsl{B}}

\def\BBB{\cc{B}}

\def\BBH{B}

\def\Bv{\bb{\BB}}

\def\CA{\mathcal{A}}

\def\CGP{w}
\def\CGPa{\CGP_{0}}
\def\CGPs{\CGP_{*}}

\def\CONST{\mathtt{C} \hspace{0.1em}}
\def\CONSTi{\mathtt{C}}

\def\CONSTIF{\CONSTi_{\IF}}

\def\const{\mathsf{c}}

\def\DP{D}

\def\DPN{D}

\def\DVL{\mathbb{D}}

\def\dimA{\mathbb{p}}
\def\dimAr{\breve{\dimA}}

\def\dimttl{\bar{\dimp}}

\def\dimp{p}

\def\dimG{\dimA_{\GP}}

\def\dimQ{\dimA_{\QP}}

\def\dimDF{\dimA_{\scalebox{0.666}{$\DFM$}}}
\def\dimtarg{\dimA}

\def\dimq{q}

\def\dimd{d}
\def\dimm{m}
\def\dimX{n \times d}

\def\dimD{\dimA_{\scalebox{0.666}{${\DPN}$}}}

\def\dimm{M}

\def\Eta{\cc{H}}

\def\err{\diamondsuit}

\def\etavs{\etav^{*}}

\def\fs{f}

\def\fn{g}

\def\fG{f_{\GP}}

\def\fGu{h}

\def\fNL{g}

\def\gp{g}

\def\GP{G}

\def\GPY{\Gamma}

\def\GPa{\GP_{0}}

\def\GPT{\mathcal{G}}

\def\GPW{\mathcal{G}}

\def\IF{\mathbbmsl{F}}

\def\IFN{\IFL}

\def\IFL{{\mathbbmss{F}}}

\def\IFtotal{F}
\def\IFT{\mathscr{\IFtotal}}

\def\IFG{\mathbbmsl{G}}
\def\IFTb{\Phi}

\def\IFL{\mathbb{F}}

\def\Kappa{\cc{K}}

\def\LT{L}
\def\LGP{\LT_{\GP}}

\def\LL{\cc{L}}

\def\LLp{\mathbbmsl{L}}

\def\loss{\wp}

\def\mm{m}

\def\mmc{\mm'}

\def\muA{\mu}

\def\Mh{M}

\def\nuov{\bb{a}}

\def\nup{\eta}
\def\nupv{\bb{\nup}}
\def\nupvs{\nupv^{*}}

\def\nui{s}

\def\nuiv{\bb{\nui}}

\def\nuo{\tau}
\def\nuov{\bb{\nuo}}

\def\pent{\operatorname{pen}}

\def\Proj{\Pi}

\def\priord{\pi}

\def\QP{Q}

\def\regrf{m}
\def\regrfs{\regrf^{*}}
\def\regrfv{\bb{\regrf}}
\def\regrfvs{\regrfv^{*}}

\def\rhot{t}

\def\rhoIF{\rho}

\def\riskt{\cc{R}}

\def\rr{\mathtt{r}}

\def\rrde{\rrdd}

\def\rrB{\rr_{\scalebox{0.666}{${\BBH}$}}}
\def\rrD{\rr_{\scalebox{0.666}{${\DPN}$}}}
\def\rrDF{\rr_{\scalebox{0.666}{$\DFM$}}}

\def\rrn{\rr}

\def\rrdd{\rr_{0}}

\def\score{\nabla}

\def\thetav{\bb{\theta}}
\def\thetavs{\thetav^{*}}

\def\Thetad{\Theta^{\circ}}

\def\Tau{T}

\def\uvd{\uv^{\circ}}

\def\ups{\upsilon}
\def\upsv{\bb{\ups}}

\def\upsvd{\upsv^{\circ}}
\def\upsvs{\upsv^{*}}

\def\upsvr{\breve{\upsv}}

\def\upsvn{\upsvd}

\def\ups{\upsilon}
\def\upsv{\bb{\ups}}

\def\UV{\mathcal{U}}

\def\UVz{\UV}

\def\Ups{\varUpsilon}
\def\Upsd{\Ups^{\circ}}

\def\VP{V}

\def\weight{w}

\def\wv{\bb{w}}

\def\WWv{W}
\def\WWvs{\WWv^{*}}

\def\xivr{\breve{\xiv}}

\def\xx{\mathtt{x}}

\def\YY{\cc{Y}}

\def\zq{z}

\def\thetitle{Estimation and inference for Deep Neuronal Networks}
\def\theruntitle{Estimation and inference for Deep Neuronal Networks}

\def\theabstract{
Nonlinear regression problem is one of the most popular and important statistical tasks.
The first methods like least squares estimation go back to Gauss and Legendre.
Recent models and developments in statistics and machine learning like Deep Neuronal Networks (DNN) or 
nonlinear PDE stimulate new research in this direction which has to address 
the important issues and challenges of modern statistical inference such as 
huge complexity and parameter dimension of the model, limited sample size, lack of convexity and identifiability,
among many others.
Classical results of nonparametric statistics in terms of rate of convergence do not really address the mentioned issues.
This paper offers a general approach to studying a nonlinear regression problem based on the notion of effective dimension.
First, a special case of models with stochastically linear structure (SLS) is studied.
The results provide finite sample expansions for the loss of the penalized maximum likelihood estimation (MLE).
The leading term of such expansions as well as the corresponding remainder are given via the effective dimension and the effective sample size. 
The obtained expansions can be used to obtain sharp risk bounds and for statistical inference.
Despite generality, all the presented bounds are nearly sharp and the classical asymptotic results can be obtained
as simple corollaries. 
Although the basic SLS assumptions are not fulfilled for nonlinear smooth regression, we explain 
how the stochastic linearity can be achieved by extending the parameter space. 
The obtained general results are specified to nonlinear smooth regression and to a DNN with one hidden layer. 
}

\def\kwdp{62F10,62E17}
\def\kwds{62J12}

\def\thekeywords{Fisher and Wilks expansions, risk bounds, nonlinear regression}

\def\thankstitle{}

\aua
{Vladimir Spokoiny}
{Spokoiny, V. }
{
    Weierstrass Institute and HU Berlin,  
    \\
    Mohrenstr. 39, 10117 Berlin, Germany
    }
{spokoiny@wias-berlin.de}
{Weierstrass-Institute and Humboldt University Berlin}
{
  Financial support by the German Research Foundation (DFG) through the Collaborative Research Center 1294 ``Data assimilation'' is gratefully acknowledged.
}

\begin{document}
\thispagestyle{empty}
{
\title{\thetitle}
\theauthors

\maketitle
\begin{abstract}
{\footnotesize \theabstract}
\end{abstract}

\ifAMS
    {\par\noindent\emph{AMS 2010 Subject Classification:} Primary \kwdp. Secondary \kwds}
    {\par\noindent\emph{JEL codes}: \kwdp}

} 

\tableofcontents


\section{Introduction}
\label{Snoninverseintro}
Nonlinear regression problem belongs to the very core of mathematical statistics and goes back at least to Gauss and Legendre. 
However, it remains an actively developed research topic in modern statistics and machine learning,
particularly due to applications to e.g. nonlinear inverse problems \cite{Nickl2018ConvergenceRF}, 
deep learning \cite{ShHi2020}, and references therein.
%
Nonlinearity of the model makes the study very involved and the cited results heavily used 
the recent advances in the theory of partial differential equations, inverse problems, empirical processes. 
We mention \cite{nickl2020} and \cite{Nickl2018ConvergenceRF} as particular illustrations 
of the major difficulties in the corresponding study.

We focus here on the problem of parameter estimation for a known nonlinear  
regression function \( \regrfv(\thetav) \) valued in some Euclidean space \( \YY \) from noisy data
\( \Yv \) satisfying 
\begin{EQA}
	\Yv 
	&=& 
	\regrfv(\thetav) + \epsv \in \YY,
\label{YAfseNI}
\end{EQA}
where \( \E \epsv = 0 \).
The standard approach to this problem is based on minimization of the fidelity 
\( \| \Yv - \regrfv(\thetav) \|^{2} \).
Usually, this problem is numerically hard and one or another regularization is used.
A typical example is given by Tikhonov regularization \( \muA \| \thetav \|^{2} \).
More generally, for a smooth signal \( \thetav \), one may consider a smooth penalty 
\( \| \GP \thetav \|^{2} \) with 
a penalizing operator \( \GP^{2} \), e.g.
\( \| \GP \thetav \|^{2} = \| \thetav \|_{H^{\alpha}}^{2} \).
A proper choice of the smoothness parameter \( \alpha \) is important for obtaining a rate 
optimal procedure over a class of smooth \( \thetav \) in Sobolev sense; see e.g. \cite{Nickl2018ConvergenceRF}.

For modern applications like Deep Neuronal Networks,
the main challenges for studying the problem of nonlinear regression are a possibly huge or even infinite dimension 
of the parameter space and a limited sample size.
The problem is not convex, even parameter identifiability is questionable.
Nonlinear inverse problems are often ill-posed.
In the linear case, the degree of ill-posedness is usually described in terms 
of eigenvalues of the corresponding operator via so called source condition; its extension on the nonlinear setup
is questionable and requires rather strong assumptions.
The classical parametric asymptotic approach hardly applies in this situation.
Typical nonparametric results provide some results in term of the convergence rate over smoothness classes; see e.g. 
\cite{Nickl2018ConvergenceRF}, \cite{ShHi2020}.
However, these rate results are usually not really informative for inference issues because 
they involve a number of hidden constants which can even explode when parameter dimension grows. 
To address such issues, new tools and ideas are called for. 

This paper offers a new approach to studying the properties of the MLE in nonlinear regression problems.
The study includes the following main steps.
First, we establish some results about MLE properties for a special class 
of so called \emph{stochastically linear smooth} (SLS) models. 
The major assumptions for the SLS setup are linearity of the stochastic component of the considered
log-likelihood w.r.t. the target parameter and concavity of the expected log-likelihood.
This allows to overcome traditional difficulties in studying the MLE or minimum contrast estimators and 
avoid the high-tech tools of the empirical process theory;
cf. \cite{Kosorok,nickl_2015}.
Instead, we only need some deviation bounds for quadratic forms; see 
\ifsqnorm{Section~\ref{SdevboundnonGauss}.}{\cite{Sp2023c}, \cite{Sp2023d}.}
Unfortunately, the main SLS assumptions fail for nonlinear regression \eqref{YAfseNI}. 
The objective function is not convex and its stochastic component is not linear in the parameter.
Later we offer a method called \emph{calming} which allows to overcome the issue 
of nonlinearity of the stochastic component by extending the parameter space and including the response in the parameter vector.
This naturally leads to a semiparametric problem in which the target parameter is 
estimated along with a high dimensional nuisance parameter.
The parameter space is enlarged, however, the problem is reduced back to the semiparametric SLS framework; see Section~\ref{SsemiMLE} for details.
Semiparametric estimation is well developed; see e.g. \cite{chen1995},
\cite{bickel1993efficient}, \cite{Kosorok} and references therein.
However, most of available results are stated in classical asymptotic setup and cannot be used for our study.
Later in Section~\ref{SsemiMLE} we revisit and reconsider the main notions and results of the semiparametric theory using the general Fisher and Wilks expansions developed for the SLS setup.
A particular focus is on the semiparametric effective dimension and on the bias arising in profile MLE estimation. 

The issue of non-concavity is even more severe.
So far, no universal method of studying the problem of non-convex optimization is available. 
We follow the standard ``localization'' idea imposing the so called ``warm start'' assumption;
see e.g. \cite{Gr2007} for the results about Gauss-Newton iterative methods in non-convex optimization.
%
Combining the mentioned ideas of calming and localization allows to state rather precise finite sample results about the properties of the pMLE 
for nonlinear regression.  

\subsection*{This paper's contributions}
This paper offers a finite-sample approach to parametric estimation and specifies it 
to nonlinear regression problem.
One can highlight few essential points of the study.

\smallskip \par \textbf{Finite sample expansions.}
The main focus of this paper is on Fisher and Wilks expansions of the pMLE with explicit error terms.
Some related results as well as an extensive literature overview can be found in \cite{OsBa2021}.
This paper's approach allows us to establish sharp results and describe the remainder of each expansion in a closed form.
Let \( L(\upsv) \) be  a random function, \( \upsv \in \Ups \subseteq \R^{\dimp} \),
\( \dimp < \infty \).
It can be viewed as log-likelihood or negative loss.
Consider
\begin{EQA}[rcccl]
	\tilde{\upsv}
	&=&
	\argmax_{\upsv} L(\upsv) ;
	\qquad
	\upsvs
	&=&
	\argmax_{\upsv} \E \, L(\upsv) .
\end{EQA}
Define the \emph{Fisher information} matrix	\( \IF(\upsv) = - \nabla^{2} \E L(\upsv) \) and write \( \IF = \IF(\upsvs) \).
Also, introduce 
the \emph{score} vector \( \nabla \zeta = \nabla \zeta(\upsvs) \) for \( \zeta(\upsv) \eqdef L(\upsv) - \E L(\upsv) \).
The \emph{effective dimension} \( \dimA \) and
the \emph{effective sample size} \( \neff \)  are given by 
\begin{EQA}[c]
	\dimA \eqdef \tr \bigl( \IF^{-1} \Var(\nabla \zeta) \bigr),
	\qquad
	\neff \eqdef \lambda_{\min}(\IF) ;
\label{dyffgwuw7ukwd6hj3kis4}
\end{EQA} 
cf.
see \cite{SP2013_rough}, \cite{SpPa2019}, or \cite{Sp2022}.
The main results can be summarized as follows:
\begin{EQA}[lccl]
	\text{\textbf{Concentration}:}
	&
	\| \IF^{1/2} (\tilde{\upsv} - \upsvs) \|
	& \lesssim &
	\sqrt{\dimA} \, ;
	\\
	\text{\textbf{Fisher expansion}:}
	&
	\| \IF^{1/2} (\tilde{\upsv} - \upsvs - \IF^{-1} \nabla \zeta) \|
	& \lesssim &
	\frac{\dimA}{\neff^{1/2}} \, ;
	\\
	\text{\textbf{Wilks expansion}:} \quad
	&
	| L(\tilde{\upsv}) - L(\upsvs) - \tfrac{1}{2} \| \IF^{-1/2} \nabla \zeta \|^{2} |
	& \lesssim &
	\frac{\dimA^{3/2}}{\neff^{1/2}} \, .
\label{f6dfhstwhdf8vhewbnde6fy}
\end{EQA}
The ambient parameter dimension does not show up in the remainder. 
This allows to apply the results to high-dimensional or even nonparametric problems. 
The results are useful in understanding the success of dimension reduction methods, manifold learning, sparse estimation,
or other overparametrized models like Deep Neuronal Networks. 
In each case, the structural model assumption helps to reduce the effective dimension of the problem.
The Fisher expansion enables us to obtain sharp risk bounds for estimation and prediction risk
and to recup the classical asymptotic parametric results like root-n consistency and normality and
asymptotic efficiency, and, at the same time, to prove rate optimality over smoothness classes for the nonparametric framework.
Moreover, Fisher and Wilks expansions can be used for inference, in particular, for studying the asymptotic behavior of the pMLE, for validation 
of resampling bootstrap procedures, testing of structural hypotheses, etc.
A recent paper \cite{GSZ2023} provides an excellent example of using a finite sample expansion for top-k ranking problem in \( \dimp \asymp n \) regime.

Section~\ref{Spriorexample} explains how the obtained risk bounds can be used for deriving the classical 
minimax rate results. 
However, it is well known that the rate results are too conservative and correspond to the worst-case function.
The bias-variance decomposition of the risk provided by Theorem~\ref{TQFiWibias} can be viewed as an upper risk bound 
within the parametric setup.
Another benefit of the parametric setup is that
lower bounds can be obtained in several ways, including Cramer-Rao inequality, \cite{van2000asymptotic}, 
van Trees inequality, \cite{GillLevit1995}, or by 
using the Laplace approximation; see \cite{katsevich2023tight,SpLaplace2022}.

\smallskip \par \textbf{Critical parameter dimension.}
As already mentioned, the obtained results are ``dimension-free'', only the effective dimension matters.
The range of applicability is given by the relation between the effective dimension and the effective sample size.
The concentration and full-dimensional Fisher results apply under the ``critical dimension condition'' \( \dimA \ll n \) meaning a sufficiently many observations per parameter
to be estimated. 
Inference on a low-dimensional parameter sub-vector requires a stronger condition \( \dimA^{2} \ll n \) or even 
\( \dimA^{3} \ll \neff \).

\smallskip \par \textbf{Sharp probability bounds under stochastic linearity.}
Any study of a maximum likelihood or minimum contrast estimator includes a bound on stochastic fluctuation of the random objective function.
Power tools of the empirical process theory can be very helpful; see e.g. \cite{Kosorok,nickl_2015}.
However, the obtained results are too rough for our purposes and do not allow to derive sharp expansions and risk bounds.
Our general results are stated under the condition of
linearity of the stochastic term of the objective function.
This significantly simplifies the study of the penalized maximum likelihood estimator.
The whole stochastic analysis can be reduced to some deviation bounds for the norm of standardized score.
In particular, one may apply the recent results from 
\ifsqnorm{Section~\ref{SdevboundnonGauss}}{\cite{Sp2023c}}
for sub-gaussian case and from 
\ifsqnorm{Section~\ref{SdevboundnonGauss}}{\cite{Sp2023d}} for the sub-exponential case. 
We, therefore, introduce a class of \emph{stochastically linear smooth} (SLS) models for which the stochastic component of the objective function
linearly depends on the parameter.
Section~\ref{SgenBounds} and Section~\ref{SsemiMLE} provide a rigorous study of parametric and
semiparametric estimation for SLS models.
The condition of stochastic linearity is automatically fulfilled for Generalized Linear Models.
But the SLS class is much larger.
We guess that any model can always be embedded into SLS framework by the so called calming trick based on extension of the parameter space.

\smallskip \par \textbf{Calming device for enforcing the stochastic linearity.}
The assumption of linearity of the stochastic term for SLS models can be viewed as quite restrictive.
However, we demonstrate on the example of nonlinear regression how this condition can be enforced by 
extending the parameter space, more precisely, by including the response vector in the parameter list.
The other examples of a calming device include error-in-operator models \cite{HoRe2008,Trabs_2018}, random design regression \cite{ChMo2022},
linear diffusion \cite{Nickl2019}, instrumental variable regression \cite{XiWa2024}, and many others. 
Such an extension of the parameter space leads to a semiparametric setup when the target parameter is to be estimated under presence of 
a high dimensional nuisance parameter. 
However - and this is a great benefit of the approach - the developed technique allows to deduce 
accurate finite sample expansions for the target parameter from the full dimensional expansions. 
The corresponding risk is given by the effective dimension of the target parameter, and, therefore, the proposed calming approach
does not lead to any reduction of estimation accuracy.

\iffourG{
\smallskip \par \textbf{Fourth order expansions.}
The use of a fourth order expansion allows to improve the critical dimension condition \( \dimA^{2} \ll n \) to \( \dimA^{3/2} \ll n \).
This can be useful e.g. for the problem of estimating a smooth functional; see \cite{Kolt2021}.
}{}

\smallskip \par \textbf{Sharp bounds for regularization bias and the quadratic risk.}
Almost any estimation procedure for a high dimensional setup includes one or another regularization technique.
A standard approach is to add a penalty term \( \pent_{\GP}(\upsv) \) to the objective function \( L(\upsv) \) 
which is responsible for the model complexity. 
Typical examples are given by ridge regression \( \pent_{\lambda}(\upsv) = \lambda \| \upsv \|^{2}/2 \),
roughness penalty \( \pent_{\GP}(\upsv) = \| \GP \upsv \|^{2}/2 \), 
sparse penalty \( \pent_{\lambda}(\upsv) = \lambda \| \upsv \|_{1} \), among many others.
Any penalization induces some bias which has to be carefully evaluated. 
In this paper, we limit ourselves to a class of smooth penalties.
A great benefit of such penalization is that the bias term can be described in closed form.
This, in turns, yields sharp risk bounds in the penalized estimation problem. 
We also explain how the choice of penalty can be used to reduced the effective dimension of the problem; see Section~\ref{Spriorexample}.
Another benefit of the approach is that it is \emph{coordinate free} and does not rely on any spectral decomposition 
and/or any spectral representation for the information and penalization operators. 
%

\smallskip \par \textbf{Linear perturbation theory.}
One of the main technical tools of the study is provided by results on perturbed optimization. 
Let \( \fs(\upsv) \) be a smooth concave function, 
\begin{EQA}
	\upsvs
	&=&
	\argmax_{\upsv} \fs(\upsv),
	\quad
	\IFN = - \nabla^{2} \fs(\upsvs) .
\label{fg5hg3gf98tkj3dciryti}
\end{EQA}
Let another function \( \fn(\upsv) \) satisfy for some vector \( \Av \)
\begin{EQA}
	\fn(\upsv) - \fn(\upsvs) 
	&=&
	\bigl\langle \upsv - \upsvs, \Av \bigr\rangle + \fs(\upsv) - \fs(\upsvs) .
\label{4hbh8njoelvt6jwgf09i}
\end{EQA}
Define
\begin{EQA}
	\upsvn
	& \eqdef &
	\argmax_{\upsv} \fn(\upsv),
	\qquad
	\fn(\upsvn)
	=
	\max_{\upsv} \fn(\upsv) .
\label{6yc63yhudf7fdy6edgehyi} 
\end{EQA}
Section~\ref{Slocalsmooth} collects some accurate bounds on \( \upsvn - \upsvs \) and \( \fn(\upsvn) - \fn(\upsvs) \)
in terms of the second, third, or fourth Gateaux derivative of the objective function
under the so called ``self-concordance'' condition introduced in \cite{NeNe1994}; see Section~\ref{Slocalsmooth}.

\smallskip \par \textbf{Sharp risk bounds in training of shallow DNN.}
The obtained results are specified to shallow networks.
It appears that overparametrization typical for DNN does not lead to ``curse of dimensionality'' issue. 
On the contrary, it is effectively used in the analysis to show local concavity of the objective function; cf.
\cite{LCEZZ2022}. 
Under mild conditions on the model, we derive root-n normality in estimation of the DNN architecture, 
the risk bounds do not involve the dimension of the parameter space.
The proposed approach offers another view on the problem of DNN training which differs substantially from 
recent findings in this area; cf. 
    \cite{ShHi2020}, 
    \cite{SYZ2019}, 
In particular, a manifold hypothesis for dimension reduction like in 
    \cite{CJLZ2019}, 
%
    \cite{KKL2019}, 
    \cite{JSLH2021}, 
    \cite{LCZL2021}, 
    \cite{cloninger2021deep}, 
    \cite{JSLH2023}, 
is not used.

  
%

%
%

\subsection*{Organization of the paper}
Sections~\ref{SgenBounds} presents the general results for SLS models.
Profile semiparametric estimation in the SLS setup is studied in Section~\ref{SsemiMLE}.
Section~\ref{Snoninverse} explains the setup and the details of the proposed calming approach
for nonlinear regression and DNN.
Useful results about linearly perturbed optimization are collected in Section~\ref{Slocalsmooth} of the appendix.
Section~\ref{Spriorexample} addresses the issue of selecting a proper penalty and  
comments how the obtained results can be used for deriving the standard asymptotic rate results.


\Chapter{Properties of the MLE \( \tilde{\upsv} \) for SLS models}
\label{SgenBounds}
This \chname collects general results about concentration and expansion of the MLE in the SLS setup
which substantially improve the bounds from \cite{SP2013_rough} and \cite{SpLaplace2022}.
%
We assume to be given a random function \( L(\upsv) \), \( \upsv \in \Ups \subseteq \R^{\dimp} \),
\( \dimp < \infty \).
This function can be viewed as log-likelihood or negative loss.
Consider in parallel two optimization problems defining 
the MLE \( \tilde{\upsv} \) and its population counterpart (the background truth) \( \upsvs \):
\begin{EQA}[rcl]
	\tilde{\upsv} 
	&=& 
	\argmax_{\upsv} L(\upsv),
	\qquad
	\upsvs 
	=
	\argmax_{\upsv} \E L(\upsv),
	\qquad
\label{tuGauLGususGE}
\end{EQA}
Define the Fisher information matrix \( \IF(\upsv) \eqdef - \nabla^{2} \E L(\upsv) \) 
and denote \( \IF = \IF(\upsvs) \). 

\Section{Basic conditions}
\label{Scondgeneric}
Now we present our major conditions.
The most important one is about linearity of the stochastic component 
\( \zeta(\upsv) = L(\upsv) - \E L(\upsv) = L(\upsv) - \E L(\upsv) \).

\medskip
\begin{description}
    \item[\label{Eref} \( \bb{(\zeta)} \)]
      \textit{The stochastic component \( \zeta(\upsv) = L(\upsv) - \E L(\upsv) \) of the process \( L(\upsv) \) is linear in 
      \( \upsv \). 
      We denote by \( \nabla \zeta \equiv \nabla \zeta(\upsv) \in \R^{\dimp} \) its gradient
      }.
\end{description}

Below we assume some concentration properties of the stochastic vector \( \nabla \zeta \).
More precisely, we require that \( \nabla \zeta \) obeys the following condition.

\begin{description}
\item[\label{EU2ref}\( \bb{(\nabla \zeta)} \)]
	\textit{There exists \( \VP^{2} \geq \Var(\nabla \zeta) \) such that  
	for all considered  \( \BBH \in \Matr_{\dimp} \) and \( \xx > 0 \)
	}
\begin{EQA}
	\P\bigl( \| \BBH^{1/2} \VP^{-1} \nabla \zeta \| \geq \zq(\BBH,\xx) \bigr)
	& \leq &
	3 \ex^{-\xx} ,
\label{2emxGPm12nz122}
	\\
	\zq^{2}(\BBH,\xx)
	& \eqdef &
	\tr \BBH + 2 \sqrt{\xx \, \tr \BBH^{2}} + 2 \xx \| \BBH \| \,  .
\label{34rtyghuioiuyhgvftid}
\end{EQA}
\end{description}

This condition can be effectively checked if the errors in the data exhibit sub-gaussian or sub-exponential behavior; see 
\ifsqnorm{Section~\ref{SdevboundnonGauss}.}{\cite{Sp2023c}, \cite{Sp2023d}.}
The important special case corresponds to \( \BBH = \IF^{-1/2} \VP^{2} \IF^{-1/2} \) 
and \( \xx \approx \log n \) leading to the bound
\begin{EQA}
	\P\bigl( \| \IF^{-1/2} \nabla \zeta \| > \zq(\BBH,\xx) \bigr)
	& \leq &
	3/n .
\label{udyvfeyejff6777dj23}
\end{EQA}
The  value \( \dimG = \tr(\IF^{-1} \VP^{2}) \) can be called the \emph{effective dimension}; see \cite{SP2013_rough}.

We also assume that the log-likelihood \( L(\upsv) \) or, equivalently, its deterministic part 
\( \E L(\upsv) \) is a concave function.
It can be relaxed using localization; see Section~\ref{Snoninverse}.

\medskip
\begin{description}
    \item[\label{LLref} \( \bb{(\mathcal{C})} \)]
      \textit{The function \( \E L(\upsv) \) is concave on \( \Ups \) which is open and convex set in \( \R^{\dimp} \).      
      }
\end{description}
\medskip

\ifLaplace{}{
In Section~\ref{Spostconcentr} we consider a stronger condition of semi-concavity of \( \E L(\upsv) \). }
%

Later we will also need some smoothness conditions on the function \( f(\upsv) = \E L(\upsv) \)
within a local vicinity of the point \( \upsvs \).
The notion of locality is given in terms of a metric tensor \( \DPN \in \Matr_{\dimp} \).
In most of the results later on, one can use \( \DPN = \IF^{1/2} \).
In general, we only assume \( \DPN^{2} \leq \dmax^{2} \IF \) for some \( \dmax > 0 \).
%
Introduce the error of the second-order Taylor approximation at a point \( \upsv \) in a direction \( \uv \) by
\begin{EQ}[rcl]
	\dltw_{3}(\upsv,\uv) 
	&=& 
	f(\upsv + \uv) - f(\upsv) - \langle \nabla f(\upsv), \uv \rangle 
	- \frac{1}{2} \langle \nabla^{2} f(\upsv), \uv^{\otimes 2} \rangle , 
	\\
	\dltwd_{3}(\upsv,\uv) 
	&=&
	\langle \nabla f(\upsv + \uv), \uv \rangle - \langle \nabla f(\upsv), \uv \rangle 
	- \langle \nabla^{2} f(\upsv), \uv^{\otimes 2} \rangle \, .
\label{dltw3vufuv12f2g}
\end{EQ}
Second order smoothness means a bound of the form
\begin{EQA}
	\dltw_{3}(\upsv,\uv) 
	& \leq & 
	\dltwb \| \DPN \uv \|^{2} \, ,
	\quad
	\dltwd_{3}(\upsv,\uv) 
	\leq 
	\dltwbd \| \DPN \uv \|^{2} \, ,
	\qquad
	\| \DPN \uv \| \leq \rr \, ,
\label{7jduusyswjdtcrtebfjvu}
\end{EQA}
for some radius \( \rr \) and small constants \( \dltwb \) and \( \dltwbd \).
These quantities can be effectively bounded under smoothness
conditions \nameref{LL3tref}, \nameref{LLsT3ref}, or \nameref{LLtS3ref} given in Section~\ref{Slocalsmooth}.
For instance, under \nameref{LL3tref}, by Lemma~\ref{LdltwLa3t}, it holds for a small constant \( \dltwu_{3} \)
\begin{EQA}
	\dltwbd
	& \leq &
	\dltwu_{3} \, \rr \, ,
	\qquad
	\dltwb
	\leq 
	\dltwu_{3} \, \rr/ 3 .
\label{swdy6fwqd6qtxcvbdyfdtw}
\end{EQA}
Also under \nameref{LLtS3ref}, the same bounds apply with \( \dltwu_{3} = \hmax_{3} \, n^{-1/2} \); see Lemma~\ref{LdltwLaGP}.

The class of models satisfying the conditions \nameref{Eref}, 
\nameref{EU2ref}, and \nameref{LLref}
with a smooth function \( f(\upsv) = \E L(\upsv) \) will be referred to as \emph{stochastically linear smooth} (SLS). 
This class includes linear regression, generalized linear models (GLM), and log-density models; 
see \cite{SpPa2019}, \cite{OsBa2021}\ifapp{ or Section~\ref{SGBvM} later.}{ or \cite{SpLaplace2022}.}
However, this class is much larger.
For instance, nonlinear regression can be adapted to the SLS framework 
by an extension of the parameter space; see 
\ifNL{Section~\ref{Snoninverse}.}{\cite{Sp2023b}.}

\Section{Concentration of the MLE \( \tilde{\upsv} \). 2S-expansions}
\label{SgenpMLE}
This section discusses properties of the MLE \( \tilde{\upsv} = \argmax_{\upsv} L(\upsv) \)
under second-order smoothness conditions.
%
Fix \( \xx > 0 \) and define with \( \VP^{2} \) from \nameref{EU2ref} and \( \BBH = \IF^{-1/2} \VP^{2} \IF^{-1/2} \) 
\begin{EQA}[c]
	\UVz
	\eqdef 
	\bigl\{ \uv \colon \| \IF^{1/2} \uv \| \leq \tfrac{4}{3} \rrB \bigr\},
	\qquad
	\rrB \eqdef \zq(\BBH,\xx) ,
\label{rm1ucDGu2r0DGu}
\end{EQA}
where \( \zq(\BBH,\xx) \) is given by \eqref{34rtyghuioiuyhgvftid}.
By \nameref{EU2ref}, on a random set \( \Omega(\xx) \) with 
\( \P(\Omega(\xx)) \geq 1 - 3 \ex^{-\xx} \), it holds \( \| \IF^{-1/2} \nabla \zeta \| \leq \rr \).
Further, for the metric tensor \( \DPN \) from \eqref{7jduusyswjdtcrtebfjvu}, define
\begin{EQA}[c]
	\dltwb
    \eqdef 
    \sup_{\uv \in \UVz}
    \frac{2 |\dltw_{3}(\upsvs,\uv)|}{\| \DPN \uv \|^{2}} \,\, ,
    \qquad
    \dltwbd
    \eqdef 
    \sup_{\uv \in \UVz} \frac{|\dltwd_{3}(\upsvs,\uv)|}{\| \DPN \uv \|^{2}} \,\, . 
\label{dtb3u1DG2d3GP}
\end{EQA}

\begin{proposition}
\label{PconcMLEgenc}
Suppose
\nameref{Eref},
\nameref{EU2ref},
and 
\nameref{LLref}.
Let also \( \DPN^{2} \leq \dmax^{2} \IF \) and \( \dltwbd \, \dmax^{2} < 1/4 \);
see \eqref{dtb3u1DG2d3GP}.
Then on \( \Omega(\xx) \), 
it holds 
\begin{EQA}
	\| \IF^{1/2} (\tilde{\upsv} - \upsvs) \|  
	& \leq &
	\frac{4}{3} \, \rrB \, ,
	\qquad
	\| \DPN (\tilde{\upsv} - \upsvs) \|  
	\leq
	\frac{4 \dmax}{3} \, \rrB
	\, . 
\label{rhDGtuGmusGU0}
\end{EQA}
\end{proposition}

\begin{proof}
Apply Proposition~\ref{Pconcgeneric} to 
\( \fs(\upsv) = \E L(\upsv) \), \( \amax = 3/4 \), and \( \Av = \nabla \zeta \).
\end{proof}

Concentration of \( \tilde{\upsv} \) around \( \upsvs \) 
can be used to establish a version of the Fisher expansion for 
the estimation error \( \tilde{\upsv} - \upsvs \) and the Wilks expansion for the excess 
\( L(\tilde{\upsv}) - L(\upsvs) \).
The result substantially improves the bounds from \cite{OsBa2021} for M-estimators and follows by Proposition~\ref{PFiWigeneric}.

\begin{theorem}
\label{TFiWititG}
Assume the conditions of Proposition~\ref{PconcMLEgenc}.
Then on \( \Omega(\xx) \!\) 
\begin{EQ}[rcl]
    2 L(\tilde{\upsv}) - 2 L(\upsvs) 
    - \bigl\| \IF^{-1/2} \nabla \zeta \bigr\|^{2}
    & \leq &
    \frac{\dltwb}{1 - \dmax^{2} \hspm \dltwb} \bigl\| \DPN \IF^{-1} \nabla \zeta \bigr\|^{2} \, ,
    \\
    2 L(\tilde{\upsv}) - 2 L(\upsvs) 
    - \bigl\| \IF^{-1/2} \nabla \zeta \bigr\|^{2}
    & \geq &
    - \frac{\dltwb}{1 + \dmax^{2} \hspm \dltwb} \bigl\| \DPN \IF^{-1} \nabla \zeta \bigr\|^{2} .
\label{3d3Af12DGttG}
\end{EQ}
Also
\begin{EQ}[rcl]
    \bigl\| \DPN \bigl( \tilde{\upsv} - \upsvs - \IF^{-1} \nabla \zeta \bigr) \bigr\|
    & \leq &
    \frac{\sqrt{3 \dltwb}}{1 - \dmax^{2} \hspm \dltwb} \, \bigl\| \DPN \IF^{-1} \nabla \zeta \bigr\| \, ,
    \\
    \bigl\| \DPN \bigl( \tilde{\upsv} - \upsvs \bigr) \bigr\|
    & \leq &
    \frac{1 + \sqrt{3 \dltwb}}{1 - \dmax^{2} \hspm \dltwb} \, \bigl\| \DPN \IF^{-1} \nabla \zeta \bigr\| \, .
\label{DGttGtsGDGm13rG}
\end{EQ}
\end{theorem}

\Section{Expansions and risk bounds under third-order smoothness}
\label{SFiWiexs3}

The results of Theorem~\ref{TFiWititG} can be refined under third-order smoothness conditions.
%
%
Namely, Proposition~\ref{Pconcgeneric2} yields the following Wilks expansion for the MLE \( \tilde{\upsv} \).

\begin{theorem}
\label{TFiWititG2}
Assume 
\nameref{Eref},
\nameref{EU2ref},
and 
\nameref{LLref}.
Let also \nameref{LL3tref} hold at \( \upsvs \)  
with a metric tensor \( \DPN \) and values \( \rr \) and \( \dltwu_{3} \) satisfying 
\begin{EQA}
	\DPN^{2} \leq \dmax^{2} \, \IF , 
	\quad \rrn \geq \frac{4\dmax}{3} \, \rrB ,
	\quad 
	\dltwu_{3} \, \dmax^{3} \, \rrB
	& < &
	\frac{1}{4} \, ,
\label{yxdhewndu7jwnjjuuMLE}
\end{EQA}
for \( \rrB \) from \eqref{rm1ucDGu2r0DGu}.
Then on \( \Omega(\xx) \), it holds 
\begin{EQA}
	\| \IF^{1/2} (\tilde{\upsv} - \upsvs) \|  
	& \leq &
	\frac{4}{3} \rrB \, ,
	\qquad
	 \| \DPN (\tilde{\upsv} - \upsvs) \|  
	\leq 
	\frac{4\dmax}{3} \, \rrB \, ,
\label{rhDGtuGmusGU0a2MLE}
\end{EQA}
and
\begin{EQ}[rcl]
    \Bigl| 2 L(\tilde{\upsv}) - 2 L(\upsvs) - \| \IF^{-1/2} \nabla \zeta \|^{2} \Bigr|
    & \leq &
    \frac{\dltwu_{3}}{2} \, \| \DPN \IF^{-1} \nabla \zeta \|^{3} 
    \, .
\label{3d3Af12DGttG2}
\end{EQ}
\end{theorem}

Under \nameref{LLsT3ref}, Proposition~\ref{PFiWigeneric2} yields an advanced Fisher expansion.
Define
\begin{EQA}[c]
	\BBH_{\DPN} = \DPN \IF^{-1} \VP^{2} \IF^{-1} \DPN,
	\\
	\dimD
	\eqdef
	\tr \BBH_{\DPN} \, ,
	\quad
	\rrD \eqdef \zq(\BBH_{\DPN},\xx)
	\leq 
	\sqrt{\tr \BBH_{\DPN}} + \sqrt{2\xx \, \| \BBH_{\DPN} \|} \, ;
\label{y7s7d7d77dfdy7fuegue3j}
\end{EQA}
cf. \eqref{34rtyghuioiuyhgvftid}.
By \nameref{EU2ref}, it holds 
\( \P(\| \DPN \, \IF^{-1} \nabla \zeta \| > \rrD) \leq 3\ex^{-\xx} \).
The result follows by limiting to the set \( \Omega(\xx) \) on which \( \| \DPN \, \IF^{-1} \nabla \zeta \| \leq \rrD \)
and by applying Proposition~\ref{PFiWigeneric2}.

\begin{theorem}
\label{TFiWititG3}
Assume 
\nameref{Eref},
\nameref{EU2ref},
and 
\nameref{LLref}.
Let \nameref{LLsT3ref} hold at \( \upsvs \) with a metric tensor 
\( \DPN \) and values \( \rr \) and \( \dltwu_{3} \) 
satisfying
\begin{EQA}[c]
	\DPN^{2} \leq \dmax^{2} \hspm \IF ,
	\quad
	\rr \geq \frac{3}{2} \, \rrD \, ,
	\quad
	\dltwu_{3} \, \dmax^{2} \hspm \rrD < \frac{4}{9} \, ,
\label{8difiyfc54wrboeMLE}
\end{EQA}
where \( \rrD \) is from \eqref{y7s7d7d77dfdy7fuegue3j}.
With \( \Omega(\xx) = \{ \| \DPN \, \IF^{-1} \nabla \zeta \| \leq \rrD \} \), 
it holds \( \P(\Omega(\xx)) \geq 1 - 3 \ex^{-\xx} \) and on \( \Omega(\xx) \)
\begin{EQA}[rcl]
    \| \DPN^{-1} \IF (\tilde{\upsv} - \upsvs - \IF^{-1} \nabla \zeta) \|
    & \leq &
    \frac{3\dltwu_{3}}{4} \| \DPN \, \IF^{-1} \nabla \zeta \|^{2} 
	\, .
\label{DGttGtsGDGm13rG22}
\end{EQA}
%
\end{theorem}

Expansion \eqref{DGttGtsGDGm13rG22} yields accurate risk bounds.
\begin{theorem}
\label{TQFiWi}
Assume \nameref{Eref},
\nameref{EU2ref},
and 
\nameref{LLref}.
Let \( \fs(\upsv) = \E L(\upsv) \) satisfy \nameref{LLsT3ref} at \( \upsvs \) with some
\( \DPN \), \( \rr \), and \( \dltwu_{3} \).
Let also
\begin{EQA}[c]
	\DPN^{2} \leq \dmax^{2} \, \IF \, ,
	\qquad 
	\rr \geq \frac{3}{2} \rrD \, ,
	\qquad
	\dmax^{2} \dltwu_{3} \, \rrD < \frac{4}{9} \, ;
\label{7deyedhfg5563w6rfygjLR}
\end{EQA}
see \eqref{y7s7d7d77dfdy7fuegue3j}.
For any linear mapping \( \QP \colon \R^{\dimp} \to \R^{\dimq} \), 
it holds on \( \Omega(\xx) \)
\begin{EQA}
	\| \QP (\tilde{\upsv} - \upsvs - \IF^{-1} \nabla \zeta) \|
	& \leq &
	\| \QP \IF^{-1} \DPN \| \, \frac{3\dltwu_{3}}{4} \, \| \DPN \IF^{-1} \nabla \zeta \|^{2}
	\, .
	\qquad
	\quad
\label{g25re9fjfregdndg}
\end{EQA}
Also, introduce 
\begin{EQA}[c]
	\riskt_{\QP} \eqdef \E \{ \| \QP \IF^{-1} \nabla \zeta \|^{2} \Ind_{\Omega(\xx)} \} 
	\leq 
	\dimQ  \, 
\label{7djhed8cjfct534etgdhdyQP}
\end{EQA}
with \( \dimQ \eqdef \E \| \QP \IF^{-1} \nabla \zeta \|^{2} = \tr \Var(\QP \IF^{-1} \nabla \zeta) \).
Then
\begin{EQA}
	\E \bigl\{ \| \QP (\tilde{\upsv} - \upsvs) \| \Ind_{\Omega(\xx)} \bigr\} 
	& \leq &
	\riskt_{\QP}^{1/2} 
	+ \| \QP \IF^{-1} \DPN \| \, \frac{3\dltwu_{3}}{4} \, \dimD \, .
\label{EtuGus11md3GQP}
\end{EQA}
Further, assume \( \E \bigl\{ \| \DPN \IF^{-1} \nabla \zeta \|^{4} \Ind_{\Omega(\xx)} \bigr\} \leq \CONSTi_{4}^{2} \, \dimD^{2} \) 
and define
\begin{EQ}[rcl]
	\alp_{\QP}
	& \eqdef & 
	\frac{\| \QP \IF^{-1} \DPN \| \, (3/4) \dltwu_{3} \, \CONSTi_{4} \, \dimD } {\sqrt{\riskt_{\QP}}}
	\, . 
\label{6dhx6whcuydsds655srew}
\end{EQ}
If \( \alp_{\QP} < 1 \) then
\begin{EQA}
	(1 - \alp_{\QP})^{2} \riskt_{\QP} 
	\leq 
	\E \bigl\{ \| \QP \, (\tilde{\upsv} - \upsvs) \|^{2} \Ind_{\Omega(\xx)} \bigr\}
	& \leq &
	(1 + \alp_{\QP})^{2} \riskt_{\QP} \, .
\label{EQtuGmstrVEQtGQP}
\end{EQA}
\end{theorem}

\Section{Effective and critical dimension in ML estimation}
\label{ScritdimMLE}
This section discusses the important question of the critical parameter dimension 
still ensuring the validity of the presented results.
To be more specific, we only consider the 3S-results of Theorem~\ref{TFiWititG3}.
Also, assume \( \dmax \equiv 1 \).
The important constant \( \dltwu_{3} \) is identified by \nameref{LLtS3ref}: \( \dltwu_{3} = \hmax_{3} /\sqrt{n} \),
where the scaling factor \( n \) means the sample size.
It can be defined as the smallest eigenvalue of the Fisher operator \( \IF \).

First, we discuss the case \( \QP = \DPN = \IF^{1/2} \).
It appears that in this full dimensional situation, all the obtained results apply and are meaningful under the condition \( \dimA \ll n \),
where \( \dimA = \tr(\BBH) \) for \( \BBH = \IF^{-1/2} \VP^{2} \IF^{-1/2} \) is the \emph{effective dimension} of the problem.
Indeed,  \( \rrD^{2} = \rrB^{2} \approx \tr(\BBH) = \dimA \), and
condition \eqref{8difiyfc54wrboeMLE} requires \( \dltwu_{3} \, \rrD \ll 1 \) which can be spelled out as \( \dimA \ll n \).
Expansion \eqref{DGttGtsGDGm13rG22} means
\begin{EQA}
	\| \IF^{1/2} (\tilde{\upsv} - \upsvs) \|
    & \leq &
    \| \IF^{-1/2} \nabla \zeta \| + \frac{3\dltwu_{3}}{4} \| \IF^{-1/2} \nabla \zeta \|^{2} \, ,
\label{Y6DF66FCC6dseeEGHBE}
\end{EQA}
and the second term on the right-hand side of this bound is smaller 
than the first one under the same condition \( \dltwu_{3} \, \rrD \ll 1 \).
Similar observations apply to bound \eqref{EQtuGmstrVEQtGQP} of Theorem~\ref{TQFiWi} 
which is meaningful only if \( \alp_{\QP} \) in \eqref{6dhx6whcuydsds655srew} is small.
As \( \riskt_{\QP} \approx \dimQ = \dimA \), the condition \( \dltwu_{3} \, \rrD \ll 1 \) implies \( \alp_{\QP} \ll 1 \)
and hence, the bound \eqref{EQtuGmstrVEQtGQP} is sharp.
We conclude that the main properties of the MLE \( \tilde{\upsv} \) 
are valid under the condition \( \dimA \ll n \) meaning sufficiently many observations 
per effective number of parameters.

The situation changes drastically if \( \QP \) is not full-dimensional as e.g. in semiparametric estimation,
when \( \QP \) projects onto a low-dimensional target component.
We will see in \Chname \ref{SsemiMLE} that in this case,
\eqref{6dhx6whcuydsds655srew} requires \( \dimA^{2} \ll n \).
\iffourG{
An interesting question about a further improvement of the error term in \eqref{g25re9fjfregdndg}
will be discussed in the next section.}{}

\iffourG{
\Section{Bounds under fourth-order smoothness}
\label{SMLE4}
This section explains how the accuracy of the expansions for MLE can be improved
and the critical dimension condition can be relaxed under 
fourth-order smoothness of \( \fs(\upsv) = \E L(\upsv) \).

Consider the third-order tensor
\( \TensU(\uv) = \frac{1}{6} \langle \nabla^{3} \fs(\upsvs), \uv^{\otimes 3} \rangle \) and its gradient
\( \nabla \TensU(\uv) = \frac{1}{2} \langle \nabla^{3} \fs(\upsvs), \uv^{\otimes 2} \rangle \). 
Define a random vector \( \svn \) by
\begin{EQ}[rcl]
	\svn 
	&=&
	\IF^{-1} \nabla \zeta + \IF^{-1} \nabla \TensU(\IF^{-1} \nabla \zeta) 
	\, .
\label{8vfjvr43f8khg54ed54}
\end{EQ}
The next result shows that the use of \( \svn \) in place of \( \IF^{-1} \nabla \zeta \)
allows to improve the accuracy of the Fisher expansion \eqref{DGttGtsGDGm13rG22} and of the Wilks expansion \eqref{3d3Af12DGttG2}.

\begin{theorem}
\label{Teffp4s}
Assume \nameref{Eref}, \nameref{LLref}, and \nameref{EU2ref}.
Let \nameref{LLsT3ref} and \nameref{LLsT4ref} hold at \( \upsvs \) and 
\begin{EQA}[c]
	\DPN^{2} \leq \dmax^{2} \, \IF \, ,
	\;\; 
	\rr \geq \frac{3}{2} \rrD \, ,
	\;\;
	\dmax^{2} \dltwu_{3} \, \rrD < \frac{4}{9} \, ,
	\;\;
	\dmax^{2} \dltwu_{4} \, \rrD^{2} < \frac{1}{3} 
	\, ,
\label{7deyedhfg5563w6rfygjLR4}
\end{EQA}
with \( \rrD \) from \eqref{y7s7d7d77dfdy7fuegue3j}.
Then \( \svn \) from \eqref{8vfjvr43f8khg54ed54} fulfills on \( \Omega(\xx) \)
\begin{EQA}[rcl]
\label{yfjvh9h4f5e53tghfugu44}
	\| \DPN^{-1} \IF ( \tilde{\upsv} - \upsvs - \svn) \|
    & \leq &
    \Bigl( \frac{\dltwu_{4}}{2} + \dmax^{2} \dltwu_{3}^{2} \Bigr) \, \| \DPN \IF^{-1} \nabla \zeta \|^{3} \, ,
	\\
	\| \DPN^{-1} \IF \, (\svn - \IF^{-1} \nabla \zeta) \|
	&=&
	\| \DPN^{-1} \nabla \TensU(\IF^{-1} \nabla \zeta) \|
	\leq 
	\frac{\dltwu_{3}}{2} \, \| \DPN \IF^{-1} \nabla \zeta \|^{2}
	\, ,
	\qquad
\label{iuvchycvf6e64rygh3224}
\end{EQA}
and
\begin{EQA}
	&& \nquad
	\bigl| 
		L(\tilde{\upsv}) - L(\upsvs) 
		- \frac{1}{2} \| \IF^{-1/2} \nabla \zeta \|^{2} - \TensU(\IF^{-1} \nabla \zeta) 
	\bigr|
	\\
	& \leq &
	\frac{\dltwu_{4} + 4 \dmax^{2} \dltwu_{3}^{2}}{8} \| \DPN \IF^{-1} \nabla \zeta \|^{4} 
    + \frac{\dmax^{2} (\dltwu_{4} + 2 \dmax^{2} \dltwu_{3}^{2})^{2}}{4} \, \| \DPN \IF^{-1} \nabla \zeta \|^{6} \, .
    \qquad
\label{87dfsjqweudsfjht7w3}
\end{EQA}
\end{theorem}

\begin{proof}
See Proposition~\ref{Pconcgeneric4} with \( \Av = \nabla \zeta \) and \( \IFN = \IF \).
\end{proof}

The obtained expansion yields the bound on the loss and risk of \( \tilde{\upsv} \).
Define 
\begin{EQ}[rcl]
	\riskt_{\QP}
	& \eqdef &
	\E \bigl\{ \| \QP \IF^{-1} \nabla \zeta \|^{2} \Ind_{\Omega(\xx)} \bigr\} \, ,
\label{7dhdrw2bfvu78u78e4ndwQ}
	\\
\label{hdvje39bug53ebfh8edx}
	\riskt_{\QP,2} 
	& \eqdef & 
	\E \bigl\{ \| \QP \svn \|^{2} \Ind_{\Omega(\xx)} \bigr\} \, .
\end{EQ}

\begin{theorem}
\label{Teff41}
Assume the conditions of Theorem~\ref{Teffp4} and let
\begin{EQA}[c]
	\E \bigl\{ \| \DPN \IF^{-1} \nabla \zeta \|^{k} \Ind_{\Omega(\xx)} \bigr\} \leq \CONSTi_{k}^{2} \, \dimD^{k/2} ,
	\qquad
	k=3,4,6 
	\, .
\label{6hjdfv8e6hyefyeheew7sk}
\end{EQA}
Then it holds for any linear mapping \( \QP \)
\begin{EQ}[rcl]
	&& \nquad
	\E \bigl\{ \| \QP \, (\tilde{\upsv} - \upsvs) \| \Ind_{\Omega(\xx)} \bigr\}
	\leq 	 
	\E \bigl\{ \| \QP \svn \| \Ind_{\Omega(\xx)} \bigr\}
	+ \| \QP \IF^{-1} \DPN \| \, 
	\Bigl( \frac{\dltwu_{4}}{2} + \dmax^{2} \dltwu_{3}^{2} \Bigr) \, \CONSTi_{3}^{2} \, \dimD^{3/2} \, ,
\label{0mkvhgjnrw3dfwe3u8gtygE1}
	\\
	&& \nquad
	\Bigl| \E \bigl\{ \| \QP \svn \| \Ind_{\Omega(\xx)} \bigr\}
	- \E \bigl\{ \| \QP \IF^{-1} \, \nabla \zeta \| \Ind_{\Omega(\xx)} \bigr\} 
	\Bigr|
	\leq 
	\| \QP \IF^{-1} \DPN \| \, \frac{\dltwu_{3}}{2} \, \dimD  \, .
\label{udtgecthwjdytdehduuc6}
\end{EQ}
With \( \riskt_{\QP,2} \) from \eqref{7dhdrw2bfvu78u78e4ndwQ}, let
\begin{EQA}
	\alp_{\QP,2}
	& \eqdef &
	\frac{\| \QP \IF^{-1} \DPN \| \, ( \dltwu_{4}/2 + \dmax^{2} \dltwu_{3}^{2}) \, \CONSTi_{6} \, \dimD^{3/2}}
		 {\sqrt{\riskt_{\QP,2}}} 
	< 1 
	\, .
\label{yfhcvhched6chejdrte}
\end{EQA}
Then
\begin{EQA}
	\bigl( 1 - \alp_{\QP,2} \bigr)^{2} \riskt_{\QP,2}
	\leq 
	\E \bigl\{ \| \QP \, (\tilde{\upsv} - \upsvs) \|^{2} \Ind_{\Omega(\xx)} \bigr\}
	& \leq &
	\bigl( 1 + \alp_{\QP,2} \bigr)^{2} \riskt_{\QP,2} \, .
\label{6shx76whnjvyehfbvyfh}
\end{EQA}
If another constant \( \alp_{\QP,1} < 1 \) ensures 
\begin{EQ}[rcl]
	&&
	\| \QP \IF^{-1} \DPN \| \, \frac{\dltwu_{3}}{2} \, \CONSTi_{4} \, \dimD 
	\leq 
	\alp_{\QP,1} \, \sqrt{\riskt_{\QP}} \, 
\label{6dhx6whcuydsds655srew4}
\end{EQ}
with \( \riskt_{\QP} \) from \eqref{7dhdrw2bfvu78u78e4ndwQ} then
\begin{EQA}
	\riskt_{\QP} (1 - \alp_{\QP,1})^{2} 
	\leq 
	\riskt_{\QP,2}
	& \leq &
	\riskt_{\QP} (1 + \alp_{\QP,1})^{2} \, .
\label{EQtuGmstrVEQtGQ2}
\end{EQA}
\end{theorem}

\begin{proof}
Rescaling of \( \DPN \) reduces the proof to \( \dmax = 1 \).
Theorem~\ref{Teffp4s} yields
\begin{EQA}
	\| \QP \, (\tilde{\upsv} - \upsvs - \svn) \|
	& \leq &
	\| \QP \IF^{-1} \DPN \| \, 
	\Bigl( \frac{\dltwu_{4}}{2} + \dltwu_{3}^{2} \Bigr) \,  
	\| \DPN \IF^{-1} \nabla \zeta \|^{3} \, ,
	\qquad
	\qquad
\label{0mkvhgjnrwwe3u8gtyg}
	\\
	\| \QP \{ \svn - \IF^{-1} \nabla \zeta \} \|
	& \leq &
	\frac{\dltwu_{3}}{2} \, \| \QP \IF^{-1} \DPN \| \, \| \DPN \IF^{-1} \nabla \zeta \|^{2} \, .
\label{vuedy766t4e3bfvyt6e}
\end{EQA}
Now \eqref{0mkvhgjnrw3dfwe3u8gtygE1} follows from \eqref{6hjdfv8e6hyefyeheew7sk} with \( k=3 \).
Next, we study the quadratic risk of \( \tilde{\upsv} \).
Define \( \epsv_{\QP} = \QP (\tilde{\upsv} - \upsvs - \svn) \). 
By \eqref{0mkvhgjnrwwe3u8gtyg} 
\begin{EQA}
	\sqrt{\E ( \| \epsv \|^{2} \Ind_{\Omega(\xx)} )}
	& \leq &
	\| \QP \IF^{-1} \DPN \| \, 
	\frac{\dltwu_{4}}{2} \,  
	\sqrt{\E \| \DPN \IF^{-1} \nabla \zeta \|^{6} \Ind_{\Omega(\xx)}} 
	\leq 
	\alp_{\QP,2} \sqrt{\riskt_{\QP,2}} \, ,
\label{d8ew3jfjhyvb6rt6543ejhvu}
\end{EQA}
and \eqref{6shx76whnjvyehfbvyfh} follows. 
Further, denote 
\begin{EQA}[rcccl]
	\loss_{\QP}
	& \eqdef &
	\QP \IF^{-1} \nabla \zeta \, ,
	\qquad
	\delta_{\QP}
	& \eqdef &
	\QP (\IF^{-1} \nabla \zeta - \svn) \, .
\label{7ejvuejfumfmrgvytweof7}
\end{EQA}
By definition, \( \riskt_{\QP} = \E \bigl\{ \| \loss_{\QP} \|^{2} \Ind_{\Omega(\xx)} \bigr\} \), 
\( \riskt_{\QP,2} = \E \bigl\{ \| \loss_{\QP} + \delta_{\QP} \|^{2} \Ind_{\Omega(\xx)} \bigr\} \), and 
\begin{EQA}
	\riskt_{\QP,2} - \riskt_{\QP}
	& = &
	\E \bigl\{ \| \delta_{\QP} \|^{2} \Ind_{\Omega(\xx)} \bigr\}
	+ 2 \E \bigl\{ \langle \loss_{\QP} , \delta_{\QP} \rangle \Ind_{\Omega(\xx)} \bigr\} .
\label{7jde7jevyreteyv8rnbe}
\end{EQA}
Also \eqref{iuvchycvf6e64rygh3224} and \eqref{6dhx6whcuydsds655srew4} imply
\begin{EQA}
	\sqrt{\E \bigl( \| \delta_{\QP} \|^{2} \Ind_{\Omega(\xx)} \bigr) }
	& \leq &
	\| \QP \IF^{-1} \DPN \| \, \frac{\dltwu_{3}}{2} \, 
	\sqrt{\E \| \DPN \IF^{-1} \nabla \zeta \|^{4} \Ind_{\Omega(\xx)} } 
	\\
	& \leq &
	\| \QP \IF^{-1} \DPN \| \, \frac{\dltwu_{3}}{2} \, \CONSTi_{4} \, \dimD 
	\leq 
	\alp_{\QP,1} \, \sqrt{\riskt_{\QP}} \, .
\label{yfhf73hjf9bryrnvbir}
\end{EQA}
This proves \eqref{EQtuGmstrVEQtGQ2}.
\end{proof}

\begin{remark}
As \( \| \DPN \, \IF^{-1} \nabla \zeta \| \leq \rrD \) on \( \Omega(\xx) \), it holds 
\begin{EQA}
	\E \bigl( \| \DPN \IF^{-1} \nabla \zeta \|^{4} \Ind_{\Omega(\xx)} \bigr)
	& \leq &
	\rrD^{2} \, \E \bigl( \| \DPN \IF^{-1} \nabla \zeta \|^{2} \Ind_{\Omega(\xx)} \bigr)
	\leq 
	\rrD^{2} \, \dimD \, .
\label{6hdnweyftegwgfywggvf636}
\end{EQA}
If \( \rrD^{2} \approx \dimD \), then  \( \CONSTi_{4} \approx 1 \) in \eqref{6dhx6whcuydsds655srew}.
\end{remark}

The results of Theorem~\ref{Teff41} enable us to improve the issue of \emph{critical dimension}.
For simplicity, let \( \QP = \DPN = \IF^{1/2} \).
Then the derived bounds are meaningful if 
\begin{EQA}
	( \dltwu_{4} + \dltwu_{3}^{2} ) \, \dimD^{3/2} 
	&=&
	o(1) .
\label{ufcjnciwf7fnrvuehj}
\end{EQA} 
Assuming \( \dltwu_{4} \asymp 1/n \) and \( \dltwu_{3}^{2} \asymp 1/n \),
we obtain the critical dimension condition
\( \dimD^{3/2} \ll n \) which is weaker than \( \dimD^{2} \ll n \).
Condition \eqref{6dhx6whcuydsds655srew4} ensuring equivalence of \( \riskt_{\QP,2} \) and \( \riskt_{\QP} \)
requires \( \dltwu_{3} \, \dimD \ll \riskt_{\QP} \) as in the 3S case.
}{}

\iffourG{
\Section{Fourth-order expansion and bias correction}
\label{SMLE4b}

Expansion \eqref{8vfjvr43f8khg54ed54} and \eqref{yfjvh9h4f5e53tghfugu44}
\begin{EQA}
	\tilde{\upsv} - \upsvs
	& \approx &
	\svn
	=
	\IF^{-1} \nabla \zeta + \IF^{-1} \nabla \TensU(\IF^{-1} \nabla \zeta)
\label{dfv8edefye3hejxs2qw}
\end{EQA}
involves the term \( \IF^{-1} \nabla \TensU(\IF^{-1} \nabla \zeta) \) which is quadratic in \( \IF^{-1} \nabla \zeta \).
If the dimension \( \dimp \) is large, one can expect that this quadratic form concentrates around its expectation.
This suggests a third-order correction of the MLE \( \tilde{\upsv} \):
\begin{EQA}
	\hat{\upsv}
	&=&
	\tilde{\upsv} - \IF^{-1} \, \E \, \nabla \TensU(\IF^{-1} \nabla \zeta),
	\\
	\tilde{\upsv} - \upsvs
	& \approx &
	\IF^{-1} \nabla \zeta + \IF^{-1} \bigl\{ \nabla \TensU(\IF^{-1} \nabla \zeta) - \E \, \nabla \TensU(\IF^{-1} \nabla \zeta) \bigr\} .
\label{tsgqhyug8g900o4kqtdje}
\end{EQA}
Define
\begin{EQA}
	\Zv
	& \eqdef &
	\nabla \TensU(\IF^{-1} \nabla \zeta),
	\\
	\xiv
	& \eqdef &
	\Zv - \E \Zv .
\label{di8d7vy7j3cciedrwsfy}
\end{EQA}

\begin{lemma}
\label{LTensmoments}
Let \( \Var(\nabla \zeta) \leq \VP^{2}  \).
For any \( \uv \in \R^{\dimp} \), it holds on \( \Omega(\xx,\uv) \) 
\begin{EQA}[c]
	| \uv^{\T} \IF^{-1} \xiv |
	\leq 
	\dltwu_{3} \| \DPN \, \IF^{-1} \uv \| \, \bigl( \sqrt{\dimD} + \sqrt{2\xx} \bigr)
\label{xdyctc4332egb8gfjhhg}
\end{EQA}

\end{lemma}

\begin{proof}
For any \( \uv \in \R^{\dimp} \), define \( \TensU[\uv] = \sum_{j=1}^{\dimp} u_{j} \TensU_{j} \).
Then 
\begin{EQA}
	\uv^{\T} \IF^{-1} \nabla \TensU(\IF^{-1} \nabla \zeta) 
	&=&
	\langle \TensU[\IF^{-1} \uv], (\IF^{-1} \nabla \zeta)^{\otimes 2} \rangle
\label{uduwdywyey8vy8ryf4ih}
\end{EQA}
and hence, \( \uv^{\T} \IF^{-1} \nabla \TensU(\IF^{-1} \nabla \zeta) \) 
is a quadratic form of \( \nabla \zeta \).
If \( \| \uv \| = 1 \) then
\begin{EQA}
	&& \nquad
	\E \bigl\{ \uv^{\T} \IF^{-1} \nabla \TensU(\IF^{-1} \nabla \zeta) \bigr\}
	=
	\E \langle \TensU[\IF^{-1} \uv], (\IF^{-1} \nabla \zeta)^{\otimes 2} \rangle
	\\
	&=&
	\tr \bigl\{ \TensU[\IF^{-1} \uv] \, \Var(\IF^{-1} \nabla \zeta) \bigr\}
	=
	\tr \bigl\{ \TensU[\IF^{-1} \uv] \, \IF^{-1} \Var(\nabla \zeta) \IF^{-1} \bigr\}
	\\
	& \leq &
	\dltwu_{3} \| \DPN \, \IF^{-1} \uv \| \tr \bigl\{ \DPN \, \IF^{-1} \Var(\nabla \zeta) \IF^{-1} \DPN \bigr\}
\label{ushduwntuyu3geerbxnc}
\end{EQA}
Further, \( \uv^{\T} \IF^{-1} \Zv \) is a quadratic form of \( \nabla \zeta \); see \eqref{uduwdywyey8vy8ryf4ih}.
Therefore, assuming \( \VP \leq \DPN \)
\begin{EQA}
	\uv^{\T} \IF^{-1} \Zv
	&=&
	\langle \TensU[\IF^{-1} \uv], (\IF^{-1} \DPN \, \DPN^{-1} \nabla \zeta)^{\otimes 2} \rangle
	=
	\langle \BBH, (\DPN^{-1} \nabla \zeta)^{\otimes 2} \rangle
\label{edfchy8uiehfi8w3erj3i}
\end{EQA}
with \( \BBH = \DPN \, \IF^{-1} \TensU[\IF^{-1} \uv] \, \IF^{-1} \DPN \).
By \nameref{EU2ref}, it holds on \( \Omega(\xx,\uv) \)
\begin{EQA}
	\bigl| \uv^{\T} \IF^{-1} \xiv \bigr|
	&=&\
	\bigl| \uv^{\T} \IF^{-1} \Zv - \E \uv^{\T} \IF^{-1} \Zv \bigr|
	\leq 
	\sqrt{\xx \tr \BBH^{2}} + 2 \xx \| \BBH \|
\label{ytd7w3y8ru3igytkhgy8}
\end{EQA}
\nameref{LLsT3ref} implies in view of \( \DPN^{2} \leq \IF \)
\begin{EQA}
	\| \BBH \|
	& \leq &
	\dltwu_{3} \| \DPN \, \IF^{-1} \uv \| \, \| \DPN \, \IF^{-1} \DPN \|^{2}
	\leq 
	\dltwu_{3} \| \DPN \, \IF^{-1} \uv \|,
	\\
	\tr \BBH
	& \leq &
	\dltwu_{3} \| \DPN \, \IF^{-1} \uv \| \, \tr (\DPN \, \IF^{-1} \DPN)^{2}
	\leq 
	\dltwu_{3} \| \DPN \, \IF^{-1} \uv \| \, \tr (\DPN \, \IF^{-1} \DPN)
	\leq 
	\dltwu_{3} \| \DPN \, \IF^{-1} \uv \| \, \dimD \, ,
	\\
	\tr \BBH^{2}
	& \leq &
	\dltwu_{3} \| \DPN \, \IF^{-1} \uv \| \, \tr \BBH
	\leq 
	\dltwu_{3}^{2} \| \DPN \, \IF^{-1} \uv \|^{2} \, \dimD \, ,
\label{f8jw3f87vjdwnjtgcnwr}
\end{EQA}
and the assertion follows.
\end{proof}
}{}

\Section{Penalization bias}
\label{Ssmoothbias}

A common approach for improving the performance of MLE is based on regularization or penalization. 
The objective function \( L(\upsv) \) is extended by including a penalty term \( \pent_{\GP}(\upsv) \) 
which is responsible for complexity (roughness) of the parameter \( \upsv \).
A typical example to keep in mind is \( \pent_{\GP}(\upsv) = \| \GP \upsv \|^{2}/2 \)
for a penalization matrix \( \GP^{2} \).
Penalization by \( \pent_{\GP}(\upsv) \) can gradially improve  stability and numerical properties of the estimator,
however, it leads to a change of the ``truth'' \( \upsvs \), and hence, to some bias.
This section describes the bias caused by a smooth penalty.
Define the penalized MLE \( \tilde{\upsv}_{\GP} \)
\begin{EQA}
	\tilde{\upsv}_{\GP}
	& \eqdef &
	\argmax_{\upsv} \LGP(\upsv)
	=
	\argmax_{\upsv} \{ L(\upsv) - \pent_{\GP}(\upsv) \} \, .
\label{vdyedbwdtf6g7h8y6ie}
\end{EQA}
Compared to \eqref{tuGauLGususGE}, consider three optimization problems
\begin{EQA}[rcl]
	\tilde{\upsv}_{\GP} 
	&=& 
	\argmax_{\upsv} \LGP(\upsv),
	\qquad
	\upsvs_{\GP} 
	=
	\argmax_{\upsv} \E \LGP(\upsv),
	\qquad
	\upsvs 
	=
	\argmax_{\upsv} \E L(\upsv) .
	\qquad
\label{tuGauLGususGEgp}
\end{EQA}
Due to 
{Proposition~\ref{PconcMLEgenc}}, 
the penalized MLE  \( \tilde{\upsv}_{\GP} \) estimates rather \( \upsvs_{\GP} \) 
than \( \upsvs \).
This section describes the bias \( \upsvs_{\GP} - \upsvs \) caused by penalization.

Define the penalized Fisher information \( \IF_{\GP}(\upsv) = - \nabla^{2} \E \LGP(\upsv) \) and
introduce \( \AvGP(\upsv) \eqdef \nabla \!\pent_{\GP}(\upsv) \).
Set \( \IF_{\GP} = \IF_{\GP}(\upsvs_{\GP}) \),
\begin{EQA}[rcccccl]
	\IF_{\GP} &=& \IF_{\GP}(\upsvs_{\GP}) \, ,
	\qquad
	\AvGP 
	& \eqdef &
	\nabla \!\pent_{\GP}(\upsvs) ,
	\quad
	\biasD
	& \eqdef &
	\| \DPN \IF_{\GP}^{-1} \AvGP \| 
	\, .
\label{ghdrd324ee4ew222}
\end{EQA}
For a quadratic penalty \( \pent_{\GP}(\upsv) = \| \GP \upsv \|^{2}/2 \), this results in
\begin{EQA}[rcccl]
	\AvGP
	&=&
	\GP^{2} \upsvs,
	\qquad
	\biasD
	& \eqdef &
	\| \DPN \IF_{\GP}^{-1} \GP^{2} \upsvs \| 
	\, .
\label{ghdrd324ee4ew222q}
\end{EQA}
Proposition~\ref{Pbiaspen} yields the following result.

\begin{proposition}
\label{Lvarusetb} 
Let \( \fG(\upsv) = \E \LGP(\upsv) \) satisfy \nameref{LLsT3ref} at \( \upsvs_{\GP} \) with some
metric tensor \( \DPN \) and values \( \rr \) and \( \dltwu_{3} \) such that
\begin{EQA}[c]
	\DPN^{2} \leq \dmax^{2} \, \IF_{\GP} \, ,
	\qquad 
	\rr \geq 3 \biasD/2 \, ,
	\qquad
	\dltwu_{3} \, \dmax^{2} \, \biasD < 4/9 ,
\label{7deyedhfg5563w6rfygj}
\end{EQA}
for \( \biasD \) from \eqref{ghdrd324ee4ew222}. 
Then 
\begin{EQA}[rcl]
\label{11ma3eaelDebQ}
	\| \DPN^{-1} \IF_{\GP} (\upsvs_{\GP} - \upsvs + \IF_{\GP}^{-1} \AvGP) \|
	& \leq &
	\frac{3\dltwu_{3}}{4} \, \biasD^{2} 
	\, .
\end{EQA}
The same bounds apply with \( \IF_{\GP}(\upsvs) \) in place of \( \IF_{\GP} = \IF_{\GP}(\upsvs_{\GP}) \).
\end{proposition}

\iffourG{
The accuracy of the expansions for the bias of a pMLE can be improved
under fourth-order smoothness of \( \fG(\upsv) = \E L_{\GP}(\upsv) \).
%
Consider the third-order tensor
\( \TensU(\uv) = \frac{1}{6} \langle \nabla^{3} \fs(\upsvs_{\GP}), \uv^{\otimes 3} \rangle \) and its gradient
\( \nabla \TensU(\uv) = \frac{1}{2} \langle \nabla^{3} \fs(\upsvs_{\GP}), \uv^{\otimes 2} \rangle \). 
Define a vector \( \bvn_{\GP} \) by
\begin{EQ}[rcl]
	\bvn_{\GP}
	&=&
	\IF_{\GP}^{-1} \AvGP + \IF_{\GP}^{-1} \nabla \Tens(\IF_{\GP}^{-1} \AvGP) 
	\, .
\label{8vfjvr43f8khg54ed5G}
\end{EQ}
Proposition~\ref{Pbiaspen} exlains the impact of using \( \bvn_{\GP} \) in place of \( \IF_{\GP}^{-1} \AvGP \).

\begin{theorem}
\label{Teffp4}
Assume \nameref{LLref}.
Let \nameref{LLsT3ref} and \nameref{LLsT4ref} hold at \( \upsvs_{\GP} \) and 
\begin{EQA}[c]
	\DPN^{2} \leq \dmax^{2} \, \IF_{\GP} \, ,
	\;\; 
	\rr \geq \frac{3}{2} \biasD \, ,
	\;\;
	\dmax^{2} \dltwu_{3} \, \biasD < \frac{4}{9} \, ,
	\;\;
	\dmax^{2} \dltwu_{4} \, \biasD^{2} < \frac{1}{3} 
	\, ,
\label{7deyedhfg5563w6rfygjLR4}
\end{EQA}
with \( \biasD = \| \DPN \IF_{\GP}^{-1} \AvGP \| \).
Then 
\begin{EQA}
\label{0mkvhgjnrw3dfwer45}
	\| \DPN^{-1} \IF_{\GP} (\upsvs_{\GP} - \upsvs + \bvn_{\GP}) \|
	& \leq &
	\Bigl( \frac{\dltwu_{4}}{2} + \dmax^{2} \dltwu_{3}^{2} \Bigr) \, \biasD^{3} \, ,
	\qquad
	\\
	\| \DPN^{-1} \IF_{\GP} (\bvn_{\GP} - \IF_{\GP}^{-1} \AvGP) \|
	& \leq &
	\frac{\dltwu_{3}}{2} \biasD^{2} \, .
\label{0mkvhgjnrw3dfwe3u8gtysg}
\end{EQA}
\end{theorem}
}{}

\Section{Loss and risk of the pMLE. Bias-variance decomposition}
\label{Slossrisksb}
Now we combine the previous results about the stochastic term \( \tilde{\upsv}_{\GP} - \upsvs_{\GP} \)
and the bias term \( \upsvs_{\GP} - \upsvs \) to obtain sharp bounds 
on the loss and risk of the pMLE \( \tilde{\upsv}_{\GP} \).

\begin{theorem}
\label{TQFiWibias}
Assume \nameref{Eref},
\nameref{EU2ref},
and 
\nameref{LLref}.
Let \( \fG(\upsv) = \E \LGP(\upsv) \) satisfy \nameref{LLsT3ref} at \( \upsvs_{\GP} \) with some
\( \DPN \), \( \rr \), and \( \dltwu_{3} \).
With \( (\rrD \vee \biasD) \eqdef \max\{ \rrD , \biasD \} \), assume
\begin{EQA}[c]
	\DPN^{2} \leq \dmax^{2} \, \IF_{\GP} \, ,
	\qquad 
	\rr \geq \frac{3}{2} (\rrD \vee \biasD) \, ,
	\qquad
	\dmax^{2} \dltwu_{3} \, (\rrD \vee \biasD) < \frac{4}{9} \, ;
\label{7deyedhfg5563w6rfygjLR}
\end{EQA}
see \eqref{y7s7d7d77dfdy7fuegue3j} and \eqref{ghdrd324ee4ew222}.
For any linear mapping \( \QP \colon \R^{\dimp} \to \R^{\dimq} \), 
it holds on \( \Omega(\xx) \)
\begin{EQA}
	\| \QP (\tilde{\upsv}_{\GP} - \upsvs - \IF_{\GP}^{-1} \nabla \zeta + \IF_{\GP}^{-1} \AvGP) \|
	& \leq &
	\| \QP \IF_{\GP}^{-1} \DPN \| \, \frac{3\dltwu_{3}}{4} \, \bigl( \| \DPN \IF_{\GP}^{-1} \nabla \zeta \|^{2} + \biasD^{2} \bigr) \,
	\, .
	\qquad
	\quad
\label{g25re9fjfregdndgb}
\end{EQA}
Also, introduce \( \dimQ \eqdef \E \| \QP \IF_{\GP}^{-1} \nabla \zeta \|^{2} = \tr \Var(\QP \IF_{\GP}^{-1} \nabla \zeta) \) and
\begin{EQA}[c]
	\riskt_{\QP} \eqdef \E \{ \| \QP \IF_{\GP}^{-1} (\nabla \zeta - \AvGP) \|^{2} \Ind_{\Omega(\xx)} \} 
	\leq 
	\dimQ + \| \QP \IF_{\GP}^{-1} \AvGP \|^{2} \, .
\label{7djhed8cjfct534etgdhdy}
\end{EQA}
Then
\begin{EQA}
	\E \bigl\{ \| \QP (\tilde{\upsv}_{\GP} - \upsvs) \| \Ind_{\Omega(\xx)} \bigr\} 
	& \leq &
	\riskt_{\QP}^{1/2} 
	+ \| \QP \IF_{\GP}^{-1} \DPN \| \, \frac{3\dltwu_{3}}{4} \bigl( \dimD + \biasD^{2} \bigr) \, .
\label{EtuGus11md3GQ}
\end{EQA}
Further, assume \( \E \bigl\{ \| \DPN \IF_{\GP}^{-1} \nabla \zeta \|^{4} \Ind_{\Omega(\xx)} \bigr\} \leq \CONSTi_{4}^{2} \, \dimD^{2} \) 
and define
\begin{EQ}[rcl]
	\alp_{\QP}
	& \eqdef & 
	\frac{\| \QP \IF_{\GP}^{-1} \DPN \| \, (3/4) \dltwu_{3} \, (\CONSTi_{4} \, \dimD + \biasD^{2})} {\sqrt{\riskt_{\QP}}}
	\, . 
\label{6dhx6whcuydsds655srewb}
\end{EQ}
If \( \alp_{\QP} < 1 \) then
\begin{EQA}
	(1 - \alp_{\QP})^{2} \riskt_{\QP} 
	\leq 
	\E \bigl\{ \| \QP \, (\tilde{\upsv}_{\GP} - \upsvs) \|^{2} \Ind_{\Omega(\xx)} \bigr\}
	& \leq &
	(1 + \alp_{\QP})^{2} \riskt_{\QP} \, .
\label{EQtuGmstrVEQtGQ}
\end{EQA}
\end{theorem}

\begin{proof}
It holds by \eqref{DGttGtsGDGm13rG22} and \eqref{11ma3eaelDebQ}
\begin{EQ}[rcl]
	\| \QP (\tilde{\upsv}_{\GP} - \upsvs_{\GP} - \IF_{\GP}^{-1} \nabla \zeta \|
	& \leq &
	\| \QP \IF_{\GP}^{-1} \DPN \| \, \frac{3\dltwu_{3}}{4} \, \| \DPN \IF_{\GP}^{-1} \nabla \zeta \|^{2} \, ,
	\\
	\| \QP (\upsvs_{\GP} - \upsvs + \IF_{\GP}^{-1} \AvGP) \|
	& \leq &
	\| \QP \IF_{\GP}^{-1} \DPN \| \, \frac{3\dltwu_{3}}{4} \, \biasD^{2} \, ,
\label{u87dcudenhedst6cbweq2q}
\end{EQ}
and hence,
\begin{EQA}
	\| \QP (\tilde{\upsv}_{\GP} - \upsvs - \IF_{\GP}^{-1} \nabla \zeta + \IF_{\GP}^{-1} \AvGP) \|
	& \leq &
	\| \QP \IF_{\GP}^{-1} \DPN \| \, \frac{3\dltwu_{3}}{4} \, \bigl( \| \DPN \IF_{\GP}^{-1} \nabla \zeta \|^{2} + \biasD^{2} \bigr)
\label{ydy6fc6fv6e43yte46fhuQ}
\end{EQA}
yielding \eqref{g25re9fjfregdndgb} and \eqref{EtuGus11md3GQ}.
Further, define 
\begin{EQA}
	\epsv_{\GP}
	& \eqdef &
	\QP \bigl\{ \tilde{\upsv}_{\GP} - \upsvs - \IF_{\GP}^{-1} (\nabla \zeta - \AvGP) \bigr\} \, .
\label{dfu8nj3eyvnhedd6ywhwhvfu7}
\end{EQA}
It holds by \eqref{u87dcudenhedst6cbweq2q}
\begin{EQA}
	\E^{1/2} \bigl\{ \| \epsv_{\QP} \|^{2} \Ind_{\Omega(\xx)} \bigr\}
	& \leq &
	\| \QP \IF_{\GP}^{-1} \DPN \| \frac{3\dltwu_{3}}{4} \, 
	\bigl\{ \E^{1/2} \bigl( \| \DPN \IF_{\GP}^{-1} \nabla \zeta \|^{4} \Ind_{\Omega(\xx)} \bigr) + \biasD^{2} \bigr\} 
	\leq 
	\alp_{\QP} \, \riskt_{\QP}^{1/2} \, ,
\label{d7hw3fujgv76444gdfu7eQ}
\end{EQA}
and therefore,
\begin{EQA}
	&& \nquad
	\E^{1/2} \bigl\{ \| \QP \, (\tilde{\upsv}_{\GP} - \upsvs) \|^{2} \Ind_{\Omega(\xx)} \bigr\}
	= 
	\E^{1/2} \bigl\{ \| \QP \IF_{\GP}^{-1} (\nabla \zeta - \AvGP) + \epsv_{\QP} \|^{2} \Ind_{\Omega(\xx)} \bigr\} 
	\\
	& \leq &
	\E^{1/2} \bigl\{ \| \QP \IF_{\GP}^{-1} (\nabla \zeta - \AvGP) \|^{2} \Ind_{\Omega(\xx)} \bigr\}
	+ \E^{1/2} \bigl\{ \| \epsv_{\QP} \|^{2} \Ind_{\Omega(\xx)} \bigr\}
	\leq 
	(1 + \alp_{\QP}) \, \riskt_{\QP}^{1/2} \, .
\label{7dnmswewe8fyhbdfc7er}
\end{EQA}
This yields \eqref{EQtuGmstrVEQtGQ}.
\end{proof}

\begin{remark}
\label{RremainderD}
The condition \( \DPN^{2} \leq \dmax^{2} \IF_{\GP} \) implies \( \| \QP \IF_{\GP}^{-1} \DPN \| \leq \dmax^{2} \| \QP \DPN^{-1} \| \)
which can be used in the remainder for all risk bounds. 
\end{remark}

\begin{remark}
Due to \eqref{EQtuGmstrVEQtGQ} 
\begin{EQA}
	\E \bigl\{ \| \QP \, (\tilde{\upsv}_{\GP} - \upsvs) \|^{2} \Ind_{\Omega(\xx)} \bigr\}
	& = &
	\bigl( \dimQ + \| \QP \IF_{\GP}^{-1} \AvGP \|^{2} \bigr) \, \bigl\{ 1 + o(1) \bigr\}.
	\qquad
\label{EQtuGuvus2tr1o1}
\end{EQA}
This relation is usually referred to as ``bias-variance decomposition''.
Our bound is sharp in the sense that for the special case of linear models, 
\eqref{EQtuGuvus2tr1o1} becomes equality.
Under the so-called ``small bias'' condition \( \| \QP \IF_{\GP}^{-1} \AvGP \|^{2} \ll \dimQ \), 
the impact of the bias induced by penalization is negligible. 
The relation \( \| \QP \IF_{\GP}^{-1} \AvGP \|^{2} \asymp \dimQ \) is called ``bias-variance trade-off'',
it leads to minimax rate of estimation;
see Section~\ref{Spriorexample}.
\end{remark}

If the constant \( \alp_{\QP} \) from \eqref{6dhx6whcuydsds655srewb} satisfies \( \alp_{\QP} \ll 1 \) then 
by \eqref{EQtuGmstrVEQtGQ}, 
\( \E \bigl\{ \| \QP \, (\tilde{\upsv}_{\GP} - \upsvs) \|^{2} \Ind_{\Omega(\xx)} \bigr\} = (1 + o(1)) \riskt_{\QP} \).
\iffourG{
Now we explain how the accuracy of the expansions for pMLE can be improved
and the critical dimension condition can be relaxed under 
fourth-order smoothness of \( \fG(\upsv) = \E L_{\GP}(\upsv) \).
%
Putting together the results on the stochastic component \( \tilde{\upsv}_{\GP} - \upsvs_{\GP} \) and 
on the bias \( \upsvs_{\GP} - \upsvs \) yields the bound on the loss and risk of the estimator \( \tilde{\upsv}_{\GP} \).
Define 
\begin{EQ}[rcl]
	\riskt_{\QP}
	& \eqdef &
	\E \bigl\{ \| \QP \IF_{\GP}^{-1} (\nabla \zeta - \AvGP) \|^{2} \Ind_{\Omega(\xx)} \bigr\} \, ,
\label{7dhdrw2bfvu78u78e4ndw}
	\\
\label{hdvje39bug53ebfh8edx}
	\riskt_{\QP,2} 
	& \eqdef & 
	\E \bigl\{ \| \QP (\svn_{\GP} - \bvn_{\GP}) \|^{2} \Ind_{\Omega(\xx)} \bigr\} \, .
\end{EQ}

\begin{theorem}
\label{Teff4s}
Assume the conditions of Theorem~\ref{Teffp4s} and Theorem~\ref{Teffp4} and let
\begin{EQA}[c]
	\E \bigl\{ \| \DPN \IF_{\GP}^{-1} \nabla \zeta \|^{k} \Ind_{\Omega(\xx)} \bigr\} \leq \CONSTi_{k}^{2} \, \dimD^{k/2} ,
	\qquad
	k=3,4,6 
	\, .
\label{6hjdfv8e6hyefyeheew7skb}
\end{EQA}
Then it holds for any linear mapping \( \QP \)
\begin{EQA}
	&& \nquad
	\E \bigl\{ \| \QP \, (\tilde{\upsv}_{\GP} - \upsvs) \| \Ind_{\Omega(\xx)} \bigr\}
	\\
	& \leq &	 
	\E \bigl\{ \| \QP (\svn_{\GP} - \bvn_{\GP}) \| \Ind_{\Omega(\xx)} \bigr\}
	+ \| \QP \IF_{\GP}^{-1} \DPN \| \, 
	\Bigl( \frac{\dltwu_{4}}{2} + \dmax^{2} \dltwu_{3}^{2} \Bigr) \, \bigl( \CONSTi_{3}^{2} \, \dimD^{3/2} + \biasD^{3} \bigr) ,
\label{0mkvhgjnrw3dfwe3u8gtygE1b}
	\\
	&& \nquad
	\Bigl| \E \bigl\{ \| \QP (\svn_{\GP} - \bvn_{\GP}) \| \Ind_{\Omega(\xx)} \bigr\}
	- \E \bigl\{ \| \QP \IF_{\GP}^{-1} \, \nabla \zeta - \QP \IF_{\GP}^{-1} \, \AvGP \| \Ind_{\Omega(\xx)} \bigr\} 
	\Bigr|
	\\
	& \leq &
	\| \QP \IF_{\GP}^{-1} \DPN \| \, \frac{\dltwu_{3}}{2} \, \bigl( \dimD + \biasD^{2} \bigr) \, .
\label{udtgecthwjdytdehduuc6}
\end{EQA}
With \( \riskt_{\QP,2} \) from \eqref{7dhdrw2bfvu78u78e4ndw}, let
\begin{EQA}
	\alp_{\QP,2}
	& \eqdef &
	\frac{\| \QP \IF_{\GP}^{-1} \DPN \| \, ( \dltwu_{4}/2 + \dmax^{2} \dltwu_{3}^{2}) \, (\CONSTi_{6} \, \dimD^{3/2} + \biasD^{3})}
		 {\sqrt{\riskt_{\QP,2}}} 
	< 1 
	\, .
\label{yfhcvhched6chejdrte}
\end{EQA}
Then
\begin{EQA}
	\bigl( 1 - \alp_{\QP,2} \bigr)^{2} \riskt_{\QP,2}
	\leq 
	\E \bigl\{ \| \QP \, (\tilde{\upsv}_{\GP} - \upsvs) \|^{2} \Ind_{\Omega(\xx)} \bigr\}
	& \leq &
	\bigl( 1 + \alp_{\QP,2} \bigr)^{2} \riskt_{\QP,2} \, .
\label{6shx76whnjvyehfbvyfhb}
\end{EQA}
If another constant \( \alp_{\QP,1} < 1 \) ensures 
\begin{EQ}[rcl]
	&&
	\| \QP \IF_{\GP}^{-1} \DPN \| \, \frac{\dltwu_{3}}{2} \, \bigl( \CONSTi_{4} \, \dimD + \biasD^{2} \bigr)
	\leq 
	\alp_{\QP,1} \, \sqrt{\riskt_{\QP}} \, 
\label{6dhx6whcuydsds655srew4b}
\end{EQ}
with \( \riskt_{\QP} \) from \eqref{7dhdrw2bfvu78u78e4ndw} then
\begin{EQA}
	\riskt_{\QP} (1 - \alp_{\QP,1})^{2} 
	\leq 
	\riskt_{\QP,2}
	& \leq &
	\riskt_{\QP} (1 + \alp_{\QP,1})^{2} \, .
\label{EQtuGmstrVEQtGQ2b}
\end{EQA}
\end{theorem}

\begin{proof}Rescaling of \( \DPN \) reduces the proof to \( \dmax = 1 \).
Theorem~\ref{Teffp4} yields
\begin{EQA}
	\| \QP \, (\tilde{\upsv}_{\GP} - \upsvs - \svn_{\GP} + \bvn_{\GP}) \|
	& \leq &
	\| \QP \IF_{\GP}^{-1} \DPN \| \, 
	\Bigl( \frac{\dltwu_{4}}{2} + \dltwu_{3}^{2} \Bigr) \,  
	\bigl( \| \DPN \IF_{\GP}^{-1} \nabla \zeta \|^{3} + \biasD^{3} \bigr) ,
	\qquad
	\qquad
\label{0mkvhgjnrwwe3u8gtygb}
	\\
	\| \QP \{ \svn_{\GP} - \bvn_{\GP} - \IF_{\GP}^{-1} (\nabla \zeta - \AvGP) \} \|
	& \leq &
	\frac{\dltwu_{3}}{2} \, \| \QP \IF_{\GP}^{-1} \DPN \| \, \bigl( \| \DPN \IF_{\GP}^{-1} \nabla \zeta \|^{2} + \biasD^{2} \bigr) \, .
\label{vuedy766t4e3bfvyt6e}
\end{EQA}
Now \eqref{0mkvhgjnrw3dfwe3u8gtygE1b} follows from \eqref{6hjdfv8e6hyefyeheew7skb} with \( k=3 \).
Next, we study the quadratic risk of \( \tilde{\upsv}_{\GP} \).
Define \( \epsv_{\QP} = \QP (\tilde{\upsv}_{\GP} - \upsvs - \svn_{\GP} + \bvn_{\GP}) \). 
By \eqref{0mkvhgjnrwwe3u8gtygb} 
\begin{EQA}
	\sqrt{\E ( \| \epsv_{\GP} \|^{2} \Ind_{\Omega(\xx)} )}
	& \leq &
	\| \QP \IF_{\GP}^{-1} \DPN \| \, 
	\Bigl( \frac{\dltwu_{4}}{2} + \dltwu_{3}^{2} \Bigr) \,  
	\Bigl( \sqrt{\E \| \DPN \IF_{\GP}^{-1} \nabla \zeta \|^{6} \Ind_{\Omega(\xx)}}  + \biasD^{3} \Bigr)
	\leq 
	\alp_{\QP,2} \sqrt{\riskt_{\QP,2}} \, ,
\label{d8ew3jfjhyvb6rt6543ejhvu}
\end{EQA}
and \eqref{6shx76whnjvyehfbvyfhb} follows. 
Further, denote 
\begin{EQA}
	\loss_{\QP}
	& \eqdef &
	\QP \IF_{\GP}^{-1} (\nabla \zeta - \AvGP) ,
	\\
	\delta_{\QP}
	& \eqdef &
	\QP (\IF_{\GP}^{-1} \nabla \zeta - \svn_{\GP}) - \QP (\IF_{\GP}^{-1} \AvGP - \bvn_{\GP}) .
\label{7ejvuejfumfmrgvytweof7}
\end{EQA}
By definition, \( \riskt_{\QP} = \E \bigl\{ \| \loss_{\QP} \|^{2} \Ind_{\Omega(\xx)} \bigr\} \), 
\( \riskt_{\QP,2} = \E \bigl\{ \| \loss_{\QP} + \delta_{\QP} \|^{2} \Ind_{\Omega(\xx)} \bigr\} \), and 
\begin{EQA}
	\riskt_{\QP,2} - \riskt_{\QP}
	& = &
	\E \bigl\{ \| \delta_{\QP} \|^{2} \Ind_{\Omega(\xx)} \bigr\}
	+ 2 \E \bigl\{ \langle \loss_{\QP} , \delta_{\QP} \rangle \Ind_{\Omega(\xx)} \bigr\} .
\label{7jde7jevyreteyv8rnbe}
\end{EQA}
Also \eqref{6hjdfv8e6hyefyeheew7skb}, \eqref{0mkvhgjnrw3dfwe3u8gtysg}, and \eqref{6dhx6whcuydsds655srew4b} imply
\begin{EQA}
	\sqrt{\E \bigl( \| \delta_{\QP} \|^{2} \Ind_{\Omega(\xx)} \bigr) }
	& \leq &
	\| \QP \IF_{\GP}^{-1} \DPN \| \, \frac{\dltwu_{3}}{2} \, 
	\Bigl( \sqrt{\E \| \DPN \IF_{\GP}^{-1} \nabla \zeta \|^{4} \Ind_{\Omega(\xx)} } + \biasD^{2} \Bigr)
	\\
	& \leq &
	\| \QP \IF_{\GP}^{-1} \DPN \| \, \frac{\dltwu_{3}}{2} \, \bigl( \CONSTi_{4} \, \dimD + \biasD^{2} \bigr)
	\leq 
	\alp_{\QP,1} \, \sqrt{\riskt_{\QP}} \, .
\label{yfhf73hjf9bryrnvbir}
\end{EQA}
This proves \eqref{EQtuGmstrVEQtGQ2b}.
\end{proof}
}{}

%


\Chapter{Profile semiparametric estimation for SLS models}
\label{SsemiMLE}
This \chname discusses the problem of the semiparametric estimation for SLS models using the profile MLE method. 
%
%
Semi-parametric statistics focuses on estimation of a low-dimensional parameter in the presence of high-dimensional or nonparametric nuisance parameters.
There exists a rich and evolving literature on this topic (e.g., 
\cite{bickel1982adaptive,pfanzagl1982lecture,MR1245941,robinson1988root,andrews1994asymptotics,robins1995semiparametric,ai2003efficient,ChCh2018,lin2023semi}.
A key concept is that of Neyman
orthogonality \cite{neyman1959optimal} of the score function; it is a property that formalizes the first-order effect
of perturbations in the nuisance terms on the target estimator. 
This notion has played an important role in semi-parametric estimation \cite{andrews1994asymptotics,newey1994asymptotic} as well as inference
for high-dimensional linear models \cite{zhang2014confidence,BeChWa2014,belloni2016post,javanmard2014confidence}. 
A popular method is ample splitting, in which different
portions of the dataset are used to estimate the non-parametric and parametric components,
(e.g., \cite{bickel1982adaptive,schick1986asymptotically,fan2012variance}).
\cite{ChCh2018} combined the notion of Neyman orthogonality with sample
splitting to construct Z-estimators that are asymptotically normal.
Other closely related procedures have been developed for 
estimating 
continuous treatment effects \cite{colangelo2020double,semenova2021debiased}; 
statistical learning with nuisance parameters \cite{foster2023orthogonal}; etc.
Most of available results have been established for an asymptotic i.i.d. setup.
The approach of this section relies on non-asymptotic full dimensional expansions for the quasi log-likelihood function
and the corresponding MLE.
It applies to a very general data generating process.

Suppose to be given a log-likelihood function \( \LL(\prmtv) = \LL(\Yv,\prmtv) \), 
where the full parameter \( \prmtv \in \R^{\dimttl} \) contains a \( \dimp \)-dimensional target of estimation 
\( \tarpv \in \Theta \subseteq \R^{\dimp} \) and a \( \dimq \)-dimensional nuisance parameter 
\( \nupv \in \Eta \subseteq \R^{\dimq} \) with \( \dimttl = \dimp + \dimq \).
We will write \( \LL(\prmtv) = \LL(\tarpv,\nupv) \).
The full dimensional MLE \( \tilde{\prmtv} = (\tilde{\tarpv},\tilde{\nupv}) \) 
is defined by the joint optimization of \( \LL(\prmtv) \)
w.r.t. the target parameter \( \tarpv \) and \( \nupv \):
\begin{EQA}
	\tilde{\prmtv}
	&=&
	(\tilde{\tarpv},\tilde{\nupv})
	=
	\argmax_{\prmtv} \LL(\prmtv)
	=
	\argmax_{(\tarpv,\nupv)} \LL(\tarpv,\nupv) .
\label{456765456789ikjms}
\end{EQA}
The \emph{profile MLE} \( \tilde{\tarpv} \) is just the \( \tarpv \)-component of \( \tilde{\prmtv} \):
\begin{EQA}
	\tilde{\tarpv}
	&=&
	\argmax_{\tarpv} \max_{\nupv} \LL(\prmtv)
	=
	\argmax_{\tarpv} \max_{\nupv} \LL(\tarpv,\nupv)
	=
	\argmax_{\tarpv} \LL(\tarpv,\tilde{\nupv}) .
\label{109oijnbvfdjhgdvfbn}
\end{EQA}
Population counterparts of \( \tilde{\prmtv} \) and \( \tilde{\tarpv} \) 
are defined by replacing the log-likelihood with its expectation:
\begin{EQ}[rcl]
\label{109oijnbvfdjhgdvfbnse}
	\prmtvs
	&=&
	\argmax_{\prmtv} \E \LL(\prmtv)
	=
	\argmax_{(\tarpv,\nupv)} \E \LL(\tarpv,\nupv) ,
	\\
	\tarpvs
	&=&
	\argmax_{\tarpv} \max_{\nupv} \E \, \LL(\prmtv)
	=
	\argmax_{\tarpv} \max_{\nupv} \E \, \LL(\tarpv,\nupv)
	=
	\argmax_{\tarpv} \E \, \LL(\tarpv,\nupvs) .
	\qquad
\end{EQ}
The background truth \( \tarpvs \) is defined as 
\begin{EQA}
	\tarpvs
	&=&
	\argmax_{\tarpv} \max_{\nupv} \E \, \LL(\prmtv)
	=
	\argmax_{\tarpv} \max_{\nupv} \E \, \LL(\tarpv,\nupv) .
\label{6hedcv8wh2fuf7r7erhfyehse}
\end{EQA}
Later we present sharp finite sample error bounds on the accuracy of estimation
\( \tilde{\tarpv} - \tarpvs \).

\Section{Full dimensional estimation}
Everywhere in this section, we assume that the basic assumptions 
\nameref{Eref}, \nameref{LLref}, and \nameref{EU2ref} are fulfilled for the full dimensional model 
given by \( \LL(\prmtv) \).
Denote \( \IFT = \IFT (\prmtvs) \).
We also fix a full-dimensional metric tensor \( \DFM \) such that \( \DFM^{2} \leq \dmax^{2} \IFT \) and 
suppose that the function \( \fs(\prmtv) = \E \LL(\prmtv) \) satisfies 
\nameref{LLsT3ref} at \( \prmtvs \) with this tensor \( \DFM \)
and some \( \rr \) and \( \dltwu_{3} \).
Due to \nameref{Eref}, the stochastic gradient
\( \nabla \zeta = \nabla \LL(\prmtv) - \nabla \E \LL(\prmtv) \) does not depend on \( \prmtv \). 
Let \( \VP^{2} \geq \Var(\nabla \zeta) \) be shown in \nameref{EU2ref}.
Then on a random set \( \Omega(\xx) \) with \( \P\bigl( \Omega(\xx) \bigr) \geq 1 - 3 \ex^{-\xx} \), it holds
\begin{EQA}
	\| \DFM \, \IFT^{-1} \nabla \zeta \|
	& \leq &
	\rrDF \, ,
\label{obkjyy643enhfgvtygyjyefs}
\end{EQA}
where with \( \BBH_{\DFM} \eqdef \DFM \IFT^{-1} \Var(\nabla \zeta) \, \IFT^{-1} \DFM \),
the radius \( \rrDF \) is given by
\begin{EQA}
	\rrDF 
	&=& 
	\zq(\BBH_{\DFM},\xx)
	\leq 
	\sqrt{\tr \BBH_{\DFM}} + \sqrt{2 \xx \, \| \BBH_{\DFM} \|} \,\, .
\label{vhyvw323sdwe4w3rdesf23}
\end{EQA}
Theorem~\ref{TFiWititG2} and Theorem~\ref{TFiWititG3} yield the following result.

\begin{theorem}
\label{Tconcsemi}
Suppose \nameref{Eref}, \nameref{LLref}, and \nameref{EU2ref} 
for the full dimensional parameter \( \prmtv \), and let 
\( \prmtvs \) be from \eqref{109oijnbvfdjhgdvfbnse}.
Assume \nameref{LLsT3ref} at \( \prmtvs \) with the tensor \( \DFM \) and
\( \rr \), \( \dltwu_{3} \) satisfying
\begin{EQA}[c]
	\DFM^{2} \leq \dmax^{2} \hspm \IFT \, ,
	\quad
	\rr \geq \frac{3}{2} \, \rrDF \, ,
	\quad
	\dmax^{2} \dltwu_{3} \, \hspm \rrDF < \frac{4}{9} \, ,
\label{8difiyfc54wrboeMLEse}
\end{EQA}
where \( \rrDF \) is defined by \eqref{vhyvw323sdwe4w3rdesf23}.
Then it holds on \( \Omega(\xx) \) 
\begin{EQA}
\label{yfjvh9h4f5e53tghfugu}
	\bigl\| \DFM^{-1} \IFT ( \tilde{\prmtv} - \prmtvs - \IFT^{-1} \nabla \zeta ) \bigr\|
    & \leq & 
    \frac{3\dltwu_{3}}{4} \, \| \DFM \, \IFT^{-1} \nabla \zeta \|^{2} \, ,
    \\
	\bigl| 2 \LL(\tilde{\prmtv}) - 2 \LL(\prmtvs) - \| \IFT^{-1/2} \nabla \zeta \|^{2} \bigr|
    & \leq &
    \frac{\dltwu_{3}}{2} \, \| \DFM \, \IFT^{-1} \nabla \zeta \|^{3} \, .
\label{kfchvyfyrdetye46t}
\end{EQA}
\end{theorem}

\medskip
This full-dimensional result yields concentration bounds for the \( \tarpv \) and \( \nupv \) components of \( \prmtv \).
Indeed, as \( \DFM^{2} \leq \dmax^{2} \hspm \IFT \), \eqref{yfjvh9h4f5e53tghfugu} implies
\begin{EQA}
	\bigl\| \DFM \bigl( \tilde{\prmtv} - \prmtvs - \IFT^{-1} \nabla \zeta \bigr) \bigr\|
    & \leq & 
    \frac{3 \dmax^{2} \dltwu_{3}}{4} \, \| \DFM \, \IFT^{-1} \nabla \zeta \|^{2} 
    \leq 
    \frac{3 \dmax^{2} \dltwu_{3} \, \rrDF}{4} \, \| \DFM \, \IFT^{-1} \nabla \zeta \|
    \, ,
    \\
    \bigl\| \DFM ( \tilde{\prmtv} - \prmtvs)) \bigr\|
    & \leq & 
    \bigl( 1 + \frac{3 \dmax^{2} \dltwu_{3} \, \rrDF}{4} \bigr) \| \DFM \, \IFT^{-1} \nabla \zeta \| \, .
\label{ergvfyde4wse4e3w4e}
\end{EQA}
For instance, if the matrix \( \IFT \) is block-diagonal, that is, 
if \( \IFT = \blk( \IFT_{\tarpv\tarpv},\IFT_{\nupv\nupv}) \) and \( \DFM = \IFT^{1/2} \),
then \eqref{ergvfyde4wse4e3w4e} implies on \( \Omega(\xx) \) under \( \dmax^{2} \dltwu_{3} \, \rrDF \leq 2/3 \) 
\begin{EQ}[rcccl]
	\| \IFT_{\tarpv\tarpv}^{1/2} (\tilde{\tarpv} - \tarpvs) \|
	& \leq &
	3 \rrDF/2 \, ,
	\qquad
	\| \IFT_{\nupv\nupv}^{1/2} (\tilde{\nupv} - \nupvs) \|
	& \leq &
	3 \rrDF/2 \, .
\label{5huu776776wewee3344juui}
\end{EQ}
%

\Section{Expansions and risk bounds for the profile MLE}
\label{SprofileMLE}
The main issue with result \eqref{5huu776776wewee3344juui} is that the radius \( \rrDF \)
is of order \( \dimDF^{1/2} \), where the value 
\( \dimDF = \tr(\B_{\DFM}) \) corresponds to the full parameter dimension 
and can be large for a high dimensional nuisance parameter \( \nupv \)   
yielding a poor estimation accuracy.
Later we discuss what can be extracted from the full dimensional expansion \eqref{yfjvh9h4f5e53tghfugu} 
concerning the target component \( \tarpv \).
%
%
Represent the negative Hessian matrix \( \IFT(\prmtv) = - \nabla^{2} \E \LL(\prmtv) \) 
of \( f(\prmtv) = \E \LL(\prmtv) \) in the block form corresponding to the components \( \tarpv \) and \( \nupv \):
\begin{EQA}
	\IFT(\prmtv)
	&=&
	- \begin{pmatrix}
		\nabla_{\tarpv\tarpv} f(\prmtv) & \nabla_{\tarpv\nupv} f(\prmtv)
		\\
		\nabla_{\nupv\tarpv} f(\prmtv) & \nabla_{\nupv\nupv} f(\prmtv)
	\end{pmatrix} 
	=
	\begin{pmatrix}
		\IFT_{\tarpv\tarpv}(\prmtv) & \IFT_{\tarpv\nupv}(\prmtv)
		\\
		\IFT_{\nupv\tarpv}(\prmtv) & \IFT_{\nupv\nupv}(\prmtv)
	\end{pmatrix}.
\label{HAuuuvvuvvTH0}
\end{EQA}
Here \( \IFT_{\nupv\tarpv}(\prmtv) = \IFT_{\tarpv\nupv}(\prmtv)^{\T} \).
Later we use the inverse of \( \IFT = \IFT(\prmtvs) \) and its block representation.
Due to Schur's complement formulas, see Section~\ref{Ssemitools}, the diagonal blocks of 
\( \IFT^{-1} \) are  \( \IFTb_{\tarpv\tarpv}^{-1} \) and \( \IFTb_{\nupv\nupv}^{-1} \), where
\begin{EQ}[rcl]
	\IFTb_{\tarpv\tarpv}
	&=&
	\IFT_{\tarpv\tarpv} - \IFT_{\tarpv\nupv} \, \IFT_{\nupv\nupv}^{-1} \, \IFT_{\nupv\tarpv} \, ,
	\\
	\IFTb_{\nupv\nupv}
	&=&
	\IFT_{\nupv\nupv} - \IFT_{\nupv\tarpv} \, \IFT_{\tarpv\tarpv}^{-1} \, \IFT_{\tarpv\nupv} \, .
\label{divjjw3d6f5f5te3gjg}
\end{EQ}
%
Now we explain how the full dimensional expansions can be used to derive the bounds for the target parameter.
We establish two results that can be viewed as finite sample analogs of the classical asymptotic results of parametric
statistics.
The Fisher expansion allows studying the properties of the profile MLE \( \tilde{\tarpv} \) with an explicit 
leading term in the error \( \tilde{\tarpv} - \tarpvs \).
It is convenient and does not restrict generality to assume that the local metric tensor \( \DFM \) has a block-diagonal structure 
\( \DFM = \blk\{ \DFM_{\tarpv\tarpv},\DFM_{\nupv\nupv} \} \).

\Subsection{Fisher and Wilks expansions}

\begin{theorem}
\label{Tsemieffp3n}
Assume the conditions of Theorem~\ref{Tconcsemi} and let \( \DFM = \blk\{ \DFM_{\tarpv\tarpv},\DFM_{\nupv\nupv} \} \).
Then on \( \Omega(\xx) \), it holds for any linear mapping \( \QP \) on \( \R^{\dimp} \)
\begin{EQA}
\label{vdud6ehry3hrf67ruyta}
	\| \QP \{ \tilde{\tarpv} - \tarpvs - (\IFT^{-1} \nabla \zeta)_{\tarpv} \} \|
	& \leq &
	\| \QP \, \DFM_{\tarpv\tarpv}^{-1} \| \, \frac{3 \dmax^{2} \dltwu_{3}}{4} \, \| \DFM \IFT^{-1} \nabla \zeta \|^{2} \, .
\end{EQA}
\end{theorem}

The proof is given in Section~\ref{Ssemitools}.
Due to Theorem~\ref{Tsemieffp3n}, the stochastic error of \( \QP (\tilde{\tarpv} - \tarpvs) \) is described by the 
vector \( \QP (\IFT^{-1} \nabla \zeta)_{\tarpv} \).
Particular examples of choosing the scaling matrix \( \QP \) include \( \QP = \sqrt{n} \Id_{\dimp} \), 
\( \QP = \DFM_{\tarpv\tarpv} \) or \( \QP = \IFTb_{\tarpv\tarpv}^{1/2} \); see \eqref{divjjw3d6f5f5te3gjg}.
By Lemma~\ref{LSchur}, it holds with the full dimensional score vector \( \nabla \zeta = (\scorem{\tarpv} \zeta,\scorem{\nupv} \zeta) \)
\begin{EQA}
	(\IFT^{-1} \nabla \zeta)_{\tarpv}
	&=&
	\IFTb_{\tarpv\tarpv}^{-1} (\scorem{\tarpv} \zeta - \IFT_{\tarpv\nupv} \, \IFT_{\nupv\nupv}^{-1} \scorem{\nupv} \zeta)
	\, .
	\qquad
\label{hf8g8hy984u3453fhfhh}
\end{EQA}
Introduce also the \emph{standardized semiparametric score} \( \xivr \) by
\begin{EQA}[rcl]
	\xivr
	& \eqdef &
	\IFTb_{\tarpv\tarpv}^{1/2} \, (\IFT^{-1} \nabla \zeta)_{\tarpv}
	=
	\IFTb_{\tarpv\tarpv}^{-1/2} (\scorem{\tarpv} \zeta - \IFT_{\tarpv\nupv} \, \IFT_{\nupv\nupv}^{-1} \scorem{\nupv} \zeta)
	\, .
\label{hf8g8hy984u3453fhfh}
\end{EQA}
Wilks expansion generalizes the prominent Wilks phenomenon \cite{Wilks1938}, \cite{FaZh2001} 
which claims for the profile MLE \( \tilde{\tarpv} \) 
in a regular parametric \( \dimp \)-dimensional model that the twice excess
\( 2 \LLp(\tilde{\tarpv}) - 2 \LLp(\tarpvs) \) for \( \LLp(\tarpv) = \sup_{\nupv} \LL(\tarpv,\nupv) \)
is asymptotically chi-squared with \( \dimp \) degrees of freedom.

\begin{theorem}
\label{TsemiWilks3}
Assume the conditions of Theorem~\ref{Tsemieffp3n}.
Then on \( \Omega(\xx) \), it holds 
with \( \LLp(\tarpv) = \sup_{\nupv} \LL(\tarpv,\nupv) \)
and \( \xivr \) from \eqref{hf8g8hy984u3453fhfh},
\begin{EQA}
	\Bigl| 
		2 \LLp(\tilde{\tarpv}) - 2 \LLp(\tarpvs) - \bigl\| \xivr \bigr\|^{2} 
	\Bigr|
	& \leq &
	\frac{\dltwu_{3}}{2} 
	\Bigl( \| \DFM \IFT^{-1} \nabla \zeta \|^{3} + \| \DFM_{\nupv\nupv} \IFT_{\nupv\nupv}^{-1} \scorem{\nupv} \zeta \|^{3} \Bigr) \, .
	\qquad
\label{vdud6ehry3hrf67ruy233}
\end{EQA}
\end{theorem}

The proof is given in Section~\ref{Ssemitools} later.

\Subsection{Identifiability and semiparametric effective/critical dimension}

With \( \xivr \) from \eqref{hf8g8hy984u3453fhfh},
the \emph{semiparametric effective dimension} \( \dimAr \) is defined by
\begin{EQA}
	\dimAr
	& \eqdef &
	\E \| \xivr \|^{2} 
	=
	\tr \Var ( \xivr ) \, .
\label{fiuvuv77he6v63hf63he}
\end{EQA}
Here we provide some evidence that \( \dimAr \) 
is of the same order as the standard efficient dimension 
\( \dimtarg = \tr \Var(\IFT_{\tarpv\tarpv}^{-1/2} \, \scorem{\tarpv} \zeta) \) 
for the target parameter \( \tarpv \) only under the so-called \emph{identifiability} condition.
Informally this condition means that the full-dimensional information matrix
\( \IFT \) can be bounded from below by a multiple of the block-diagonal matrix 
\( \blk\{ \IFT_{\tarpv\tarpv},\IFT_{\nupv\nupv} \} \).
%
%

\begin{lemma}
\label{Leffdimrho}
Assume for some matrix \( \BBH_{\tarpv} \in \Matr_{\dimp} \) and some \( \rhoIF < 1 \)
\begin{EQA}[rcl]
	\IFT_{\tarpv\tarpv}^{-1/2} \, \Var(\scorem{\tarpv} \zeta) \, \IFT_{\tarpv\tarpv}^{-1/2} 
	& \leq &
	\BBH_{\tarpv} 
	\, ,
	\\
	\IFT_{\tarpv\tarpv}^{-1/2} \, \IFT_{\tarpv\nupv} \, \IFT_{\nupv\nupv}^{-1} \, \Var(\scorem{\nupv} \zeta) \IFT_{\nupv\nupv}^{-1} \, \IFT_{\nupv\tarpv} \,\IFT_{\tarpv\tarpv}^{-1/2} 
	& \leq &
	\rhoIF^{2} \BBH_{\tarpv} \, .
\label{dywehw237dfcuw2jwuyfhje3}
\end{EQA}
Then
\begin{EQA}[rcccl]
	\Var(\xivr)
	& \leq &
	\frac{1 + \rhoIF}{1 - \rhoIF} \, \BBH_{\tarpv} \, ,
\label{7djdcjhdyhdy56ewghdythL}
	\qquad
	\dimAr
	=
	\E \| \xivr \|^{2}
	& \leq &
	\frac{1 + \rhoIF}{1 - \rhoIF} \tr (\BBH_{\tarpv}) .
\end{EQA}
\end{lemma}
The proof is given in Section~\ref{Ssemitools}.

%

\medskip
Now we discuss the issue of \emph{critical dimension} in profile semiparametric estimation.
The Fisher expansion \eqref{vdud6ehry3hrf67ruyta} is only meaningful 
if the error term on the right-hand side is sufficiently small.
Namely, with \( \QP = \DFM_{\tarpv\tarpv} \),
the quantity \( \dimQ = \E \| \QP (\IFT^{-1} \nabla \zeta)_{\tarpv} \|^{2} \) 
of order the effective target dimension, 
while \( \| \DFM \IFT^{-1} \nabla \zeta \|^{2} \)
corresponds to the full effective parameter dimension 
\( \dimDF = \tr \bigl\{ \DFM \IFT^{-1} \Var(\nabla \zeta) \, \IFT^{-1} \DFM \bigr\} \);
cf. \eqref{vhyvw323sdwe4w3rdesf23}.
Hence, we need the condition \( \dltwu_{3} \, \dimDF \ll \dimQ^{1/2} \).
Under self-concordance condition \nameref{LLtS3ref}, the value \( \dltwu_{3} \) is of order \( n^{-1/2} \).
For a small target effective dimension, the related \emph{critical dimension} condition for 
the Fisher expansion \eqref{vdud6ehry3hrf67ruyta} 
reads as \( \dimDF \ll n^{1/2} \).
The Wilks expansion of Theorem~\ref{TsemiWilks3} is even more demanding. 
The right-hand side of \eqref{vdud6ehry3hrf67ruy233} is relatively small under 
\( \dltwu_{3} \, \dimDF^{3/2} \ll \dimAr \).
For \( \dimAr \) fixed, this requires \( \dimDF \ll n^{1/3} \) and can be quite restrictive.
\iffourG{
The next section explains, how the condition can be improved under fourth-order smoothness assumptions. 

\Section{Profile estimation. 4S bounds}
\label{Sfourthsemi}
This section explains how the advanced expansions from Section~\ref{SMLE4} based on Proposition~\ref{Pconcgeneric4} 
under \nameref{LLsT3ref} and \nameref{LLsT4ref} can be used to substantially improve
the error terms in the expansions for the profile MLE \( \tilde{\tarpv} \).
We follow the line of Section~\ref{SMLE4}.
For \( \fs(\prmtv) = \E \LL(\prmtv) \), consider the third-order tensor
\( \TensU(\uv) = \frac{1}{6} \langle \nabla^{3} \fs(\prmtvs), \uv^{\otimes 3} \rangle \)
and its gradient \( \frac{1}{2} \langle \nabla^{3} \fs(\prmtvs), \uv^{\otimes 2} \rangle \).
With \( \IFT = \IFT(\prmtvs) \), define 
\begin{EQA}[rcccl]
	\svn 
	&=&
	\IFT^{-1} \bigl\{ \nabla \zeta + \nabla \TensU(\IFT^{-1} \nabla \zeta) \bigr\} 
	&=&
	(\svn_{\tarpv},\svn_{\nupv}) \, .
\label{8vfjvr43223efryfuweefgse}
\end{EQA}
Also remind \( \dimDF = \E \| \DFM \IFT^{-1} \nabla \zeta \|^{2} \).
We apply Theorem~\ref{Teffp4}.

\begin{theorem}
\label{Tsemieffp4}
Suppose \nameref{Eref}, \nameref{LLref}, and \nameref{EU2ref} 
for the full dimensional parameter \( \prmtv \).
Let \nameref{LLsT3ref} and \nameref{LLsT4ref} hold at \( \prmtvs \) with the metric tensor 
\( \DFM = \blk\{ \DFM_{\tarpv\tarpv},\DFM_{\nupv\nupv} \} \) and values \( \rr \), \( \dltwu_{3} \), and \( \dltwu_{4} \) satisfying
\begin{EQA}[c]
	\DFM^{2} \leq \dmax^{2} \, \IFT \, ,
	\;\; 
	\rr \geq \frac{3}{2} \rrDF \, ,
	\;\;
	\dmax^{2} \dltwu_{3} \, \rrDF < \frac{4}{9} \, ,
	\;\;
	\dmax^{2} \dltwu_{4} \, \rrDF^{2} < \frac{1}{3} 
	\, ,
\label{7deyedhfg5563w6rfygjLR4se}
\end{EQA}
for \( \rrDF \) from \eqref{vhyvw323sdwe4w3rdesf23}.
Then on \( \Omega(\xx) \), the estimate \( \tilde{\prmtv} \) satisfies concentration bound \eqref{ergvfyde4wse4e3w4e}.
For any linear mapping \( \QP \) of \( \tarpv \)
\begin{EQ}[rcl]
	\| \QP \, ( \tilde{\tarpv} - \tarpvs - \svn_{\tarpv} ) \|
	& \leq &
	\| \QP \, \DFM_{\tarpv\tarpv}^{-1} \| \, \dmax^{2} 
	\Bigl( \frac{\dltwu_{4}}{2} + \dmax^{2} \dltwu_{3}^{2} \Bigr) \, \| \DFM \IFT^{-1} \nabla \zeta \|^{3} 
	\, .
	\qquad
	\qquad
\label{0mkvhgjnrwwe3u8gtygse}
\end{EQ}
and
\begin{EQ}[rcl]
	\| \QP (\svn - \IFT^{-1} \, \nabla \zeta)_{\tarpv} \|
	& \leq &
	\| \QP \, \DFM_{\tarpv\tarpv}^{-1} \| \, \frac{\dmax^{2} \dltwu_{3}}{2} \, 
	\| \DFM \IFT^{-1} \nabla \zeta \|^{2} 
	\, .
	\qquad
\label{vuedy766t4e3bfvyt6ese}
\end{EQ}
\end{theorem}

Now we study the risk of \( \tilde{\tarpv} \) using Theorem~\ref{Teff41}.
Define 
\begin{EQA}[c]
	\riskt_{\QP}
	\eqdef
	\E \bigl\{ \| \QP (\IFT^{-1} \, \nabla \zeta)_{\tarpv} \|^{2} \Ind_{\Omega(\xx)} \bigr\} \, .
\label{7dhdrw2bfvu78u78e4ndwse}
\end{EQA}
Assuming \( \E \{ \IFT^{-1} \, \nabla \zeta \, \Ind_{\Omega(\xx)} \} \approx 0 \), we derive
\begin{EQA}
	\riskt_{\QP}
	& \approx &
	\tr \Var\bigl( \QP (\IFT^{-1} \, \nabla \zeta)_{\tarpv} \bigr) \, .
\label{6jdfctewbfybneftvbej}
\end{EQA}
Also, define
\begin{EQA}
\label{hdvje39bug53ebfh8edxse}
	\riskt_{\QP,2} 
	& \eqdef & 
	\E \bigl\{ \| \QP \svn_{\tarpv} \|^{2} \Ind_{\Omega(\xx)} \bigr\} \, .
\label{yfhcvhched6chejdrtese}
\end{EQA}

\begin{theorem}
\label{Tseff42}
Assume the conditions of Theorem~\ref{Tsemieffp4} and let
\begin{EQA}[c]
	\E \bigl\{ \| \DFM \IFT^{-1} \nabla \zeta \|^{k} \Ind_{\Omega(\xx)} \bigr\} \leq \CONSTi_{k}^{2} \, \dimDF^{k/2} ,
	\qquad
	k=3,4,6 \, .
\label{6hjdfv8e6hyefyeheew7skGPT}
\end{EQA}
For a linear mapping \( \QP \) and \( \riskt_{\QP,2} \) from \eqref{hdvje39bug53ebfh8edxse}, it holds
\begin{EQA}
	&& \nquad
	\E \bigl\{ \| \QP \, (\tilde{\tarpv} - \tarpvs) \| \Ind_{\Omega(\xx)} \bigr\}
	\leq 	 
	\E \bigl\{ \| \QP (\svn_{\tarpv} - \bvn_{\tarpv}) \| \Ind_{\Omega(\xx)} \bigr\}
	\\
	&&
	+ \, \| \QP \, \DFM_{\tarpv\tarpv}^{-1} \| \, \dmax^{2} 
	\Bigl( \frac{\dltwu_{4}}{2} + \dmax^{2} \dltwu_{3}^{2} \Bigr) \, \CONSTi_{3}^{2} \, \dimDF^{3/2} 
\label{0mkvhgjnrw3dfwe3u8gtygE1se}
\end{EQA}
and 
\begin{EQA}
	&& \nquad
	\Bigl| \E \bigl\{ \| \QP \svn_{\tarpv} \| \Ind_{\Omega(\xx)} \bigr\}
	- \E \bigl\{ \| \QP \, (\IFT^{-1} \, \nabla \zeta)_{\tarpv} \| \Ind_{\Omega(\xx)} \bigr\} 
	\Bigr|
	\leq 
	\| \QP \, \DFM_{\tarpv\tarpv}^{-1} \| \, \frac{\dmax^{2} \dltwu_{3}}{2} \, \dimDF \, .
\label{udtgecthwjdytdehduuc6se}
\end{EQA}
Furthermore,
\begin{EQA}
	\bigl( 1 - \alp_{\QP,2} \bigr)^{2} \riskt_{\QP,2}
	\leq 
	\E \bigl\{ \| \QP \, (\tilde{\tarpv} - \tarpvs) \|^{2} \Ind_{\Omega(\xx)} \bigr\}
	& \leq &
	\bigl( 1 + \alp_{\QP,2} \bigr)^{2} \riskt_{\QP,2} \, 
\label{6shx76whnjvyehfbvyfhse}
\end{EQA}
provided that
\begin{EQA}
	\alp_{\QP,2}
	& \eqdef &
	\frac{\| \QP \, \DFM_{\tarpv\tarpv}^{-1} \| \, \dmax^{2} ( \dltwu_{4}/2 + \dmax^{2} \dltwu_{3}^{2} ) \, \CONSTi_{6} \, \dimDF^{3/2} }
	 	 {\sqrt{\riskt_{\QP,2}}} 
	< 1 \, .
\label{cuyhe4f8jfdebnuise}
\end{EQA}
If another value \( \alp_{\QP,1} < 1 \) is such that  
\begin{EQ}[rcl]
	&&
	\| \QP \, \DFM_{\tarpv\tarpv}^{-1} \| \, \frac{\dmax^{2} \dltwu_{3}}{2} \, \CONSTi_{4} \, \dimDF 
	\leq 
	\alp_{\QP,1} \, \sqrt{\riskt_{\QP}} \, 
\label{6dhx6whcuydsds655srew4se}
\end{EQ}
with \( \riskt_{\QP} \) from \eqref{7dhdrw2bfvu78u78e4ndwse} then
\begin{EQA}
	\riskt_{\QP} (1 - \alp_{\QP,1})^{2} 
	\leq 
	\riskt_{\QP,2}
	& \leq &
	\riskt_{\QP} (1 + \alp_{\QP,1})^{2} \, .
\label{EQtuGmstrVEQtGQse}
\end{EQA}
\end{theorem}

The presented results can be spelled out as follows:
\begin{EQA}
	\tilde{\tarpv} - \tarpvs
	& \approx &
	\svn_{\tarpv} 
\label{y6fcckfcvuy7uujrf7yfh}
\end{EQA}
with high accuracy.
Moreover, it holds for the quadratic risk of \( \tilde{\tarpv} \)
\begin{EQA}
	\E \bigl\{ \| \QP \, (\tilde{\tarpv} - \tarpvs) \|^{2} \Ind_{\Omega(\xx)} \bigr\}
	& \approx &
	\E \bigl\{ \svn_{\tarpv}^{2} \Ind_{\Omega(\xx)} \bigr\}
	=
	\riskt_{\QP,2} \, .
\label{uykjf83whed6vhy4534dxk}
\end{EQA}
If \( \alp_{\QP,1} \) is small then \( \riskt_{\QP,2} \approx \riskt_{\QP} \); see \eqref{7dhdrw2bfvu78u78e4ndwse}.
Therefore, in this case, a third-order correction is not necessary, one can use 
\( \tilde{\tarpv} - \tarpvs \approx (\IFT^{-1} \, \nabla \zeta)_{\tarpv} \).

\medskip
The next result improves the Wilks expansion \eqref{vdud6ehry3hrf67ruy233} from Theorem~\ref{Tsemieffp3n} by using fourth-order smoothness
condition.

\begin{theorem}
\label{TsemiWilks4}
Assume the conditions of Theorem~\ref{Tsemieffp4}.
With \( \LLp(\tarpv) = \sup_{\nupv} \LL(\tarpv,\nupv) \), it holds on \( \Omega(\xx) \)
\begin{EQA}
	&& \nquad
	\bigl| \LLp(\tilde{\tarpv}) - \LLp(\tarpvs) - \frac{1}{2} \| \xivr \|^{2}
	- \TensU(\IFT^{-1} \nabla \zeta) + \TensU(\IFT_{\nupv\nupv}^{-1} \scorem{\nupv} \zeta) \bigr|
	\\
	& \leq &
	\frac{\dltwu_{4} + 4 \dmax^{2} \dltwu_{3}^{2}}{8} 
	\bigl( \| \DFM \IFT^{-1} \nabla \zeta \|^{4} + \| \DFM_{\nupv\nupv} \, \IFT_{\nupv\nupv}^{-1} \scorem{\nupv} \zeta \|^{4} \bigr) 
    \\
    && 
    + \, \frac{\dmax^{2} (\dltwu_{4} + 2 \dmax^{2} \dltwu_{3}^{2})^{2}}{4} \, 
    \bigl( \| \DFM \IFT^{-1} \nabla \zeta \|^{6} + \| \DFM_{\nupv\nupv} \, \IFT_{\nupv\nupv}^{-1} \scorem{\nupv} \zeta \|^{6} \bigr) .
\label{iu8dcjhfc767y63whf6hdfu44}
\end{EQA}
\end{theorem}

The improved results enable us to reconsider the critical dimension issue.
We assume \( \dltwu_{3} \asymp n^{-1/2} \), \( \dltwu_{4} \asymp n^{-1} \), and
\( \rrDF \approx \dimDF^{1/2} \).
Bounds \eqref{0mkvhgjnrwwe3u8gtygse} and \eqref{0mkvhgjnrw3dfwe3u8gtygE1se} are meaningful if 
\( \dimDF^{3/2} \ll n \),
which improves the relation \( \dimDF^{2} \ll n \dimtarg \)
required for Theorem~\ref{Tsemieffp3n}.
However, without third-order correction, by \eqref{vuedy766t4e3bfvyt6ese} 
\begin{EQ}[rcl]
	\| \QP \{ \svn - \IFT^{-1} \, \nabla \zeta \}_{\tarpv} \|
	& \leq &
	\| \QP \, \DFM_{\tarpv\tarpv}^{-1} \| \, \frac{\dmax^{2} \dltwu_{3}}{2} \, \| \IFT^{-1/2} \nabla \zeta \|^{2} \, 
	\qquad
\label{vuedy766t4e3bfvyt6ecr}
\end{EQ}
and for \( \| \QP \, \DFM_{\tarpv\tarpv}^{-1} \| \asymp 1 \), the right-hand side 
is of order \( \dltwu_{3} \, \dimDF \asymp n^{-1/2} \, \dimDF \) as in 3G case. 
For the Wilks expansion of Theorem~\ref{TsemiWilks4}, the right-hand side of \eqref{iu8dcjhfc767y63whf6hdfu44} is relatively small 
provided that \( \dimDF^{2} \ll n \, \dimtarg \), where \( \dimtarg = \tr \Var(\IFT_{\tarpv\tarpv}^{-1/2} \, \scorem{\tarpv} \zeta) \).
This improves the condition \( \dimDF^{3} \ll n \, \dimtarg^{2} \) from the 3S case.
}{}

\Section{Penalization in profile MLE}
Suppose a smooth penalty function 
\( \pent_{\GPT}(\prmtv) \) to be given.
A popular example is a quadratic penalty \( \| \GPT \prmtv \|^{2}/2 \).
We accept a general penalty function.
The full dimensional penalized MLE \( \tilde{\prmtv}_{\GPT} = (\tilde{\tarpv}_{\GPT},\tilde{\nupv}_{\GPT}) \) 
is defined by the joint optimization of \( \LL_{\GPT}(\prmtv) = \LL(\prmtv) - \pent_{\GPT}(\prmtv) \)
w.r.t. the target parameter \( \tarpv \) and \( \nupv \):
\begin{EQA}
	\tilde{\prmtv}_{\GPT}
	&=&
	(\tilde{\tarpv}_{\GPT},\tilde{\nupv}_{\GPT})
	=
	\argmax_{\prmtv} \LL_{\GPT}(\prmtv)
	=
	\argmax_{(\tarpv,\nupv)} \LL_{\GPT}(\tarpv,\nupv) .
\label{456765456789ikjms}
\end{EQA}
The \emph{profile MLE} \( \tilde{\tarpv}_{\GPT} \) is just the \( \tarpv \)-component of \( \tilde{\prmtv}_{\GPT} \):
\begin{EQA}
	\tilde{\tarpv}_{\GPT}
	&=&
	\argmax_{\tarpv} \max_{\nupv} \LL_{\GPT}(\prmtv)
	=
	\argmax_{\tarpv} \max_{\nupv} \LL_{\GPT}(\tarpv,\nupv)
	=
	\argmax_{\tarpv} \LL_{\GPT}(\tarpv,\tilde{\nupv}_{\GPT}) .
\label{109oijnbvfdjhgdvfbn}
\end{EQA}
Population counterparts of \( \tilde{\prmtv}_{\GPT} \) and \( \tilde{\tarpv}_{\GPT} \) 
are defined by replacing the log-likelihood with its expectation:
\begin{EQ}[rcl]
\label{109oijnbvfdjhgdvfbns}
	\prmtvs_{\GPT}
	&=&
	\argmax_{\prmtv} \E \LL_{\GPT}(\prmtv)
	=
	\argmax_{(\tarpv,\nupv)} \E \LL_{\GPT}(\tarpv,\nupv) ,
	\\
	\tarpvs_{\GPT}
	&=&
	\argmax_{\tarpv} \max_{\nupv} \E \, \LL_{\GPT}(\prmtv)
	=
	\argmax_{\tarpv} \max_{\nupv} \E \, \LL_{\GPT}(\tarpv,\nupv)
	=
	\argmax_{\tarpv} \E \, \LL_{\GPT}(\tarpv,\nupvs_{\GPT}) .
	\qquad
\end{EQ}
The corresponding Fisher information matrix is \( \IFT_{\GPT} = \IFT_{\GPT}(\prmtvs_{\GPT}) \)
with the blocks \( \IFT_{\GPT,\tarpv\tarpv} \), \( \IFT_{\GPT,\tarpv\nupv} \), \( \IFT_{\GPT,\nupv\nupv} \).
The background truth \( \tarpvs \) is defined as 
\begin{EQA}
	\tarpvs
	&=&
	\argmax_{\tarpv} \max_{\nupv} \E \, \LL(\prmtv)
	=
	\argmax_{\tarpv} \max_{\nupv} \E \, \LL(\tarpv,\nupv) .
\label{6hedcv8wh2fuf7r7erhfyeh}
\end{EQA}
Penalization by \( \pent_{\GPT}(\tarpv,\nupv) \) yields some bias in the sense that \( \tilde{\tarpv}_{\GPT} \)
effectively estimates \( \tarpvs_{\GPT} \) rather than \( \tarpvs \).
Proposition~\ref{Lvarusetb} provides
a bound on the norm of \( \QF (\prmtvs_{\GPT} - \prmtvs) \) 
for a linear mapping \( \QF \) on the full dimensional parameter space. 
%
Similarly to \eqref{11ma3eaelDebQ}, we state the following bound. 

\begin{proposition}
\label{Pvarusetbse}
Define  
\begin{EQA}[rcccl]
	\AvGPT 
	& \eqdef &
	\nabla \!\pent_{\GP}(\upsvs) ,
	\quad
	\biasDF
	& \eqdef &
	\| \DFM \IFT_{\GPT}^{-1} \AvGPT \| 
	\, .
\label{56shdwsjw76e27fur4hyf}
\end{EQA}
For the function \( \fs_{\GPT}(\prmtv) \), assume \nameref{LLsT3ref}  at \( \prmtvs_{\GPT} \) with 
\( \DFM = \blk\{ \DFM_{\tarpv\tarpv},\DFM_{\nupv\nupv} \} \) and the values \( \rr \), \( \dltwu_{3} \)
satisfying
\begin{EQA}[c]
	\DFM^{2} \leq \dmax^{2} \hspm \IFT_{\GPT} ,
	\quad
	\rr \geq \frac{3}{2} \, \biasDF \, ,
	\quad
	\dltwu_{3} \, \dmax^{2} \hspm \biasDF < \frac{4}{9} \, .
\label{8difiyfc54wrboebias}
\end{EQA}
For any linear operator \( \QP \) on \( \R^{\dimp} \), it holds 
\begin{EQA}
	\bigl\| \QP \, \bigl\{ \tarpvs_{\GPT} - \tarpvs + (\IFT_{\GPT}^{-1} \, \AvGPT)_{\tarpv} \bigr\} \bigr\| 
	& \leq & 
	\| \QP \, \DFM_{\tarpv\tarpv}^{-1} \| \, \frac{3 \dmax^{2} \dltwu_{3}}{4} \, \biasDF^{2} \, .
	\qquad
\label{0mkvhgjnrw3dfwe3u8gtyse}
\end{EQA}
If \( \Avn_{\GPT} = (\Avn_{\GPT,\tarpv},\Avn_{\GPT,\nupv}) \) then with 
\( \IFTb_{\GPT,\tarpv\tarpv} = 
\IFT_{\GPT,\tarpv\tarpv} - \IFT_{\GPT,\tarpv\nupv} \, \IFT_{\GPT,\nupv\nupv}^{-1} \, \IFT_{\GPT,\nupv\tarpv} \) 
and \( \IFTb_{\GPT,\nupv\nupv} = 
\IFT_{\GPT,\nupv\nupv} - \IFT_{\GPT,\nupv\tarpv} \, \IFT_{\GPT,\tarpv\tarpv}^{-1} \, \IFT_{\GPT,\tarpv\nupv} \); 
cf. \eqref{divjjw3d6f5f5te3gjg},
\begin{EQA}[rcl]
	(\IFT_{\GPT}^{-1} \, \Avn_{\GPT})_{\tarpv} 
	&=&
	\IFTb_{\GPT,\tarpv\tarpv}^{-1} (\Avn_{\GPT,\tarpv} - \IFT_{\tarpv\nupv} \, \IFT_{\GPT,\nupv\nupv}^{-1} \, \Avn_{\GPT,\nupv}) 
	\\
	&=&
	\IFTb_{\GPT,\tarpv\tarpv}^{-1} \, \Avn_{\GPT,\tarpv} 
	- \IFT_{\GPT,\tarpv\tarpv}^{-1} \, \IFT_{\tarpv\nupv} \, \IFTb_{\GPT,\nupv\nupv}^{-1} \, \Avn_{\GPT,\nupv} \, .
	\qquad
\label{jufyedhned78j3ibh87i4g}
\end{EQA}
\end{proposition}


\Subsection{Separable penalty}
In many situations, it is quite natural to assume an additive structure 
\( \pent_{\GPT}(\tarpv,\nupv) = \pent_{\GP}(\tarpv) + \pent_{\GPY}(\nupv) \).
This includes the case when the penalty only depends on one component \( \tarpv \) or \( \nupv \).
In the quadratic case, this means a block-diagonal structure 
\( \GPT^{2} = \blk\{ \GP^{2},\GPY^{2} \} \) yielding \( \| \GPT \prmtv \|^{2} = \| \GP \tarpv \|^{2} + \| \GPY \nupv \|^{2} \).

%
\begin{lemma}
Let \( \IFT \) be block-diagonal with \( \IFT = \blk\bigl( \IFT_{\tarpv\tarpv}, \IFT_{\nupv\nupv} \bigr) \)
and \( \pent_{\GPT}(\tarpv,\nupv) = \pent_{\GP}(\tarpv) + \pent_{\GPY}(\nupv) \).
Then 
\begin{EQA}
	(\IFT_{\GPT}^{-1} \, \AvGPT)_{\tarpv} 
	&=&
	\IFT_{\GPT,\tarpv\tarpv}^{-1} \, \Avn_{\GPT,\tarpv} \, .
\label{8vnvcte3434gffujdey}
\end{EQA}
\end{lemma}

Therefore, in the block-diagonal case, penalization of the nuisance component 
does not change the leading term in expansion \eqref{0mkvhgjnrw3dfwe3u8gtyse}.
However, this does not apply in the general situation and
penalization of the nuisance component may result in some estimation bias for the target component.
Schur's complement formula \eqref{8jdkvtwfe6xhejdcedvscy} for \( \IFT_{\GPT}^{-1} \) yields a closed-form representation.

\begin{lemma}
If \( \pent_{\GPT}(\tarpv,\nupv) = \pent_{\GP}(\tarpv) + \pent_{\GPY}(\nupv) \), then 
\begin{EQA}
	(\IFT_{\GPT}^{-1} \, \AvGPT)_{\tarpv} 
	&=&
	\IFTb_{\GPT,\tarpv\tarpv}^{-1} (\Avn_{\GPT,\tarpv} - \IFT_{\tarpv\nupv} \, \IFT_{\GPT,\nupv\nupv}^{-1} \, \Avn_{\GPT,\nupv}) 
	\\
	&=&
	\IFTb_{\GPT,\tarpv\tarpv}^{-1} \, \Avn_{\GPT,\tarpv} 
	- \IFT_{\GPT,\tarpv\tarpv}^{-1} \, \IFT_{\tarpv\nupv} \, \IFTb_{\GPT,\nupv\nupv}^{-1} \, \Avn_{\GPT,\nupv} .
\label{jufyedhned78j3ibh87i4}
\end{EQA}
\end{lemma}

\Subsection{Risk of the profile penalized MLE}
Putting together the results of Theorem~\ref{Tsemieffp3n} and Proposition~\ref{Pvarusetbse} yields the local risk bound as in 
Theorem~\ref{TQFiWibias}.

\begin{theorem}
\label{CTFiWibiasse}
Assume the conditions of Theorem~\ref{Tsemieffp3n} and Proposition~\ref{Pvarusetbse}.
It holds on \( \Omega(\xx) \) for any linear mapping \( \QP \) on \( \R^{\dimp} \)
\begin{EQA}
	\bigl\| 
		\QP \bigl\{ \tilde{\tarpv}_{\GPT} - \tarpvs - (\IFT_{\GPT}^{-1} \nabla \zeta)_{\tarpv} + (\IFT_{\GPT}^{-1} \Avn_{\GPT})_{\tarpv} \bigr\} 
	\bigr\|
	& \leq &
	\| \QP \, \DFM_{\tarpv\tarpv}^{-1} \| \, \frac{3 \dmax^{2} \dltwu_{3}}{4} \, 
	(\| \DFM \IFT_{\GPT}^{-1} \nabla \zeta \|^{2} + \biasDF^{2}) \, .
\label{rG1r1d3GbvGbseQP}
\end{EQA}
Also, define 
\begin{EQA}[c]
	\BBH_{\QP} \eqdef \Var\bigl\{ \QP (\IFT_{\GPT}^{-1} \nabla \zeta)_{\tarpv} \bigr\} ,
	\qquad
	\dimQ \eqdef \tr \BBH_{\QP} , 
	\qquad 
	\\
	\riskt_{\QP} 
	\eqdef 
	\E \bigl\{ 
		\bigl\| \QP (\IFT_{\GPT}^{-1} \nabla \zeta)_{\tarpv} - \QP (\IFT_{\GPT}^{-1} \Avn_{\GPT})_{\tarpv} \bigr\|^{2} \Ind_{\Omega(\xx)} 
	\bigr\} 
	\leq 
	\dimA_{\QP} + \| \QP \, (\IFT_{\GPT}^{-1} \, \GPT^{2} \prmtvs)_{\tarpv} \|^{2} \, .
	\qquad
\label{t6skjid7ehfcr56tegfdse}
\end{EQA}
Then with \( \dimDF = \E \| \DFM \IFT_{\GPT}^{-1} \nabla \zeta \|^{2} \), it holds
\begin{EQA}
	&& \nquad
	\E \, \bigl\{ \| \QP (\tilde{\tarpv}_{\GPT} - \tarpvs) \| \Ind_{\Omega(\xx)} \bigr\} 
	\leq 
	\riskt_{\QP}^{1/2}
	+ \| \QP \, \DFM_{\tarpv\tarpv}^{-1} \| \, \frac{3 \dmax^{2} \dltwu_{3}}{4} \, \bigl( \dimDF + \biasDF^{2} \bigr) 
	\, .
\label{EtuGus11md3Gse}
\end{EQA}
If 
\begin{EQA}[c]
	\E \bigl\{ \| \DFM \IFT_{\GP}^{-1} \nabla \zeta \|^{4} \Ind_{\Omega(\xx)} \bigr\} \leq \CONSTi_{4}^{2} \, \dimDF^{2} \, ,
\label{6hjdfv8e6hyefyeheew7skGPT}
\end{EQA}
and a constant \( \alp_{\QP} \) ensures 
\begin{EQ}[rcl]
	\alp_{\QP}
	& \eqdef & 
	\frac{\| \QP \, \DFM_{\tarpv\tarpv}^{-1} \| \, (3/4) \dmax^{2} \dltwu_{3} \, (\CONSTi_{4} \, \dimDF + \biasDF^{2})} {\sqrt{\riskt_{\QP}}} 
	< 1 \, ,
\label{6dhx6whcuydsds655srse}
\end{EQ}
then 
\begin{EQA}
	(1 - \alp_{\QP})^{2} \riskt_{\QP}
	\leq 
	\E \bigl\{ \| \QP (\tilde{\tarpv}_{\GPT} - \tarpvs) \|^{2} \Ind_{\Omega(\xx)} \bigr\}
	& \leq &
	(1 + \alp_{\QP})^{2} \riskt_{\QP} \, .
\label{EQtuGmstrVEQtGse}
\end{EQA}
\end{theorem}


\def\deta{N}
\Chapter{Nonlinear regression}
\label{Snoninverse}

Let the data \( \Yv \in \R^{n} \) following the model equation
\begin{EQA}
	\Yv
	&=&
	\regrfv(\thetav) + \epsv \in \R^{n} \, ,
\label{YvKSsepYS}
\end{EQA}
where \( \regrfv(\thetav) = \bigl( \regrf_{i}(\thetav), \, i \leq n \bigr) \in \R^{n} \) is a nonlinear 
mapping (operator) of the source signal \( \thetav \in \R^{\dimp} \) to the target space \( \R^{n} \).
Later we consider the problem of inverting the relation \( \E \Yv = \regrfv(\thetav) \) from noisy observations \( \Yv \).
A classical example is provided by nonlinear parametric regression \( \E Y_{i} = \regrf(X_{i},\thetav) \)
with a deterministic design \( X_{i} \in \R^{\dimd} \) and \( \regrf_{i}(\thetav) = \regrf(X_{i},\thetav) \),
\( i \leq n \).
More recent examples include deep neuronal networks where \( \thetav \) codes the whole DNN architecture. 
The SLS approach from 
\ifMLE{\Chname \ref{SgenBounds}}{\cite{Sp2022}} 
does not apply to this model because 
the major SLS assumptions \nameref{Eref} and \nameref{LLref} are not fulfilled in this setup due to nonlinearity of the regression function.
However, we explain below how the problem can be transformed back to the SLS framework by
extending the parameter space and localization.

\Section{Nonlinear least squares estimator}

MLE and penalized MLE procedures are often used for recovering the parameter \( \thetav \).
Assuming a zero mean nearly standardized noise \( \epsv \), 
the MLE approach leads to the nonlinear least squares problem of maximizing the random function   
\begin{EQA}
	L(\thetav) 
	&=& 
	- \frac{1}{2} \| \Yv - \regrfv(\thetav) \|^{2} .
\label{LLt2s2YkKPtt12}
\end{EQA}
The background truth \( \thetavs \) for the original model \eqref{YvKSsepYS} can be defined as 
\begin{EQA}
	\thetavs
	&=&
	\argmin_{\thetav} \| \E \Yv - \regrfv(\thetav) \|^{2} . 
\label{juctyeebftw26efu37w}
\end{EQA}
This definition allows to incorporate the case of model misspecification when \( \E \Yv \neq \regrfv(\thetav) \)
for all \( \thetav \in \Theta \).
The nonlinear function \( \regrfv(\thetav) \) in the data fidelity term \( \| \Yv - \regrfv(\thetav) \|^{2} \) creates fundamental problems 
for studying the behavior of the pMLE.
In particular, the stochastic component \( \zeta(\thetav) \) of \( L(\thetav) \) reads 
\begin{EQA}
	\zeta(\thetav)
	&=&
	L(\thetav) - \E L(\thetav)
	=
	\regrfv(\thetav) (\Yv - \E \Yv) 
	=
	\regrfv(\thetav) \, \epsv
\label{suyw2u1jsytq1qqjkdgq1qg}
\end{EQA}
and it is not linear in \( \thetav \) unless \( \regrfv(\thetav) \) is linear.
A standard sufficient condition for concavity of \( \E L(\thetav) \) reads 
\( - \nabla^{2} \E L(\thetav) \geq 0 \) and its global check for all \( \thetav \)
and a general nonlinear \( \regrfv(\cdot) \) is tricky.
Existing approaches to analysis of \( \tilde{\thetav} \) utilise deep tools from empirical processes;
see e.g. \cite[Chapter 13]{wainwright2019high}, \cite{nickl_2015}, 
\cite{Nickl2018ConvergenceRF}, \cite{nickl2022bayesian} and references therein. 
The obtained results mainly describe the rate of estimation and have been stated in the 
asymptotic minimax sense.
Our objective is to establish sharp finite sample results under realistic and mild conditions.
Unfortunately, the well-developed SLS approach does not apply to the model \eqref{YvKSsepYS}.
Both fundamental conditions \nameref{Eref} about linearity of the stochastic component and \nameref{LLref} 
about concavity of the expected log-likelihood are not fulfilled for nonlinear regression functions
\( \regrfv(\thetav) \). 
In what follows we present some ideas which enable us to reduce the study to the SLS case.
Concavity issue is addressed by using the ideas of ``warm start'' and ``localization''; 
see Section~\ref{Swarmstart} and Section~\ref{Swarmstartnl}.
The proposed ``calming'' approach enforces condition \nameref{Eref} by extending the parameter space
and relaxing the structural relation; see Section~\ref{Scalming}.

\Subsection{Noiseless case and local concavity}
\label{Swarmstart}
First, we discuss the noiseless case with deterministic observations 
\( \Yv = \E \Yv = \regrfvs = (\regrf_{i}(\thetavs)) \) leading to minimization of the fidelity 
\( \| \regrfv(\thetav) - \regrfvs \|^{2} \).
Concavity condition \nameref{LLref} is usually too restrictive for a nonlinear regression model
described by the log-likelihood \( L(\thetav) \) from \eqref{LLt2s2YkKPtt12}. 
This section explains how this condition can be relaxed using a ``warm start'' assumption.
It holds with \( L(\thetav) = - \| \regrfv(\thetav) - \regrfvs \|^{2}/2 \)
\begin{EQA}
	\IF(\thetav)
	=
	- \nabla^{2} L(\thetav)
	&=&
	\sumi \nabla \regrf_{i}(\thetav) \, \nabla \regrf_{i}(\thetav)^{\T}
	+ \sumi \bigl\{ \regrf_{i}(\thetav) - \regrfs_{i} \bigr\} \nabla^{2} \regrf_{i}(\thetav) .
	\qquad
\label{dcrf5rd5er5rde3w5w3wefdyegf}
\end{EQA}
A sufficient condition for concavity of \( L(\thetav) \) is \( \IF(\thetav) \geq 0 \), \( \thetav \in \Theta \).
Weak concavity means that \( \IF(\thetav) + \GP_{0}^{2} \geq 0 \), \( \thetav \in \Theta \),
for some \( \GP_{0}^{2} \).
Even weak concavity can be quite restrictive on the whole domain \( \Theta \).
Below we try to show how this condition can be relaxed by localization to a subset \( \Thetad \subseteq \Theta \)
containing the truth \( \thetavs \). 
Define
\begin{EQA}
	\DFN^{2}(\thetav)
	&=&
	\nabla \regrfv(\thetav) \, \nabla \regrfv(\thetav)^{\T}
	=
	\sumi \nabla \regrf_{i}(\thetav) \, \nabla \regrf_{i}(\thetav)^{\T} \in \Matr_{\dimp} \, . 
\label{tsfdc6frdw6red6wt7qagdudwo}
\end{EQA}
Injectivity of \( \regrfv(\cdot) \) means that \( \DFN^{2}(\thetav) \) is positive definite and well posed. 
If \( \sumi \nabla^{2} \regrf_{i}(\thetav) \leq \CONST \sumi \nabla \regrf_{i}(\thetav) \, \nabla \regrf_{i}(\thetav)^{\T} \)
and \( \max_{i \leq  n} |\regrf_{i}(\thetav) - \regrfs_{i}| \) is small for all \( \thetav \in \Thetad \), 
then the desired local concavity follows easily from \eqref{dcrf5rd5er5rde3w5w3wefdyegf}.
Usually, for a reasonable starting guess \( \thetav_{0} \), the local set \( \Thetad \) is taken in the form
\begin{EQA}[rcccl]
\label{ocnfhyetetgdgrwebgdgtdgt}
	\Thetad
	&=&
	\bigl\{ \thetav \colon \| \DFN_{0} (\thetav - \thetav_{0}) \| \leq \rrdd \bigr\} ,
	\qquad
	\DFN_{0}^{2} 
	&=&
	\nabla \regrfv(\thetav_{0}) \, \nabla \regrfv(\thetav_{0})^{\T} ,
\label{ocnfhyetetgdgrwebgdgtdgtD}
\end{EQA}
with a proper radius \( \rrdd \).
Localization is naturally combined with Gauss-Newton approximation: the regression function 
\( \regrfv(\thetav) \) is approximated by a linear function 
\( \regrfv(\thetav_{0}) + \langle \nabla \regrfv(\thetav_{0}), \thetav - \thetav_{0} \rangle \)
leading to the Gauss-Newton iteration 
\begin{EQA}
	\thetav_{1}
	&=&
	\argmin_{\thetav} 
	\| \regrfvs - \regrfv(\thetav_{0}) + \nabla \regrfv(\thetav_{0}) (\thetav - \thetav_{0}) \|^{2}
	\\
	&=&
	\thetav_{0} + \DFN_{0}^{-2} \, \nabla \regrfv(\thetav_{0}) \{ \regrfvs - \regrfv(\thetav_{0}) \} ;
\label{uviuwe3jewuvw3efki}
\end{EQA}
see e.g. \cite{Gr2007} for a detailed analysis and further references.
The quality of a starting guess is important, and in practice, several steps are necessary to achieve a desirable
accuracy.
Extending the approach to the case of noisy observations \( \Yv \) is not a simple task.
In particular, the most important step of showing 
\( \nabla^{2} L(\thetav) \leq 0 \) does not work for \( \Yv \) random 
because \( \| \Yv - \regrfv(\thetav) \|_{\infty} \) is not small whatever \( \thetav \) is considered.
We, however, show that the proposed approach performs essentially as a Gauss-Newton iteration with a very good 
starting guess \( \thetav_{0} = \thetavs \).  

\Subsection{Data smoothing and calming}
\label{Scalming}
This section explains the main ideas of the calming approach which allows to address nonlinearity of the stochastic 
component \( \zeta(\thetav) \); see \eqref{suyw2u1jsytq1qqjkdgq1qg}.
The basic idea is to extend the parameter space by introducing the additional parameter \( \etav \)
representing the image \( \regrfv(\thetav) \) and relaxing the structural relation \( \etav = \regrfv(\thetav) \).
This also allows us to address the issue of model misspecification.
To cope with a possibly large observation noise, we also
introduce an additional smoothing \( \Zv = \SmOp \, \Yv \) in the image space by a linear smoothing 
operator \( \SmOp \colon \R^{n} \to \R^{\dimq} \).
Further, define \( \Regrfv(\thetav) \eqdef \SmOp \, \regrfv(\thetav) \) and 
represent \eqref{YvKSsepYS} by two relations
\( \SmOp \, \Yv \approx \etav + \epsv \) and \( \etav \approx \Regrfv(\thetav) \).
Then maximization of \( L(\thetav) = - \frac{1}{2} \| \Yv - \regrfv(\thetav) \|^{2} \) is replaced by
maximization of
\begin{EQA}
	\LL(\thetav, \etav)
	&=&
	- \frac{1}{2}\| \SmOp \Yv - \etav \|^{2} - \frac{\lambda}{2} \| \SmOp \, \regrfv(\thetav) - \etav \|^{2} 
	=
	- \frac{1}{2}\| \Zv - \etav \|^{2} - \frac{\lambda}{2} \| \Regrfv(\thetav) - \etav \|^{2} 
	\qquad
\label{sm2Yg2lkKf22}
\end{EQA}
with a Lagrange multiplier \( \lambda \).
Later we fix a natural choice \( \lambda = 1 \) which is sufficient for our setup.
Now we proceed with a couple of parameters \( \upsv = (\thetav, \etav) \).
The profile MLE \( \tilde{\thetav} \) is defined as the component of \( \tilde{\upsv} \):
\begin{EQA}
	\tilde{\thetav}
	&=&
	\argmax_{\thetav} \max_{\etav} \LL(\upsv). 
\label{gcvtyehdfnvyt3wbnbkjg}
\end{EQA}
Expression \eqref{sm2Yg2lkKf22} is quadratic in \( \etav \) for a fixed \( \thetav \).
This helps to get a closed-form solution for the partial optimization problem w.r.t. \( \etav \) for 
\( \thetav \) fixed:
\begin{EQA}
	\tilde{\etav}(\thetav)
	&=&
	\argmin_{\etav} 
	\left\{ 
		\| \Zv - \etav \|^{2} + \| \Regrfv(\thetav) - \etav \|^{2} 
	\right\} 
	=
	\frac{1}{2} 
		\bigl\{ \Zv + \Regrfv(\thetav) 
		\bigr\} \, .
\label{djf53rfyrfeghtteat}
\end{EQA}
Moreover, plugging this \( \etav \) in \eqref{gcvtyehdfnvyt3wbnbkjg} yields
\begin{EQA}
	\tilde{\thetav}
	&=&
	\argmin_{\thetav} \bigl\{ \| \Zv - \Regrfv(\thetav) \|^{2} \, .
\label{kvujbvug5tfdv67g7hu8y}
\end{EQA}
Therefore, 
calming does not change the usual least squares procedure. 

The main benefit of representation \eqref{sm2Yg2lkKf22} is in possibility of applying the general SLS theory.
The stochastic data only enter in the quadratic term \( \| \Zv - \etav \|^{2} \), 
this incredibly simplifies the stochastic analysis. 
A dependence on \( \thetav \) is a bit more complicated due to the structural term
\( \| \Regrfv(\thetav) - \etav \|^{2} \) which penalizes for deviations from the forward non-linear structural relation 
\( \etav = \Regrfv(\thetav) = \SmOp \, \regrfv(\thetav) \).
However, this structural term is now deterministic and smooth.
Another benefit of calming is in introducing the image/response \( \etav \) as an additional parameter which
is also estimated by the procedure with a possibility of inference and uncertainty quantification.

The calming approach transforms the original problem into a SLS setup by extending the parameter space.
This enables us to apply the general results from \Chname \ref{SgenBounds} and \Chname \ref{SsemiMLE}
to the estimator \( \tilde{\thetav} \).
First, we present the sufficient conditions on the regression function \( \regrfv(\thetav) \) and then state the results.

\Section{Main definitions and conditions}
\label{ScondNL}
%
%
The target of estimation \( \upsvs = (\thetavs,\etavs) \) for the extended model \eqref{sm2Yg2lkKf22}
is defined by maximizing the expected log-likelihood: with \( \regrfvs = \E \Yv \) and 
\( \Regrfvs = \SmOp \regrfvs \)
\begin{EQA}[rclcl]
	\upsvs
	&=&
	\argmax_{\upsv = (\thetav,\etav)} \E \LL(\upsv) 
	&=&
	\argmin_{\upsv = (\thetav,\etav) \in \Ups} 
	\left\{ 
		\| \Regrfvs - \etav \|^{2} + \| \Regrfv(\thetav) - \etav \|^{2} 
	\right\} .
	\qquad
\label{djf53rfyrf6r6fr6fegh}
\end{EQA}
The \( \thetav \)-component \( \thetavs \) of \( \upsvs \)  
solves the original problem in which the smoothed response
\( \Zv = \SmOp \, \Yv \) is replaced by the auxiliary parameter \( \etavs \):
\begin{EQA}
	\thetavs
	&=&
	\argmin_{\thetav \in \Theta} \| \Regrfv(\thetav) - \etavs \|^{2} 
	\, .
\label{djf53rfyrf6r6fr6feghto}
\end{EQA}
\Subsection{Smoothing operator \( \SmOp \) }
\label{SopPhiv}
The proposed calming approach relies on a proper choice of the linear operator 
\( \SmOp \colon \R^{n} \to \R^{\dimq} \) given by a \( \dimq \times n \) matrix \( (\smop_{j,i}) \).
Denote its rows by 
\( \smopv_{j}^{\T} = (\smop_{j,1},\ldots,\smop_{j,n}) \) with \( \smopv_{j} \in \R^{n} \), \( j \leq \dimq \).
We mention two natural ways of choosing the operator \( \SmOp \).
The first one is the most general:
consider the row vectors \( \smopv_{j} \) of the matrix \( \SmOp \) as basis/feature vectors in 
\( \R^{n} \), \( j \leq  \dimq \).
A proper basis choice should provide \( \SmOp \, \SmOp^{\T} \asymp \Id_{\dimq} \).
A simple example is given by \( \dimq = n \) and \( \SmOp = \Id_{n} \).
However, it is desirable to ensure some additional smoothing effect by applying \( \SmOp \) in the image space
in the sense \( \| \SmOp \epsv \|_{\infty} = o(1) \).
Using the ideas of compressed sensing, one can randomly generate \( \SmOp \) with i.i.d. entries \( \smop_{j,i} \)
satisfying \( \E \smop_{j,i} = 0 \) and \( \E \smop_{j,i}^{2} = 1/n \). 
One more natural choice is given by tangent space approximation at \( \thetav_{0} \) yielding 
\( \dimq = \dimp \) and
\( \SmOp = \nabla \regrfv(\thetav_{0}) \).

It is important to ensure that the use of the calming device does not lead to a significant loss of information in the data.
Multiplication with \( \SmOp \) informally yields a kind of projection of the data \( \Yv \) on the subspace in \( \R^{n} \)
spanned by the rows \( \smopv_{j}^{\T} \) of \( \SmOp \).
In the case of linear regression \( \regrfv(\thetav) = \Psiv^{\T} \thetav \), 
the related condition of ``no information loss'' means \( \Psiv \, \Proj_{\SmOp} \, \Psiv^{\T} \approx \Psiv \, \Psiv^{\T} \),
where \( \Proj_{\SmOp} = \SmOp^{\T} (\SmOp \, \SmOp^{\T})^{-1} \SmOp \) is the projector in \( \R^{n} \) on the image of \( \SmOp \).
In the general case, we replace \( \Psiv \) with the gradient \( \nabla \regrfv(\thetav) \).

\begin{description}
    \item[\label{Phiref} \( \bb{(\SmOp)} \)]
      \emph{ With \( \Proj_{\SmOp} = \SmOp^{\T} (\SmOp \, \SmOp^{\T})^{-1} \SmOp \), it holds 
      for 
      some constant \( \CONSTSO \geq 1 \)
      }
\begin{EQA}
	\nabla \regrfv(\thetavs) \, \nabla \regrfv(\thetavs)^{\T}
	& \leq &
	\CONSTSO \nabla \regrfv(\thetavs) \, \Proj_{\SmOp} \, \nabla \regrfv(\thetavs)^{\T} \, .
\label{wfegybkf53hruygdyfys}
\end{EQA}
\end{description}

If this condition is fulfilled with \( \CONSTSO \) close to one, 
the use of \( \SmOp \)-mapping does not lead to any substantial loss of information. 

\Subsection{Local conditions on \( \Regrfv(\thetav) \) and warm start}

For the \( \dimq \)-vector \( \Regrfv(\thetav) = \SmOp \, \regrfv(\thetav) \), its gradient 
\( \nabla \Regrfv(\thetav) = \nabla \regrfv(\thetav) \, \SmOp^{\T} \) is a \( \dimp \times \dimq \)-matrix  with columns 
\( \nabla \Regrf_{j}(\thetav) = \nabla \regrfv(\thetav) \, \smopv_{j} \),
where the \( \smopv_{j}^{\T} \)'s are rows of \( \SmOp \).
Define
\begin{EQA}
	\DPNL^{2}(\thetav)
	&=&
	\frac{1}{2} \nabla \Regrfv(\thetav) \,\, \nabla \Regrfv(\thetav)^{\T}
	=
	\frac{1}{2} \sum_{j=1}^{\dimq} \nabla \Regrf_{j}(\thetav) \,\, \nabla \Regrf_{j}(\thetav)^{\T}
	\in \Matr_{\dimp}
	 \, . 
\label{tsfdc6frdw6red6wt7qagdudwg}
\end{EQA}
In what follows, similarly to the noiseless case, we limit ourselves to a local elliptic set 
\( \Thetad = \bigl\{ \thetav \colon \| \DPNLc (\thetav - \thetav_{0}) \| \leq \rrdd \bigr\} \);
cf. \eqref{ocnfhyetetgdgrwebgdgtdgt}; where \( \thetav_{0} \) 
is an initial guess, \( \DPNLc^{2} = \DPNL^{2}(\thetav_{0}) \), and \( \rrdd \) is a properly selected radius.
%
%
An important ``warm start'' condition means that the starting guess \( \thetav_{0} \)
is reasonable and the targets \( \thetavs \) from \eqref{djf53rfyrf6r6fr6feghto} 
is within \( \Thetad \).

\begin{description}

\item[\( \bb{(\thetavs)} \)\label{thetaMref}]
	\emph{It holds \( \thetavs \in \Thetad \).	
	}
\end{description}
Conditions of this kind are often applied in nonlinear optimization for studying, e.g. Gauss-Newton iterations;
see e.g. \cite{Gr2007}.

Later we assume the following regularity and smoothness conditions.

\begin{description}
    \item[\label{LLM2ref} \( \bb{(\nabla \Regrfv)} \)]
      \emph{ 
      For some \( \dltwbD \leq 1/3 \)
      and any \( \thetav \in \Thetad \), it holds
\begin{EQA}
	(1 - \dltwbD) \, \DPNLc^{2}
	& \leq &
	\DPNL^{2}(\thetav)
	\leq 
	(1 + \dltwbD) \, \DPNLc^{2} .
\label{fchghiu87686574e5rtyyu}
\end{EQA}
}
%
%
%
%
%
\def\deta{N}
    \item[\label{LLMkref} \( \bb{(\nabla^{k} \Regrf)} \)]
      \emph{ For \( k \in \{ 2,3,4 \} \) and some \( \deta \geq 1 \), 
      uniformly over \( \thetav \in \Thetad \) and \( \uv \in \R^{\dimp} \)
\begin{EQA}
	\sum_{j=1}^{\dimq} \langle \nabla^{k} \Regrf_{j}(\thetav), \uv^{\otimes k} \rangle^{2}
	& \leq &
	\frac{1}{\deta^{k-1}} \, \| \DPNL(\thetav) \uv \|^{2k} \, .
\label{wfegy7r5qrw35edfhgdyfys}
\end{EQA}
}
\end{description}

\begin{remark}
The constant \( \deta \) in \eqref{wfegy7r5qrw35edfhgdyfys} may depend on \( k \).
We use the same \( \deta \) for ease of notation.
Smoothness \nameref{LLtS3ref} of \nameref{LLtS4ref} of \( \Regrfv \) yields \( \deta^{-1} \asymp \| \DPNLc^{-2} \| = n^{-1} \)%
; see Section~\ref{Slocalsmooth}.
\end{remark}
The radius \( \rrde \) of the local set \( \Thetad \) 
should be sufficiently large to ensure that the full dimensional estimator \( \tilde{\thetav} \) 
concentrates on this set and, at the same time, sufficiently small to ensure a proper localization;
see Theorem~\ref{TtifvsGPeff} later.
%
%


\Subsection{Full dimensional information matrix and identifiability}
Localization is an important tool in establishing local identifiability for
the full dimensional parameter \( \upsv \). 
Define \( \upsv_{0} = (\thetav_{0},\etav_{0}) \) with \( \etav_{0} = \Regrfv(\thetav_{0}) \) so that
the structural relation \( \etav \approx \Regrfv(\thetav) \) is precisely 
fulfilled at the starting point \( \upsv_{0} = (\thetav_{0},\etav_{0}) \).
It is convenient to consider local sets of product structure in \( \thetav \) and \( \etav \) in the form
\begin{EQA}[rcl]
	\Upsd
	&=&
	\Thetad \times \Etad
	=
	\bigl\{ \upsv = (\thetav,\etav) \colon 
		\| \DPNLc (\thetav - \thetav_{0}) \|^{2} \leq \rrde^{2}, \,\,
		\| \etav - \etav_{0} \|^{2} \leq \rrde^{2} 
	\bigr\} .
	\qquad
\label{cujduyfdt5c544retd5td54e4}
\end{EQA}
%
%
The gradient and Hessian of \( \LL(\upsv) \) read as follow: 
with \( \nabla \Regrfv(\thetav) = \nabla \regrfv(\thetav) \, \SmOp^{\T} \in \R^{\dimp \times \dimq} \)
\begin{EQA}
	\frac{\partial}{\partial\thetav} \LL(\thetav,\etav)
	&=& 
	- \nabla \Regrfv(\thetav) \bigl\{ \Regrfv(\thetav) - \etav \bigr\} \, ,
	\\
	\frac{\partial}{\partial\etav} \LL(\thetav,\etav)
	&=& 
	(\Zv - \etav) 
	+ \bigl\{ \Regrfv(\thetav) - \etav \bigr\} \, ,
\label{Lth122Y1ppp}
\end{EQA}
and
\begin{EQA}
	\IFT(\upsv)
	\eqdef
	- \nabla^{2} \LL(\upsv)
	&=&
	\begin{pmatrix}
	\IF(\upsv) & - \nabla \Regrfv(\thetav) 
	\\
	- \nabla \Regrfv(\thetav)^{\T} & 2 \, \Id_{\dimq} 
	\end{pmatrix}
	\, 
\label{G2Kzsm2PY}
\end{EQA}
with the upper left diagonal block
\begin{EQ}[rcl]
	\IF(\upsv)
	& \eqdef &
	\nabla \Regrfv(\thetav) \, \nabla \Regrfv(\thetav)^{\T} 
		+ \sum_{j=1}^{\dimq} \{ \Regrf_{j}(\thetav) - \eta_{j} \} \, \nabla^{2} \Regrf_{j}(\thetav) 
	\, .
	\qquad
\label{DPGP2HPGP22}
\end{EQ}
For our results we need that the matrix 
\( \IFT(\upsv) \) is well posed with a reasonable 
conditional number for all \( \upsv \in \Upsd \).
As in Section~\ref{Swarmstart} for the noiseless case, we use the ideas of warm start and localization. 
Lemma~\ref{Sconccalm} states 
\begin{EQA}
	\IFT(\upsv)
	& \geq &
	\dmax^{-2} \blk\bigl\{ \DPNL^{2}(\thetav), \, \Id_{\dimq} \bigr\} 
\label{d78dwifeire7f73udbcjw}
\end{EQA} 
for \( \dmax^{2} = 2 \) and any point \( \upsv = (\thetav,\etav) \) in a local vicinity \( \Upsd \)
of \( \upsv_{0} \).  
This particularly implies that the function \( \LL(\upsv) \) is strongly concave on \( \Upsd \).
Denote \( \IFT = \IFT(\upsvs) \).
The inverse matrix \( \IFT^{-1} \) will be used in our results. 
By Schur's complement formula, 
\begin{EQA}
	\IFT^{-1}
	&=&
	\begin{pmatrix}
		\IFTb_{\thetav\thetav}^{-1} & \frac{1}{2} \IFTb_{\thetav\thetav}^{-1} \, \nabla \Regrfv  
		\\
		\frac{1}{2} (\IFTb_{\thetav\thetav}^{-1} \, \nabla \Regrfv)^{\T} & \IFTb_{\etav\etav}^{-1}
	\end{pmatrix} \,
\label{y67ytef8y8wuj3i4gf5rgt34h}
\end{EQA}
with
\begin{EQ}[rcl]
	\IFTb_{\thetav\thetav}
	&=&
	\IF - \frac{1}{2} \nabla \Regrfv \, \nabla \Regrfv^{\T} 
	=
	\IF - \DPNL^{2} \, ,
	\\
	\IFTb_{\etav\etav}
	&=&
	2 \Id_{\dimq} - \nabla \Regrfv^{\T} \, \IF^{-1} \, \nabla \Regrfv \, ;
\label{ydsydfyfyedtw3e23w4er}
\end{EQ}
see \eqref{8jdkvtwfe6xhejdcedvscy} of Lemma~\ref{LSchur}.

%
\Subsection{Stochastic term} 
\label{Scondgllo} 
Now we check the general conditions \nameref{Eref} and \nameref{EU2ref}
from Section~\ref{Scondgeneric} for the considered case
with the extended log-likelihood \( \LL(\upsv) \) from \eqref{sm2Yg2lkKf22}.
We heavily use that the data \( \Yv \) only 
enter in the fidelity term \( \| \Zv - \etav \|^{2}/2 \) with \( \Zv = \SmOp \, \Yv \)
and the stochastic term linearly depends on \( \etav \) and is free of
\( \thetav \).
By \eqref{sm2Yg2lkKf22}, it holds for \( \zeta(\upsv) = \LL(\upsv) - \E \, \LL(\upsv) \)
\begin{EQA}
	\nabla \zeta
	&=&
	\begin{pmatrix}
		0
		\\
		\nabla_{\etav} \zeta
	\end{pmatrix}
	=
	\begin{pmatrix}
		0
		\\
		\SmOp \epsv
	\end{pmatrix} \, ,
\label{nanafvgvsm2e}
\end{EQA}
and condition \nameref{Eref} is fulfilled.
Bounding \( \nabla \zeta \) can be easily reduced to a similar question for \( \SmOp \epsv \).
Later we assume the following condition.

\begin{description}
\item[\( \bb{(\SmOp \epsv)} \)\label{ED0ref}]
  \emph{The vector \( \SmOp \epsv \) satisfies 
  for all considered \( \xx > 0 \) 
  }
\begin{EQA}[c]
	\P\bigl( \| \SmOp \epsv \| > \zq(\HVB^{2},\xx) \bigr)
	\leq 
	3 \ex^{-\xx} ,
\label{expzetaclocGP}
\end{EQA}
\emph{where} 
\begin{EQA}
	\HVB^{2} 
	& \eqdef &
	\Var(\SmOp \epsv) = \SmOp \, \Var(\epsv) \, \SmOp^{\T} ,
	\\
	\zq(\HVB^{2},\xx)
	& \eqdef & 
	\sqrt{\tr \HVB^{2}} + \sqrt{2 \xx \, \| \HVB^{2} \|} \,\, .
	\qquad
\label{vcdguhwedu82quhwwet543}
\end{EQA}
\end{description} 

\noindent
General results from  
\ifsqnorm{Section~\ref{SdevboundnonGauss}}{\cite{Sp2023c}, \cite{Sp2023d}} 
ensure such deviation bounds under exponential moment conditions on \( \SmOp \epsv \).
If \( \HVB^{2} \leq \sigma^{2} \, \Id_{\dimq} \) then
\begin{EQA}
	\zq(\HVB^{2},\xx)
	& \leq &
	\sigma (\sqrt{\dimq} + \sqrt{2\xx}) .
\label{ojjuy6664fhjkjklolo}
\end{EQA}
With \( \nabla \zeta = (0,\SmOp \epsv) \), \( \IFT = \IFT(\upsvs) \), and
\( \DFM^{2} = \blk\{ \DPNL^{2}, \Id_{\dimq} \} \), it holds by Lemma~\ref{LDFIFuv}
\begin{EQA}
	\| \DFM \, \IFT^{-1} \nabla \zeta \|
	& \leq &
	2 \| \SmOp \epsv \| \, .
\label{jvbhgdftet6ed56hh}
\end{EQA}
By condition \nameref{ED0ref},
on a set \( \Omega(\xx) \) with \( \P(\Omega(\xx)) \geq 1 - 3 \ex^{-\xx} \)
\begin{EQA}
	\| \DFM \, \IFT^{-1} \nabla \zeta \|
	\leq 
	\rrDF 
	& \eqdef &
	2 \zq(\HVB^{2},\xx) \, .
	\qquad
\label{tggtfrew3345456dwjkj}
\end{EQA}
By \eqref{y67ytef8y8wuj3i4gf5rgt34h}
\begin{EQA}
	\bigl( \IFT^{-1} \nabla \zeta \bigr)_{\thetav}
	&=&
	\frac{1}{2} \IFTb_{\thetav\thetav}^{-1} \, \nabla \Regrfv \, \SmOp \epsv \, .
\label{dufhiwr93rujgkjl65y3}
\end{EQA}

\Section{Expansions and accuracy guarantees for nonlinear regression}
%
%
Restricting the parameter set \( \Ups \) to the local set \( \Upsd \) from \eqref{cujduyfdt5c544retd5td54e4}
and the calming device bring the original problem back to the SLS setup 
with a linear stochastic component and a smooth and concave expected log-likelihood.
This allows us to apply the well-developed general results from 
\ifNL{\Chname \ref{SgenBounds} and \Chname \ref{SsemiMLE}}{Section~2 of \cite{Sp2022}}.

\begin{theorem}
\label{TtifvsGPeff}
Assume \nameref{LLM2ref}, 
\nameref{LLMkref} for \( k=2,3 \), 
\nameref{thetaMref}, 
\nameref{ED0ref}. 
Let 
\begin{EQA}[c]
	2 \, \rrde  
	<  
	\frac{\sqrt{\deta}}{2} \,\, ,
	\quad
	\rrdd \geq \frac{3}{2} \, \rr_{\DFM} \, ,
	\quad
	\const_{3} \, \rr_{\DFM} < \frac{2 \sqrt{\deta}}{9} \, ,
\label{8difiyfc54wrboeNL}
\end{EQA}
where \( \rr_{\DFM} \) is from \eqref{tggtfrew3345456dwjkj}, \( \rrdd \) from \eqref{cujduyfdt5c544retd5td54e4}, and the constant
\( \const_{3} \) depends on \( \dltwbD \) from \nameref{LLM2ref} only.
Then 
on a set \(\Omega(\xx) \) with \( \P\bigl( \Omega(\xx) \bigr) \geq 1 - 3 \ex^{-\xx} \), the estimate 
\( \tilde{\upsv} \) satisfies
\begin{EQA}[rcl]
    \| \DFM^{-1} \IFT (\tilde{\upsv} - \upsvs - \IFT^{-1} \nabla \zeta) \|
    & \leq &
    \frac{3 \const_{3}}{4 \sqrt{\deta}} \| 2 \SmOp \epsv \|^{2} 
	\, .
\label{DGttGtsGDGm13rG22NL}
\end{EQA}
\end{theorem}

\begin{proof}
We apply Theorem~\ref{Tconcsemi} after restricting the parameter set to \( \Upsd \) from \eqref{cujduyfdt5c544retd5td54e4}.
Condition \nameref{Eref} is fulfilled by construction, \nameref{LLref} follows by \eqref{d78dwifeire7f73udbcjw} and \nameref{LLM2ref}.
Further, \nameref{EU2ref} follows by \nameref{ED0ref} and \eqref{jvbhgdftet6ed56hh}.
Lemma~\ref{Lupssmooth} ensures \nameref{LLsT3ref} with \( \dltwu_{3} = \const_{3} / \sqrt{\deta} \).
Condition \eqref{8difiyfc54wrboeMLEse} with \( \dmax^{2} = 2 \) follow from \eqref{8difiyfc54wrboeNL}.
Therefore, all the conditions of Theorem~\ref{Tconcsemi} are fulfilled and also
by \eqref{jvbhgdftet6ed56hh} \( \| \DFM \, \IFT^{-1} \nabla \zeta \| \leq 2 \| \SmOp \epsv \| \).
Now \eqref{DGttGtsGDGm13rG22NL} follows from \eqref{yfjvh9h4f5e53tghfugu}.
\end{proof}

\begin{remark}
Lemma~\ref{Lupssmooth} describes explicitly \( \const_{3} \) and a similar quantity \( \const_{4} \).
\end{remark}

\begin{remark}
Bound \eqref{DGttGtsGDGm13rG22NL} and \( \DFM^{2} \leq 2 \IFT \) imply on \( \Omega(\xx) \)
\begin{EQA}
    \| \DFM (\tilde{\upsv} - \upsvs - \IFT^{-1} \nabla \zeta) \|
	& \leq &
	2 \| \DFM^{-1} \IFT (\tilde{\upsv} - \upsvs - \IFT^{-1} \nabla \zeta) \|
    \leq 
    \frac{3 \const_{3}}{2 \sqrt{\deta}} \| 2 \SmOp \epsv \|^{2}
\label{duhyvf83nheri3f6rt6}
\end{EQA} 
yielding by \eqref{jvbhgdftet6ed56hh} and \eqref{8difiyfc54wrboeNL}
\begin{EQA}
	\| \DFM (\tilde{\upsv} - \upsvs) \|
    & \leq &
    \| \DFM \IFT^{-1} \nabla \zeta \| + \frac{6 \const_{3}}{\sqrt{\deta}} \| \SmOp \epsv \|^{2} \, 
    \leq 
    \| 2 \SmOp \epsv \| \bigl( 1 + \frac{3 \const_{3} \, \rr_{\DFM}}{\sqrt{\deta}} \bigr) 
    \leq 
    \frac{5}{3} \| 2 \SmOp \epsv \| \, .
\label{ycduydjwe37bvhjri934}
\end{EQA}
The main problem with this bound is that the value \( \| \SmOp \epsv \| \) is of order \( \zq(\HVB^{2},\xx) \) and it
corresponds to the full parameter dimension \( \dimq \) and can be quite large.
A great benefit of expansion \eqref{vdud6ehry3hrf67ruytaNL} is that it allows to improve the leading term of the error in the Fisher expansion 
by projecting on the target direction.
\end{remark}

The next result specifies Theorem~\ref{Tsemieffp3n} and Theorem~\ref{TsemiWilks3} to nonlinear regression.

\begin{theorem}
\label{TtifvsGPeffse}
Under the conditions of Theorem~\ref{TtifvsGPeff}, for any linear mapping \( \QP \) on \( \R^{\dimp} \)
\begin{EQA}
\label{vdud6ehry3hrf67ruytaNL}
	\| \QP \{ \tilde{\thetav} - \thetavs - (\IFT^{-1} \nabla \zeta)_{\thetav} \} \|
	& \leq &
	\| \QP \, \DPNL^{-1} \| \, \frac{3 \, \const_{3}}{2 \sqrt{\deta}} \, \| 2 \SmOp \epsv \|^{2} \, .
\end{EQA}
Moreover, it holds on \(\Omega(\xx) \) with \( \LLp(\thetav) = \sup_{\etav} \LL(\thetav,\etav) \)
\begin{EQA}[rcl]
	\Bigl| 
		2 \LLp(\tilde{\thetav}) - 2\LLp(\thetavs)
		- \| \xivr \|^{2} 
	\Bigr|
	& \leq &
	\frac{\const_{3}}{\sqrt{\deta}} \, \| 2 \SmOp \epsv \|^{3} \, ,
\label{222bb23G2L2rGsonlr}
\end{EQA}
where with \( \IFTb_{\thetav\thetav} \) and \( \IFTb_{\etav\etav} \) from \eqref{ydsydfyfyedtw3e23w4er}
\begin{EQA}
	\xivr
	& \eqdef &
	\IFTb_{\thetav\thetav}^{1/2} \, (\IFT^{-1} \nabla \zeta)_{\thetav}
	=
	\IFTb_{\thetav\thetav}^{-1/2} \, \nabla \Regrfv \, \SmOp \epsv \, .
\label{d7fje3gv7hre3je8gb53jii}
\end{EQA}
\end{theorem}

\begin{proof}
Statements \eqref{vdud6ehry3hrf67ruytaNL} follows from \eqref{vdud6ehry3hrf67ruyta}.
Further, \eqref{vdud6ehry3hrf67ruy233}
and \( \| \DFM_{\nupv\nupv} \IFT_{\nupv\nupv}^{-1} \scorem{\nupv} \zeta \| \leq 2 \| \SmOp \epsv \| \)
yield \eqref{222bb23G2L2rGsonlr}.
\end{proof}

The next result specifies Theorem~\ref{CTFiWibiasse} to our regression setup.
For a linear mapping \( \QP \) on \( \thetav \), define
\begin{EQA}[rcl]
\label{t6skjid7ehfcr56tegfd}
	\dimQ
	& \eqdef &
	\tr \Var\bigl\{ \QP (\IFT^{-1} \nabla \zeta)_{\thetav} \bigr\} \, .
	\qquad 
	\\
	\riskt_{\QP} 
	& \eqdef &
	\E \bigl\{ 
		\bigl\| \QP (\IFT^{-1} \nabla \zeta)_{\thetav} \bigr\|^{2} \Ind_{\Omega(\xx)} 
	\bigr\} 
	\leq 
	\dimQ \, .
	\qquad
\end{EQA}
The next result provides upper bounds for the loss and risk of \( \tilde{\thetav} \).

\begin{theorem}
\label{TNLriskQ}
Assume the conditions of Theorem~\ref{TtifvsGPeff}.
Then it holds on \( \Omega(\xx) \) 
\begin{EQA}
	\bigl\| \QP \{ \tilde{\thetav} - \thetavs 
			- (\IFT^{-1} \nabla \zeta)_{\thetav} \}
	\bigr\|
	& \leq &
	\| \QP \, \DPNL^{-1} \| \,
	\frac{3 \, \const_{3}}{2 \sqrt{\deta}} \, \| 2 \SmOp \epsv \|^{2}  \, ;
	\qquad 
\label{rG1r1d3GbvGbNLQ}
\end{EQA}
see \eqref{dufhiwr93rujgkjl65y3}.
Further, with 
\( \dimDF = \E \| \DFM \, \IFT^{-1} \nabla \zeta \|^{2} \leq 4 \E \| \SmOp \epsv \|^{2} \), 
it holds
\begin{EQA}
	&& \nquad
	\E \, \bigl\{ \| \QP (\tilde{\thetav} - \thetavs) \| \Ind_{\Omega(\xx)} \bigr\} 
	\leq 
	\dimQ^{1/2}
	+ \| \QP \, \DPNL^{-1} \| \, \frac{3 \, \const_{3}}{2 \sqrt{\deta}} \, \dimDF  
	\, .
\label{EtuGus11md3Gre}
\end{EQA}
If \( 4 \E^{1/2} \bigl\{ \| \SmOp \epsv \|^{4} \Ind_{\Omega(\xx)} \bigr\} \leq \CONSTi_{4} \, \dimDF \) and
a constant \( \alp_{\QP} \) ensures 
\begin{EQ}[rcl]
	\alp_{\QP}
	& \eqdef & 
	\frac{\| \QP \, \DPNL^{-1} \| \, (3/4)\const_{3} \, \CONSTi_{4} \, \dimDF} {\sqrt{\deta \, \dimQ}} 
	< 1 \, ,
\label{6dhx6whcuydsds655srseNL}
\end{EQ}
then 
\begin{EQA}
	(1 - \alp_{\QP})^{2} \riskt_{\QP}
	\leq 
	\E \bigl\{ \| \QP (\tilde{\thetav} - \thetavs) \|^{2} \Ind_{\Omega(\xx)} \bigr\}
	& \leq &
	(1 + \alp_{\QP})^{2} \riskt_{\QP} \, .
\label{EQtuGmstrVEQtGseNL}
\end{EQA}
\end{theorem}

\subsection*{Critical dimension condition}
Now we discuss the issue of \emph{critical dimension}.
The \emph{full effective dimension} \( \dimDF \) of the extended model satisfies by Lemma~\ref{LDFIFuv}
\begin{EQA}
	\dimDF
	&=&
	\E \| \DFM \, \IFT^{-1} \nabla \zeta \|^{2}
	\leq 
	\E \| 2 \SmOp \epsv \|^{2}
	\leq 
	4 \tr \HVB^{2}
\label{hcy7y6ehyudf76fbrjj}
\end{EQA}
and for the homogeneous noise \( \HVB^{2} \leq \sigma^{2} \, \Id_{\dimq} \), it holds
\( \dimDF \leq 4 \sigma^{2} \, \dimq \).
One can see that complexity of the full dimensional model measured by the effective dimension \( \dimDF \)
can be controlled via the dimension of the image space \( \dimq \).
The \emph{target effective dimension} \( \dimQ \) is given by \eqref{t6skjid7ehfcr56tegfd} and it can be significantly
smaller than \( \dimDF \).
To be specific assume \( \QP = \IFTb_{\thetav\thetav}^{1/2} = (\IF - \DPNL^{2})^{1/2} \).
Under homogeneous noise \( \HVB^{2} \leq \sigma^{2} \Id_{\dimq} \), it holds 
by \eqref{d7fje3gv7hre3je8gb53jii} and Lemma~\ref{LDFblk2}
\begin{EQA}
	\dimQ
	&=&
	\Var\bigl\{ \IFTb_{\thetav\thetav}^{1/2} (\IFT^{-1} \nabla \zeta)_{\thetav} \bigr\}
	=
	\Var\bigl\{ \IFTb_{\thetav\thetav}^{-1/2} \, \nabla \Regrfv \, \SmOp \epsv \bigr\}
	\\
	& \leq &
	\sigma^{2} \tr \bigl( \IFTb_{\thetav\thetav}^{-1} \, \nabla \Regrfv \, \nabla \Regrfv^{\T} \bigr)
	=
	2 \sigma^{2} \tr \bigl( \IFTb_{\thetav\thetav}^{-1} \, \DPNL^{2} \bigr)
	\, .
\label{7jdcd73wjhdifct63whd}
\end{EQA}
This value corresponds to the effective dimension of the target parameter \( \thetav \).

The \emph{effective sample size} \( n \) is defined via the constant \( \deta \) from \nameref{LLMkref}.
We use \( \deta \asymp n \).
Theorem~\ref{TtifvsGPeff} and Theorem~\ref{TNLriskQ} require that
the error terms in all the expansions are sufficiently small.
In particular, bound \eqref{EtuGus11md3Gre} assumes that
\( \deta^{-1/2} \, \dimDF \ll \dimQ \) or 
\( \dimDF^{2} \ll n \, \dimQ^{2} \).
These conditions also ensure that \( \alp_{\QP} \ll 1 \); see \eqref{6dhx6whcuydsds655srseNL}.
However, in the case when the finite effective target dimension \( \dimQ \)
is much smaller than \( \dimDF \), the condition \( \dimDF^{2} \ll n \) can be quite restrictive.
Wilks expansion \eqref{222bb23G2L2rGsonlr} is even more demanding.
For the leading term of the expansion, it holds \( \| \xivr \|^{2} \approx \dimQ \), and 
\eqref{222bb23G2L2rGsonlr} is only meaningful if \( n^{-1/2} \, \dimDF^{3/2} \ll \dimQ \). 


\iffourG{
\Section{Profile MLE. 4G bounds}
\label{SNLfourth}
This section presents advanced risk bounds for the profile MLE \( \tilde{\thetav} \) based on fourth-order expansions.
We follow the line of Section~\ref{Sfourthsemi}.
Introduce the third-order tensor
\( \TensU(\wv) \eqdef \frac{1}{6} \bigl\langle \nabla^{3} \fs(\upsvs),\wv^{\otimes 3} \bigr\rangle \),
where \( - 2 \fs(\upsv) = - 2 \fs(\thetav,\etav) = \| \Regrfvs - \etav \|^{2} + \| \Regrfv(\thetav) - \etav \|^{2} \); see \eqref{djf53rfyrf6r6fr6fegh}.
With \( \nabla \zeta = (0,\SmOp \epsv) \), define 
\begin{EQ}[rcccl]
	\svn 
	&=&
	\IFT^{-1} \bigl\{ \nabla \zeta + \nabla \TensU(\IFT^{-1} \nabla \zeta) \bigr\} 
	&=&
	(\svn_{\thetav},\svn_{\etav}) \, .
\label{8vfjvr43223efryfuweefgNL}
\end{EQ}
Theorem~\ref{Tsemieffp4} yields the following result.

\begin{theorem}
\label{TNLrisk4a}
Assume the conditions of Theorem~\ref{TtifvsGPeff} and, in addition, \nameref{LLMkref} holds for \( k=4 \)
and \( \deta \, \const_{4} \, \rrDF^{2} < 1/6 \),
where \( \rrDF \) is from \eqref{tggtfrew3345456dwjkj}, 
\( \rrdd \) from \eqref{cujduyfdt5c544retd5td54e4}, and the constant
\( \const_{4} \) depends on \( \dltwbD \) from \nameref{LLM2ref} only.
Then on \( \Omega(\xx) \), 
for any linear mapping \( \QP \) of \( \thetav \)
\begin{EQ}[rcl]
	\| \QP \, ( \tilde{\thetav} - \thetavs - \svn_{\thetav} ) \|
	& \leq &
	\| \QP \, \DPNL^{-1} \| \, 
	\deta^{2} (\const_{4} + 4 \const_{3}^{2}) \, \| 2 \SmOp \epsv \|^{3} \, ,
	\qquad
\label{0mkvhgjnrwwe3u8gtygNL}
\end{EQ}
and
\begin{EQ}[rcl]
	\| \QP (\svn - \IFT^{-1} \, \nabla \zeta)_{\thetav} \|
	& \leq &
	\| \QP \, \DPNL^{-1} \| \, \deta \, \const_{3} \, \| 2 \SmOp \epsv \|^{2} \, .
	\qquad
\label{vuedy766t4e3bfvyt6ere}
\end{EQ}
Also, let \( \E \bigl\{ \| 2 \SmOp \epsv \|^{3} \Ind_{\Omega(\xx)} \bigr\} \leq \CONSTi_{3}^{2} \, \dimDF^{3/2} \) 
with \( \dimDF = \E \| 2 \SmOp \epsv \|^{2} \).
Then
\begin{EQA}
	&& \nquad
	\E \bigl\{ \| \QP \, (\tilde{\thetav} - \thetavs) \| \Ind_{\Omega(\xx)} \bigr\}
	\leq 
	\E \bigl\{ \| \QP \svn_{\thetav}  \| \Ind_{\Omega(\xx)} \bigr\}
	+ \, \| \QP \, \DPNL^{-1} \| \, 
	\deta^{2} (\const_{4} + 4 \const_{3}^{2}) \, \CONSTi_{3} \, \dimDF^{3/2}  
	\, ,
\label{0mkvhgjnrw3dfwe3u8gtygE1NL}
	\\
	&& \nquad
	\Bigl| \E \bigl\{ \| \QP \svn_{\thetav} \| \Ind_{\Omega(\xx)} \bigr\}
	- \E \bigl\{ \| \QP \, (\IFT^{-1} \, \nabla \zeta)_{\thetav} \| \Ind_{\Omega(\xx)} \bigr\} 
	\Bigr|
	\leq 
	\| \QP \, \DPNL^{-1} \| \, \deta \, \const_{3} \, \dimDF
	\, .
\label{udtgecthwjdytdehduuc6NL}
\end{EQA}
\end{theorem}

\begin{proof}
Lemma~\ref{Lupssmooth} checks that \( f(\upsv) = \E \LL(\upsv) \) follows \nameref{LLsT3ref} and \nameref{LLsT4ref} for all \( \upsv \in \Upsd \) with 
\( \dltwu_{3} =	\const_{3} \, \deta \) and \( \dltwu_{4} =	\const_{4} \, \deta^{2} \),
where \( \const_{3} \) and \( \const_{4} \) depend \( \dltwbD \) from \nameref{LLM2ref} only.
Now the results follow from Theorem~\ref{Tsemieffp4} with \( \dmax^{2} = 2 \).
\end{proof}

\noindent
The squared risk can be bounded in a similar way using Theorem~\ref{Tseff42}.
Define 
\begin{EQA}
	\riskt_{\QP}
	& \eqdef &
	\E \bigl\{ \| \QP \, (\IFT^{-1} \, \nabla \zeta)_{\thetav} \|^{2} \Ind_{\Omega(\xx)} 
	\bigr\} \, ,
\label{7dhdrw2bfvu78u78e4ndwNL}
	\\
	\riskt_{\QP,2} 
	& \eqdef & 
	\E \bigl\{ \| \QP \svn_{\thetav} \|^{2} \Ind_{\Omega(\xx)} \bigr\} \, .
\label{yfhcvhched6chejdrteNL}
\end{EQA}

\begin{theorem}
\label{Tseff42NL}
Assume the conditions of Theorem~\ref{TNLrisk4a} and let
\begin{EQA}[c]
	\E \bigl\{ \| 2 \SmOp \epsv \|^{k} \Ind_{\Omega(\xx)} \bigr\} \leq \CONSTi_{k}^{2} \, \dimDF^{k/2} ,
	\qquad
	k=3,4,6 \, .
\label{6hjdfv8e6hyefyeheew7skNL}
\end{EQA}
For a linear mapping \( \QP \) and \( \riskt_{\QP,2} \) from \eqref{yfhcvhched6chejdrteNL}, it holds
\begin{EQA}
	\bigl( 1 - \alp_{\QP,2} \bigr)^{2} \riskt_{\QP,2}
	\leq 
	\E \bigl\{ \| \QP \, (\tilde{\thetav} - \thetavs) \|^{2} \Ind_{\Omega(\xx)} \bigr\}
	& \leq &
	\bigl( 1 + \alp_{\QP,2} \bigr)^{2} \riskt_{\QP,2} \, 
\label{6shx76whnjvyehfbvyfhseNL}
\end{EQA}
provided that
\begin{EQA}
	\alp_{\QP,2}
	& \eqdef &
	\| \QP \, \DPNL^{-1} \| \, 
	\frac{\deta^{2} ( \const_{4}/3 + \const_{3}^{2} ) \, \CONSTi_{6} \, \dimDF^{3/2}}{\sqrt{\riskt_{\QP,2}}} 
	< 1 \, .
\label{cuyhe4f8jfdebnuiseNL}
\end{EQA}
If another value \( \alp_{\QP,1} < 1 \) is such that  
\begin{EQ}[rcl]
	&&
	\| \QP \, \DPNL^{-1} \| \, \frac{\deta \, \const_{3}}{2} \, \CONSTi_{6} \, \dimDF^{3/2} 
	\leq 
	\alp_{\QP,1} \, \sqrt{\riskt_{\QP}} \, 
\label{6dhx6whcuydsds655srew4seNL}
\end{EQ}
with \( \riskt_{\QP} \) from \eqref{7dhdrw2bfvu78u78e4ndwNL} then
\begin{EQA}
	\riskt_{\QP} (1 - \alp_{\QP,1})^{2} 
	\leq 
	\riskt_{\QP,2}
	& \leq &
	\riskt_{\QP} (1 + \alp_{\QP,1})^{2} \, .
\label{EQtuGmstrVEQtGQseNL}
\end{EQA}
\end{theorem}

\subsection*{Critical dimension condition}

Now we discuss how the fourth-order expansions improve the issue of critical dimension.
We again suppose that \( \| \QP \DPNL^{-1} \| = 1 \) and \( \deta = n^{-1/2} \).
In view of \( \| 2 \SmOp \epsv \|^{2} \approx \dimDF \) on \( \Omega(\xx) \),
the error term in \eqref{0mkvhgjnrwwe3u8gtygNL} is of order \( n^{-1} \, \dimDF^{3/2} \) 
and it is small provided that \( \dimDF^{3/2} \ll n \).
This is a substantial relaxation of \( n^{-1/2} \, \dimDF \ll 1 \) as in \eqref{rG1r1d3GbvGbNLQ}.
This improvement is due to the fact that the full dimensional error term of the expansion becomes smaller with the degree of expansion,
while the leading term \( \svn \) can be just projected on the target direction \( \thetav \).
}{}


\Section{Deep Neuronal Network with one hidden layer}
This section specifies the general results to the case of a Deep Neuronal Network (DNN).
There is a vast literature on this problem, recent papers include
    \cite{ShHi2020}, 
    \cite{SYZ2019}, 
    \cite{FGZ2022}, 
    \cite{KK2022}, 
    \cite{SJLH2022}, 
    \cite{LCEZZ2022}, 
    \cite{Si2022}. 
Most of the results are stated in terms of asymptotic minimax rate of convergence 
for the regression function. 
The common problem in these results is that the involved function classes
are very large leading to extremely poor rate results. 
Manifold learning approach is used to reduce the curse of dimensionality issue; see
    \cite{CJLZ2019}, 
%
    \cite{KKL2019}, 
    \cite{JSLH2021}, 
    \cite{LCZL2021}, 
    \cite{cloninger2021deep}, 
    \cite{ZCWLZ2023}, 
    \cite{SJLH2023}, 
    \cite{JSLH2023}. 
However, even under manifold assumption, the rates are not particularly informative. 
The approach of this paper attempts to restate the results in terms of efficient dimension \( \dimA \).
In addition, this paper focuses on understanding the DNN structure rather 
than estimation of the response function.
The main results bound the squared error of estimation of the network parameter for 
a limited sample size \( n \) 
in a highly overparametrized regime by a multiple of \( \dimA/n \).
In our results, we frequently exploit the ``blessing of dimension'' phenomenon: 
overparametrization appears to be very useful and helps to simplify the study.


\Subsection{Regression with DNN}
Given data \( (Y_{i},\Xv_{i}) \), \( i \leq n \), with response \( Y_{i} \) and regressors \( \Xv_{i} \in \R^{\dimq} \), we aim to 
build a network providing a small prediction error of the response \( \E Y_{i} - \regrf(\Xv_{i}) \).
A DNN to be constructed transforms the input regressor \( \Xv_{i} \) into \( \XXv_{i} \),
that is, \( \XXv_{i} = \DNN(\Xv_{i}) \) for all \( i \leq n \).
This transformed vector is used for linear prediction.
Most of our analysis will be done for a DNN with one hidden layer.
Given an input vector \( \Xv = (\Xv_{i}) \in \R^{\dimd} \), 
the DNN transformation and the corresponding response are defined by
\begin{EQA}
	\XXv
	&=&
	\sigma(\WWv \Xv) 
	\in \R^{\Mh},
	\qquad
	\regrf_{i} = \XXv_{i}^{\T} \betav,
\label{c83w8l88y2ygce8ec8eh3}
\end{EQA}
where \( \Mh \) is the number of nodes in the hidden layer,
\( \betav \in \R^{\Mh} \),
\( \WWv = (\weight_{\mm j}) \) is a weighting \( \Mh \times \dimd \)-matrix with rows \( \WWv_{\mm} \),
and \( \sigma \) is an entry-wise activating function, e.g.
\begin{EQA}
	\sigma(t)
	&=&
	\alp^{-1} \log(1 + \ex^{\alp t}) .
\label{gx72jdcciw9d79idhieg6g}
\end{EQA}
This is a smoothed version of the popular ReLU activation \( \sigma(t) = t_{+} \).
The structure of this neuronal network is described by the structural parameter \( \WWv \) defining for each
input feature vector \( \Xv_{i} \)  the transformed 
feature vector \( \XXv_{i} = \sigma(\WWv \Xv_{i}) \), or \( \XXv_{\mm i} = \sigma(\WWv_{\mm} \Xv_{i}) \),
\( m=1,\ldots,\Mh \). 
Here \( \Xv \) is a \( \dimd \times n \) matrix while \( \XXv \) is a \( \Mh \times n \) matrix, and \( \XXv = \sigma(\WWv \Xv) \) means the entrywise operation.
The transformed vector \( \XXv_{i} = \sigma(\WWv \Xv_{i}) \in \R^{\Mh} \) 
is viewed as a new feature vector used for linear regression: 
\( \regrf_{i} = \XXv_{i}^{\T} \betav \).
To ensure ``small noise'' condition, we perform a presmoothing operation 
\( \SmOp \colon \R^{n} \to \R^{\dimq} \) 
on data \( \Yv \) leading to the set of observables
\begin{EQA}
	\Zv
	& \eqdef &
	\SmOp \Yv \in \R^{\dimq} \, .
\label{udswjcyehdcyendcyeuu}
\end{EQA}
Later we assume that \( \SmOpT = \SmOp^{\T} \SmOp \) satisfies to \( \| \SmOpT \|_{\oper} \leq 1 \).
The DNN structure leads to 
\begin{EQA}[c]
	\E \Zv = \SmOp \sigma(\WWv \Xv)^{\T} \betav 
	\, .
\label{dcu3ehvc7hejf7fucgue}
\end{EQA}
Following the calming idea, we 
introduce an auxiliary ``response'' vector \( \zv \), 
decouple the structural relation 
\( \zv = \SmOp \XXv^{\T} \betav \), and treat \( \zv \)
as a nuisance parameter from \( \R^{\dimq} \).
The DNN structure will be enforced by the structural penalty 
\( \frac{\lambda^{*}}{2} \| \XXv - \sigma(\WWv \Xv) \|^{2} \).
Later, we apply \( \lambda^{*} = 1 \).
To ensure stability and identifiability, we also introduce a ridge penalty 
on each variable \( \WWv \), \( \XXv \), \( \betav \), and \( \zv \):
\begin{EQA}[rcllrcllrcllrcl]
	\pent(\WWv)
	&=&
	\frac{\lambda_{\wwv}}{2} \| \WWv \|^{2}
	\, ,
	\;
	&
	\pent(\XXv)
	&=&
	\frac{\lambda_{\XXv}}{2} \| \XXv \|^{2}
	\, ,
	\;
	&
	\pent(\betav)
	&=&
	\frac{\lambda_{\betav}}{2} \| \betav \|^{2}
	\, ,
	\;
	&
	\pent(\zv)
	&=&
	\frac{\lambda_{\zv}}{2} \| \zv \|^{2}
	\, .
\label{sywywdg7w7gccwwhduwscj}
\end{EQA}
The choice of the coefficients \( \lambda_{\wwv} \), \( \lambda_{\xxv} \), \( \lambda_{\betav} \), and \( \lambda_{\zv} \) 
is important and will be discussed in this section.
In particular, the choice \( \lambda_{\zv} = 3 \), \( \lambda_{\xxv} \geq 4 \), and
\( \lambda_{\wwv} = \lambda_{\betav} = 1 \) can be recommended. 
Everywhere, the collection of weights \( \WWv \) and the transformed design \( \XXv \) are considered as vectors 
from the corresponding spaces \( \R^{\dimp \times \dimq} \) and \( \R^{\dimp \times n} \), their \( \ell_{2} \)-norms
means the usual Euclidean norm in this space.  
The whole approach is based on joint maximization of the penalized log-likelihood:
for \( \prmtv = (\WWv,\XXv,\betav,\zv) \)
\begin{EQA}[rcl]
	\LL(\prmtv)
	&=&
	- \frac{1}{2} \| \Zv - \zv \|^{2}
	- \frac{1}{2} \| \zv - \SmOp \XXv^{\T} \betav \|^{2}
	- \frac{1}{2} \| \XXv - \sigma(\WWv \Xv) \|^{2}
	\\
	&&
	- \, \frac{\lambda_{\wwv}}{2} \| \WWv \|^{2}
	- \frac{\lambda_{\xxv}}{2} \| \XXv \|^{2}
	- \frac{\lambda_{\betav}}{2} \| \betav \|^{2}
	- \frac{\lambda_{\zv}}{2} \| \zv \|^{2}
	\, .
\label{dtg6wtwyc7wctdc7wh5f}
\end{EQA}
%
%
The estimator \( \tilde{\prmtv} = (\tilde{\WWv},\tilde{\XXv},\tilde{\betav},\tilde{\zv}) \) is defined by maximization 
of this log-likelihood, while the background truth \( \prmtvs = (\WWvs,\XXvs,\betavs,\zvs) \) maximizes its expectation:
\begin{EQA}[rcccl]
	\tilde{\prmtv}
	&=&
	\argmax_{\prmtv} \LL(\prmtv),
	\qquad
	\prmtvs
	&=&
	\argmax_{\prmtv} \E \LL(\prmtv).
\label{cv6hnir54433hfg8i43}
\end{EQA}
The main results describe the error of estimation \( \tilde{\prmtv} - \prmtvs \).

\Subsection{Fisher information and smoothness conditions}
An important step of the analysis concerns the Fisher information matrix
\( \IFT(\prmtv) = - \nabla^{2} \LL(\prmtv) \).
Lemma~\ref{LHessDNN} states that \( \IFT(\prmtv) > 0 \) and bounds this matrix from below by 
a simple block-diagonal matrix.

We now check local smoothness condition \nameref{LLsT3ref} for \( \LL(\prmtv) \) from \eqref{dtg6wtwyc7wctdc7wh5f} 
in the vicinity of \( \prmtvs = (\WWvs,\XXvs,\betavs,\zvs) \).
Define for \( \prmtv = (\WWv,\XXv,\betav,\zv) \) 
\begin{EQ}[ccl]
	\IF_{\mm}(\prmtv)
	&=&
	\IF_{\mm}(\WWv_{\mm})
	\eqdef 
	\sumi \sigma'\bigl(\WWv_{\mm} \Xv_{i} \bigr)^{2} \, \Xv_{i} \, \Xv_{i}^{\T} 
	+ \frac{4}{3} \lambda_{\wwv} \Id_{\dimd}
	\, ,
	\quad
	m=1,\ldots, \Mh
	\, ,
	\qquad
\label{df7duyd5tw3hdy7dgywehf}
	\\
	\IFT_{\wwv\wwv}(\prmtv)
	&=&
	\IFT_{\wwv\wwv}(\WWv)
	=
	\blk\{ \IF_{1}(\WWv_{1}), \ldots, \IF_{\Mh}(\WWv_{\Mh}) \} 
	\, .
\label{c73gth9n9retdf54df3rt5g}
\end{EQ}
Introduce a full dimensional metric tensor \( \DFM \) by
\begin{EQA}[rcl]
	\DFM^{2}
	& \eqdef &
	\blk\bigl\{ \DPM_{\wwv}^{2} \, ,
		\; \Id_{\xxv} \, ,
		\; \DPM_{\betav}^{2} \, ,
		\; \Id_{\zv}
	\bigr\} ,
\label{u8cduv9i23wioj4276dgtxuwq}
\end{EQA}
where
\begin{EQA}[rcccl]
	\DPM_{\wwv}^{2}
	& \eqdef &
	\frac{3}{4} \, \blk\{ \IF_{1}, \ldots , \IF_{\Mh} \}
	\, ,
	\qquad
	\DPM_{\betav}^{2}
	& \eqdef &
	\frac{1}{2} \XXvs \SmOpT \, {\XXvs}^{\T} + \lambda_{\betav} \Id_{\betav} \, ,
\label{dsxycfydfdf6dy35t6ert}
\end{EQA}
and \( \IF_{\mm} = \IF_{\mm}(\WWvs_{\mm}) \); see \eqref{df7duyd5tw3hdy7dgywehf}.
Lemma~\ref{LHessDNN} implies \( \DFM^{2} \leq \IFT \).

With \( \sigma(t) = \alp^{-1} \log(1 + \ex^{\alp t}) \), define the values \( \deta \) and \( \detab \) by
\begin{EQ}[rcl]
	\deta^{-1/2}
	& = &
	\max_{\mm \leq \Mh, \, i \leq n} \Bigl\{ 
		\frac{\sigma''(\WWvs_{\mm} \Xv_{i})}{\sigma'(\WWvs_{\mm} \Xv_{i})} \, \| \IF_{\mm}^{-1/2} \Xv_{i} \| 
	\Bigr\} 
	\, ,
	\\
	\detab^{-1/2}
	& = &
	\max_{\mm \leq \Mh, \, i \leq n} \| \IF_{\mm}^{-1/2} \Xv_{i} \|
	\, .
	\qquad
\label{6dft7wduwr56twrvbfjd}
\end{EQ}
Also, with \( \DPM_{\betav} \) from \eqref{dsxycfydfdf6dy35t6ert}, define \( \deta_{1} \) by
\begin{EQA}[rcl]
	\deta_{1}^{-1/2}
	& = &
	\| \DPM_{\betav}^{-1} \| \, .
\end{EQA}
Lemma~\ref{Ltau3DNN} shows that under \( \| \betav \| \leq b_{0} \, , \, \XXv \SmOpT \XXv^{\T} \leq 2 \XXvs \SmOpT {\XXvs}^{\T} \), 
\nameref{LLsT3ref} is fulfilled with 
\begin{EQA}[rcl]
	\dltwu_{3} 
	& \eqdef &
	\frac{6 \sqrt{\ex}}{\deta^{1/2}} + \frac{6 (b_{0} + 1)}{\deta_{1}^{1/2}} 
	\, .
\label{cjvcuydhyer4yer5ewdxcd}
\end{EQA}

\Subsection{Effective dimension and accuracy of DNN training}
This section discusses the full effective dimension of the model \eqref{dtg6wtwyc7wctdc7wh5f} and states the main results.
The data \( \Yv \) only enter in the fidelity term \( \| \Zv - \zv \|^{2} \) with \( \Zv = \SmOp \Yv \).
The  stochastic component \( \zeta(\prmtv) = \langle \Zv - \E \Zv, \zv \rangle \) is linear in \( \prmtv \).
For the score vector \( \nabla \zeta \), it holds \( \nabla \zeta = (0,0,0,	\scorem{\zv} \zeta) \) with
\begin{EQA}[rcl]
	\scorem{\zv} \zeta
	=
	\Zv - \E \Zv 
	\, ,
	\quad
	\Var(\scorem{\zv} \zeta)
	&=&
	\Var(\SmOp \Yv)
	=
	\SmOp \Var(\Yv) \SmOp^{\T}
	\, .
	\qquad
\label{xcydtc5445tgyjubhi988r}
\end{EQA}
Later we assume that \( \SmOp \epsv \eqdef \SmOp (\Yv - \E \Yv) \) satisfies \nameref{ED0ref} from Section~\ref{ScondNL}.
With \( \IFT = \IFT(\prmtvs) = - \nabla^{2} \LL(\prmtvs) \),
the full effective dimension \( \dimA = \dimA(\prmtvs) \) is given by
\begin{EQA}
	\dimA
	& \eqdef &
	\tr\bigl\{ \IFT^{-1} \Var(\nabla \zeta) \bigr\} 
	\, .
	\qquad
\label{fdvud56h3wiv55433egy}
\end{EQA}
Lemma~\ref{LHessDNN} yields \( \IFT^{-1} \leq \DFM^{-2} \) for \( \DFM^{2} \) from \eqref{u8cduv9i23wioj4276dgtxuwq} and
\begin{EQA}[rcl]
	\| \DFM \IFT^{-1} \nabla \zeta \|
	& \leq &
	\| \DFM^{-1} \nabla \zeta \|
	=
	\| \SmOp \epsv \| 
	\, ,
	\\
	\dimA
	& \leq &
	\tr\bigl\{ \DFM^{-2} \Var(\nabla \zeta) \bigr\}
	\leq 
	\E \| \SmOp \epsv \|^{2}
	\, .
\label{ydu2hd5trggts4eqfwxg}
\end{EQA}
Moreover, the effective radius \( \rr_{\DFM} \) from \eqref{y7s7d7d77dfdy7fuegue3j} fulfills with \( \HVB^{2} \eqdef \Var(\SmOp \epsv) \)
\begin{EQA}[c]
	\rr_{\DFM}
	\leq 
	\zq(\HVB^{2},\xx)
	\leq 
	\sqrt{\tr \HVB^{2}} + \sqrt{2 \xx \| \HVB^{2} \|}
	\, .
\end{EQA}
For the simplest case with homogeneous independent noise \( \Var(\Yv) = \sigma^{2} \Id_{n} \), a \( \dimq \)-dimensional projector \( \SmOp \),
it follows \( \dimA \leq \dimq \sigma^{2} \).

Before stating the main result, we have to resolve an identifiability issue. 
The transformed feature vector \( \XXv \) and the vector of coefficients \( \betav \) enters in the structural penalty 
\( \| \zv - \SmOp \XXv^{\T} \betav \|^{2} \) via the product \( \XXv^{\T} \betav \).
Ridge penalties \( \lambda_{\xxv} \| \XXv \|^{2} \) and \( \lambda_{\betav} \| \betav \|^{2} \) help to identify 
each component of this product after restricting to a bounded set 
\begin{EQA}[rcl]
	\Upsd
	& \eqdef &
	\bigl\{ \prmtv = (\WWv,\XXv,\betav,\zv) \colon \| \betav \| \leq b_{0} \, , \, \XXv \SmOpT \XXv^{\T} \leq 2 \XXvs \SmOpT {\XXvs}^{\T} \bigr\}
	\, .
\label{cvhv5654433ergfg7gtueded}
\end{EQA}
Here \( \SmOpT = \SmOp^{\T} \SmOp \) and \( b_{0} \) is a fixed constant, e.g. \( b_{0} = 2 \).
Our main result is the Fisher expansion for the error of estimation \( \tilde{\prmtv} - \prmtvs \approx \IFT^{-1} \nabla \zeta \)
after restricting \( \prmtv \) to \( \Upsd \).
It follows from Theorem~\ref{TQFiWi}.

\begin{theorem}
\label{TDNN}
Assume \nameref{ED0ref} and \( \| \SmOp^{\T} \SmOp \|_{\oper} \leq 1 \).
For \( \Upsd \) from \eqref{cvhv5654433ergfg7gtueded} and \( \LL(\prmtv) \) from \eqref{dtg6wtwyc7wctdc7wh5f} with
\( \lambda_{\zv} \geq 3 \) and \( \lambda_{\xxv} \geq 4 + 3 b_{0}^{2}/2 \),
let the estimator \( \tilde{\prmtv} \) and the background truth \( \prmtvs \) be defined by \eqref{cv6hnir54433hfg8i43} 
after limiting \( \prmtv \) to \( \Upsd \).
Let \( \rr_{\DFM} \, \dltwu_{3} \leq 4/9 \) for \( \dltwu_{3} \) from \eqref{cjvcuydhyer4yer5ewdxcd}.
Then with \( \DFM \) from \eqref{u8cduv9i23wioj4276dgtxuwq}, 
for any \( \QP \), it holds on a random set \( \Omega(\xx) \) with \( \P(\Omega(\xx)) \geq 1 - \ex^{-\xx} \)
\begin{EQA}[rcl]
	\| \QP (\tilde{\prmtv} - \prmtvs - \IFT^{-1} \nabla \zeta) \|
	& \leq &
	\| \QP \DFM^{-1} \| \frac{3 \dltwu_{3}}{4} \| \SmOp \epsv \|^{2} \, .
\label{duhs8iu3wfc0wlvfuvc}
\end{EQA}
Also, introduce 
\begin{EQA}[c]
	\dimQ \eqdef \E \{ \| \QP \IFT^{-1} \nabla \zeta \|^{2} \Ind_{\Omega(\xx)} \} \, ,
\label{7djhed8cjfct534etgdhdyQPD}
\end{EQA}
assume \( \E \bigl\{ \| \SmOp \epsv \|^{4} \Ind_{\Omega(\xx)} \bigr\} \leq \CONSTi_{4}^{2} \, \dimA^{2} \), 
and define
\begin{EQ}[rcl]
	\alp_{\QP}
	& \eqdef & 
	\frac{\| \QP \DFM^{-1} \| \, (3/4) \dltwu_{3} \, \CONSTi_{4} \, \dimA } {\sqrt{\dimQ}}
	\, . 
\label{6dhx6whcuydsds655srewD}
\end{EQ}
If \( \alp_{\QP} < 1 \) then
\begin{EQA}
	(1 - \alp_{\QP})^{2} \dimQ 
	\leq 
	\E \bigl\{ \| \QP \, (\tilde{\prmtv} - \prmtvs) \|^{2} \Ind_{\Omega(\xx)} \bigr\}
	& \leq &
	(1 + \alp_{\QP})^{2} \dimQ \, .
\label{EQtuGmstrVEQtGQPD}
\end{EQA}
\end{theorem}

This result can be used for each component of \( \prmtvs \).
Denote \( \IFTinv = \IFT^{-1} \).
Then
\begin{EQA}
	\IFT^{-1} \nabla \zeta 
	&=&
	(\IFTinv_{\wwv\zv} \SmOp \epsv, \, \IFTinv_{\xxv\zv} \SmOp \epsv, \, \IFTinv_{\betav\zv} \SmOp \epsv, \, \IFTinv_{\zv\zv} \SmOp \epsv) 
	\, .
\label{bvu7654rfgy6y6rdfmjoie}
\end{EQA}
and, in particular, with \( \DPM_{\wwv} \) and \( \DPM_{\betav} \) from \eqref{dsxycfydfdf6dy35t6ert}, it holds on \( \Omega(\xx) \)
\begin{EQA}
\label{vdud6ehry3hrf67ruytaD}
	\| \DPM_{\wwv} ( \tilde{\WWv} - \WWvs - \IFTinv_{\wwv\zv} \SmOp \epsv ) \|
	& \leq &
	(3/4) \dltwu_{3} \, \| \SmOp \epsv \|^{2} 
	\, ,
	\\
	\| \tilde{\XXv} - \XXvs - \IFTinv_{\xxv\zv} \SmOp \epsv ) \|
	& \leq &
	(3/4) \dltwu_{3} \, \| \SmOp \epsv \|^{2} 
	\, ,
	\\
	\| \DPM_{\betav} ( \tilde{\betav} - \betavs - \IFTinv_{\betav\zv} \SmOp \epsv ) \|
	& \leq &
	(3/4) \dltwu_{3} \, \| \SmOp \epsv \|^{2} 
	\, ,
	\\
	\| \tilde{\zv} - \zvs - \IFTinv_{\zv\zv} \SmOp \epsv ) \|
	& \leq &
	(3/4) \dltwu_{3} \, \| \SmOp \epsv \|^{2} 
	\, .
\end{EQA}
Under a proper choice of the presmoothing operator \( \SmOp \), the leading term \( \IFT^{-1} \nabla \zeta \) of the expansion \eqref{duhs8iu3wfc0wlvfuvc} 
is nearly normal as well as all its components.
In particular, if each matrix \( \IF_{\mm} = \IF_{\mm}(\WWvs) \) from \eqref{df7duyd5tw3hdy7dgywehf} satisfies \( \lambda_{\min}(\IF_{\mm}) \geq \CONST n \) then
the estimator \( \tilde{\WWv} \) of the DNN architecture is root-n normal.
The result \eqref{duhs8iu3wfc0wlvfuvc} can be effectively used for inference purposes such as model check, model selection, pruning, etc.

\Section{DNN tools}
This section collects some technical statements for the DNN study.

\Subsection{Fisher information matrix}
Denote \( \prmtv = (\WWv,\XXv,\betav,\zv) \).
This section discusses some geometric properties of the function \( \LL(\prmtv) \).
In particular, we show that it is strongly concave in a vicinity of any point 
\( \prmtv = (\WWv,\XXv,\betav,\zv) \) with \( \XXv = \sigma(\WWv \Xv) \) and \( \zv = \SmOp \XXv^{\T} \betav \).
An important step of the study is to show that the Fisher information matrix 
\( \IFT(\prmtv) = - \nabla^{2} \LL(\prmtv) \) is well-posed:
\begin{EQA}[rcl]
	\IFT(\prmtv) 
	& \eqdef & 
	- \nabla^{2} \LL(\prmtv)
	>
	0
\label{dsufgwu876hy2wqjd}
\end{EQA}
for all considered \( \prmtv \).
%
Difficulties for the analysis come from the nonlinear transformation of the product \( \WWv \Xv \in \R^{\Mh \times n} \)
in the structural penalty \( \| \XXv - \sigma(\WWv \Xv) \|^{2}/2 \) yielding a complex non-convex behavior of this term.
However, overparametrization enables us to limit the parameter set  
\( (\WWv,\XXv,\betav,\zv) \) with \( \XXv = \sigma(\WWv \Xv) \).
At any such point, the second derivative of the nonlinear function \( \sigma(\WWv \Xv) \) enters 
in the expansion for the Hessian of the structural penalty with the factor \( \XXv - \sigma(\WWv \Xv) \) and hence, vanishes.
We will use this fact repeatedly in our derivations. 
The same applies to the product structure \( \SmOp \XXv^{\T} \betav \) in \( \| \zv - \SmOp \XXv^{\T} \betav \|^{2} \).
Under the structural equation \( \zv = \SmOp \XXv^{\T} \betav \), the second derivative term in the Hessian can be ignored.

We start by describing a directional second-order derivative of \( \LL(\prmtv) \).

\begin{lemma}
\label{LpartderDNN}
Fix \( \prmtv = (\WWv,\XXv,\betav,\zv) \) with \( \XXv = \sigma(\WWv \Xv) \) and \( \zv = \SmOp \XXv^{\T} \betav \). 
Then for any \( \uv = (\wwv,\xxv,\betav,\zv) \) with 
\( \wwv = (\wwv_{\mm}) \in \R^{\Mh \times \dimd} \), \( \xxv = (\xxv_{i}) \in \R^{\Mh \times n} \), 
\( \bv \in \R^{\dimp} \), and \( \hv \in \R^{\dimq} \)
\begin{EQA}
	&& \nquad
	- \frac{d^{2}}{dt^{2}} \LL(\prmtv + t \uv) \bigg|_{t=0}
	=
	(1 + \lambda_{\hv}) \| \hv \|^{2} + \lambda_{\xxv} \| \xxv \|^{2} + \lambda_{\betav} \| \bv \|^{2} + \lambda_{\wwv} \| \wwv \|^{2}
	\\
	&&  
	+ \sumi \sum_{\mm=1}^{\Mh}  
	\bigl\{ \sigma'(\WWv_{\mm} \Xv_{i}) \, \wwv_{\mm} \Xv_{i} - \xxv_{\mm i} \bigr\}^{2}	
	+ \| \hv - \SmOp (\XXv^{\T} \bv + \xxv^{\T} \betav) \|^{2}
	\, .
	\qquad
\label{d7vyu7e3hirjo3gbi8r}
\end{EQA} 
\end{lemma}

\begin{proof}
The structural equation \( \XXv = \sigma(\WWv \Xv) \) enables us to ignore the term in the decomposition of 
the second derivative of the structural penalty \( \| \XXv - \sigma(\WWv \Xv) \|^{2}/2 \) that contains the second derivative of the nonlinear function 
\( \sigma(\WWv \Xv) \).
The same applies to the penalty \( \| \zv - \SmOp \XXv^{\T} \betav \|^{2} \).
This yields the result.
\end{proof}

The next step is to bound from below the matrix \( \IFT(\prmtv) \) by a block-diagonal matrix.

\begin{lemma}
\label{LHessDNN}
Assume \( \lambda_{\zv} \geq 3 \).
Given \( b_{0} > 0 \), let \( \lambda_{\xxv} \geq 4 + \frac{3 b_{0}^{2}}{2} \).
For \( \prmtv = (\WWv,\XXv,\betav,\zv) \) with \( \XXv = \sigma(\WWv \Xv) \), \( \zv = \SmOp \XXv^{\T} \betav \), 
and \( \| \betav \| \leq b_{0} \)
\begin{EQA}[rcl]
	\IFT(\prmtv)
	& \geq &
	\blk\Bigl\{ 
		\frac{3}{4} \, \IFT_{\wwv\wwv} \, , \; 
		\Bigl( \lambda_{\xxv} - 3 - \frac{3 b_{0}^{2}}{2} \Bigr) \Id_{\xxv} \, , \;
		\frac{1}{2} \, \XXv \SmOpT \XXv^{\T} + \lambda_{\betav} \Id_{\betav} \, , \;
		(\lambda_{\zv} - 2) \Id_{\zv}
	\Bigr\}
	\, .
\end{EQA}
\end{lemma}

\begin{proof}
For \( \uv = (\wwv,\xxv,\betav,\zv) \), we apply \eqref{d7vyu7e3hirjo3gbi8r} and note that with \( \rho = 3 \)
\begin{EQA}[rcl]
	&& \nquad
	\rho \| \xxv \|^{2} + \lambda_{\wwv} \Id_{\wwv} + \sumi \sum_{\mm=1}^{\Mh}  
	\bigl\{ \xxv_{\mm i} - \sigma'(\WWv_{\mm} \Xv_{i}) \, \wwv_{\mm} \Xv_{i} \bigr\}^{2}
	\\
	& = &
	\sumi \sum_{\mm=1}^{\Mh}  
	\biggl\{ 
		\sqrt{1 + \rho} \,\, \xxv_{\mm i} 
		- \frac{\sigma'(\WWv_{\mm} \Xv_{i})}{\sqrt{1 + \rho}} \,\, \wwv_{\mm} \Xv_{i} 
	\biggr\}^{2}
	+ \frac{\rho}{1+\rho} \sum_{\mm=1}^{\Mh} \wwv_{\mm} \IF_{\mm}(\WWv_{\mm}) \wwv_{\mm}^{\T} 
	\\
	& \geq &
	\frac{3}{4} \sum_{\mm=1}^{\Mh} \wwv_{\mm} \IF_{\mm}(\WWv_{\mm}) \wwv_{\mm}^{\T} 
	=
	\frac{3}{4} \bigl\langle \IFT_{\wwv\wwv}(\WWv), \wwv^{\otimes 2} \bigr\rangle
	\, .
\end{EQA}
Similarly,
the inequality \( (a + b + c)^{2} + 3a^{2} + 3b^{2}/2 \geq c^{2}/2 \) for any \( a,b,c \) implies
\begin{EQA}[c]
	\| \hv - \SmOp \xxv^{\T} \betav - \SmOp \XXv^{\T} \bv \|^{2} + 3 \| \hv \|^{2} + \frac{3}{2} \| \SmOp \xxv^{\T} \betav \|^{2}
	\geq 
	\frac{1}{2} \, \bigl\| \SmOp \XXv^{\T} \bv \bigr\|^{2} 
	\, .
\end{EQA}
By \( \SmOpT = \SmOp^{\T} \SmOp \) and \( \| \betav \| \leq b_{0} \) 
\begin{EQA}[c]
	\frac{1}{2} \, \bigl\| \SmOp \XXv^{\T} \bv \bigr\|^{2} 
	=
	\frac{1}{2} \, \bv^{\T} \XXv \SmOpT \XXv^{\T} \bv \, .
	\qquad
	\| \SmOp \xxv^{\T} \betav \|^{2} 
	=
	\betav^{\T} \xxv \SmOpT \xxv^{\T} \betav
	\leq 
	b_{0}^{2} \, \| \xxv \SmOpT \xxv^{\T} \|_{\oper}
\end{EQA}
and
\begin{EQA}[rcl]
	\frac{3}{2} \, b_{0}^{2} \| \xxv \SmOpT \xxv^{\T} \|_{\oper} + 3 \| \hv \|^{2} +  \| \hv - \SmOp \XXv^{\T} \bv - \SmOp \xxv^{\T} \betav \|^{2}
	& \geq &
	\frac{1}{2} \, \bv^{\T} \XXv \SmOpT \XXv^{\T} \bv
	\, .
\end{EQA}
Now by \eqref{d7vyu7e3hirjo3gbi8r}
\begin{EQA}[rcl]
	&& \nquad
	\bigl\langle \IFT(\prmtv), \uv^{\otimes 2} \bigr\rangle
	\geq 
	\frac{3}{4} \bigl\langle \IFT_{\wwv\wwv}(\WWv), \wwv^{\otimes 2} \bigr\rangle
	+ ( \lambda_{\xxv} - 3 ) \| \xxv \|^{2} 
	- \frac{3}{2} \, b_{0}^{2} \, \| \xxv \SmOpT \xxv^{\T} \|_{\oper}
	\\
	&&
	+ \, \lambda_{\betav} \| \bv \|^{2} 
	+ \frac{1}{2} \, \bv^{\T} \XXv \SmOpT \XXv^{\T} \bv
	+ (\lambda_{\zv} - 2) \| \hv \|^{2} 
	\, .
\end{EQA}
This yields the assertion by \( \| \SmOpT \|_{\oper} = \| \SmOp^{\T} \SmOp \|_{\oper} \leq 1 \).
\end{proof}

\Subsection{Smoothness conditions}
The first step is in evaluating the variability of each matrix 
\( \IF_{\mm}(\WWv_{\mm}) \) in an elliptic vicinity around \( \WWvs_{\mm} \).

\begin{lemma}
\label{LFmDNN}
Let \( \rr \) satisfy \( 2 (\ex \vee \alp) \, \rr < \deta^{1/2} \); see \eqref{6dft7wduwr56twrvbfjd}.
Fix any \( (\WWv,\XXv) \) satisfying \( \XXv = \sigma(\WWv \Xv) \) and
\begin{EQA}[rcl]
	\| (\WWv_{\mm} - \WWvs_{\mm}) \IF_{\mm}^{1/2} \|
	& \leq &
	\rr \, ,
	\qquad
	\mm =1,\ldots,\Mh
	\, .
\label{xcyjervb7gwe3dfc5fctr5}
\end{EQA}
Then it holds for all \( \mm \leq \Mh \) 
\begin{EQA}[rcl]
	(1 - \delta) \IF_{\mm} 
	& \leq &
	\IF_{\mm}(\WWv_{\mm})
	\leq 
	(1 + \delta) \IF_{\mm} \, ,
	\qquad
	\delta \eqdef 2 \ex \, \rr \, \deta^{-1/2} .
\label{vduf786grey8hrjvy6weyt6h}
\end{EQA}
\end{lemma}

\begin{proof}
For \( \sigma(t) = \alp^{-1} \log(1 + \ex^{\alp t}) \), we can use
\begin{EQA}[rcl]
	\sigma'(t)
	&=&
	\frac{\ex^{\alp t}}{1 + \ex^{\alp t}} \, ,
	\qquad
	\sigma''(t)
	=
	\frac{\alp \ex^{\alp t}}{(1 + \ex^{\alp t})^{2}}
	\, ,
\label{dgceeo3f6d7dh3ecywyewh}
\end{EQA}
yielding \( 0 < \sigma'(t) < 1 \) and \( 0 < \sigma''(t)/\sigma'(t) \leq \alp \).
Moreover, for any \( u \) with \( \alp |u| \leq 1/2 \)
\begin{EQA}[rcl]
	\frac{\sigma'(t+u)}{\sigma'(t)} 
	& \leq &
	\sqrt{\ex} \, ,
	\qquad
	\frac{\sigma'(t+u) \, \sigma''(t+u)}{\sigma'(t) \, \sigma''(t)} 
	\leq 
	\ex \, .
\label{cdud5654544duid723g7}
\end{EQA}
Fix \( \mm \) and \( \WWv_{\mm} \) with \( \| (\WWv_{\mm} - \WWvs_{\mm,\GPW}) \IF_{\mm}^{1/2} \| \leq \rr \).
By definition, for any \( \bv \in \R^{\dimd} \),
\begin{EQA}[rcl]
	\bv^{\T} \bigl\{ \IF_{\mm}(\WWv_{\mm}) - \IF_{\mm} \bigr\} \bv
	&=&
	\sumi \bigl\{ \sigma'(\WWv_{\mm} \Xv_{i})^{2} - \sigma'(\WWvs_{\mm} \Xv_{i})^{2} \bigr\} \, (\bv^{\T} \Xv_{i})^{2} \, .
\end{EQA}
Further, with \( \wwv_{\mm} = \WWv_{\mm} - \WWvs_{\mm} \), it holds for some \( t \in [0,1] \)
\begin{EQA}[rcl]
	\sigma'(\WWv_{\mm} \Xv_{i})^{2} - \sigma'(\WWvs_{\mm} \Xv_{i})^{2}
	&=&
	2 \sigma'(\WWvs_{\mm} \Xv_{i} + t \wwv_{\mm} \Xv_{i}) \, \sigma''(\WWvs_{\mm} \Xv_{i} + t \wwv_{\mm} \Xv_{i}) \, \wwv_{\mm} \Xv_{i} \, .
\end{EQA}
Now we show that \( \alp |\wwv_{\mm} \Xv_{i}| \leq 1/2 \) for all \( i \leq n \).
Indeed, by \eqref{6dft7wduwr56twrvbfjd}
\begin{EQA}[c]
	\alp |\wwv_{\mm} \Xv_{i}| 
	\leq 
	\alp \| \wwv_{\mm} \IF_{\mm}^{1/2} \| \,\, \| \IF_{\mm}^{-1/2} \Xv_{i} \|
	\leq 
	\alp \, \rr \, \max_{i \leq n} \| \IF_{\mm}^{-1/2} \Xv_{i} \|
	\leq 
	\alp \, \rr \, \detab^{-1/2}
	\leq 
	1/2 .
\end{EQA}
This and \eqref{cdud5654544duid723g7} imply
\begin{EQA}[rcl]
	\bigl| \sigma'(\WWv_{\mm} \Xv_{i})^{2} - \sigma'(\WWvs_{\mm} \Xv_{i})^{2} \bigr|
	& \leq &
	2 \ex \, \sigma'(\WWvs_{\mm} \Xv_{i}) \, \sigma''(\WWvs_{\mm} \Xv_{i}) \, \bigl| \wwv_{\mm} \Xv_{i} \bigr|
	\\
	& \leq &
	2 \ex \, \max_{i \leq n} \Bigl\{ 
		\frac{\sigma''(\WWvs_{\mm} \Xv_{i})}{\sigma'(\WWvs_{\mm} \Xv_{i})} \, |\wwv_{\mm} \Xv_{i}| 
	\Bigr\} \,\, \sigma'(\WWvs_{\mm} \Xv_{i})^{2}
\end{EQA}
and again by \eqref{6dft7wduwr56twrvbfjd}
\begin{EQA}[rcl]
	&& \nquad
	\bv^{\T} \bigl\{ \IF_{\mm}(\WWv_{\mm}) - \IF_{\mm} \bigr\} \bv
	\leq 
	2 \ex \sumi \sigma'(\WWvs_{\mm} \Xv_{i}) \, \sigma''(\WWvs_{\mm} \Xv_{i}) \, (\bv^{\T} \Xv_{i})^{2}
	\\
	& \leq &
	2 \ex \, \max_{i \leq n} \Bigl\{ 
		\frac{\sigma''(\WWvs_{\mm} \Xv_{i})}{\sigma'(\WWvs_{\mm} \Xv_{i})} \, |\wwv_{\mm} \Xv_{i}| 
	\Bigr\} \,\, \bv^{\T} \IF_{\mm} \bv
	\\
	& \leq &
	2 \ex \, \max_{i \leq n} \Bigl\{ 
		\frac{\sigma''(\WWvs_{\mm} \Xv_{i})}{\sigma'(\WWvs_{\mm} \Xv_{i})} \, \| \IF_{\mm}^{-1/2} \Xv_{i} \| 
	\Bigr\} \, \bigl\| \wwv_{\mm} \IF_{\mm}^{1/2} \bigr\| \,\, \bv^{\T} \IF_{\mm} \bv
	\leq 
	\frac{2 \ex \, \rr}{\deta^{1/2}} \,\, \bv^{\T} \IF_{\mm} \bv
\end{EQA}
as claimed.
\end{proof}

Fix \( \prmtv = (\WWv, \XXv, \betav, \zv) \) from a local vicinity of \( \prmtvs \). 
The third derivative of the fidelity \( \| \Zv - \zv \|^{2} \) and of the quadratic penalties 
in \( \LL(\prmtv) \) vanishes.
Therefore, we need to bound the third derivatives of the structural penalties 
\( \| \XXv - \sigma(\WWv \Xv) \|^{2}/2 \) and \( \| \zv - \SmOp \XXv^{\T} \betav \|^{2} \).

\begin{lemma}
\label{LtauDNN}
Let \( 2 (\ex \vee \alp) \, \rr < \deta^{1/2} \).
Fix any \( (\WWv,\XXv) \) satisfying \( \XXv = \sigma(\WWv \Xv) \) and \eqref{xcyjervb7gwe3dfc5fctr5}.
Also, fix any \( \wwv \in \R^{\dimp \times \dimd} \) and \( \xxv \in \R^{\dimp \times n} \) and define
\begin{EQA}[rcl]
	h(t)
	& \eqdef &
	\frac{1}{2} \, \bigl\| \XXv + t \xxv - \sigma\bigl( (\WWv + t \wwv) \Xv \bigr) \bigr\|^{2} \, ,
	\qquad
	t \in \R
	\, .
\end{EQA}
Then
\begin{EQA}[c]
	\bigl| h'''(0) \bigr|
	\leq 
	\frac{6 \sqrt{\ex}}{\sqrt{\deta}}  
	\bigl( \| \xxv \|^{2} + \| \DPM_{\wwv} \wwv \|^{2} \bigr)^{3/2} 
	\, .	
	\qquad
\label{dvuyguwc6hywjverg4hbl}
\end{EQA}
\end{lemma}

\begin{proof}
Fix \( (\WWv,\XXv,\betav,\zv) \) satisfying \( \XXv = \sigma(\WWv \Xv) \) and \eqref{xcyjervb7gwe3dfc5fctr5}.
It holds  
\begin{EQA}
	&& \nquad
	- \frac{d^{3}}{dt^{3}} \sumi \frac{1}{2} \| \XXv_{i} + t \xxv_{i} - \sigma(\WWv \Xv_{i} + t \wwv \Xv_{i}) \|^{2} \bigg|_{t=0}
	\\
	&=&
	- \, 3 \sum_{\mm=1}^{\Mh} \sumi \bigl\{ \sigma'(\WWv_{\mm} \Xv_{i}) \, \wwv_{\mm} \Xv_{i} - \xxv_{\mm i} \bigr\} \, 
		\sigma''(\WWv_{\mm} \Xv_{i}) \, (\wwv_{\mm} \Xv_{i})^{2}
	\, .
\label{d7vyu7e3hirjo3gbi8r3}
\end{EQA}
Further, \eqref{df7duyd5tw3hdy7dgywehf} and \eqref{cdud5654544duid723g7} imply for every \( \mm \leq \Mh \)
\begin{EQA}
	&& \nquad
	\biggl| \sumi \sigma'(\WWv_{\mm} \Xv_{i}) \,  
	\sigma''(\WWv_{\mm} \Xv_{i}) \, (\wwv_{\mm} \Xv_{i})^{3} \biggr|
	\leq 
	2 \ex \, \biggl| \sumi \sigma'(\WWvs_{\mm} \Xv_{i}) \,  
	\sigma''(\WWvs_{\mm} \Xv_{i}) \, (\wwv_{\mm} \Xv_{i})^{3} \biggr|
	\\
	& \leq &
	2 \ex \, \max_{i \leq n} 
	\Bigl\{ 
		\frac{\sigma''(\WWvs_{\mm} \Xv_{i})}{\sigma'(\WWvs_{\mm} \Xv_{i})} \, |\wwv_{\mm} \Xv_{i}| 
	\Bigr\} \,\,
	\sumi \sigma'(\WWvs_{\mm} \Xv_{i})^{2} \, (\wwv_{\mm} \Xv_{i})^{2}
	\\
	& \leq &
	2 \ex \, \max_{i \leq n} 
	\Bigl\{ 
		\frac{\sigma''(\WWvs_{\mm} \Xv_{i})}{\sigma'(\WWvs_{\mm} \Xv_{i})} \, |\wwv_{\mm} \Xv_{i}| 
	\Bigr\} \, \wwv_{\mm} \IF_{\mm} \wwv_{\mm}^{\T}
	\\
	& \leq &
	2 \ex \, \deta^{-1/2} \bigl\| \wwv_{\mm} \IF_{\mm}^{1/2} \bigr\| \,\,
	\wwv_{\mm} \IF_{\mm} \wwv_{\mm}^{\T}
	\, .
\label{tusd253ugucgeyc5e4wext2}
\end{EQA}
Similarly, by Cauchy-Schwarz inequality and \eqref{cdud5654544duid723g7}
\begin{EQA}
	&& \nquad
	\Biggl( 
		\sumi 
		\xxv_{\mm i} \, \sigma''(\WWv_{\mm} \Xv_{i}) \, (\wwv_{\mm} \Xv_{i})^{2} 
	\Biggr)^{2}
	\leq 
	\| \xxv_{\mm} \|^{2} \,\,
	\sumi 
	\sigma''(\WWv_{\mm} \Xv_{i})^{2} \, (\wwv_{\mm} \Xv_{i})^{4}
	\\
	& \leq &
	\| \xxv_{\mm} \|^{2} \,\,
	2 \ex \, \sumi 
	\sigma''(\WWvs_{\mm} \Xv_{i})^{2} \, (\wwv_{\mm} \Xv_{i})^{4}
	\\
	& \leq &
	2 \ex \, \| \xxv_{\mm} \|^{2} 
	\max_{i \leq n} \Bigl\{ 
		\frac{\sigma''(\WWvs_{\mm} \Xv_{i})}{\sigma'(\WWvs_{\mm} \Xv_{i})} \, |\wwv_{\mm} \Xv_{i}| 
	\Bigr\}^{2}
	\sumi 
	\sigma'(\WWvs_{\mm} \Xv_{i})^{2} \, (\wwv_{\mm} \Xv_{i})^{2}
	\\
	& \leq &
	2 \ex \, \| \xxv_{\mm} \|^{2} 
	\deta^{-1} 
	\bigl\| \wwv_{\mm} \IF_{\mm}^{1/2} \bigr\|^{2} \,\, \wwv_{\mm} \IF_{\mm} \wwv_{\mm}^{\T} 
	=
	2 \ex \, \deta^{-1} \| \xxv_{\mm} \|^{2} \bigl( \wwv_{\mm} \IF_{\mm} \wwv_{\mm}^{\T} \bigr)^{2} \, .
\label{8sd9w9rwoob77e3ijfoe}
\end{EQA}
Putting together all the bounds yields
\begin{EQA}[rcl]
	&& \nquad
	\frac{d^{3}}{dt^{3}} \, \frac{1}{2} \, \bigl\| \XXv + t \xxv - \sigma(\WWv \Xv + t \wwv \Xv) \bigr\|^{2} \bigg|_{t=0}
	\leq 
	\frac{3 \sqrt{2\ex}}{\sqrt{\deta}}  
	\sum_{\mm=1}^{\Mh} \bigl( \| \xxv_{\mm} \| + \| \wwv_{\mm} \IF_{\mm}^{1/2} \| \bigr) \,\, \wwv_{\mm} \IF_{\mm} \wwv_{\mm}^{\T}  
	\\
	& \leq &
	\frac{6 \sqrt{\ex}}{\sqrt{\deta}}  
	\sum_{\mm=1}^{\Mh} \bigl( \| \xxv_{\mm} \|^{2} + \wwv_{\mm} \IF_{\mm} \wwv_{\mm}^{\T} \bigr)^{1/2} \,\, \wwv_{\mm} \IF_{\mm} \wwv_{\mm}^{\T}
	\\
	& \leq &
	\frac{6 \sqrt{\ex}}{\sqrt{\deta}}  
	\biggl( \sum_{\mm=1}^{\Mh} \| \xxv_{\mm} \|^{2} + \wwv_{\mm} \IF_{\mm} \wwv_{\mm}^{\T} 
	\biggr)^{3/2} 
	\, ,
\end{EQA}
and \eqref{dvuyguwc6hywjverg4hbl} follows.
\end{proof}

Now we consider the second structural penalty \( \| \zv - \SmOp \XXv^{\T} \betav \|^{2} \).

\begin{lemma}
\label{LthdSmOp}
Fix \( \prmtv = (\WWv,\XXv,\betav,\zv) \) with \( \| \betav \| \leq b_{0} \) and \( \XXv \SmOpT \XXv^{\T} \leq 2 \XXvs \SmOpT {\XXvs}^{\T} \),
and any directions \( \xxv \in \R^{\dimp \times n} \), \( \bv \in \R^{\dimp} \), \( \hv \in \R^{\dimq} \).
Then the function 
\begin{EQA}[c]
	h(t)
	\eqdef
	\frac{1}{2} \bigl\| \SmOp (\XXv + t \, \xxv)^{\T} (\betav + t \, \bv) - (\zv + t \, \hv) \bigr\|^{2}
\end{EQA}
satisfies
\begin{EQA}[c]
	\bigl| h'''(0) \bigr|
	\leq 
	\Bigl\{ 6 (b_{0} + 1) \, \| \xxv \SmOpT \xxv^{\T} \|_{\oper} + 6 \| \DPM_{\betav} \bv \|^{2} + 3 \| \hv \|^{2}
	\Bigr\} \| \bv \|
	\, .
\label{ic98ddhuer5rtg7h78en}
\end{EQA}
\end{lemma}

\begin{proof}
For any smooth vector function \( g(t) \in \R^{\dimq} \), it holds
\begin{EQA}[lcl]
	\bigl( \| g(t) \|^{2} \bigr)'
	&=&
	2 \bigl\langle g(t),g'(t) \bigr\rangle \, ,
	\\
	\bigl( \| g(t) \|^{2} \bigr)''
	&=&
	2 \bigl\langle g'(t),g'(t) \bigr\rangle + 2 \bigl\langle g(t),g''(t) \bigr\rangle \, , 
	\\
	\bigl( \| g(t) \|^{2} \bigr)'''
	&=&
	6 \bigl\langle g'(t),g''(t) \bigr\rangle + 2 \bigl\langle g(t),g'''(t) \bigr\rangle \
	\, .
\end{EQA}
For \( g(t) = \SmOp (\XXv + t \, \xxv)^{\T} (\betav + t \, \bv) - (\zv + t \, \hv) \), this implies by \( g'''(t) \equiv 0 \)
\begin{EQA}[c]
	\frac{1}{2} \bigl( \| g(0) \|^{2} \bigr)'''
	=
	3 \bigl\langle g'(0),g''(0) \bigr\rangle
	=
	6 \bigl\langle \SmOp \xxv^{\T} \betav + \SmOp \XXv^{\T} \bv - \hv, \, \SmOp \xxv^{\T} \bv \bigr\rangle .
\end{EQA}
As \( \| \betav \| \leq b_{0} \) and \( \XXv \SmOpT \XXv^{\T} \leq 2 \XXvs \SmOpT {\XXvs}^{\T} \)
\begin{EQA}[rcccl]
	\bigl| \bigl\langle \SmOp \xxv^{\T} \betav , \, \SmOp \xxv^{\T} \bv \bigr\rangle \bigr|
	&=& 
	\bigl| \betav^{\T} \xxv \SmOpT \xxv^{\T} \bv \bigr|
	& \leq &
	b_{0} \, \| \xxv \SmOpT \xxv^{\T} \|_{\oper} \, \| \bv \| 
	\, ,
	\\
	\bigl| \bigl\langle \SmOp \XXv^{\T} \bv , \, \SmOp \xxv^{\T} \bv \bigr\rangle \bigr|
	& \leq &
	\| \SmOp \XXv^{\T} \bv \| \, \| \xxv \SmOpT \xxv^{\T} \|_{\oper}^{1/2} \,\, \| \bv \|
	& \leq & 
	\frac{1}{2} \bigl( \bv^{\T} \XXv \SmOpT \XXv^{\T} \bv + \| \xxv \SmOpT \xxv^{\T} \|_{\oper} \bigr) \| \bv \|
	\, ,
	\\
	\bigl| \bigl\langle \hv , \, \SmOp \xxv^{\T} \bv \bigr\rangle \bigr|
	& \leq & 
	\| \hv \| \, \| \xxv \SmOpT \xxv^{\T} \|_{\oper}^{1/2} \, \| \bv \| 
	& \leq & 
	\frac{1}{2} \bigl( \| \xxv \SmOpT \xxv^{\T} \|_{\oper} + \| \hv \|^{2} \bigr) \, \| \bv \|
	\, ,
\end{EQA}
and \eqref{ic98ddhuer5rtg7h78en} follows.
\end{proof}

Lemma~\ref{LtauDNN} and Lemma~\ref{LthdSmOp} imply the following bound.

\begin{lemma}
\label{Ltau3DNN}
For \( \prmtv = (\WWv,\XXv,\betav,\zv) \), assume \eqref{xcyjervb7gwe3dfc5fctr5}, 
\( 2 (\ex \vee \alp) \, \rr < \deta^{1/2} \),
\( \| \betav \| \leq b_{0} \),
\( \XXv \SmOpT \XXv^{\T} \leq 2 \XXvs \SmOpT {\XXvs}^{\T} \).
For any \( \uv = (\wwv,\xxv,\bv,\hv) \) with 
\( \xxv \in \R^{\Mh \times n} \), \( \wwv \in \R^{\Mh \times \dimd} \),
\( \bv \in \R^{\dimp} \), and \( \hv \in \R^{\dimq} \),
\begin{EQA}[rcl]
	\| \DFM \uv \|^{2}
	&=&
	\| \DPM_{\wwv} \wwv \|^{2} + \| \xxv \|^{2} + \| \DPM_{\betav} \, \bv \|^{2} + \| \hv \|^{2} 
	\, .
\label{fyeijowf6677eufje4}
\end{EQA}
Then \( h(t) = \LL(\prmtv + t \uv) \) fulfills
\begin{EQA}[c]
	|h'''(0)|
	\leq 
	\Bigl( \frac{6 \sqrt{\ex}}{\deta^{1/2}} + \frac{6 (b_{0} + 1)}{\deta_{1}^{1/2}} \Bigr) \, \| \DFM \uv \|^{3} 
	\, .
\end{EQA}
\end{lemma}


\Section{A multilayer DNN}

Construction \eqref{dtg6wtwyc7wctdc7wh5f} can be extended to a \( K \)-layer network using recurrence 
\begin{EQA}
	\XXv^{(k)}
	&=&
	\sigma(\WWv^{(k)} \XXv^{(k-1)})
\label{7huuct44gwuchee6bccnwR}
\end{EQA}
for \( k = 1,\ldots,K \) and \( \XXv^{(0)} = \Xv \).
Let \( \Mh_{k} \) be the number of neurons at \( k \)th layer, 
so that \( \XXv^{(k)} \in \R^{\Mh_{k} \times n} \).
Also set \( \Mh_{0} = \dimd \).
Usually, one also introduces a shift parameter vector \( \bv^{(k)} \) in \eqref{7huuct44gwuchee6bccnwR}: 
\( \XXv^{(k)} = \sigma(\bv^{(k)} + \WWv^{(k)} \XXv^{(k-1)}) \).
This, however, can be incorporated in \eqref{7huuct44gwuchee6bccnwR} by extending each 
vector \( \XXv_{i}^{(k)} \) with \( \XXv_{i,0}^{(k)} = 1 \).
%
The full parameter set \( \prmtv \) includes for each layer \( k \) a \( \Mh_{k} \times \Mh_{k-1} \)-matrix \( \WWv^{(k)} \) of weights, 
a \( \Mh_{k} \times n \) matrix \( \XXv^{(k)} \) of transformed feature vectors \( \XXv_{i}^{(k)} \).
In total, we have \( n_{\wwv} = \sum_{k=1}^{K} \Mh_{k} \, \Mh_{k-1} \) weights \( \WWv_{\mm\mmc}^{(k)} \)
and \( \dimX = n \sum_{k=1}^{K} \Mh_{k} \) values \( \XXv_{i\mm}^{(k)} \).
The vector of coefficients \( \betav \) and the response vector \( \zv \) also enter in \( \prmtv \).
With \( \prmtv = (\WWv^{(1)}, \XXv^{(1)}, \ldots,\WWv^{(K)},\XXv^{(K)},\betav,\zv) \), this leads to the log-likelihood
\begin{EQA}
	\LL(\prmtv)
	&=&
	- \frac{1}{2} \| \Zv - \zv \|^{2} 
	- \frac{1}{2} \| \zv - \SmOp {\XXv^{(K)}}^{\T} \betav \|^{2}
	- \frac{\lambda_{\xxv,K}}{2} \| \xxv^{(K)} \|^{2}	
	- \frac{\lambda_{\betav}}{2} \| \betav \|^{2}	
	- \frac{\lambda_{\zv}}{2} \| \zv \|^{2}	
	\\
	&& 
	- \, \frac{\lambda_{\wwv}}{2} \| \wwv \|^{2}	
	- \frac{\lambda_{\xxv}}{2} \| \xxv \|^{2}	
	- \frac{1}{2} \sum_{k=1}^{K} \| \XXv^{(k)} - \sigma(\WWv^{(k)} \XXv^{(k-1)}) \|^{2}  
	\, .
\label{dtg6wtwyc7wctdc7wh5fkR}
\end{EQA}
Overparametrization enables us to limit the parameter set to those  
\( \prmtv = (\WWv,\XXv,\betav,\zv) \) for which \( \XXv^{(k)} = \sigma(\WWv^{(k)} \XXv^{(k-1)}) \) and 
\( \zv = \SmOp {\XXv^{(K)}}^{\T} \betav \).
Then it holds for any \( \uv = (\wwv,\xxv,\bv,\hv) \) with
\( \wwv = (\wwv_{m}^{(k)}) \in \R^{\dimW} \), 
\( \xxv = (\xxv_{i}^{(k)}) \in \R^{\dimX} \), 
\( \bv \in \R^{\dimp} \), and \( \hv \in \R^{\dimq} \) similarly to \eqref{d7vyu7e3hirjo3gbi8r}
\begin{EQA}
	&& \nquad
	- \frac{d^{2}}{dt^{2}} \LL(\prmtv + t \uv ) \bigg|_{t=0}
	=
	\lambda_{\xxv} \| \xxv \|^{2} + \lambda_{\wwv} \| \wwv \|^{2}
	\\
	&+& \sum_{k=1}^{K} \sum_{\mm=1}^{\Mh_{k}} \sumi
	\Bigl\{ \sigma'(\WWv_{\mm}^{(k)} \XXv_{i}^{(k-1)}) \, 
	\bigl( \wwv_{\mm}^{(k)} \XXv_{i}^{(k-1)} + \WWv_{\mm}^{(k)} \xxv_{i}^{(k-1)} \bigr) - \xxv_{\mm i}^{(k)} \Bigr\}^{2}	
	\\
	&+& \| \hv - \SmOp ({\XXv^{(K)}}^{\T} \bv + {\xxv^{(K)}}^{\T} \betav) \|^{2}
	+ \lambda_{\xxv,K} \| \xxv^{(K)} \|^{2}
	+ (1 + \lambda_{\hv}) \| \hv \|^{2} 
	+ \lambda_{\betav} \| \bv \|^{2}
	\, .
	\qquad
\label{d7vyu7e3hirjo3gbi8rKR}
\end{EQA}
With 
\( \prmtv^{(k)} = (\WWv^{(k)},\XXv^{(k-1)}) \),
define for \( k \leq K \) 
\begin{EQA}[rcl]
	\IF_{\mm}^{(k)}(\prmtv^{(k)})
	& \eqdef &
	\sumi \sigma'\bigl(\WWv_{\mm}^{(k)} \XXv_{i}^{(k-1)} \bigr)^{2} \, \XXv_{i}^{(k-1)} \, {\XXv_{i}^{(k-1)}}^{\T} 
	+ 2 \lambda_{\wwv} \Id_{\Mh_{k}} \, ,
	\quad
	m = 1,\ldots,\Mh_{k} \, ,
	\\
	\IF^{(k)}(\prmtv^{(k)})
	& \eqdef &
	\blk\bigl\{ \IF_{1}^{(k)}(\prmtv^{(k)}), \ldots, \IF_{\Mh_{k}}^{(k)}(\prmtv^{(k)}) \bigr\} \, .
\label{df7duyd5tw3hdy7dgywehfK}
\end{EQA}
Similarly define for \( k = 2, \ldots, K \) and \( i =1, \ldots, n \)
\begin{EQA}
	\IFG_{i}^{(k)}(\prmtv^{(k)})
	& \eqdef &
	\sum_{\mm=1}^{\Mh_{k}}  
	\sigma'\bigl( \WWv_{\mm}^{(k)} \XXv_{i}^{(k-1)} \bigr)^{2} \, \,
	{\WWv_{\mm}^{(k)}}^{\T} \WWv_{\mm}^{(k)} 
	\, ,
	\\
	\IFG^{(k)}(\prmtv^{(k)})
	& \eqdef &
	\blk\bigl\{ \IFG_{1}^{(k)}(\prmtv^{(k)}), \ldots, \IFG_{n}^{(k)}(\prmtv^{(k)}) \bigr\} 
	\, .
\label{df7duyd5twvyehders}
\end{EQA}
Also set \( \IFG^{(1)}(\prmtv^{(1)}) \equiv 0 \).

\begin{lemma}
For \( \prmtv = (\WWv,\XXv, \betav, \zv) \), let 
\( \XXv^{(k)} \equiv \sigma(\WWv^{(k)} \XXv^{(k-1)}) \), \( k=1,\ldots,K \), and \( \zv = \SmOp {\XXv^{(K)}}^{\T} \betav \).
Also, assume \( \| \betav \| \leq b_{0} \) and
\begin{EQA}[rcl]
	\max_{i =1,\ldots, n} \,\, \max_{k = 1,\ldots,K} \, \frac{3}{2} \| \IFG_{i}^{(k)}(\prmtv^{(k)}) \|_{\oper}
	& \leq &
	G_{0}
	\leq 
	\lambda_{\xxv} - 4
	\, .
\label{cufcyc6cf55df5t63gh9}
\end{EQA}
Then it holds
\begin{EQA}[rcl]
	\IFT(\prmtv) 
	& \geq &
	\blk\Bigl\{ 
		\frac{1}{2} \IF^{(1)}(\prmtv^{(1)}), \,\,
		\Id_{\XXv^{(1)}} 
		\, , \ldots,
		\frac{1}{2} \IF^{(K)}(\prmtv^{(K)}), \,\,
		\Id_{\XXv^{(K)}} , \,\,
	\\
	&&
	\frac{1}{2} \, \XXv^{(K)} \SmOpT {\XXv^{(K)}}^{\T} + \lambda_{\betav} \Id_{\betav} \, , \;
	(\lambda_{\zv} - 2) \Id_{\zv}
	\Bigr\}
	\, .
\label{duvduewhctsudkwfwft}
\end{EQA}
\end{lemma}
\begin{proof}
The use of \( (a + b + c)^{2} + 3a^{2} + 3b^{2}/2 \geq c^{2}/2 \) implies for \( k=1,\ldots,K \)
\begin{EQA}[rcl]
	&& \nquad
	3 \| \xxv^{(k)} \|^{2} 
	+ \sum_{\mm=1}^{\Mh_{k}} \sumi \Bigl\{ \sigma'(\WWv_{\mm}^{(k)} \XXv_{i}^{(k-1)}) \, 
	\bigl( \wwv_{\mm}^{(k)} \XXv_{i}^{(k-1)} + \WWv_{\mm}^{(k)} \xxv_{i}^{(k-1)} \bigr) - \xxv_{\mm i}^{(k)} \Bigr\}^{2}
	+ \lambda_{\wwv} \| \wwv^{(k)} \|^{2}
	\\
	& \geq &
	- \frac{3}{2} \sumi \sum_{\mm=1}^{\Mh_{k}} \bigl| \sigma'(\WWv_{\mm}^{(k)} \XXv_{i}^{(k-1)}) \, \WWv_{\mm}^{(k)} \xxv_{i}^{(k-1)} \bigr|^{2}
	+ \frac{1}{2} \sum_{\mm=1}^{\Mh_{k}} \sumi \bigl| \sigma'(\WWv_{\mm}^{(k)} \XXv_{i}^{(k-1)}) \, \wwv_{\mm}^{(k)} \XXv_{i}^{(k-1)} \bigr|^{2}
	\\
	&=&
	- \frac{3}{2} \sumi {\xxv_{i}^{(k-1)}}^{\T} \IFG_{i}^{(k)}(\prmtv^{(k)}) \xxv_{i}^{(k-1)}
	+ \frac{1}{2} \sum_{\mm=1}^{\Mh_{k}} \wwv_{\mm}^{(k)} \IF_{\mm}^{(k)}(\prmtv^{(k)}) {\wwv_{\mm}^{(k)}}^{\T}
	\\
	& \geq &
	- \frac{3 G_{0}}{2} \sumi \| \xxv_{i}^{(k-1)} \|^{2}
	+ \frac{1}{2} \sum_{\mm=1}^{\Mh_{k}} \wwv_{\mm}^{(k)} \IF_{\mm}^{(k)}(\prmtv^{(k)}) {\wwv_{\mm}^{(k)}}^{\T}
	\, .
\end{EQA}
Similarly, by \( \SmOp^{\T} \SmOp = \SmOpT \)
\begin{EQA}[rcl]
	&& \nquad
	\| \hv - \SmOp ({\XXv^{(K)}}^{\T} \bv + {\xxv^{(K)}}^{\T} \betav) \|^{2} + 3 \| \hv \|^{2} + \frac{3}{2} \| \SmOp {\xxv^{(K)}}^{\T} \betav \|^{2}
	\\
	& \geq &
	\frac{1}{2} \| \SmOp {\XXv^{(K)}}^{\T} \bv \|^{2}
	=
	\frac{1}{2} \bv^{\T} \XXv^{(K)} \SmOpT \, {\XXv^{(K)}}^{\T} \bv 
	\, .
\end{EQA}
Also, \( \| \betav \| \leq b_{0} \) and \( \| \SmOpT \| \leq 1 \) imply
\begin{EQA}[c]
	\| \SmOp {\xxv^{(K)}}^{\T} \betav \|^{2}
	\leq 
	b_{0}^{2} \, \| \xxv^{(K)} \SmOpT \, {\xxv^{(K)}}^{\T} \|_{\oper}
	\leq 
	b_{0}^{2} \, \| \xxv^{(K)} \|^{2} 
	\, .
\end{EQA}
The assertion follows from \eqref{d7vyu7e3hirjo3gbi8rKR} by \eqref{cufcyc6cf55df5t63gh9} and \( \lambda_{\xxv} \geq 4 + G_{0} \).
\end{proof}

\newpage
\appendix


\Chapter{Local smoothness and linearly perturbed optimization}
\label{Slocalsmooth}
This section discusses the problem of linearly and quadratically perturbed optimization
of a smooth and concave function \( \fs(\upsv) \), \( \upsv \in \R^{\dimp} \).

\Section{Gateaux smoothness and self-concordance}
Below we assume 
the function \( \fs(\upsv) \) to be strongly concave with the negative Hessian 
\( \IFN(\upsv) \eqdef - \nabla^{2} \fs(\upsv) \in \Matr_{\dimp} \) positive definite. 
Also, assume \( \fs(\upsv) \) three or sometimes even four times Gateaux differentiable in \( \upsv \in \Ups \).
For any particular direction \( \uv \in \R^{\dimp} \), we consider the univariate function 
\( \fs(\upsv + t \uv) \) and measure its smoothness in \( t \).
Local smoothness of \( \fs \) will be described by the relative error of the Taylor expansion 
of the third or fourth order.
Namely, define
\begin{EQ}[rcl]
	\dltw_{3}(\upsv,\uv) 
	&=& 
	\fs(\upsv + \uv) - \fs(\upsv) - \langle \nabla \fs(\upsv), \uv \rangle 
	- \frac{1}{2} \langle \nabla^{2} \fs(\upsv), \uv^{\otimes 2} \rangle , 
	\\
	\dltwd_{3}(\upsv,\uv) 
	&=&
	\langle \nabla \fs(\upsv + \uv), \uv \rangle - \langle \nabla \fs(\upsv), \uv \rangle 
	- \langle \nabla^{2} \fs(\upsv), \uv^{\otimes 2} \rangle \, ,
\label{dltw3vufuv12f2ga}
\end{EQ}
and
\begin{EQA}
	\dltw_{4}(\upsv,\uv)
	& \eqdef &
	\fs(\upsv + \uv) - \fs(\upsv) - \langle \nabla \fs(\upsv), \uv \rangle 
	- \frac{1}{2} \langle \nabla^{2} \fs(\upsv), \uv^{\otimes 2} \rangle
	- \frac{1}{6} \langle \nabla^{3} \fs(\upsv), \uv^{\otimes 3} \rangle \, .
\label{hvcduywgedfuyg2y1y35e3wweg}
\end{EQA}
Now, for each \( \upsv \), suppose to be given a positive symmetric operator 
\( \DFN(\upsv) \in \Matr_{\dimp} \) 
defining a local metric and a local vicinity around \( \upsv \):
\begin{EQA}
	\UVz_{\rr}(\upsv)
	&=&
	\bigl\{ \uv \in \R^{\dimp} \colon \| \DFN(\upsv) \uv \| \leq \rr \bigr\}
\label{ed7sycf7wedwgedq2ftwdfgtv}
\end{EQA}
for some radius \( \rr \).

Local smoothness properties of \( \fs \) at \( \upsv \) are given via the quantities
\begin{EQA}[rcccl]
    \dltwb(\upsv)
    & \eqdef &
    \sup_{\uv \colon \| \DFN(\upsv) \uv \| \leq \rr} \,
    \frac{2|\dltw_{3}(\upsv,\uv)|}{\| \DFN(\upsv) \uv \|^{2}} 
    \,\, ,
    \quad
    \dltwbd(\upsv)
    & \eqdef &
    \sup_{\uv \colon \| \DFN(\upsv) \uv \| \leq \rr} \, \frac{|\dltwd_{3}(\upsv,\uv)|}{\| \DFN(\upsv) \uv \|^{2}} \,\, . 
    \qquad
\label{dtb3u1DG2d3GPg}
\end{EQA}
The definition yields for any \( \uv \) with \( \| \DFN(\upsv) \uv \| \leq \rr \)
\begin{EQ}[rcccl]
	\bigl| \dltw_{3}(\upsv,\uv) \rangle \bigr|
	& \leq &
	\frac{\dltwb(\upsv)}{2} \| \DFN(\upsv) \uv \|^{2} 
	\, ,
	\qquad
	\bigl| \dltwd_{3}(\upsv,\uv) \bigr|
	& \leq &
	\dltwbd(\upsv) \| \DFN(\upsv) \uv \|^{2}
	\, .
	\qquad
\label{dta3u1DG2d3GPa1g}
\end{EQ}
%
The approximation results can be improved 
provided a third order upper bound on the error of Taylor expansion. 

\begin{description}
    \item[\label{LL3tref} \( \bb{(\mathcal{T}_{3})} \)]
      \textit{For some \( \dltwu_{3} \)}
\begin{EQA}
	\bigl| \dltw_{3}(\upsv,\uv) \bigr|
	& \leq &
	\frac{\dltwu_{3}}{6} \| \DFN(\upsv) \, \uv \|^{3} \, ,
	\quad
	\bigl| \dltwd_{3}(\upsv,\uv) \bigr|
	\leq 
	\frac{\dltwu_{3}}{2} \| \DFN(\upsv) \, \uv \|^{3} \, ,
	\quad
	\uv \in \UVz_{\rr}(\upsv).
\label{bd3xu16f3uo3st}
\end{EQA}
\end{description}
 
\begin{description}
    \item[\label{LL4tref} \( \bb{(\mathcal{T}_{4})} \)]
      \textit{For some \( \dltwu_{4} \)}
\begin{EQA}
	\bigl| \dltw_{4}(\upsv,\uv) \bigr|
	& \leq &
	\frac{\dltwu_{4}}{24} \| \DFN(\upsv) \, \uv \|^{4} \, ,
	\qquad
	\uv \in \UVz_{\rr}(\upsv).
\label{1mffmxum5st}
\end{EQA}
\end{description}

We also present a version of \nameref{LL3tref} resp. \nameref{LL4tref} in terms of the third (resp. fourth) derivative of \( \fs \).
\begin{description}
    \item[\label{LLsT3ref} \( \bb{(\mathcal{T}_{3}^{*})} \)]
    \emph{\( \fs(\upsv) \) is three times differentiable and 
	}
\begin{EQA}
    \sup_{\uv \colon \| \DFN(\upsv) \uv \| \leq \rr} \,\, \sup_{\zv \in \R^{\dimp}} \,\, 
    \frac{\bigl| \langle \nabla^{3} \fs(\upsv + \uv), \zv^{\otimes 3} \rangle \bigr|}
		 {\| \DFN(\upsv) \zv \|^{3}} 
	& \leq &
	\dltwu_{3} \, .
\label{jcxhydtferyu9j3d6vhew}
\end{EQA}

    \item[\label{LLsT4ref} \( \bb{(\mathcal{T}_{4}^{*})} \)]
    \emph{\( \fs(\upsv) \) is four times differentiable and 
	}
\begin{EQA}
    \sup_{\uv \colon \| \DFN(\upsv) \uv \| \leq \rr} \,\, \sup_{\zv \in \R^{\dimp}} \,\, 
    \frac{\bigl| \langle \nabla^{4} \fs(\upsv + \uv), \zv^{\otimes 4} \rangle \bigr|}
		 {\| \DFN(\upsv) \zv \|^{4}} 
	& \leq &
	\dltwu_{4} \, .
\label{jcxhydtferyu9j3d6vhew4}
\end{EQA}

\end{description}

%
\noindent
By Banach's characterization \cite{Banach1938}, \nameref{LLsT3ref} implies
\begin{EQA}
	\bigl| \langle \nabla^{3} \fs(\upsv + \uv), \zv_{1} \otimes \zv_{2} \otimes \zv_{3} \rangle \bigr|
	& \leq &	 
	\dltwu_{3} \| \DFN(\upsv) \zv_{1} \| \, \| \DFN(\upsv) \zv_{2} \| \, \| \DFN(\upsv) \zv_{3} \| \, 
\label{jbuyfg773jgion94euyyfg}
\end{EQA}
for any \( \uv \) with \( \| \DFN(\upsv) \uv \| \leq \rr \) and all \( \zv_{1} , \zv_{2}, \zv_{3} \in \R^{\dimp} \).
Similarly under \nameref{LLsT4ref}
\begin{EQA}
	\bigl| \langle \nabla^{4} \fs(\upsv + \uv), \zv_{1} \otimes \zv_{2} \otimes \zv_{3} \otimes \zv_{4} \rangle \bigr|
	& \leq &	 
	\dltwu_{4} \prod_{k=1}^{4} \| \DFN(\upsv) \zv_{k} \| \, ,
	\quad 
	\zv_{1} , \zv_{2}, \zv_{3}, \zv_{4} \in \R^{\dimp} \, .
	\qquad
\label{jbuyfg773jgion94euyyfg4}
\end{EQA}

\begin{lemma}
\label{LdltwLa3t}
Under \nameref{LL3tref} or \nameref{LLsT3ref},
the values \( \dltwb(\upsv) \) and \( \dltwbd(\upsv) \) from \eqref{dtb3u1DG2d3GPg} satisfy
\begin{EQA}[rcccl]
\label{gtcdsftdffrvsewsea}
	\dltwb(\upsv)
	& \leq &
	\frac{\dltwu_{3} \, \rr}{3 } \, ,
	\qquad
	\dltwbd(\upsv)
	& \leq &
	\frac{\dltwu_{3} \, \rr}{2} \, ,
	\qquad
	\upsv \in \Upsd .
\label{gtcdsftdfvtwdsefhfdvfrvsewseG}
\end{EQA}
\end{lemma}

\begin{proof}
For any \( \uv \in \UVz_{\rr}(\upsv) \) with \( \| \DFN(\upsv) \uv \| \leq \rr \)
\begin{EQA}
	\bigl| \dltw_{3}(\upsv,\uv) \bigr|
	& \leq &
	\frac{\dltwu_{3}}{6} \, \| \DFN(\upsv) \uv \|^{3} 
	\leq 
	\frac{\dltwu_{3} \, \rr}{6} \, \| \DFN(\upsv) \uv \|^{2},
\label{jrgeteteer2234587654}
\end{EQA}
and the bound for \( \dltwb(\upsv) \) follows.
The proof for \( \dltwbd(\upsv) \) is similar.
\end{proof}

The values \( \dltwu_{3} \) and \( \dltwu_{4} \) are usually very small.
Some quantitative bounds are given later in this section
under the assumption that the function \( \fs(\upsv) \) can be written in the form \( - \fs(\upsv) = n \hL(\upsv) \) 
for a fixed smooth function \( h(\upsv) \) with the Hessian \( \nabla^{2} \hL(\upsv) \). 
The factor \( n \) has meaning of the sample size%
\ifapp{; see \Chname \ref{ScritdimMLE} or \Chname \ref{SGBvM}.}{.}

\begin{description}
    \item[\label{LLtS3ref} \( \bb{(\mathcal{S}_{3}^{*})} \)]
      \emph{ \( - \fs(\upsv) = n \hL(\upsv) \) for \( \hL(\upsv) \) three times differentiable and
\begin{EQA}
	\sup_{\uv \colon \| \HL(\upsv) \uv \| \leq \rr/\sqrt{n}} 
	\frac{\bigl| \langle \nabla^{3} \hL(\upsv + \uv), \uv^{\otimes 3} \rangle \bigr|}{\| \HL(\upsv) \uv \|^{3}}
	& \leq &
	\hmax_{3} \, .
\end{EQA}
}
    \item[\label{LLtS4ref} \( \bb{(\mathcal{S}_{4}^{*})} \)]
      \emph{ the function \( \hL(\cdot) \) satisfies \nameref{LLtS3ref} and  
\begin{EQA}
	\sup_{\uv \colon \| \HL(\upsv) \uv \| \leq \rr/\sqrt{n}}
	\frac{\bigl| \langle \nabla^{4} \hL(\upsv + \uv), \uv^{\otimes 4} \rangle \bigr|}{\| \HL(\upsv) \uv \|^{4}}
	& \leq &
	\hmax_{4} \, .
\end{EQA}
}
\end{description}

\noindent
\nameref{LLtS3ref} and \nameref{LLtS4ref}
are local versions of the so-called self-concordance condition; see \cite{Ne1988} and \cite{OsBa2021}.
In fact, they require that each univariate function \( \hL(\upsv + t \uv) \) of \( t \in \R \)
is self-concordant with some universal constants \( \hmax_{3} \) and \( \hmax_{4} \).
Under \nameref{LLtS3ref} and \nameref{LLtS4ref}, with \( \DFN^{2}(\upsv) = n \, \HL^{2}(\upsv) \), the values 
\( \dltw_{3}(\upsv,\uv) \), \( \dltw_{4}(\upsv,\uv) \), and \( \dltwb(\upsv) \), \( \dltwbd(\upsv) \) can be
bounded.

\begin{lemma}
\label{LdltwLaGP}
Suppose \nameref{LLtS3ref}.
Then 
\nameref{LL3tref} follows with \( \dltwu_{3} = \hmax_{3} n^{-1/2} \).
Moreover, for \( \dltwb(\upsv) \) and \( \dltwbd(\upsv) \) from \eqref{dtb3u1DG2d3GPg}, it holds
\begin{EQA}[rcccl]
	\dltwb(\upsv)
	& \leq &
	\frac{\hmax_{3} \, \rr}{3 n^{1/2}} \, ,
	\qquad
	\dltwbd(\upsv)
	& \leq &
	\frac{\hmax_{3} \, \rr}{2 n^{1/2}} \, .
\label{gtcdsftdfvtwdsefhfdvfrvsewseGP}
\end{EQA}
Also \nameref{LL4tref} follows from \nameref{LLtS4ref} with \( \dltwu_{4} = \hmax_{4} n^{-1} \).
\end{lemma}

\begin{proof}
For any \( \uv \in \UVz_{\rr}(\upsv) \) and \( t \in [0,1] \), by the Taylor expansion of the third order
\begin{EQA}
	|\dltw(\upsv,\uv)|
	& \leq &
	\frac{1}{6} \bigl| \langle \nabla^{3} \fs(\upsv + t \uv), \uv^{\otimes 3} \rangle \bigr|
	=
	\frac{n}{6} \, \bigl| \langle \nabla^{3} \hL(\upsv + t \uv), \uv^{\otimes 3} \rangle \bigr|
	\leq 
	\frac{n \, \hmax_{3}}{6} \, \| \HL(\upsv) \uv \|^{3} 
	\\
	&=&
	\frac{n^{-1/2} \, \hmax_{3}}{6} \, \| \DFN(\upsv) \uv \|^{3}
	\leq 
	\frac{n^{-1/2} \, \hmax_{3} \, \rr}{6} \, \| \DFN(\upsv) \uv \|^{2} \, .
\label{jrgeteteer2234587654}
\end{EQA}
This implies \nameref{LL3tref} as well as \eqref{gtcdsftdfvtwdsefhfdvfrvsewseGP}; see \eqref{dta3u1DG2d3GPa1g}.
The statement about \nameref{LL4tref} is similar.
\end{proof}

\def\AFN{\mathbbmsl{U}}
\def\Avm{\bb{M}}

\Section{Optimization after linear perturbation. A basic lemma}
\label{Squadnquad}

Let \( \fs(\upsv) \) be a smooth concave function, 
\begin{EQA}
	\upsvs
	&=&
	\argmax_{\upsv} \fs(\upsv),
\label{fg5hg3gf98tkj3dciryt}
\end{EQA}
and \( \IFN = - \nabla^{2} \fs(\upsvs) \).
Later we study the question of how the point of maximum and the value of maximum of \( \fs \) change if we add a linear or quadratic 
component to \( \fs \).
More precisely, let another function \( \fn(\upsv) \) satisfy for some vector \( \Av \)
\begin{EQA}
	\fn(\upsv) - \fn(\upsvs) 
	&=&
	\bigl\langle \upsv - \upsvs, \Av \bigr\rangle + \fs(\upsv) - \fs(\upsvs) .
\label{4hbh8njoelvt6jwgf09}
\end{EQA}
A typical example corresponds to \( \fs(\upsv) = \E L(\upsv) \) and \( \fn(\upsv) = L(\upsv) \) 
for a random function \( L(\upsv) \) with a linear stochastic component \( \zeta(\upsv) = L(\upsv) - \E L(\upsv) \)%
\ifapp{; see \nameref{Eref}.}{.}
Then \eqref{4hbh8njoelvt6jwgf09} is satisfied with \( \Av =	\nabla \zeta \).
Define
\begin{EQA}
	\upsvn
	& \eqdef &
	\argmax_{\upsv} \fn(\upsv),
	\qquad
	\fn(\upsvn)
	=
	\max_{\upsv} \fn(\upsv) .
\label{6yc63yhudf7fdy6edgehy} 
\end{EQA}
The aim of the analysis is to evaluate the quantities \( \upsvn - \upsvs \) and
\( \fn(\upsvn) - \fn(\upsvs) \).
First, we consider the case of a quadratic function \( \fs \).

\begin{lemma}
\label{Pquadquad}
Let \( \fs(\upsv) \) be quadratic with \( \nabla^{2} \fs(\upsv) \equiv - \IFN \) and
\( \fn(\upsv) \) satisfy \eqref{4hbh8njoelvt6jwgf09}. 
Then 
\begin{EQA}
	\upsvn - \upsvs
	&=&
	\IFN^{-1} \Av,
	\qquad
	\fn(\upsvn) - \fn(\upsvs)
	=
	\frac{1}{2} \| \IFN^{-1/2} \Av \|^{2} .
\label{kjcjhchdgehydgtdtte35}
\end{EQA}
\end{lemma}

\begin{proof}
If \( \fs(\upsv) \) is quadratic, then, of course, under \eqref{4hbh8njoelvt6jwgf09}, \( \fn(\upsv) \) is quadratic as well
with \( - \nabla^{2} \fn(\upsv) \equiv \IFN \).
This implies
\begin{EQA}
	\nabla \fn(\upsvs) - \nabla \fn(\upsvn)
	&=&
	\IFN (\upsvn - \upsvs) .
\label{dcudydye67e6dy3wujhds7}
\end{EQA}
Further, \eqref{4hbh8njoelvt6jwgf09} and \( \nabla \fs(\upsvs) = 0 \) yield \( \nabla \fn(\upsvs) = \Av \).
Together with \( \nabla \fn(\upsvn) = 0 \), this implies
\( \upsvn - \upsvs = \IFN^{-1} \Av \).
The Taylor expansion of \( \fn \) at \( \upsvn \) yields by \( \nabla \fn(\upsvn) = 0 \)
\begin{EQA}
	\fn(\upsvs) - \fn(\upsvn)
	&=&
	- \frac{1}{2} \| \IFN^{1/2} (\upsvn - \upsvs) \|^{2}
	=
	- \frac{1}{2} \| \IFN^{-1/2} \Av \|^{2} 
\label{8chuctc44wckvcuedje}
\end{EQA}
and the assertion follows.
\end{proof}
  
The next result describes the concentration properties of \( \upsvn \) from \eqref{6yc63yhudf7fdy6edgehy} in a local elliptic set
of the form
\begin{EQA}
	\CA(\rr)
	& \eqdef &
	\{ \upsv \colon \| \IFN^{1/2} (\upsv - \upsvs) \| \leq \rr \} ,
\label{0cudc7e3jfuyvct6eyhgwe}
\end{EQA}
where \( \rr \) is slightly larger than \( \| \IFN^{-1/2} \Av \| \).

\begin{proposition}
\label{Pconcgeneric}
Let \( \fs(\upsv) \) be a strongly concave function with \( \fs(\upsvs) = \max_{\upsv} \fs(\upsv) \)  
and \( \IFN = - \nabla^{2} \fs(\upsvs) \).
Let further \( \fn(\upsv) \) and \( \fs(\upsv) \) be related by \eqref{4hbh8njoelvt6jwgf09} with some vector \( \Av \).
Fix \( \amax < 1 \) and \( \rrn \) such that \( \| \IFN^{-1/2} \Av \| \leq \amax \, \rrn \).
Suppose now that \( \fs(\upsv) \) satisfy \eqref{dtb3u1DG2d3GPg} for \( \upsv = \upsvs \), 
\( \DFN(\upsvs) = \DFN \leq \dmax \, \IFN^{1/2} \) with some \( \dmax > 0 \) 
and \( \dltwbd \) such that 
\begin{EQA}
	1 - \amax - \dltwbd \dmax^{2}
	& > &
	0 .
\label{rrm23r0ut3ua}
\end{EQA}
Then for \( \upsvn \) from \eqref{6yc63yhudf7fdy6edgehy}, it holds 
\begin{EQA}
	\| \IFN^{1/2} (\upsvn - \upsvs) \|  
	& \leq &
	\rrn \, 
	\quad
	\text{ and }
	\quad
	\| \DFN (\upsvn - \upsvs) \|
	\leq 
	\dmax \, \rrn \, . 
\label{rhDGtuGmusGU0a}
\end{EQA}
\end{proposition}

\begin{proof}
Rescaling \( \DFN \) by \( \dmax^{-1} \) reduces the proof to \( \dmax = 1 \).
The bound \( \| \IFN^{-1/2} \Av \| \leq \amax \, \rrn \) implies for any \( \uv \)
\begin{EQA}
	\bigl| \langle \Av, \uv \rangle \bigr|
	& = &
	\bigl| \langle \IFN^{-1/2} \Av, \IFN^{1/2} \uv \rangle \bigr|
	\leq 
	\amax \, \rrn \| \IFN^{1/2} \uv \| \, .
\label{LLoDGm1nzua}
\end{EQA}
%
Let \( \upsv \) be a point on the boundary of the set \( \CA(\rrn) \) from \eqref{0cudc7e3jfuyvct6eyhgwe}.
We also write \( \uv = \upsv - \upsvs \).
The idea is to show that the derivative  \( \frac{d}{dt} \fn(\upsvs + t \uv) < 0 \) 
is negative for \( t > 1 \).
Then all the extreme points of \( \fn(\upsv) \) are within \( \CA(\rrn) \).
We use the decomposition
\begin{EQA}
	\fn(\upsvs + \rhot \uv) - \fn(\upsvs)
	&=&
	\langle \Av, \uv \rangle \, \rhot 
	+ \fs(\upsvs + \rhot \uv) - \fs(\upsvs) .
\label{LGtsGtuLGtsa}
\end{EQA}
With \( \fGu(t) = \fs(\upsvs + \rhot \uv) - \fs(\upsvs) + \langle \Av, \uv \rangle \, \rhot \), it holds
\begin{EQA}
	\frac{d}{d \rhot} \fs(\upsvs + \rhot \uv)
	&=&
	- \langle \Av, \uv \rangle + \fGu'(\rhot) .
\label{frddtLtGstua}
\end{EQA}
By definition of \( \upsvs \), it also holds \( \fGu'(0) = \langle \Av, \uv \rangle \).
The identity \( - \nabla^{2} \fs(\upsvs) = \IFN \) yields \( - \fGu''(0) = \| \IFN^{1/2} \uv \|^{2} \).
Bound \eqref{dta3u1DG2d3GPa1g} implies for \( | \rhot | \leq 1 \)
\begin{EQA}
	\bigl| \fGu'(\rhot) - \fGu'(0) - \rhot \fGu''(0) \bigr|
	& \leq &
	\rhot \, \| \DFN \uv \|^{2} \, \dltwbd \, .
\label{fptfp0fpttfpp13a}
\end{EQA}
For \( \rhot = 1 \), we obtain by \eqref{rrm23r0ut3ua} 
\begin{EQA}
	\fGu'(1) 
	& \leq &
	- \langle \Av, \uv \rangle - \| \IFN^{1/2} \uv \|^{2} + \| \DFN \uv \|^{2} \, \dltwbd
	\leq 
	- \| \IFN^{1/2} \uv \|^{2} (1 -  \dltwbd - \amax)
	< 0 .
\label{fp1fpp13d3rGa}
\end{EQA}
Moreover, concavity of \( \fGu(\rhot) \) imply that \( \fGu'(\rhot) - \fGu'(0) \) decreases in 
\( \rhot \) for \( \rhot > 1 \).
Further, summing up the above derivation yields 
\begin{EQA}
	\frac{d}{dt} \fn(\upsvs + \rhot \uv) \Big|_{\rhot=1}
	& \leq &
	- \| \IFN^{1/2} \uv \|^{2} (1 - \amax - \dltwbd)
	< 0 .
\label{ddtLGtstu33a}
\end{EQA}
As \( \frac{d}{d \rhot} \fn(\upsvs + \rhot \uv) \) decreases with \( \rhot \geq 1 \) together with 
\( \fGu'(\rhot) \) due to \eqref{frddtLtGstua}, the same applies to all such \( \rhot \).
This implies the assertion.
\end{proof}

The result of Proposition~\ref{Pconcgeneric} allows to localize the point \( \upsvn = \argmax_{\upsv} \fn(\upsv) \)
in the local vicinity \( \CA(\rrn) \) of \( \upsvs \).
The use of smoothness properties of \( \fn \) or, equivalently, of \( \fs \), in this vicinity helps to obtain
rather sharp expansions for \( \upsvn - \upsvs \) and for \( \fn(\upsvn) - \fn(\upsvs) \).

\begin{proposition}
\label{PFiWigeneric}
Under the conditions of Proposition~\ref{Pconcgeneric}, 
\begin{EQ}[rcccl]
    - \frac{\dltwb}{1 + \dmax^{2} \dltwb} \| \DFN \, \IFN^{-1} \Av \|^{2}
    & \leq &
    2 \fn(\upsvn) - 2 \fn(\upsvs) 
    - \| \IFN^{-1/2} \Av \|^{2}
    & \leq &
    \frac{\dltwb}{1 - \dmax^{2} \hspm \dltwb} \| \DFN \, \IFN^{-1} \Av \|^{2} \, .
    \qquad \quad
\label{3d3Af12DGttGa}
\end{EQ}
Also
\begin{EQ}[rcl]
    \| \DFN (\upsvn - \upsvs - \IFN^{-1} \Av) \|
    & \leq &
    \frac{\sqrt{3 \dltwb}}{1 - \dmax^{2} \hspm \dltwb} \, \| \DFN \, \IFN^{-1} \Av \| \, ,
    \\
    \| \DFN (\upsvn - \upsvs) \|
    & \leq &
    \frac{1 + \sqrt{3 \dltwb}}{1 - \dmax^{2} \hspm \dltwb} \, \| \DFN \, \IFN^{-1} \Av \| \, .
\label{DGttGtsGDGm13rGa}
\end{EQ}
\end{proposition}

\begin{proof}
As in the proof of Proposition~\ref{Pconcgeneric}, rescaling \( \DFN \) by \( \dmax^{-1} \) reduces the statement to \( \dmax = 1 \).
By \eqref{dtb3u1DG2d3GPg}, for any \( \upsv \in \CA(\rrn) \)
\begin{EQA}
	\Bigl| 
		\fs(\upsvs) - \fs(\upsv) - \frac{1}{2} \| \IFN^{1/2} (\upsv - \upsvs) \|^{2} 
	\Bigr|
	& \leq &
	\frac{\dltwb}{2} \| \DFN (\upsv - \upsvs) \|^{2} .
\label{d3GrGELGtsG12}
\end{EQA}
Further, 
\begin{EQA}[rcl]
	&& \nquad
	\fn(\upsv) - \fn(\upsvs) - \frac{1}{2} \| \IFN^{-1/2} \Av \|^{2}
	\\
	&=&
	\bigl\langle \upsv - \upsvs, \Av \bigr\rangle
	+ \fs(\upsv) - \fs(\upsvs) - \frac{1}{2} \| \IFN^{-1/2} \Av \|^{2} 
	\\
	&=&
	- \frac{1}{2} \bigl\| \IFN^{1/2} (\upsv - \upsvs) - \IFN^{-1/2} \Av \bigr\|^{2}
	+ \fs(\upsv) - \fs(\upsvs) + \frac{1}{2} \| \IFN^{1/2} (\upsv - \upsvs) \|^{2} .
	\qquad 
\label{12ELGuELusG}
\end{EQA}
As \( \upsvn \in \CA(\rrn) \) and it maximizes \( \fn(\upsv) \), we derive by \eqref{d3GrGELGtsG12}
\begin{EQA}
	&& \nquad
	\fn(\upsvn) - \fn(\upsvs) - \frac{1}{2} \| \IFN^{-1/2} \Av \|^{2}
	=
	\max_{\upsv \in \CA(\rrn)} 
	\Bigl\{ 
		\fn(\upsv) - \fn(\upsvs) - \frac{1}{2} \| \IFN^{-1/2} \Av \|^{2} 
	\Bigr\}
	\\
	& \leq &
	\max_{\upsv \in \CA(\rrn)} 
	\Bigl\{ 
		- \frac{1}{2} \bigl\| \IFN^{1/2} (\upsv - \upsvs) - \IFN^{-1/2} \Av \bigr\|^{2} 
		+ \frac{\dltwb}{2} \| \DFN (\upsv - \upsvs) \|^{2}
	\Bigr\} .
\label{d3G1212222B} 
\end{EQA}
Denote \( \uv = \IFN^{1/2} (\upsv - \upsvs) \), \( \xiv = \IFN^{-1/2} \Av \), and 
\( \BFN = \IFN^{-1/2} \, \DFN^{2} \, \IFN^{-1/2} \).
As \( \DFN^{2} \leq \IFN \) and \( \dltwb < 1 \), it holds \( \| \BFN \| \leq 1 \) and 
\begin{EQA}
	&& \nquad
	\max_{\upsv \in \CA(\rrn)} \bigl\{ - \bigl\| \IFN^{1/2} (\upsv - \upsvs) - \IFN^{-1/2} \Av \bigr\|^{2} 
		+ \dltwb \| \DFN (\upsv - \upsvs) \|^{2} \bigr\}
	\\
	&=&
	\max_{\| \uv \| \leq \rr} \bigl\{ - \| \uv - \xiv \|^{2} + \dltwb \, \uv^{\T} \BFN \uv \bigr\}
	=
	\xiv^{\T} \bigl\{ (\Id - \dltwb \, \BFN)^{-1} - \Id \bigr\} \xiv
	\leq 
	\frac{\dltwb}{1 - \dltwb} \xiv^{\T} \BFN \, \xiv
\label{d7eneyf6g53geygywn}
\end{EQA}
yielding
\begin{EQA}
	\fn(\upsvn) - \fn(\upsvs) - \frac{1}{2} \| \IFN^{-1/2} \Av \|^{2}
	& \leq &
	\frac{\dltwb}{2(1 - \dltwb)} \| \DFN \, \IFN^{-1} \Av \|^{2} . 
\label{fd3G122B2} 
\end{EQA}
Similarly 
\begin{EQA}
	&& \nquad 
	\fn(\upsvn) - \fn(\upsvs) - \frac{1}{2} \| \IFN^{-1/2} \Av \|^{2}
	\\
	& \geq &
	\max_{\upsv \in \CA(\rrn)} 
	\Bigl\{ 
		- \frac{1}{2} \bigl\| \IFN^{1/2} (\upsv - \upsvs) - \IFN^{-1/2} \Av \bigr\|^{2} 
		- \frac{\dltwb}{2} \| \DFN (\upsv - \upsvs) \|^{2}
	\Bigr\}
	\\
	& \geq &
	\frac{1}{2} \xiv^{\T} \bigl\{ (\Id + \dltwb \, \BFN)^{-1} - \Id \bigr\} \xiv
	\geq 
	- \frac{\dltwb }{2(1 + \dltwb)} \, \| \DFN \, \IFN^{-1} \Av \|^{2} . 
	\quad 
\label{fd3G122B2m} 
\end{EQA}
These bounds imply 
\eqref{3d3Af12DGttGa}.

Now we derive similarly to \eqref{12ELGuELusG} that for \( \upsv \in \CA(\rrn) \)
\begin{EQA}
	\fn(\upsv) - \fn(\upsvs) 
	& \leq &
	\bigl\langle \upsv - \upsvs, \Av \bigr\rangle
	- \frac{1}{2} \| \IFN^{1/2} (\upsv - \upsvs) \|^{2}
	+ \frac{\dltwb}{2} \| \DFN (\upsv - \upsvs) \|^{2} .
\label{LGvLGvsGf1d3G2}
\end{EQA}
A particular choice \( \upsv = \upsvn \) yields
\begin{EQA}
	\fn(\upsvn) - \fn(\upsvs) 
	& \leq &
	\bigl\langle \upsvn - \upsvs, \Av \bigr\rangle
	- \frac{1}{2} \| \IFN^{1/2} (\upsvn - \upsvs) \|^{2}
	+ \frac{\dltwb}{2} \| \DFN (\upsvn - \upsvs) \|^{2} .
\label{21GsvtvGDG3G2}
\end{EQA}
Combining this inequality with \eqref{fd3G122B2m} allows to bound
\begin{EQA}
	\bigl\langle \upsvn - \upsvs, \Av \bigr\rangle
	- \frac{1}{2} \| \IFN^{1/2} (\upsvn - \upsvs) \|^{2}
	+ \frac{\dltwb}{2} \| \DFN (\upsvn - \upsvs) \|^{2} 
	& \geq &
	\frac{1}{2} \xiv^{\T} (\Id + \dltwb \, \BFN)^{-1} \xiv .
\label{2m1DGd3G123G}
\end{EQA}
With 
\( \uvd = \IFN^{1/2} (\upsvn - \upsvs) \), this implies
\begin{EQA}
	2 \bigl\langle \uvd, \xiv \bigr\rangle - {\uvd}^{\T}(1 - \dltwb \BFN) \uvd 
	& \geq &
	\xiv^{\T} (\Id + \dltwb \, \BFN)^{-1} \xiv \, ,
\label{dtxi2fd1d22}
\end{EQA}
and hence,
\begin{EQA}
	&& \nquad
	\bigl\{ \uvd - (\Id - \dltwb \BFN)^{-1} \xiv \bigr\}^{\T} (\Id - \dltwb \BFN) \bigl\{ \uvd - (\Id - \dltwb \BFN)^{-1} \xiv \bigr\}
	\\
	& \leq &
	\xiv^{\T} \bigl\{ (\Id - \dltwb \, \BFN)^{-1} - (\Id + \dltwb \, \BFN)^{-1} \bigr\} \xiv
	\leq 
	\frac{2 \dltwb}{(1 + \dltwb) (1 - \dltwb)} \, \xiv^{\T} \BFN \, \xiv \, .
\label{uv11wxi22w1w}
\end{EQA}
Introduce \( \| \cdot \|_{\afn} \) by \( \| \xv \|_{\afn}^{2} \eqdef \xv^{\T} (\Id - \dltwb \BFN) \xv \),
\( \xv \in \R^{\dimp} \).
Then
\begin{EQA}
	\| \uvd - (\Id - \dltwb \BFN)^{-1} \xiv \|_{\afn}^{2}
	& \leq &
	\frac{2 \dltwb}{1 - \dltwb^{2}} \, \xiv^{\T} \BFN \, \xiv \, .
\label{dcumwf6vhehe6fbwhfr}
\end{EQA}
As 
\begin{EQA}
	\| \xiv - (\Id - \dltwb \BFN)^{-1} \xiv \|_{\afn}^{2}
	& = &
	\dltwb^{2} (\BFN \xiv)^{\T} (\Id - \dltwb \BFN)^{-1} \BFN \xiv
	\leq 
	\frac{\dltwb^{2}}{1 - \dltwb} \, \| \BFN \xiv \|^{2}
	\leq 
	\frac{\dltwb^{2}}{1 - \dltwb} \, \xiv^{\T} \BFN \, \xiv
\label{c8jjkie74he3tftdy3fy}
\end{EQA}
we conclude for \( \dltwb \leq 1/3 \) by the triangle inequality
\begin{EQA}
	\| \uvd - \xiv \|_{\afn}
	& \leq &
	\biggl( \sqrt{\frac{\dltwb^{2}}{1 - \dltwb}} + \sqrt{\frac{2 \dltwb}{1 - \dltwb^{2}}} \biggr)
	\sqrt{\xiv^{\T} \BFN \, \xiv}
	\leq 
	\sqrt{\frac{3 \dltwb}{1 - \dltwb}} \,\, \sqrt{\xiv^{\T} \BFN \, \xiv} \, ,
\label{uxiBws2w1w31w}
\end{EQA}
and \eqref{DGttGtsGDGm13rGa} follows by \( \Id - \dltwb \BFN \geq (1 - \dltwb) \Id \).
\end{proof}

\begin{remark}
\label{Rfsfnlinp}
The roles of the functions \( \fs \) and \( \fn \) are exchangeable.
In particular, the results from \eqref{DGttGtsGDGm13rGa} apply with
\( \IFN = - \nabla^{2} \fn(\upsvn) = - \nabla^{2} \fs(\upsvn) \) provided that 
\eqref{dtb3u1DG2d3GPg} is fulfilled at \( \upsv = \upsvn \).
\end{remark}

\Subsection{Basic lemma under third order smoothness}
The results of Proposition~\ref{PFiWigeneric} can be refined if
\( \fs \) satisfies condition \nameref{LL3tref}.

\begin{proposition}
\label{Pconcgeneric2}
Let \( \fs(\upsv) \) be a strongly concave function with \( \fs(\upsvs) = \max_{\upsv} \fs(\upsv) \)  
and \( \IFN = - \nabla^{2} \fs(\upsvs) \).
Let \( \fn(\upsv) \) fulfill \eqref{4hbh8njoelvt6jwgf09} with some vector \( \Av \).
Suppose that \( \fs(\upsv) \) follows \nameref{LL3tref} at \( \upsvs \) with 
\( \DFN^{2} \), \( \rrn \), and \( \dltwu_{3} \) such that 
\begin{EQA}
	\DFN^{2} \leq \dmax^{2} \, \IFN , 
	\quad \rrn \geq \frac{4\dmax}{3} \, \| \IFN^{-1/2} \Av \| ,
	\quad 
	\dmax^{3} \dltwu_{3} \, \| \IFN^{-1/2} \Av \|
	& < &
	\frac{1}{4} \, .
\label{yxdhewndu7jwnjjuu}
\end{EQA}
Then \( \upsvn = \argmax_{\upsv} \fn(\upsv) \) satisfies 
\begin{EQA}
	\| \IFN^{1/2} (\upsvn - \upsvs) \|  
	& \leq &
	\frac{4}{3} \| \IFN^{-1/2} \Av \| \, ,
	\qquad
	 \| \DPN (\upsvn - \upsvs) \|  
	\leq 
	\frac{4\dmax}{3} \, \| \IFN^{-1/2} \Av \| \,. 
\label{rhDGtuGmusGU0a2}
\end{EQA}
Moreover, 
\begin{EQA}
    \Bigl| 2 \fn(\upsvn) - 2 \fn(\upsvs) - \| \IFN^{-1/2} \Av \|^{2} \Bigr|
    & \leq &
    \frac{\dltwu_{3}}{2} \, \| \DFN \, \IFN^{-1} \Av \|^{3} \, .
    \qquad
\label{3d3Af12DGttGa2}
\end{EQA}
\end{proposition}

\begin{proof}
W.l.o.g. assume \( \dmax = 1 \).
The first statement follows from Proposition~\ref{Pconcgeneric} with \( \amax = 3/4 \) because 
\nameref{LL3tref} ensures \eqref{dtb3u1DG2d3GPg} with \( \dltwbd(\upsv) = \dltwu_{3} \, \rr/2 \)
and \eqref{yxdhewndu7jwnjjuu} implies \eqref{rrm23r0ut3ua}.

As \( \nabla \fs(\upsvs) = 0 \), \nameref{LL3tref} implies for any \( \upsv \in \CA(\rrn) \)
\begin{EQA}
	\Bigl| 
		\fs(\upsvs) - \fs(\upsv) - \frac{1}{2} \| \IFN^{1/2} (\upsv - \upsvs) \|^{2} 
	\Bigr|
	& \leq &
	\frac{\dltwu_{3}}{6} \| \DFN (\upsv - \upsvs) \|^{3}
	\, .
	\qquad \quad
\label{d3GrGELGtsG122}
\end{EQA}
Further, 
\begin{EQA}[rcl]
	&& \nquad
	\fn(\upsv) - \fn(\upsvs) - \frac{1}{2} \| \IFN^{-1/2} \Av \|^{2}
	\\
	&=&
	\bigl\langle \upsv - \upsvs, \Av \bigr\rangle
	+ \fs(\upsv) - \fs(\upsvs) - \frac{1}{2} \| \IFN^{-1/2} \Av \|^{2} 
	\\
	&=&
	- \frac{1}{2} \bigl\| \IFN^{1/2} (\upsv - \upsvs) - \IFN^{-1/2} \Av \bigr\|^{2}
	+ \fs(\upsv) - \fs(\upsvs) + \frac{1}{2} \| \IFN^{1/2} (\upsv - \upsvs) \|^{2} .
	\qquad 
\label{12ELGuELusG2}
\end{EQA}
As \( \upsvn \in \CA(\rrn) \) and it maximizes \( \fn(\upsv) \), we derive by \eqref{d3GrGELGtsG122} and Lemma~\ref{Llin23}
with \( \AFN = \IFN^{1/2} \DFN^{-1} \) and \( \afv = \DFN \, \IFN^{-1} \Av \)
\begin{EQA}
	&& \nquad
	2 \fn(\upsvn) - 2 \fn(\upsvs) - \| \IFN^{-1/2} \Av \|^{2}
	=
	\max_{\upsv \in \CA(\rrn)} 
	\Bigl\{ 
		2 \fn(\upsv) - 2 \fn(\upsvs) - \| \IFN^{-1/2} \Av \|^{2} 
	\Bigr\}
	\\
	& \leq &
	\max_{\upsv \in \CA(\rrn)} 
	\Bigl\{ 
		- \bigl\| \IFN^{1/2} (\upsv - \upsvs) - \IFN^{-1/2} \Av \bigr\|^{2} 
		+ \frac{\dltwu_{3}}{3} \| \DFN (\upsv - \upsvs) \|^{3}
	\Bigr\} 
	\leq 
	\frac{\dltwu_{3}}{2} \, \| \DFN \, \IFN^{-1} \Av \|^{3} \, .
\label{d3G1212222B} 
\end{EQA}
Similarly 
\begin{EQA}
	&&  
	2 \fn(\upsvn) - 2 \fn(\upsvs) - \| \IFN^{-1/2} \Av \|^{2}
	\\
	&& \quad
	\geq 
	\max_{\upsv \in \CA(\rrn)} 
	\Bigl\{ 
		- \bigl\| \IFN^{1/2} (\upsv - \upsvs) - \IFN^{-1/2} \Av \bigr\|^{2} 
		- \frac{\dltwu_{3}}{3} \| \DFN (\upsv - \upsvs) \|^{3}
	\Bigr\}
	\geq 
	- \frac{\dltwu_{3} }{2} \, \| \DFN \, \IFN^{-1} \Av \|^{3} \, . 
	\qquad \quad
\label{fd3G122B2m2} 
\end{EQA}
This implies \eqref{3d3Af12DGttGa2}.
\end{proof}

\begin{lemma}
\label{Llin23}
Let \( \AFN \geq \Id \).
Fix some \( \rr \) and let \( \afv \in \R^{\dimp} \) satisfy \( (3/4) \rr \leq \| \afv \| \leq \rr \).
If \( \dltwu \, \rr \leq 1/3 \), then
\begin{EQA}
\label{hdcf6wyheuv76e34r35eycv}
	\max_{\| \uv \| \leq \rr} \Bigl( \frac{\dltwu}{3} \| \uv \|^{3} - (\uv - \afv)^{\T} \AFN (\uv - \afv) \Bigr)
	& \leq & 
	\frac{\dltwu}{2} \, \| \afv \|^{3} \, ,
	\\
	\min_{\| \uv \| \leq \rr} \Bigl( \frac{\dltwu}{3} \| \uv \|^{3} + (\uv - \afv)^{\T} \AFN (\uv - \afv) \Bigr)
	& \leq & 
	\frac{\dltwu}{2} \, \| \afv \|^{3} \, .
\label{hdcf6wyheuv76e34r35eycvm}
\end{EQA}
\end{lemma}

\begin{proof}
Replacing \( \| \uv \|^{3} \) with \( \rr \| \uv \|^{2} \) reduces 
the problem to quadratic programming. 
It holds with \( \rho \eqdef \dltwu \rr/3 \) and \( \afv_{\rho} \eqdef (\AFN - \rho \Id)^{-1} \AFN \afv \)
\begin{EQA}
	&& \nquad
	\frac{\dltwu}{3} \| \uv \|^{3} - (\uv - \afv)^{\T} \AFN (\uv - \afv)
	\leq 
	\frac{\dltwu \rr}{3} \| \uv \|^{2} - (\uv - \afv)^{\T} \AFN (\uv - \afv)
	\\
	&=&
	- \uv^{\T} \bigl( \AFN - \rho \Id \bigr) \uv + 2 \uv^{\T} \AFN \afv - \afv^{\T} \AFN \afv
	\\
	&=&
	- (\uv - \afv_{\rho})^{\T} (\AFN - \rho \Id) (\uv - \afv_{\rho}) 
	+ \afv_{\rho}^{\T} (\AFN - \rho \Id) \afv_{\rho} - \afv^{\T} \AFN \afv
	\\
	& \leq &
	\afv^{\T} \bigl\{ \AFN (\AFN - \rho \Id)^{-1} \AFN - \AFN \bigr\} \afv 
	=
	\rho \afv^{\T} \AFN (\AFN - \rho \Id)^{-1} \afv .
\label{f9kht446iffrtednftehy}
\end{EQA}
This implies in view of \( \AFN \geq \Id \), \( \rr \leq (4/3) \| \afv \| \), and \( \rho \leq 1/9 \)
\begin{EQA}
	&& \nquad
	\max_{\| \uv \| \leq \rr} \Bigl( \frac{\dltwu}{3} \| \uv \|^{3} - (\uv - \afv)^{\T} \AFN (\uv - \afv) \Bigr)
	\\
	& \leq &
	\frac{\rho}{1-\rho} \| \afv \|^{2}
	\leq 
	\frac{\dltwu \rr}{3(1-\rho)} \| \afv \|^{2}
	\leq 
	\frac{4\dltwu}{9 (1-\rho)} \| \afv \|^{3}
	\leq 
	\frac{\dltwu}{2} \| \afv \|^{3} \, ,
\label{gydw7guywbudvrgte7yruw}
\end{EQA}
and \eqref{hdcf6wyheuv76e34r35eycv} follows.
For \eqref{hdcf6wyheuv76e34r35eycvm} note that
\begin{EQA}
	&& \nquad
	\min_{\| \uv \| \leq \rr} \Bigl( \frac{\dltwu}{3} \| \uv \|^{3} + (\uv - \afv)^{\T} \AFN (\uv - \afv) \Bigr)
	\leq  
	\min_{\uv} \Bigl( \frac{\dltwu \rr}{3} \| \uv \|^{2} + (\uv - \afv)^{\T} \AFN (\uv - \afv) \Bigr)
	\\
	& \leq &
	\afv^{\T} \bigl\{ \AFN - \AFN (\AFN + \rho \Id)^{-1} \AFN \bigr\} \afv 
	=
	\rho \afv^{\T} \AFN (\AFN + \rho \Id)^{-1} \afv 
	\leq 
	\frac{\dltwu \rr}{3}  \| \afv \|^{2}
	\leq 
	\frac{4 \dltwu}{9}  \| \afv \|^{3} \, ,
\label{7jdfvy5433feugywdjw7krfu}
\end{EQA}
and \eqref{hdcf6wyheuv76e34r35eycvm} follows as well.
\end{proof}


\Subsection{Advanced approximation under locally uniform smoothness}
The bounds of Proposition~\ref{Pconcgeneric2} can be made more accurate if \( \fs \) follows
\nameref{LLsT3ref} and \nameref{LLsT4ref} and one can apply the Taylor expansion around any point close to \( \upsvs \).
In particular, the improved results do not involve the value \( \| \IFN^{-1/2} \Av \| \) which can be large or even 
infinite in some situation; see Section~\ref{Slinquadr}.


\begin{proposition}
\label{PFiWigeneric2}
Let \( \fs(\upsv) \) be a strongly concave function with \( \fs(\upsvs) = \max_{\upsv} \fs(\upsv) \)  
and \( \IFN = - \nabla^{2} \fs(\upsvs) \).
Assume \nameref{LLsT3ref} at \( \upsvs \) with \( \DFN^{2} \), \( \rrn \), and \( \dltwu_{3} \) such that 
\begin{EQA}[c]
	\DFN^{2} \leq \dmax^{2} \, \IFN ,
	\quad
	\rrn \geq \frac{3}{2} \| \DFN \, \IFN^{-1} \Av \| \, ,
	\quad
	\dmax^{2} \dltwu_{3} \| \DFN \, \IFN^{-1} \Av \| < \frac{4}{9} \, .
\label{8difiyfc54wrboer7bjfr}
\end{EQA}
Then \( \| \DFN (\upsvn - \upsvs) \| \leq (3/2) \| \DFN \, \IFN^{-1} \Av \| \) and moreover,
\begin{EQA}[rcl]
    \| \DFN^{-1} \IFN (\upsvn - \upsvs - \IFN^{-1} \Av) \|
    & \leq &
    \frac{3\dltwu_{3}}{4} \| \DFN \, \IFN^{-1} \Av \|^{2} 
	\, .
\label{DGttGtsGDGm13rGa2}
\end{EQA}
\end{proposition}

\begin{proof}
W.l.o.g. assume \( \dmax = 1 \).
If the function \( \fs \) is quadratic and concave with the maximum at \( \upsvs \) then the linearly perturbed function
\( \fn \) is also quadratic and concave with the maximum at \( \upsvr = \upsvs + \IFN^{-1} \Av \).
In general, the point \( \upsvr \) is not the maximizer of \( \fn \), however, it is very close to \( \upsvn \).
We use that \( \nabla \fs(\upsvs) = 0 \) 
and \( - \nabla^{2} \fs(\upsvs) = \IFN \).
The main step of the proof is given by the next lemma.

\begin{lemma}
\label{Ldltw4s}
Assume \nameref{LLsT3ref} at \( \upsv \).
With \( \DFN = \DFN(\upsv) \), let \( \UVz_{\rr} = \{ \uv \colon \| \DFN \uv \| \leq \rr \} \).
Then 
\begin{EQA}[ccl]
	\bigl\| \DFN^{-1} \bigl\{ \nabla \fs(\upsv + \uv) - \nabla \fs(\upsv) 
	- \langle \nabla^{2} \fs(\upsv), \uv \rangle \bigr\} 
	\bigr\|
	& \leq &
	\frac{\dltwu_{3}}{2} \, \| \DFN \uv \|^{2} \, ,
	\quad
	\uv \in \UVz_{\rr} \, .
\label{y6sdjsdy7erwmcuecuid}
\label{y6sdjsdy7erwmcuecuid2}
\end{EQA}
Also for all \( \uv, \uv_{1} \in \UVz_{\rr} \) 
\begin{EQ}[rcl]
	\bigl\| \DFN^{-1} \bigl\{ \nabla^{2} \fs(\upsv + \uv_{1}) - \nabla^{2} \fs(\upsv + \uv) \bigr\} \DFN^{-1}	\bigr\|
	& \leq &
	\dltwu_{3} \, \| \DFN (\uv_{1} - \uv) \|
\label{jhfuy7f7dfyedye663eh}
\end{EQ}
and
\begin{EQA}
	\bigl\| \DFN^{-1} \bigl\{ \nabla \fs(\upsv + \uv_{1}) - \nabla \fs(\upsv + \uv) - \nabla^{2} \fs(\upsv) (\uv_{1} - \uv) \bigr\}
	\bigr\|
	& \leq &
	\frac{3\dltwu_{3}}{2} \, \| \DFN (\uv_{1} - \uv) \|^{2} \, .
	\qquad
\label{6dcfujcu8ed8edsudyf5tre35}
\end{EQA}
Moreover, under \nameref{LLsT4ref}, for any \( \uv \in \UVz_{\rr} \),
\begin{EQA}
	\bigl\| \DFN^{-1} \bigl\{ \nabla \fs(\upsv + \uv) - \nabla \fs(\upsv) 
	- \langle \nabla^{2} \fs(\upsv), \uv \rangle 
	- \frac{1}{2} \langle \nabla^{3} \fs(\upsv), \uv^{\otimes 2} \rangle \bigr\} \bigr\|
	& \leq &
	\frac{\dltwu_{4}}{6} \, \| \DFN \uv \|^{3} \, .
	\qquad
\label{y6sdjsdy7erwmcuecuid4}
\end{EQA}
\end{lemma}

\begin{proof}
Denote
\begin{EQA}
	\Av
	& \eqdef &
	\nabla \fs(\upsv + \uv) - \nabla \fs(\upsv) 
	- \langle \nabla^{2} \fs(\upsv), \uv \rangle 
	\, .
\label{7dcjef6chjfer7v54etghf}
\end{EQA}
For any vector \( \wv \in \R^{\dimp} \), \nameref{LLsT3ref} and \eqref{jbuyfg773jgion94euyyfg} imply
\begin{EQA}
	\bigl| \langle \Av, \wv \rangle \bigr|
	& \leq &
	\frac{\dltwu_{3}}{2} \, \| \DFN \uv \|^{2} \, \| \DFN \wv \| .
\label{dfuedru8cm34e7edmvcghlo}
\end{EQA}
Therefore,
\begin{EQA}
	\| \DFN^{-1} \Av \|
	&=&
	\sup_{\| \wv \| = 1} \bigl| \langle \DFN^{-1} \Av,\wv \rangle \bigr|
	=
	\sup_{\| \wv \| = 1} \bigl| \langle \Av,\DFN^{-1} \wv \rangle \bigr|
	\leq 
	\frac{\dltwu_{3}}{2} \, \| \DFN \uv \|^{2} 
\label{udfjd6vhfwe36vneo}
\end{EQA}
which yields the first statement. 
For \eqref{y6sdjsdy7erwmcuecuid4}, apply
\begin{EQA}
	\Av
	& \eqdef &
	\nabla \fs(\upsv + \uv) - \nabla \fs(\upsv) 
	- \langle \nabla^{2} \fs(\upsv), \uv \rangle 
	- \frac{1}{2} \langle \nabla^{3} \fs(\upsv), \uv^{\otimes 2} \rangle
\label{t6dtwsghwesyfyghe322w2w}
\end{EQA} 
and use \nameref{LLsT4ref} and \eqref{jbuyfg773jgion94euyyfg4} instead of \nameref{LLsT3ref} and \eqref{jbuyfg773jgion94euyyfg}.
Further, with \( \Bv_{1} \eqdef \nabla^{2} \fs(\upsv + \uv_{1}) - \nabla^{2} \fs(\upsv + \uv) \) and \( \Delta = \uv_{1} - \uv \), 
by \nameref{LLsT3ref}, for any \( \wv \in \R^{\dimp} \) and some \( t \in [0,1] \),
\begin{EQA}
	&& \nquad
	\bigl| \langle \DFN^{-1} \bigl\{ \nabla^{2} \fs(\upsv + \uv_{1}) - \nabla^{2} \fs(\upsv + \uv) \bigr\} \, \DFN^{-1}, \wv^{\otimes 2} \rangle \bigr|
	=
	\bigl| \langle \Bv_{1}, (\DFN^{-1} \wv)^{\otimes 2} \rangle \bigr|
	\\
	& = &
	\bigl| \bigl\langle \nabla^{3} \fs(\upsv + \uv + t \Delta), \Delta \otimes (\DFN^{-1} \wv)^{\otimes 2} \bigr\rangle \bigr|
	\leq 
	\dltwu_{3} \| \DFN \Delta \| \, \| \wv \|^{2}
	\, .
\label{d76jnwef8j2egtftftyjr5}
\end{EQA}
This proves \eqref{jhfuy7f7dfyedye663eh}.
Similarly, 
for some \( t \in [0,1] \)
\begin{EQA}
	&& \nquad
	\bigl| \bigl\langle 
		\DFN^{-1} \bigl\{ \nabla \fs(\upsv + \uv_{1}) - \nabla \fs(\upsv + \uv) \bigr\} - \nabla^{2} \fs(\upsv + \uv) \Delta \bigr\},\wv  
	\bigr\rangle \bigr|
	\\
	&=&
	\frac{1}{2} \bigl| \bigl\langle \nabla^{3} \fs(\upsv + \uv + t \Delta), \Delta \otimes \Delta \otimes \DFN^{-1} \wv \bigr\rangle
	\bigr|
	\leq 
	\frac{\dltwu_{3}}{2} \| \DFN \Delta \|^{2} \, \| \wv \|
\label{ufhfyt3hbgvbigy4jfiu}
\end{EQA}
and with \( \Bv = \nabla^{2} \fs(\upsv + \uv) - \nabla^{2} \fs(\upsv) \), by \eqref{jhfuy7f7dfyedye663eh},
\begin{EQA}
	\bigl\| \DFN^{-1} \Bv \Delta \bigr\| 
	& \leq &
	\| \DFN^{-1} \Bv \, \DFN^{-1} \| \,\, \| \DFN \Delta \| 
	\leq 
	\dltwu_{3} \| \DFN \Delta \|^{2} \, .
\label{u8ifke0gfjw23gfsd4gy}
\end{EQA}
This completes the proof of \eqref{6dcfujcu8ed8edsudyf5tre35}.
\end{proof}

Now \eqref{y6sdjsdy7erwmcuecuid} of Lemma~\ref{Ldltw4s} 
yields
\begin{EQA}
	\bigl\| \DFN^{-1} \nabla \fn(\upsvr) \bigr\|
	&=&
	\bigl\| \DFN^{-1} \{ \nabla \fs(\upsvs + \IFN^{-1} \Av) - \nabla \fs(\upsvs) + \Av \} \bigr\| 
	\leq
	\frac{\dltwu_{3}}{2} \| \DFN \, \IFN^{-1} \Av \|^{2} \, .
	\qquad
\label{stedsyteuhyfnmgu3}
\end{EQA}
As \( \| \DFN \, \IFN^{-1} \Av \| \leq 2 \rrn/3 \), condition \nameref{LLsT3ref} can be applied 
in the \( \rrn/3 \)-vicinity of \( \upsvr \).
Fix any \( \upsv \) with \( \| \DFN (\upsv - \upsvr) \| \leq \rrn/3 \) and define 
\( \Delta = \upsv - \upsvr \).
By \eqref{6dcfujcu8ed8edsudyf5tre35} of Lemma~\ref{Ldltw4s}
\begin{EQA}
	&& \nquad
	\bigl\| \DFN^{-1} \{ \nabla \fn(\upsv) - \nabla \fn(\upsvr) + \IFN \Delta \} \bigr\|
	=
	\bigl\| \DFN^{-1} \{ \nabla \fs(\upsv) - \nabla \fs(\upsvr) + \IFN \Delta \} \bigr\|
	\leq 
	\frac{3 \dltwu_{3}}{2} \| \DFN \Delta \|^{2} \, .
\label{7djdxhes5ewtwee6e6eed}
\end{EQA}
In particular, this and \eqref{stedsyteuhyfnmgu3} yield
\begin{EQA}
	\bigl\| \DFN^{-1} \{ \nabla \fn(\upsvr + \Delta) + \IFN \Delta \} \bigr\|
	& \leq &
	2 \dltwu_{3} \| \DFN \Delta \|^{2} \, .
\label{7hewjfuv7bh7tur4jwsfycn}
\end{EQA}
For any \( \uv \) with \( \| \uv \| = 1 \), this implies
\begin{EQA}
	\bigl| \bigl\langle \nabla \fn(\upsvr + \Delta) + \IFN \Delta , \DFN^{-1} \uv \bigr\rangle \bigr|
	& \leq &
	2 \dltwu_{3} \| \DFN \Delta \|^{2} \, .
\label{hf9jmw2f7vhehe46fghfnwh}
\end{EQA}
Suppose now that \( \| \DFN \Delta \| = \rrn/3 \)
and consider the function \( h(t) = \fn(\upsvr + t \Delta) \).
Then  
\( h'(t) = \langle \nabla \fn(\upsvr + t \Delta), \Delta \rangle \)
and \eqref{hf9jmw2f7vhehe46fghfnwh} implies with \( \uv = \DFN \Delta/\| \DFN \Delta \| \)
\begin{EQA}
	\bigl| \langle \nabla \fn(\upsvr + \Delta), \Delta \rangle + \| \IFN^{1/2} \Delta \|^{2} \bigr|
	& \leq &
	2 \dltwu_{3} \| \DFN \Delta \|^{3} \, .
\label{ufudhdyf5d53egru7e4u3}
\end{EQA}
As \( \IFN \geq \DFN^{2} \), this yields
\begin{EQA}
	h'(1)
	& \leq &
	2 \dltwu_{3} \| \DFN \Delta \|^{3} - \| \DFN \Delta \|^{2} .
\label{67vhyfdhcyeyghevy2hegtie3}
\end{EQA}
Similarly, \eqref{stedsyteuhyfnmgu3} yields by \( \| \DFN \, \IFN^{-1} \Av \| = 2\rrn/3 \)
\begin{EQA}
	|h'(0)|
	=
	\bigl| \langle \nabla \fn(\upsvr), \Delta \rangle \bigr|
	& \leq &
	\frac{\dltwu_{3}}{2} \| \DFN \, \IFN^{-1} \Av \|^{2} \, \| \DFN \Delta \| 
	=
	\frac{2\dltwu_{3}}{9} \, \rr^{2} \, \| \DFN \Delta \|\, .
\label{yhev7h3bfgreewewweddh2}
\end{EQA}
Concavity of \( \fn(\cdot) \) ensures that \( t^{*} = \argmax_{t} h(t) \) satisfies \( |t^{*}| \leq 1 \)
provided that 
\begin{EQA}[c]
	h'(1) < - |h'(0)| ,
	\qquad
	h'(-1) < |h'(0)|.
\label{8gafh8fe5w8gbfuvyteyvjv}
\end{EQA}
Due to \eqref{67vhyfdhcyeyghevy2hegtie3}, \eqref{yhev7h3bfgreewewweddh2}, and
\( \| \DFN \Delta \| = \rrn/3 \),  the latter condition reads
\begin{EQA}[c]
	\frac{2\dltwu_{3}}{9} \rr^{2} \, \| \DFN \Delta \| + 2 \dltwu_{3} \| \DFN \Delta \|^{3} - \| \DFN \Delta \|^{2}
	=
	\| \DFN \Delta \| \, \rrn \Bigl( \frac{2\dltwu_{3} \, \rrn}{9} + \frac{2 \dltwu_{3} \, \rrn}{9} - \frac{1}{3} \Bigr)
	< 
	0 .
\label{thewuvf7ehehctdebenvjyw}
\end{EQA}
which is fulfilled because of \( \dltwu_{3} \| \DFN \, \IFN^{-1} \Av \| \leq 4/9 \) 
and \( \| \DFN \, \IFN^{-1} \Av \| = 2\rrn/3 \).
We summarize that \( \upsvn = \argmax_{\upsv} \fn(\upsv) \) satisfies \( \| \DFN \, (\upsvn - \upsvr) \| \leq \rrn/3 \)
while \( \| \DFN (\upsvr - \upsvs) \| = \| \DFN \, \IFN^{-1} \Av \| = 2 \rrn/3 \).
Therefore, \( \| \DFN (\upsvn - \upsvs) \| \leq \rrn \).
This allows us to use \nameref{LLsT3ref} at this point for establishing \eqref{DGttGtsGDGm13rGa2}.
By definition \( \nabla \fn(\upsvn) = 0 \) and hence,
\begin{EQA}
	\| \DFN^{-1} \{ \nabla \fn(\upsvs + \IFN^{-1} \Av) - \nabla \fn(\upsvn) \} \|
	& \leq &
	\frac{\dltwu_{3}}{2} \| \DFN \, \IFN^{-1} \Av \|^{2} \, .
\label{fhy5345etvty46dgw3}
\end{EQA}
By \eqref{6dcfujcu8ed8edsudyf5tre35} of Lemma~\ref{Ldltw4s}, it holds with \( \Delta = \upsvs + \IFN^{-1} \Av - \upsvn \)
\begin{EQA}
	\bigl\| \DFN^{-1} \{ \nabla \fn(\upsvs + \IFN^{-1} \Av) - \nabla \fn(\upsvn) - \nabla^{2} \fn(\upsvs) \Delta \} \bigr\| 
	& \leq &
	\frac{3 \dltwu_{3}}{2} \| \DFN \Delta \|^{2} \, .
\label{ygtdtdt55636ytv2wg3}
\end{EQA}
Combining with \eqref{fhy5345etvty46dgw3} yields 
\begin{EQA}
	\| \DFN^{-1} \IFN \Delta \|
	& \leq &
	\frac{3 \dltwu_{3}}{2} \| \DFN \Delta \|^{2} + \frac{\dltwu_{3}}{2} \| \DFN \, \IFN^{-1} \Av \|^{2} 
	\leq 
	\frac{3 \dltwu_{3}}{2} \| \DFN^{-1} \IFN \Delta \|^{2} + \frac{\dltwu_{3}}{2} \| \DFN \, \IFN^{-1} \Av \|^{2} \, . 
\label{f7f7fv7f66e6ehb3wyc3}
\end{EQA}
As \( 2x \leq \alpha x^{2} + \beta \) with \( \alpha = 3 \dltwu_{3} \), \( \beta = \dltwu_{3} \| \DFN \, \IFN^{-1} \Av \|^{2} \), and
\( x = \| \DFN^{-1} \IFN \Delta \| \in (0,1/\alpha) \) implies \( x \leq \beta/(2 - \alpha\beta) \),
this yields
\begin{EQA}
	\| \DFN^{-1} \IFN (\upsvn - \upsvs - \IFN^{-1} \Av) \|
	& \leq &
	\frac{\dltwu_{3}}{2 - 3 \dltwu_{3}^{2} \| \DFN \, \IFN^{-1} \Av \|^{2}} \| \DFN \, \IFN^{-1} \Av \|^{2} 
\label{fiufu7df56rgyhvnghbvt3}
\end{EQA}
and \eqref{DGttGtsGDGm13rGa2} follows by 
\( \dltwu_{3} \| \DFN \, \IFN^{-1} \Av \| \leq 4/9 \).
\end{proof}

\begin{remark}
\label{Rbiasgeneric}
As in Remark~\ref{dtb3u1DG2d3GPg}, the roles of \( \fs \) and \( \fn \) can be exchanged.
In particular, \eqref{DGttGtsGDGm13rGa2} applies with \( \IFN = \IFN(\upsvn) \) provided that 
\nameref{LLsT3ref} is also fulfilled at \( \upsvn \).
\end{remark}

\noindent
If \( \fs \) is fourth-order smooth and \nameref{LLsT4ref} holds then expansion \eqref{DGttGtsGDGm13rGa2} can further be refined.  

\begin{proposition}
\label{Pconcgeneric4}
Let \( \fs(\upsv) \) be a strongly concave function with \( \fs(\upsvs) = \max_{\upsv} \fs(\upsv) \)  
and \( \IFN = - \nabla^{2} \fs(\upsvs) \), and let
\( \fs(\upsv) \) follow \nameref{LLsT3ref} and \nameref{LLsT4ref} with some 
\( \DFN^{2} \), \( \dltwu_{3} \), \( \dltwu_{4} \), and \( \rrn \) satisfying 
\begin{EQA}[c]
	\DFN^{2} \leq \dmax^{2} \, \IFN \, ,
	\;\;
	\rrn = \frac{3}{2} \| \DFN \, \IFN^{-1} \Av \| \, ,
	\;\;
	\dmax^{2} \dltwu_{3} \| \DFN \, \IFN^{-1} \Av \| < \frac{4}{9} \, ,
	\;\;
	\dmax^{2} \dltwu_{4} \| \DFN \, \IFN^{-1} \Av \|^{2} < \frac{1}{3} \, .
	\qquad
\label{8difiyfc54wrbosT4}
\end{EQA}
Let \( \fn(\upsv) \) fulfill \eqref{4hbh8njoelvt6jwgf09} with some vector \( \Av \) and \( \fn(\upsvn) = \max_{\upsv} \fn(\upsv) \).
Then \( \| \DFN (\upsvn - \upsvs) \| \leq (3/2) \| \DFN \, \IFN^{-1} \Av \| \). 
Further, define
\begin{EQA}
	\avn 
	&=&
	\IFN^{-1} \{ \Av + \nabla \Tens(\IFN^{-1} \Av) \} \, ,
\label{8vfjvr43223efryfuweef}
\end{EQA}
where \( \Tens(\uv) = \frac{1}{6} \langle \nabla^{3} \fs(\upsvs), \uv^{\otimes 3} \rangle \) for \( \uv \in \R^{\dimp} \).
Then
\begin{EQ}[rcl]
    \| \DFN^{-1} \IFN (\upsvn - \upsvs - \avn) \|
    & \leq &
    (\dltwu_{4}/2 + \dmax^{2} \dltwu_{3}^{2}) \, \| \DFN \, \IFN^{-1} \Av \|^{3} \, .
\label{DGttGtsGDGm13rGa4}
\end{EQ}
Also
\begin{EQA}
    && \nquad
    \Bigl| \fn(\upsvn) - \fn(\upsvs) - \frac{1}{2} \| \IFN^{-1/2} \Av \|^{2} - \Tens(\IFN^{-1} \Av) \Bigr|
    \\
    & \leq &
    \frac{\dltwu_{4} + 4 \dmax^{2} \dltwu_{3}^{2}}{8} \| \DFN \, \IFN^{-1} \Av \|^{4} 
    + \frac{\dmax^{2} \, (\dltwu_{4} + 2 \dmax^{2} \dltwu_{3}^{2})^{2} }{4} \, \| \DFN \, \IFN^{-1} \Av \|^{6} \, 
    \qquad
\label{3d3Af12DGttGa4}
\end{EQA}
and 
\begin{EQA}
	\bigl| \Tens(\IFN^{-1} \Av) \bigr|
	& \leq &
	\frac{\dltwu_{3}}{6} \| \DFN \, \IFN^{-1} \Av \|^{3} \, .
\label{u7jcc7e45hfiobvioeye6yhy}
\end{EQA}
\end{proposition}

\begin{proof}
W.l.o.g. assume \( \dmax = 1 \) and \( \upsvs = 0 \).
Proposition~\ref{PFiWigeneric2} yields \eqref{DGttGtsGDGm13rGa2}.
By \nameref{LLsT3ref} 
\begin{EQA}
	&& \nquad
	\| \DFN^{-1} \, \IFN (\avn - \IFN^{-1} \Av) \|
	=
	\| \DFN^{-1} \, \nabla \Tens(\IFN^{-1} \Av) \|
	\\
	&=&
	\sup_{\| \uv \| = 1} 3 \bigl| \langle \Tens, \IFN^{-1} \Av \otimes \IFN^{-1} \Av \otimes \DFN^{-1} \uv \rangle \bigr|
	\leq 
	\frac{\dltwu_{3}}{2} \| \DFN \, \IFN^{-1} \Av \|^{2} \, .
	\qquad
\label{bhvfwfdsdxexsdwsvwe33a}
\end{EQA}
As \( \DFN^{-1} \, \IFN \geq \IFN^{1/2} \geq \DFN \), this implies by \( \dltwu_{3} \| \DFN \, \IFN^{-1} \Av \| \leq 4/9 \)
\begin{EQA}[rcl]
	\| \DFN \avn \|
	& \leq &
	\| \DFN \, \IFN^{-1} \Av \| + \| \DFN \, \IFN^{-1}  \, \nabla \Tens(\IFN^{-1} \Av) \|
	\\
	& \leq &
	\Bigl( 1 + \frac{\dltwu_{3}}{2} \| \DFN \, \IFN^{-1} \Av \| \Bigr) \| \DFN \, \IFN^{-1} \Av \|
	\leq 
	\frac{11}{9} \, \| \DFN \, \IFN^{-1} \Av \| 
	\qquad
\label{iuvchycvf6e64rygh322}
\end{EQA}
and
\begin{EQA}
	\| \IFN^{1/2} \avn - \IFN^{-1/2} \Av \|
	& \leq &
	\frac{\dltwu_{3}}{2} \| \DFN \, \IFN^{-1} \Av \|^{2} \, .
\label{iuvchycvf6e64ryghF}
\end{EQA}
Next, again by \nameref{LLsT3ref}, for any \( \wv \)
\begin{EQA}
	\| \DFN^{-1} \, \nabla^{2} \Tens(\wv) \, \DFN^{-1} \|
	&=&
	\sup_{\| \uv \| = 1} 6 \bigl| \langle \Tens, \wv \otimes (\DFN^{-1} \uv)^{\otimes 2} \rangle \bigr|
	\leq 
	\dltwu_{3} \| \DFN \wv \| \, .
\label{dunj3df7y76e4hf743j}
\end{EQA}
The tensor \( \nabla^{2} \Tens(\uv) \) is linear in \( \uv \), hence
\begin{EQA}
	&& \nquad
	\sup_{t \in [0,1]} \| \DFN^{-1} \, \nabla^{2} \Tens(t \avn + (1-t) \IFN^{-1} \Av) \, \DFN^{-1} \| 
	\\
	&=&
	\max\{ \| \DFN^{-1} \, \nabla^{2} \Tens(\IFN^{-1} \Av) \, \DFN^{-1} \|, \| \DFN^{-1} \nabla^{2} \Tens(\avn) \DFN^{-1} \| \}
	\leq 
	\dltwu_{3} \, \max\{ \| \DFN \, \IFN^{-1} \Av \|, \| \DFN \avn \| \} \, .
\label{huyd76hj3fyt7yfj4e}
\end{EQA}
Based on \eqref{iuvchycvf6e64rygh322}, assume \( \| \DFN \, \IFN^{-1} \Av \| \leq \| \DFN \avn \| \leq (11/9) \| \DFN \, \IFN^{-1} \Av \| \).
Then \eqref{bhvfwfdsdxexsdwsvwe33a} yield
\begin{EQA}
	&& \nquad
	\| \DFN^{-1} \nabla \Tens(\avn) - \DFN^{-1} \nabla \Tens(\IFN^{-1} \Av) \|
	\\
	& \leq &
	\sup_{t \in [0,1]} \| \DFN^{-1} \, \nabla^{2} \Tens(t \avn + (1-t) \IFN^{-1} \Av) \, \DFN^{-1} \| \,\, 
	\| \DFN \, \IFN^{-1} (\avn - \IFN^{-1} \Av) \|
	\\
	& \leq &
	\frac{\dltwu_{3}^{2}}{2} \, \| \DFN \, \IFN^{-1} \Av \|^{2} \, \| \DFN \avn \|
	\leq 
	\frac{2\dltwu_{3}^{2}}{3} \, \| \DFN \, \IFN^{-1} \Av \|^{3}\, .
\label{ukjikio3278eu7grt64rsa}
\end{EQA}
Further, \( - \nabla^{2} \fs(0) = \IFN \), 
\( \nabla \Tens(\avn) = \frac{1}{2} \langle \nabla^{3} \fs(0),\avn \otimes \avn \rangle \).
By \eqref{y6sdjsdy7erwmcuecuid4} of Lemma~\ref{Ldltw4s} and \eqref{iuvchycvf6e64rygh322} 
\begin{EQA}
	\bigl\| \DFN^{-1} \{ \nabla \fs(\avn) + \IFN \avn - \nabla \Tens(\avn) \} \bigr\| 
	& \leq &
	\frac{\dltwu_{4}}{6} \| \DFN \avn \|^{3} 
	\leq 
	\frac{(11/9)^{3}\dltwu_{4}}{6} \| \DFN \, \IFN^{-1} \Av \|^{3}
	\leq 
	\frac{\dltwu_{4}}{3} \| \DFN \, \IFN^{-1} \Av \|^{3} \, .
\label{stedsyteuhyfnmgu4}
\end{EQA}
Next we bound \( \bigl\| \DFN^{-1} \{ \nabla \fn(\avn) - \nabla \fn(\upsvn) \} \bigr\| \).
As \( \nabla \fn(\upsvn) = 0 \), \eqref{4hbh8njoelvt6jwgf09} and \eqref{8vfjvr43223efryfuweef} imply 
\begin{EQA}[rcl]
	&& \nquad
	\bigl\| \DFN^{-1} \{ \nabla \fn(\avn) - \nabla \fn(\upsvn) \} \bigr\|
	=
	\bigl\| \DFN^{-1} \nabla \fn(\avn) \bigr\|
	=
	\bigl\| \DFN^{-1} \{ \nabla \fn(\avn) + \IFN \avn - \nabla \Tens(\IFN^{-1} \Av) - \Av \} \bigr\|
	\\
	& \leq &
	\bigl\| \DFN^{-1} \{ \nabla \fs(\avn) + \IFN \avn - \nabla \Tens(\avn) \} \bigr\| 
	+ \| \DFN^{-1} \{ \nabla \Tens(\avn) - \nabla \Tens(\IFN^{-1} \Av) \} \|
	\leq 
	\err_{1} \, ,
\label{fhy5345etvty46dgw35w3a}
\end{EQA}
where 
\begin{EQA}[c]
	\err_{1}
	\eqdef 
	\frac{\dltwu_{4} + 2 \dltwu_{3}^{2}}{3} \| \DFN \, \IFN^{-1} \Av \|^{3} 
\label{fhjvmvgt4judjewtrwgd}
\end{EQA}
and by \eqref{8difiyfc54wrbosT4} 
\begin{EQA}
	3 \dltwu_{3} \, \err_{1}
	&=&
	\dltwu_{3} \| \DFN \, \IFN^{-1} \Av \| \; \dltwu_{4} \| \DFN \, \IFN^{-1} \Av \|^{2} 
	+ 2 \dltwu_{3}^{3} \| \DFN \, \IFN^{-1} \Av \|^{3}
	< 
	\frac{1}{3} \, .
\label{fhjvmvgt4judjewtrwg}
\end{EQA}
Further, \( \nabla^{2} \fn(0) = \nabla^{2} \fs(0) = - \IFN \), 
and \eqref{6dcfujcu8ed8edsudyf5tre35} of Lemma~\ref{Ldltw4s} implies 
\begin{EQA}
	&& \nquad
	\bigl\| \DFN^{-1} \{ \nabla \fn(\avn) - \nabla \fn(\upsvn) + \IFN (\avn - \upsvn) \} \bigr\|
	\\
	&=&
	\bigl\| \DFN^{-1} \{ \nabla \fs(\avn) - \nabla \fs(\upsvn) + \IFN (\avn - \upsvn) \} \bigr\|
	\leq 
	\frac{3\dltwu_{3}}{2} \| \DFN (\avn - \upsvn) \|^{2} .
\label{ygtdtdt55636ytv2wgdfa}
\end{EQA}
Combining with \eqref{fhy5345etvty46dgw35w3a} yields in view of \( \DFN^{2} \leq \IFN \)
\begin{EQA}
	\| \DFN^{-1} \IFN (\avn - \upsvn) \|
	& \leq &
	\frac{3\dltwu_{3}}{2} \| \DFN (\avn - \upsvn) \|^{2} + \err_{1}
	\leq 
	\frac{3\dltwu_{3}}{2} \| \DFN^{-1} \IFN (\avn - \upsvn) \|^{2} + \err_{1} \, .
\label{ufjvfchyeghdftdf67dejh}
\end{EQA}
As \( 2x \leq \alpha x^{2} + \beta \) with \( \alpha = 3 \dltwu_{3} \), \( \beta = 2 \err_{1} \), and
\( x \in (0,1/\alpha) \) implies \( x \leq \beta/(2 - \alpha\beta) \), 
we conclude by \eqref{fhjvmvgt4judjewtrwg}
\begin{EQA}
	\| \DFN^{-1} \IFN (\avn - \upsvn) \|
	& \leq &
	\frac{\err_{1}}{1 - 3 \dltwu_{3} \,\err_{1}}  
	\leq 
	\frac{\dltwu_{4} + 2 \dltwu_{3}^{2}}{2} \| \DFN \, \IFN^{-1} \Av \|^{3} \, ,
\label{fiufu7df56rgyhvnghbvtwea}
\end{EQA}
and \eqref{DGttGtsGDGm13rGa4} follows.

Next we bound \( \fn(\upsvn) - \fn(0) = \fn(\upsvn) - \fn(\avn) + \fn(\avn) - \fn(0) \).
By \eqref{iuvchycvf6e64ryghF} and \( \DFN^{2} \leq \IFN \)
\begin{EQA}
	\frac{1}{2} \| \IFN^{-1/2} \Av \|^{2} - \langle \Av, \avn \rangle + \frac{1}{2} \| \IFN^{1/2} \avn \|^{2}
	&=&
	\frac{1}{2} \| \IFN^{1/2} \avn - \IFN^{-1/2} \Av \|^{2}
	\leq 
	\frac{\dltwu_{3}^{2}}{8} \| \DFN \, \IFN^{-1} \Av \|^{4} \, .
\label{7djdfuyv7gh7tur5tuy}
\end{EQA}
This together with \( \nabla \fs(0) = 0 \), \( - \nabla^{2} \fs(0) = \IFN \geq \DFN^{2} \), \nameref{LLsT4ref}, and \eqref{iuvchycvf6e64rygh322} implies
\begin{EQA}
	&& \nquad
	\Bigl| \fn(\avn) - \fn(0) - \frac{1}{2} \| \IFN^{-1/2} \Av \|^{2} - \Tens(\avn) \Bigr|
	\\
	&=&
	\Bigl| \fs(\avn) - \fs(0) + \langle \Av, \avn \rangle - \frac{1}{2} \| \IFN^{-1/2} \Av \|^{2} - \Tens(\avn) \Bigr|
	\\
	& \leq &
	\Bigl| \fs(\avn) - \fs(0) + \frac{1}{2} \| \IFN^{1/2} \avn \|^{2} - \Tens(\avn) \Bigr| 
	+ \frac{\dltwu_{3}^{2}}{8} \| \DFN \, \IFN^{-1} \Av \|^{4}
	\\
	& \leq &
	\frac{\dltwu_{4}}{24} \| \DFN \, \avn \|^{4} + \frac{\dltwu_{3}^{2}}{8} \| \DFN \, \IFN^{-1} \Av \|^{4} 
	\leq 
	\Bigl( \frac{\dltwu_{4}}{10} + \frac{\dltwu_{3}^{2}}{8} \Bigr) \| \DFN \, \IFN^{-1} \Av \|^{4} \, .
\label{ghd6w2hehfyttet2hbdgfwa}
\end{EQA}
Further, by \( \nabla \fn(\upsvn) = 0 \) and \( \nabla^{2} \fn(\cdot) \equiv \nabla^{2} \fs(\cdot) \), it holds for some 
\( \upsv \in [\avn,\upsvn] \) 
\begin{EQA}
	2 \bigl| \fn(\avn) - \fn(\upsvn) \bigr|
	& = &
	\bigl| \langle \nabla^{2} \fs(\upsv) , (\avn - \upsvn)^{\otimes 2} \rangle \bigr| \, .
\label{dy6eh3hft636yhg3ffa}
\end{EQA}
The use of \( - \nabla^{2} \fs(0) = \IFN \geq \DFN^{2} \) and \eqref{jhfuy7f7dfyedye663eh} of Lemma~\ref{Ldltw4s} yields by 
\( \| \DFN \upsv \| \leq \rrn = \frac{3}{2} \| \DFN \, \IFN^{-1} \Av \| \), \( \dltwu_{3} \| \DFN \, \IFN^{-1} \Av \| < \frac{4}{9} \), and \eqref{fiufu7df56rgyhvnghbvtwea}
\begin{EQA}
	&& \nquad
	2 \bigl| \fn(\avn) - \fn(\upsvn) \bigr|
	\leq 
	\| \IFN^{1/2} (\avn - \upsvn) \|^{2}
	+ \bigl| \bigl\langle \nabla^{2} \fs(\upsv) - \nabla^{2} \fs(0), (\avn - \upsvn)^{\otimes 2} \bigr\rangle \bigr|	
	\\
	& \leq &
	(1 + \dltwu_{3} \rr) \| \IFN^{1/2} (\avn - \upsvn) \|^{2} 
	\leq 
	\frac{(5/3) (\dltwu_{4} + 2 \dltwu_{3}^{2})^{2}}{4} \, \| \DFN \, \IFN^{-1} \Av \|^{6}
	\, .
\label{f8mn4erf6ffyruyn4e3u8fhw3g}
\end{EQA}
Moreover, it holds with \( \Delta \eqdef \IFN^{-1} \nabla \Tens(\IFN^{-1} \Av) \)
for some \( t \in [0,1] \)
\begin{EQA}
	&& \nquad
	\bigl| \Tens(\avn) - \Tens(\IFN^{-1} \Av) \bigr|
	=
	\bigl| \Tens(\IFN^{-1} \Av + \Delta) - \Tens(\IFN^{-1} \Av) \bigr|
	=
	\bigl| \bigl\langle \nabla \Tens(\IFN^{-1} \Av + t \Delta ), \Delta \bigr\rangle \bigr| \,\, 
	\\
	& \leq &
	\frac{\dltwu_{3}}{2} \| \DFN (\IFN^{-1} \Av + t \Delta) \|^{2} \, \| \DFN \Delta \|
	=
	\frac{\dltwu_{3}}{2} \| \DFN \, \IFN^{-1} \Av + t \, \DFN \Delta) \|^{2} \, \| \DFN \Delta \|
	\, .
\label{7jc6hw3f6hw3fuvne8dgwx}
\end{EQA}
As in \eqref{bhvfwfdsdxexsdwsvwe33a} 
\( \| \DFN \Delta \| \leq \| \DFN^{-1} \nabla \Tens(\IFN^{-1} \Av) \| \leq (\dltwu_{3}/2) \| \DFN \, \IFN^{-1} \Av \|^{2} \), and
by \( \dltwu_{3} \| \DFN \, \IFN^{-1} \Av \| \leq 1/2 \)
\begin{EQA}
	&& \nquad
	\bigl| \Tens(\avn) - \Tens(\IFN^{-1} \Av) \bigr|
	\leq 
	\frac{(5/4)^{2} \dltwu_{3}^{2} }{4} \| \DFN \, \IFN^{-1} \Av \|^{4} \, .
\label{dc7hhbejrfweugdf7weuhduw}
\end{EQA}
Summing up the obtained bounds yields \eqref{3d3Af12DGttGa4}.
\eqref{u7jcc7e45hfiobvioeye6yhy} follows from \nameref{LLsT3ref}.
\end{proof}

\Subsection{Quadratic penalization}
\label{Slinquadr}
Here we discuss the case when \( \fn(\upsv) - \fs(\upsv) \) is quadratic.
The general case can be reduced to the situation with \( \fn(\upsv) = \fs(\upsv) - \| \GP \upsv \|^{2}/2 \).
To make the dependence of \( \GP \) more explicit, denote 
\( \fG(\upsv) \eqdef \fs(\upsv) - \| \GP \upsv \|^{2}/2 \),
\begin{EQA}
	\upsvs 
	&=& 
	\argmax_{\upsv} \fs(\upsv),
	\quad
	\upsvs_{\GP} = \argmax_{\upsv} \fG(\upsv) 
	=
	\argmax_{\upsv} \bigl\{ \fs(\upsv) - \| \GP \upsv \|^{2}/2 \bigr\}.
	\qquad
\label{8cvkfc9fujf6jnmcer4cd}
\end{EQA}
We study the bias \( \upsvs_{\GP} - \upsvs \) induced by this penalization.
To get some intuition, consider first the case of a quadratic function \( \fs(\upsv) \).

\begin{lemma}
\label{Lbiasquadgen}
Let \( \fs(\upsv) \) be quadratic with \( \IFN \equiv - \nabla^{2} \fs(\upsv) \) and \( \IFN_{\GP} = \IFN + \GP^{2} \).
Then 
\begin{EQA}[rcl]
	\upsvs_{\GP} - \upsvs
	&=&
	- \IFN_{\GP}^{-1} \GP^{2} \upsvs,
	\\
	\fG(\upsvs_{\GP}) - \fG(\upsvs)
	&=&
	\frac{1}{2} \| \IFN_{\GP}^{-1/2} \GP^{2} \upsvs \|^{2} \, .
\label{kjfydf554dfwertdgwdf}
\end{EQA}
\end{lemma} 

\begin{proof}
Quadraticity of \( \fs(\upsv) \) implies quadraticity of \( \fG(\upsv) \) with \( \nabla^{2} \fG(\upsv) \equiv - \IFN_{\GP} \) and
\begin{EQA}
	\nabla \fG(\upsvs_{\GP}) - \nabla \fG(\upsvs)
	&=&
	- \IFN_{\GP} \, (\upsvs_{\GP} - \upsvs) .
\label{dcudydye67e6dy3wujhdqu}
\end{EQA}
Further, \( \nabla \fs(\upsvs) = 0 \) yielding \( \nabla \fG(\upsvs) = - \GP^{2} \upsvs \).
Together with \( \nabla \fG(\upsvs_{\GP}) = 0 \), this implies
\( \upsvs_{\GP} - \upsvs = - \IFN_{\GP}^{-1} \GP^{2} \upsvs \).
The Taylor expansion of \( \fG \) at \( \upsvs_{\GP} \) yields 
\begin{EQA}
	\fG(\upsvs) - \fG(\upsvs_{\GP})
	&=&
	- \frac{1}{2} \| \IFN_{\GP}^{1/2} (\upsvs - \upsvs_{\GP}) \|^{2}
	=
	- \frac{1}{2} \| \IFN_{\GP}^{-1/2} \GP^{2} \upsvs \|^{2} 
\label{8chuctc44wckvcuedjequ}
\end{EQA}
and the assertion follows.
\end{proof}

Now we turn to the general case with \( \fs \) satisfying \nameref{LLsT3ref}.

\begin{proposition}
\label{Pbiasgeneric} 
Let \( \fG(\upsv) = \fs(\upsv) - \| \GP \upsv \|^{2}/2 \) be concave
and follow \nameref{LLsT3ref} with some \( \DFN^{2} \), \( \dltwu_{3} \), and \( \rrn \) satisfying for \( \dmax > 0 \)
\begin{EQA}[c]
	\DFN^{2} \leq \dmax^{2} \, \IFN_{\GP} \, ,
	\qquad 
	\rrn \geq 3 \bias_{\GP}/2 \, ,
	\qquad
	\dmax^{2} \dltwu_{3} \, \bias_{\GP} < 4/9 ,
\label{7fjgjgvuvy44erd52f}
\end{EQA}
where 
\begin{EQA}
	\bias_{\GP}
	&=& 
	\| \DFN \, \IFN_{\GP}^{-1} \GP^{2} \upsvs \| \, .
\label{fd9dfhy4ye6fuydfrerf}
\end{EQA} 
Then 
\begin{EQA}[c]
	\| \DFN (\upsvs_{\GP} - \upsvs) \| \leq 3 \bias_{\GP}/2 .
\label{odf6fdyr6e4deuewjug}
\end{EQA}
Moreover, 
\begin{EQA}[rcl]
	\bigl\| \DFN^{-1} \IFN_{\GP} (\upsvs_{\GP} - \upsvs + \IFN_{\GP}^{-1} \GP^{2} \upsvs) \bigr\|
	& \leq &
	\frac{3\dltwu_{3}}{4} \, \bias_{\GP}^{2}
	\, ,
\label{11ma3eaelDebbgen}
	\\
	\Bigl| 2 \fG(\upsvs_{\GP}) - 2 \fG(\upsvs) - \frac{1}{2} \| \IFN_{\GP}^{-1/2} \GP^{2} \upsvs \|^{2} \Bigr|
	& \leq &
	\frac{\dltwu_{3}}{2} \, \bias_{\GP}^{3}
	\, .
\label{7ywsjhd7wjhjdbiui84kje}
\end{EQA}
\end{proposition}

\begin{proof}
Define \( \fn_{\GP}(\upsv) \) by
\begin{EQA}
	\fn_{\GP}(\upsv) - \fn_{\GP}(\upsvs_{\GP})
	& = &
	\fG(\upsv) - \fG(\upsvs_{\GP}) + \langle \GP^{2} \upsvs, \upsv - \upsvs_{\GP} \rangle .
\label{f7nhw3dghydfte5w35}
\end{EQA}
The function \( \fG \) is concave, the same holds for \( \fn_{\GP} \) from \eqref{f7nhw3dghydfte5w35}.
Hence, \( \nabla \fn_{\GP}(\upsvs) = 0 \) implies \( \upsvs = \argmax \fn_{\GP}(\upsv) \).
By definition, \( \nabla \fG(\upsvs) = - \GP^{2} \upsvs \) yielding \( \nabla \fn_{\GP}(\upsvs) = - \GP^{2} \upsvs + \GP^{2} \upsvs = 0 \).
Now the results follow from Propositions~\ref{PFiWigeneric2} and \ref{Pconcgeneric2} applied 
with \( \fs(\upsv) = \fn_{\GP}(\upsv) = \fG(\upsv) - \langle \Av,\upsv \rangle \),
\( \fn(\upsv) = \fG(\upsv) \), and \( \Av = - \GP^{2} \upsvs \).
\end{proof}


The bound on the bias can be further improved under fourth-order smoothness of \( \fs \) 
using the results of Proposition~\ref{Pconcgeneric4}.

\begin{proposition}
\label{Pbiasgeneric4} 
Let \( \fs \) be concave and \( \upsvs = \argmax_{\upsv} \fs(\upsv) \).
Let also
\( \fs(\upsv) \) follow \nameref{LLsT3ref} and \nameref{LLsT4ref} with \( \IFN_{\GP} = - \nabla^{2} \fs(\upsvs) + \GP^{2} \) and some 
\( \DFN^{2} \), \( \dltwu_{3} \), \( \dltwu_{4} \), and \( \rrn \) satisfying 
\begin{EQA}[c]
	\DFN^{2} \leq \dmax^{2} \, \IFN_{\GP} \, ,
	\quad
	\rrn = \frac{3}{2} \bias_{\GP} \, ,
	\quad
	\dmax^{2} \dltwu_{3} \, \bias_{\GP} < \frac{4}{9} \, ,
	\quad
	\dmax^{2} \dltwu_{4} \, \bias_{\GP}^{2} < \frac{1}{3} \, 
	\qquad
\label{8difiyfc54wrbosT4G}
\end{EQA}
for \( \bias_{\GP} \) from \eqref{fd9dfhy4ye6fuydfrerf}.
Then \eqref{odf6fdyr6e4deuewjug} holds.
Furthermore, 
define
\begin{EQA}
	\bvn_{\GP}
	&=&
	\IFN_{\GP}^{-1} \{ \GP^{2} \upsvs + \nabla \Tens(\IFN_{\GP}^{-1} \GP^{2} \upsvs) \} \, 
\label{8vfjvr43223efryfuweefG}
\end{EQA}
with \( \Tens(\uv) = \frac{1}{6} \langle \nabla^{3} \fs(\upsvs), \uv^{\otimes 3} \rangle \) and
\( \nabla \Tens = \frac{1}{2} \langle \nabla^{3} \fs(\upsvs), \uv^{\otimes 2} \rangle \).
Then 
\begin{EQA}[c]
	\| \DFN (\bvn_{\GP} - \IFN_{\GP}^{-1} \GP^{2} \upsvs) \|
	\leq 
	\frac{\dltwu_{3}}{2} \bias_{\GP}^{2} 
	\leq 
	\frac{\dltwu_{3} \, \rr_{\GP}}{3} \, \bias_{\GP} \, ,
\label{iuvchycvf6e64rygh322b}
\end{EQA}
and
\begin{EQA}
    \| \DFN^{-1} \IFN_{\GP} (\upsvs_{\GP} - \upsvs + \bvn_{\GP}) \|
    & \leq &
    \frac{\dltwu_{4} + 2 \dmax^{2} \dltwu_{3}^{2}}{2} \, \bias_{\GP}^{3} \, ,
\label{DGttGtsGDGm13rGa4b}
    \\
    \Bigl| 
    	\fG(\upsvs_{\GP}) - \fG(\upsvs) 
		- \frac{1}{2} \| \IFN_{\GP}^{-1/2} \GP^{2} \upsvs \|^{2} 
		- \Tens(\IFN_{\GP}^{-1} \GP^{2} \upsvs) 
    \Bigr|
    & \leq &
    \frac{\dltwu_{4} + 4 \dmax^{2} \dltwu_{3}^{2}}{8} \bias_{\GP}^{4} 
    + \frac{\dmax^{2} \, (\dltwu_{4} + 2 \dmax^{2} \dltwu_{3}^{2})^{2} }{4} \, \bias_{\GP}^{6} \, .
    \qquad
\label{3d3Af12DGttGa4b}
\end{EQA}
\end{proposition}

\Subsection{A smooth penalty}
\label{Slinsmooth}
The case of a general smooth penalty \( \pent_{\GP}(\upsv) \) can be studied similarly to the quadratic case.
Denote 
\( \fG(\upsv) \eqdef \fs(\upsv) - \pent_{\GP}(\upsv) \),
\begin{EQA}
	\upsvs 
	&=& 
	\argmax_{\upsv} \fs(\upsv),
	\quad
	\upsvs_{\GP} = \argmax_{\upsv} \fG(\upsv) 
	=
	\argmax_{\upsv} \bigl\{ \fs(\upsv) - \pent_{\GP}(\upsv) \bigr\}.
	\qquad
\label{8cvkfc9fujf6jnmcer4pen}
\end{EQA}
We study the bias \( \upsvs_{\GP} - \upsvs \) induced by this penalization.
The statement of Proposition~\ref{Pbiasgeneric} and its proof can be extended to this situation 
by redefining \( \Avm \eqdef \nabla \pent_{\GP}(\upsvs) \).

\begin{proposition}
\label{Pbiaspen} 
Let \( \fG(\upsv) = \fs(\upsv) - \pent_{\GP}(\upsv) \) be concave
and follow \nameref{LLsT3ref} with some \( \DFN^{2} \), \( \dltwu_{3} \), and \( \rrn \) satisfying for \( \dmax > 0 \)
\begin{EQA}[c]
	\DFN^{2} \leq \dmax^{2} \, \IFN_{\GP} \, ,
	\qquad 
	\rrn \geq 3 \bias_{\GP}/2 \, ,
	\qquad
	\dmax^{2} \dltwu_{3} \, \bias_{\GP} < 4/9 ,
\label{7fjgjgvuvy44erd52fpen}
\end{EQA}
where 
\begin{EQA}
	\bias_{\GP}
	& \eqdef & 
	\| \DFN \, \IFN_{\GP}^{-1} \Avm \| \, ,
	\qquad
	\Avm \eqdef \nabla \pent_{\GP}(\upsvs) \, .
\label{fd9dfhy4ye6fuydfrerfM}
\end{EQA} 
Then 
\begin{EQA}[c]
	\| \DFN (\upsvs_{\GP} - \upsvs) \| \leq 3 \bias_{\GP}/2 .
\label{odf6fdyr6e4deuewjugM}
\end{EQA}
Moreover, 
\begin{EQA}[rcl]
	\bigl\| \DFN^{-1} \IFN_{\GP} (\upsvs_{\GP} - \upsvs + \IFN_{\GP}^{-1} \Avm) \bigr\|
	& \leq &
	\frac{3\dltwu_{3}}{4} \, \bias_{\GP}^{2}
	\, ,
\label{11ma3eaelDebbpen}
	\\
	\Bigl| 2 \fG(\upsvs_{\GP}) - 2 \fG(\upsvs) - \frac{1}{2} \| \IFN_{\GP}^{-1/2} \Avm \|^{2} \Bigr|
	& \leq &
	\frac{\dltwu_{3}}{2} \, \bias_{\GP}^{3}
	\, .
\label{7ywsjhd7wjhjdbiui84pen}
\end{EQA}
If, in addition, \( \fG(\upsv) \) satisfies \nameref{LLsT4ref} and \( \dmax^{2} \dltwu_{4} \, \bias_{\GP}^{2} < \frac{1}{3} \), 
then with \( \Tens_{\GP}(\uv) = \frac{1}{6} \langle \nabla^{3} \fG(\upsvs), \uv^{\otimes 3} \rangle \),
\( \nabla \Tens_{\GP} = \frac{1}{2} \langle \nabla^{3} \fG(\upsvs), \uv^{\otimes 2} \rangle \), and
\begin{EQA}
	\bvn_{\GP}
	&=&
	\IFN_{\GP}^{-1} \{ \Avm + \nabla \Tens_{\GP}(\IFN_{\GP}^{-1} \Avm) \} \, ,
\label{8vfjvr43223efryfuweefG}
\end{EQA}
it holds
\begin{EQA}
    \| \DFN^{-1} \IFN_{\GP} (\upsvs - \upsvs_{\GP} - \bvn_{\GP}) \|
    & \leq &
    \frac{\dltwu_{4} + 2 \dmax^{2} \dltwu_{3}^{2}}{2} \, \bias_{\GP}^{3} \, ,
\label{DGttGtsGDGm13rGa4bpen}
    \\
    \Bigl| \fG(\upsvs_{\GP}) - \fG(\upsvs) - \frac{1}{2} \| \IFN_{\GP}^{-1/2} \Avm \|^{2} - \Tens_{\GP}(\IFN_{\GP}^{-1} \Avm) \Bigr|
    & \leq &
    \frac{\dltwu_{4} + 4 \dmax^{2} \dltwu_{3}^{2}}{8} \bias_{\GP}^{4} 
    + \frac{\dmax^{2} \, (\dltwu_{4} + 2 \dmax^{2} \dltwu_{3}^{2})^{2} }{4} \, \bias_{\GP}^{6} \, .
    \qquad
\label{3d3Af12DGttGa4bpen}
\end{EQA}
\end{proposition}


\Chapter{Schur complement}
Consider a symmetric \( \dimttl \times \dimttl \) matrix \( \F \) with block representation 
\begin{EQA}[c]
	\F
	=
	\begin{pmatrix}
	\F_{\tarpv\tarpv} & \F_{\tarpv\nupv} \\
	\F_{\nupv\tarpv} & \F_{\nupv\nupv}
	\end{pmatrix} .
\label{HAuuuvvuvvTHS}
\end{EQA} 

\begin{lemma}
\label{LSchur}
Let the diagonal blocks \( \F_{\tarpv\tarpv} , \F_{\nupv\nupv} \) of \( \F \) be positive definite.
Define
\begin{EQA}[rcccl]
	\Fb_{\tarpv\tarpv} 
	& \eqdef & 
	\F_{\tarpv\tarpv} - \F_{\tarpv\nupv} \, \F_{\nupv\nupv}^{-1} \, \F_{\nupv\tarpv} ,
	\qquad
	\Fb_{\nupv\nupv}
	& \eqdef &
	\F_{\nupv\nupv} - \F_{\nupv\tarpv} \, \F_{\tarpv\tarpv}^{-1} \, \F_{\tarpv\nupv} 
	\, .
\label{y7djw38vjer4yw3jfclweus}
\end{EQA}
If \( \Fb_{\tarpv\tarpv} \) or \( \Fb_{\nupv\nupv} \) is also positive definite then \( \F \)
is positive definite as well.
It holds 
\begin{EQA}
	\begin{pmatrix}
		\F_{\tarpv\tarpv} & \F_{\tarpv\nupv} \\
		\F_{\nupv\tarpv} & \F_{\nupv\nupv}
	\end{pmatrix}^{-1}
	&=&
	\begin{pmatrix}
		\Id_{\dimp} & 0 \\
		- \F_{\nupv\nupv}^{-1} \, \F_{\nupv\tarpv} &\Id_{\dimq}
	\end{pmatrix}
	\,\,
	\begin{pmatrix}
		\Fb_{\tarpv\tarpv}^{-1} & 0 \\
		0 	& \F_{\nupv\nupv}^{-1}
	\end{pmatrix}
	\,\,
	\begin{pmatrix}
		\Id_{\dimp} &  - \F_{\tarpv\nupv} \, \F_{\nupv\nupv}^{-1} \\
		0 &	\Id_{\dimq}
	\end{pmatrix}
	\\
	&=&
	\begin{pmatrix}
		\Fb_{\tarpv\tarpv}^{-1} & - \Fb_{\tarpv\tarpv}^{-1} \, \F_{\tarpv\nupv} \, \F_{\nupv\nupv}^{-1} \\
		- \F_{\nupv\nupv}^{-1} \, \F_{\nupv\tarpv} \, \Fb_{\tarpv\tarpv}^{-1} & 
		\F_{\nupv\nupv}^{-1} + \F_{\nupv\nupv}^{-1} \, \F_{\nupv\tarpv} \, \Fb_{\tarpv\tarpv}^{-1} \, \F_{\tarpv\nupv} \, \F_{\nupv\nupv}^{-1} 
	\end{pmatrix}	
	\, 
\label{8jdkvtwfe6xhejdcedvscy}
\end{EQA}
and 
\begin{EQA}
	\begin{pmatrix}
		\F_{\tarpv\tarpv} & \F_{\tarpv\nupv} \\
		\F_{\nupv\tarpv} & \F_{\nupv\nupv}
	\end{pmatrix}^{-1}
	&=&
	\begin{pmatrix}
		\F_{\tarpv\tarpv}^{-1} + \F_{\tarpv\tarpv}^{-1} \, \F_{\tarpv\nupv} \, \Fb_{\nupv\nupv}^{-1} \, \F_{\nupv\tarpv} \, \F_{\tarpv\tarpv}^{-1} & - \F_{\tarpv\tarpv}^{-1} \, \F_{\tarpv\nupv} \, \Fb_{\nupv\nupv}^{-1} \\
		- \Fb_{\nupv\nupv}^{-1} \, \F_{\nupv\tarpv} \, \F_{\tarpv\tarpv}^{-1} & \Fb_{\nupv\nupv}^{-1} 
	\end{pmatrix}	
	\, .
\label{8jdkvtwfe6xhejdcedvscye}
\end{EQA}
In particular, this implies 
\( \Fb_{\tarpv\tarpv}^{-1} \, \F_{\tarpv\nupv} \, \F_{\nupv\nupv}^{-1} 
\equiv \F_{\tarpv\tarpv}^{-1} \, \F_{\tarpv\nupv} \, \Fb_{\nupv\nupv}^{-1} \),
\begin{EQA}
\label{t6rygfujio789uiew45ygi}
	\F_{\tarpv\tarpv}^{-1} + \F_{\tarpv\tarpv}^{-1} \, \F_{\tarpv\nupv} \, \Fb_{\nupv\nupv}^{-1} \, \F_{\nupv\tarpv} \, \F_{\tarpv\tarpv}^{-1}
	& \equiv &
	\Fb_{\tarpv\tarpv}^{-1}
	\, ,
	\\
	\F_{\nupv\nupv}^{-1} + \F_{\nupv\nupv}^{-1} \, \F_{\nupv\tarpv} \, \Fb_{\tarpv\tarpv}^{-1} \, \F_{\tarpv\nupv} \, \F_{\nupv\nupv}^{-1}
	& \equiv &
	\Fb_{\nupv\nupv}^{-1} 
	\, .
\end{EQA}
Moreover, for any \( \wv = (\tarpv,\nupv) \in \R^{\dimttl} \), it holds \( \| \F^{1/2} \, \wv \| \geq \| \Fb_{\nupv\nupv}^{1/2} \, \tarpv \| \) and
\begin{EQA}
	\| \F^{1/2} \, \wv \|^{2} 
	&=& 
	\| \Fb_{\nupv\nupv}^{1/2} \, \tarpv \|^{2} 
	+ \| \F_{\nupv\nupv}^{1/2} (\nupv - \F_{\nupv\nupv}^{-1} \, \F_{\nupv\tarpv} \, \tarpv) \|^{2} 
\label{jf76fugfyf5te345e4ghvhiu}
	\\
	\| \F^{-1/2} \, \wv \|^{2}
	&=&
	\| \Fb_{\tarpv\tarpv}^{-1/2} \, (\tarpv - \F_{\tarpv\nupv} \, \F_{\nupv\nupv}^{-1} \nupv)  \|^{2} 
	+ \| \F_{\nupv\nupv}^{-1/2} \nupv \|^{2} ;
\label{7ycdkjic8e38dd98le9dlo}
	\\
	\bigl( \F^{-1} \wv \bigr)_{\tarpv}  
	&=& 
	\Fb_{\tarpv\tarpv}^{-1} (\tarpv - \F_{\tarpv\nupv} \, \F_{\nupv\nupv}^{-1} \nupv) 
	=
	\Fb_{\tarpv\tarpv}^{-1} \tarpv - \F_{\tarpv\tarpv}^{-1} \, \F_{\tarpv\nupv} \, \Fb_{\nupv\nupv}^{-1} \nupv 
	\, .
\label{poryerjnjvyt65e6yjer}
\end{EQA}
Furthermore, suppose 
\begin{EQA}
	\| 
	\F_{\tarpv\tarpv}^{-1/2} \, \F_{\tarpv\nupv} \, \F_{\nupv\nupv}^{-1} \, \F_{\nupv\tarpv} \, \F_{\tarpv\tarpv}^{-1/2} 
	\|
	& \leq &
	\rhoIF^{2}
	< 
	1 .
\label{vcjcfvedtesqgghwqLco}
\end{EQA}
Then it holds for \( \F_{0} \eqdef \blk\{ \F_{\tarpv\tarpv},\F_{\nupv\nupv} \} \)
\begin{EQA}
	(1 - \rhoIF) \F_{0}
	\leq 
	\F
	& \leq &
	(1 + \rhoIF) \F_{0}
\label{f6eh3ewfd6ehev65r43tre}
\end{EQA}
and also
\begin{EQA}
	(1 - \rhoIF^{2}) \, \F_{\tarpv\tarpv}
	\leq 
	\Fb_{\tarpv\tarpv}
	& \leq &
	\F_{\tarpv\tarpv} \, ,
	\qquad
	(1 - \rhoIF^{2}) \, \F_{\nupv\nupv}
	\leq 
	\Fb_{\nupv\nupv}
	\leq 
	\F_{\nupv\nupv} \, .
\label{yg3w5dffctvyry4r7er7dgfe}
\end{EQA}
\end{lemma}

\begin{proof}
The block inversion follows by Schur's complement formula; see e.g. \cite{Boyd2004}[Appendix A.5.5].
Minimizing \( \| \F^{1/2} \, \wv \|^{2} = \tarpv^{\T} \, \F_{\tarpv\tarpv} \, \tarpv + 2 \tarpv^{\T} \, \F_{\tarpv\nupv} \, \nupv + \nupv^{\T} \, \F_{\nupv\nupv} \, \nupv \) 
w.r.t. \( \nupv \) leads to \( \nupv_{0} = - \F_{\nupv\nupv}^{-1} \, \F_{\nupv\tarpv} \tarpv \) 
and by quadraticity of \( \| \F^{1/2} \, \wv \|^{2} \) in \( \nupv \)
\begin{EQA}
	\| \F^{1/2} \, \wv \|^{2} 
	&=& 
	\tarpv^{\T} \, \F_{\tarpv\tarpv} \, \tarpv + 2 \tarpv^{\T} \, \F_{\tarpv\nupv} \, \nupv_{0} + \nupv_{0}^{\T} \, \F_{\nupv\nupv} \, \nupv_{0} 
	+ \| \F_{\nupv\nupv}^{1/2} (\nupv - \nupv_{0}) \|^{2} 
	\\
	&=&
	\| \Fb_{\nupv\nupv}^{1/2} \, \tarpv \|^{2} + \| \F_{\nupv\nupv}^{1/2} (\nupv - \nupv_{0}) \|^{2} .
\label{fji43324efgedygffhbu}
\end{EQA} 
This proves \eqref{jf76fugfyf5te345e4ghvhiu}.
Further, represent \( \F^{-1} \) using Gauss elimination \eqref{8jdkvtwfe6xhejdcedvscy}:
\begin{EQA}
	\F^{-1}
	&=&
	\begin{pmatrix}
		\Id_{\dimp} & 0 \\
		- \F_{\nupv\nupv}^{-1} \, \F_{\nupv\tarpv} &	\Id_{\dimq}
	\end{pmatrix}
	\,\,
	\begin{pmatrix}
		\Fb_{\tarpv\tarpv}^{-1} & 0 \\
		0 	& \F_{\nupv\nupv}^{-1}
	\end{pmatrix}
	\,\,
	\begin{pmatrix}
		\Id_{\dimp} &  - \F_{\tarpv\nupv} \, \F_{\nupv\nupv}^{-1} \\
		0 &	\Id_{\dimq}
	\end{pmatrix} \, .
\label{ucjhdcytew3hjf64ehfkesu}
\end{EQA}
Then 
\begin{EQA}
	\wv^{\T} \F^{-1} \wv 
	&=& 
	\begin{pmatrix}
		\tarpv - \F_{\tarpv\nupv} \, \F_{\nupv\nupv}^{-1} \nupv \\
		\nupv 
	\end{pmatrix}^{\T} \, 
	\begin{pmatrix}
		\Fb_{\tarpv\tarpv}^{-1} & 0 \\
		0 	& \F_{\nupv\nupv}^{-1}
	\end{pmatrix}
	\begin{pmatrix}
		\tarpv - \F_{\tarpv\nupv} \, \F_{\nupv\nupv}^{-1} \nupv \\
		\nupv
	\end{pmatrix} ,
\label{ucjd78renjg8we3vjfe}
\end{EQA}
and \eqref{7ycdkjic8e38dd98le9dlo} follows.
Also \eqref{8jdkvtwfe6xhejdcedvscy} implies \eqref{poryerjnjvyt65e6yjer}.

Next, define \( \F_{0} = \blk\{ \F_{\tarpv\tarpv},\F_{\nupv\nupv} \} \), 
\( U = \F_{\tarpv\tarpv}^{-1/2} \, \F_{\tarpv\nupv} \, \F_{\nupv\nupv}^{-1/2} \), and consider the matrix 
\begin{EQA}
	\F_{0}^{-1/2} \, \F \, \F_{0}^{-1/2} - \Id_{\dimttl}
	&=&
	\begin{pmatrix}
		0 & \F_{\tarpv\tarpv}^{-1/2} \, \F_{\tarpv\nupv} \, \F_{\nupv\nupv}^{-1/2}
		\\
		\F_{\nupv\nupv}^{-1/2} \, \F_{\nupv\tarpv} \, \F_{\tarpv\tarpv}^{-1/2} & 0
	\end{pmatrix}
	=
	\begin{pmatrix}
		0 & U
		\\
		U^{\T} & 0
	\end{pmatrix} .
\label{gswrew35e35e35e5txzgtqw}
\end{EQA}
Condition \eqref{vcjcfvedtesqgghwqLco} implies \( \| U U^{\T} \| \leq \rhoIF^{2} \) and hence,
\begin{EQA}
	- \rhoIF \, \Id_{\dimttl} 
	& \leq &
	\F_{0}^{-1/2} \, \F \, \F_{0}^{-1/2} - \Id_{\dimttl} 
	\leq 
	\rhoIF \, \Id_{\dimttl} \, .
\label{ye376evhgder52wedsytg}
\end{EQA}
Moreover,
\begin{EQA}
	\Fb_{\tarpv\tarpv} 
	& = &
	\F_{\tarpv\tarpv} - \F_{\tarpv\nupv} \, \F_{\nupv\nupv}^{-1} \, \F_{\nupv\tarpv} 
	=
	\F_{\tarpv\tarpv}^{1/2} (\Id_{\dimp} - U U^{\T}) \F_{\tarpv\tarpv}^{1/2} 
	\geq 
	(1 - \rhoIF^{2}) \F_{\tarpv\tarpv} \, ,
\label{ljhy6furf4jfu8rdfcweerdd}
\end{EQA}
and similarly for \( \Fb_{\nupv\nupv} \).
\end{proof}

For some situations, the nuisance variable \( \nuiv \) is by itself a composition of a few other subvectors.
We only consider the case of two variables \( \nuiv = (\zv,\nuov) \).
Denote by \( \F_{\targv\targv} \), \( \F_{\targv\zv} \), \( \F_{\targv\nuov} \), 
\( \F_{\zv\zv} \), \( \F_{\zv\nuov} \), \( \F_{\nuov\nuov} \) the corresponding blocks of \( \F \), that is,
\begin{EQA}
	\F
	&=&
	\begin{pmatrix}
		\F_{\targv\targv} & \F_{\targv\zv} & \F_{\targv\nuov} \\
		\F_{\zv\targv} & \F_{\zv\zv} & \F_{\zv\nuov} \\
		\F_{\nuov\targv} & \F_{\nuov\zv} & \F_{\nuov\nuov}
	\end{pmatrix} .
\label{576vner6734nf9h4ede}
\end{EQA}

\begin{lemma}
\label{LblocksIFT}
For the matrix \( \F \) from \eqref{576vner6734nf9h4ede}, suppose that 
\begin{EQ}[rcl]
	\| \F_{\targv\targv}^{-1/2} \F_{\targv\zv} \, \F_{\zv\zv}^{-1/2} \|
	& \leq &
	\rhoIF_{\targv\zv} \, ,
	\\
	\| \F_{\targv\targv}^{-1/2} \F_{\targv\nuov} \, \F_{\nuov\nuov}^{-1/2} \|
	& \leq &
	\rhoIF_{\targv\nuov} \, ,
	\\
	\| \F_{\nuov\nuov}^{-1/2} \F_{\nuov\zv} \, \F_{\zv\zv}^{-1/2} \|
	& \leq &
	\rhoIF_{\zv\nuov} \, ,
\label{vcjcfvedt7wdwhesqgLcompn}
\end{EQ}
and \( \max\{ \rhoIF_{\targv\zv} + \rhoIF_{\targv\nuov} \, , \rhoIF_{\targv\zv} + \rhoIF_{\zv\nuov} \, , 
\rhoIF_{\targv\nuov} + \rhoIF_{\zv\nuov} \} \leq 1 \).
Then 
\begin{EQA}
	\F
	& \geq &
	\begin{pmatrix}
		(1 - \rhoIF_{\targv\zv} - \rhoIF_{\targv\nuov}) \F_{\targv\targv} & 0 & 0 \\ 
		0 & (1 - \rhoIF_{\targv\zv} - \rhoIF_{\zv\nuov}) \F_{\zv\zv} & 0 \\
		0 & 0 & (1 - \rhoIF_{\targv\nuov} - \rhoIF_{\zv\nuov}) \F_{\nuov\nuov}
	\end{pmatrix}
\label{kjygy5t55gt2fc7w2qdxkwqiu}
\end{EQA}
\end{lemma}

\begin{proof}
By block-normalization we can reduce the proof to the case \( \F_{\targv\targv} = \Id \), \( \F_{\zv\zv} = \Id \),
\( \F_{\nuov\nuov} = \Id \).
Then it suffices to check positive semi-definiteness of the matrix
\begin{EQA}
	\BBH
	&=&
	\begin{pmatrix}
		(\rhoIF_{\targv\zv} + \rhoIF_{\targv\nuov}) \Id & \F_{\targv\zv} & \F_{\targv\nuov} \\
		\F_{\zv\targv} & (\rhoIF_{\targv\zv} + \rhoIF_{\zv\nuov}) \Id & \F_{\zv\nuov} \\
		\F_{\nuov\targv} & \F_{\nuov\zv} & (\rhoIF_{\targv\nuov} + \rhoIF_{\zv\nuov}) \Id
	\end{pmatrix}
\label{hythgwncuyuewuywhxjkqg}
\end{EQA}
with \( \| \F_{\targv\zv} \F_{\zv\targv} \| \leq \rhoIF_{\targv\zv}^{2} \),
\( \| \F_{\targv\nuov} \F_{\nuov\targv} \| \leq \rhoIF_{\targv\nuov}^{2} \),
\( \| \F_{\zv\nuov} \F_{\nuov\zv} \| \leq \rhoIF_{\zv\nuov}^{2} \).
Such a matrix is positive semi-definite by Gershgorin theorem or general results for diagonal dominant matrices.
We, however, present a simple proof.
For any \( \prmtv = (\targv,\zv,\nuov) \), we can use
\begin{EQA}
	2 |\targv^{\T} \F_{\targv\zv} \zv|
	& \leq &
	2 \| \targv \| \, \| \F_{\targv\zv} \zv \|
	\leq 
	2 \rhoIF_{\targv\zv} \| \targv \| \, \| \zv \|
	\leq 
	\rhoIF_{\targv\zv} (\| \targv \|^{2} + \| \zv \|^{2})
\label{gbtythhxct6wwnhxcd76}
\end{EQA}
and similarly for \( 2 \targv^{\T} \F_{\targv\nuov} \nuov \) and \( 2 \zv^{\T} \F_{\zv\nuov} \nuov \).
Therefore,
\begin{EQA}
	\prmtv^{\T} \BBH \prmtv
	&=&
	(\rhoIF_{\targv\zv} + \rhoIF_{\targv\nuov}) \| \targv \|^{2}
	+ (\rhoIF_{\targv\zv} + \rhoIF_{\zv\nuov}) \| \zv \|^{2}
	+ (\rhoIF_{\targv\nuov} + \rhoIF_{\zv\nuov}) \| \nuov \|^{2}
	\\
	&&
	+ \, 2 \targv^{\T} \F_{\targv\zv} \zv
	+ 2 \targv^{\T} \F_{\targv\nuov} \nuov
	+ 2 \zv^{\T} \F_{\zv\nuov} \nuov
	\\
	& \geq &
	(\rhoIF_{\targv\zv} + \rhoIF_{\targv\nuov}) \| \targv \|^{2}
	+ (\rhoIF_{\targv\zv} + \rhoIF_{\zv\nuov}) \| \zv \|^{2}
	+ (\rhoIF_{\targv\nuov} + \rhoIF_{\zv\nuov}) \| \nuov \|^{2}
	\\
	&&
	- \, \rhoIF_{\targv\zv} (\| \targv \|^{2} + \| \zv \|^{2})
	- \rhoIF_{\targv\nuov} (\| \targv \|^{2} + \| \nuov \|^{2})
	- \rhoIF_{\zv\nuov} (\| \zv \|^{2} + \| \nuov \|^{2})
	= 
	0
\label{bhxc4rt2yer6r67t478ds}
\end{EQA}
and the assertion follows.
\end{proof}

In some cases, it is more convenient to apply another block-diagonal matrix \( \DFN = (\DFN_{\targv\targv},\DFN_{\zv\zv},\DFN_{\nuov\nuov}) \).
\begin{lemma}
\label{LblocksDFN}
For the matrix \( \F \) from \eqref{576vner6734nf9h4ede}, suppose that 
\begin{EQ}[rcccccl]
	\| \DFN_{\targv\targv}^{-1} \, \F_{\targv\zv} \, \DFN_{\zv\zv}^{-1} \|
	& \leq &
	\aIF_{\targv\zv} \, ,
	\quad
	\| \DFN_{\targv\targv}^{-1} \, \F_{\targv\nuov} \, \DFN_{\nuov\nuov}^{-1} \|
	& \leq &
	\aIF_{\targv\nuov} \, ,
	\quad
	\| \DFN_{\nuov\nuov}^{-1} \, \F_{\nuov\zv} \, \DFN_{\zv\zv}^{-1} \|
	& \leq &
	\aIF_{\zv\nuov} \, .
	\qquad
\label{vcjcfvedt7wdwhesqgLcompnD}
\end{EQ}
Let also
\begin{EQ}[rcccccl]
	\| \DFN_{\targv\targv}^{-1} \, \F_{\targv\targv} \, \DFN_{\targv\targv}^{-1} \|
	& \geq &
	\bIF_{\targv\targv}^{2} \, ,
	\quad
	\| \DFN_{\zv\zv}^{-1} \, \F_{\zv\zv} \, \DFN_{\zv\zv}^{-1} \|
	& \geq &
	\bIF_{\zv\zv}^{2} \, ,
	\quad
	\| \DFN_{\nuov\nuov}^{-1} \, \F_{\nuov\nuov} \, \DFN_{\nuov\nuov}^{-1} \|
	& \geq &
	\bIF_{\nuov\nuov}^{2} \, ,
	\qquad
\label{vcjcfvedt7wdwhesqgLcompnd}
\end{EQ}
and 
\begin{EQ}[rcl]
	\bIF_{\targv\targv}^{2}  
	- \frac{\bIF_{\targv\targv} \, \aIF_{\targv\zv}}{\bIF_{\zv\zv}} - \frac{\bIF_{\targv\targv} \, \aIF_{\targv\nuov}}{\bIF_{\nuov\nuov}}
	& \geq &
	\dmax^{-2} \, ,
	\\
	\bIF_{\zv\zv}^{2}  
	- \frac{\bIF_{\zv\zv} \, \aIF_{\targv\zv}}{\bIF_{\targv\targv}} - \frac{\bIF_{\zv\zv} \, \aIF_{\zv\nuov}}{\bIF_{\nuov\nuov}}
	& \geq &
	\dmax^{-2} \, ,
	\\
	\bIF_{\nuov\nuov}^{2}  
	- \frac{\bIF_{\nuov\nuov} \, \aIF_{\targv\nuov}}{\bIF_{\zv\zv}} - \frac{\bIF_{\nuov\nuov} \, \aIF_{\zv\nuov}}{\bIF_{\zv\zv}}
	& \geq &
	\dmax^{-2} \, .
\label{dyc7f7yddfydrhye35wsss}
\end{EQ}
Then 
\begin{EQA}
	\F
	& \geq &
	\dmax^{-2} \DFM^{2} \, .
\label{kjygy5t55gt2fc7w2qDFM}
\end{EQA}
\end{lemma}

\begin{proof}
Conditions of the lemma imply \eqref{vcjcfvedt7wdwhesqgLcompn} with 
\begin{EQA}[c]
	\rhoIF_{\targv\zv}
	=
	\frac{\aIF_{\targv\zv}}{\bIF_{\targv\targv} \, \bIF_{\zv\zv}} \, , 
	\quad
	\rhoIF_{\targv\nuov}
	=
	\frac{\aIF_{\targv\nuov}}{\bIF_{\targv\targv} \, \bIF_{\nuov\nuov}} \, , 
	\quad
	\rhoIF_{\nuov\zv}
	=
	\frac{\aIF_{\nuov\zv}}{\bIF_{\nuov\nuov} \, \bIF_{\zv\zv}} \, , \quad
\label{d8dcjed7yw3md7ujdsdrdlk}
\end{EQA}
Now we apply Lemma~\ref{LblocksIFT} and note that by \eqref{dyc7f7yddfydrhye35wsss}
\begin{EQA}
	(1 - \rhoIF_{\targv\zv} - \rhoIF_{\targv\nuov}) \F_{\targv\targv} - \dmax^{-2} \DFM_{\targv\targv}^{2}
	& \geq &
	(1 - \rhoIF_{\targv\zv} - \rhoIF_{\targv\nuov}) \DFM_{\targv\targv}^{2} \bIF_{\targv\targv}^{2} - \dmax^{-2} \DFM_{\targv\targv}^{2}
	\geq 
	0 \, ,
	\\
	(1 - \rhoIF_{\targv\zv} - \rhoIF_{\zv\nuov}) \F_{\zv\zv} - \dmax^{-2} \DFM_{\zv\zv}^{2}
	& \geq &
	(1 - \rhoIF_{\targv\zv} - \rhoIF_{\zv\nuov}) \DFM_{\zv\zv}^{2} \bIF_{\zv\zv}^{2} - \dmax^{-2} \DFM_{\zv\zv}^{2}
	\geq 
	0 \, ,
	\\
	(1 - \rhoIF_{\targv\nuov} - \rhoIF_{\zv\nuov}) \F_{\nuov\nuov} - \dmax^{-2} \DFM_{\nuov\nuov}^{2}
	& \geq &
	(1 - \rhoIF_{\targv\nuov} - \rhoIF_{\zv\nuov}) \DFM_{\nuov\nuov}^{2} \bIF_{\nuov\nuov}^{2} - \dmax^{-2} \DFM_{\nuov\nuov}^{2}
	\geq 
	0 \, .
\label{ychywy7huyr78wjduy7djws}
\end{EQA}
This implies the assertion.
\end{proof}

\Chapter{Tools and proofs}
\label{Stoolssemi}

\Section{Semiparametric estimation. Tools and proofs}
\label{Ssemitools}

This section collects some technical tools and proofs of the main results on profile MLE.

\Subsection{Proof of Theorem~\ref{Tsemieffp3n}}
We apply Theorem~\ref{TFiWititG3}.
It holds for any linear mapping \( \QF \) on \( \R^{\dimttl} \), 
\begin{EQA}
	\bigl\| \QF (\tilde{\prmtv} - \prmtvs - \IFT^{-1} \nabla \zeta) \bigr\|
	&=&
	\bigl\| \QF \, \IFT^{-1} \DFM \, \DFM^{-1} \IFT (\tilde{\prmtv} - \prmtvs - \IFT^{-1} \nabla \zeta) \bigr\|
	\\
	& \leq &
	\| \QF \, \IFT^{-1} \DFM \| \, 
	\bigl\| \DFM^{-1} \IFT (\tilde{\prmtv} - \prmtvs - \IFT^{-1} \nabla \zeta) \bigr\| .
\label{8fcjfctrt54fdfhwufgnh}
\end{EQA}
It remains to note that \( \QF \prmtv = \QP \tarpv \) and \( \DFM^{2} \leq \dmax^{2} \IFT \) imply
\begin{EQA}[c]
	\| \QF \, \IFT^{-1} \DFM \| 
	\leq 
	\dmax^{2} \| \QF \DFM^{-1} \| 
	\leq 
	\| \QP \, \DFM_{\tarpv\tarpv}^{-1} \| ;
\label{y7edrdjvci823jgujy6gtuyr}
\end{EQA} 
cf. Remark~\ref{RremainderD}.

\Subsection{Proof of Theorem~\ref{TsemiWilks3}}
Represent
\begin{EQA}
	&& \nquad
	\LLp(\tilde{\tarpv}) - \LLp(\tarpvs)
	=
	\LL(\tilde{\tarpv},\tilde{\nupv}) - \LL(\tarpvs,\nupvs)
	- \bigl\{ \LL(\tarpvs,\tilde{\nupv}(\tarpvs)) - \LL(\tarpvs,\nupvs) \bigr\} ,
	\qquad
\label{iu8dcjhfc767y63whf6hdfu}
\end{EQA}
where
\begin{EQA}
	\tilde{\nupv}(\tarpvs)
	& \eqdef &
	\sup_{\nupv} \LL(\tarpvs,\nupv) .
\label{h9odf0f0kw3t5yikgh}
\end{EQA}
Application of \eqref{kfchvyfyrdetye46t} of Theorem~\ref{Tconcsemi} to each of two parts of the decomposition yields
\begin{EQ}[rcl]
	\bigl| 2 \LL(\tilde{\tarpv},\tilde{\nupv}) - 2 \LL(\tarpvs,\nupvs) 
		- \| \IFT^{-1/2} \nabla \zeta \|^{2}
	\bigr|
	& \leq &
	\dltwu_{3} \| \DFM \, \IFT^{-1} \nabla \zeta \|^{3} ,
	\\
	\bigl| 2 \LL(\tarpvs,\tilde{\nupv}(\tarpvs)) - 2 \LL(\tarpvs,\nupvs) 
		- \| \IFT_{\nupv\nupv}^{-1/2} \scorem{\nupv} \zeta \|^{2}
	\bigr|
	& \leq &
	\dltwu_{3} \| \DFM_{\nupv\nupv} \, \IFT_{\nupv\nupv}^{-1} \scorem{\nupv} \zeta \|^{3} .
\label{9fvjvft53w3wghbuinh4}
\end{EQ}
These two bounds imply the assertion.
To see this, consider first the case when the matrix \( \IFT \) is block-diagonal, that is,
\( \IFT = \blk\{ \IFT_{\tarpv\tarpv},\IFT_{\nupv\nupv} \} \).
Then 
\begin{EQA}
	 \| \IFT^{-1/2} \nabla \zeta \|^{2} 
	 &=& 
	\| \IFT_{\tarpv\tarpv}^{-1/2} \scorem{\tarpv} \zeta \|^{2} + \| \IFT_{\nupv\nupv}^{-1/2} \scorem{\nupv} \zeta \|^{2} 
\label{ivfiv8e35ftg6y674hdjw}
\end{EQA}
and \eqref{vdud6ehry3hrf67ruy233} follows from bound \eqref{kfchvyfyrdetye46t} 
and decomposition \eqref{iu8dcjhfc767y63whf6hdfu}.
In the case of a general matrix \( \IFT \), \eqref{7ycdkjic8e38dd98le9dlo} of Lemma~\ref{LSchur} implies
with \( \IFTb_{\tarpv\tarpv} 
= \IFT_{\tarpv\tarpv} - \IFT_{\tarpv\nupv} \, \IFT_{\nupv\nupv}^{-1} \IFT_{\nupv\tarpv} \)
\begin{EQA}
	\| \IFT^{-1/2} \, \nabla \zeta \|^{2}
	&=&
	\| \IFTb_{\tarpv\tarpv}^{-1/2} \, (\scorem{\tarpv} \zeta - \IFT_{\tarpv\nupv} \IFT_{\nupv\nupv}^{-1} \scorem{\nupv} \zeta)  \|^{2} 
	+ \| \IFT_{\nupv\nupv}^{-1/2} \scorem{\nupv} \zeta \|^{2} \, .
\label{7ycdkjic8e346gghrehjsi98de}
\end{EQA}
This identity yields \eqref{vdud6ehry3hrf67ruy233} as in the block-diagonal case
due to \eqref{hf8g8hy984u3453fhfh}.


\ifshort{
\Subsection{Proof of Theorem~\ref{TsemiWilks4}}

As in the proof of \eqref{vdud6ehry3hrf67ruy233}, let us decompose 
\begin{EQA}
	&& \nquad
	\LLp(\tilde{\tarpv}) - \LLp(\tarpvs)
	= 
	\LL(\tilde{\tarpv},\tilde{\nupv}) - \LL(\tarpvs,\nupvs)
	- \bigl\{ \LL(\tarpvs,\tilde{\nupv}(\tarpvs)) - \LL(\tarpvs,\nupvs) \bigr\} .
\label{iu8dcjhfc767y63whf6hdfu4}
\end{EQA}
Applying \eqref{87dfsjqweudsfjht7w3} of Theorem~\ref{Teffp4} to each part of the decomposition yields
on \( \Omega(\xx) \)
\begin{EQ}[rcl]
	\bigl| 2 \LL(\tilde{\tarpv},\tilde{\nupv}) - 2 \LL(\tarpvs,\nupvs) 
		- \| \IFT^{-1/2} \nabla \zeta \|^{2} - \TensU(\IFT^{-1} \nabla \zeta)
	\bigr|
	& \leq &
	\err_{\prmtv} \, ,
	\\
	\bigl| 2 \LL(\tarpvs,\tilde{\nupv}(\tarpvs)) - 2 \LL(\tarpvs,\nupvs) 
		- \| \IFT_{\nupv\nupv}^{-1/2} \scorem{\nupv} \zeta \|^{2}
		- \TensU(\IFT_{\nupv\nupv}^{-1} \scorem{\nupv} \zeta)
	\bigr|
	& \leq &
	\err_{\nupv} \, ,
\label{9fvjvft53w3wghbuinh4}
\end{EQ}
with
\begin{EQA}
	\err_{\prmtv}
	& \leq &
	\frac{\dltwu_{4} + 4 \dltwu_{3}^{2}}{8} \| \DFM \IFT^{-1} \nabla \zeta \|^{4} 
    + \frac{(\dltwu_{4} + 3 \dltwu_{3}^{2})^{2}}{9} \, \| \DFM \IFT^{-1} \nabla \zeta \|^{6} \, ,
    \qquad
    \\
	\err_{\nupv}
	& \leq &
	\frac{\dltwu_{4} + 4 \dltwu_{3}^{2}}{8} \| \HPM \IFT_{\nupv\nupv}^{-1} \scorem{\nupv} \zeta \|^{4} 
    + \frac{(\dltwu_{4} + 3 \dltwu_{3}^{2})^{2}}{9} \, \| \HPM \IFT_{\nupv\nupv}^{-1} \scorem{\nupv} \zeta \|^{6} \, .
    \qquad
\label{87dfsjqweudsfjht7w3se}
\end{EQA}
Here we assume a natural embedding of \( \HPM^{-1} \scorem{\nupv} \zeta \) in \( \R^{\dimttl} \)
with the \( \tarpv \)-component equal to zero.
It remains to mention that 
\( \| \IFT^{-1/2} \nabla \zeta \|^{2} = \| \IFT_{\nupv\nupv}^{-1/2} \scorem{\nupv} \zeta \|^{2} 
+ \| \DPMb \xivr \|^{2} \); see \eqref{7ycdkjic8e346gghrehjsi98de}.
}{}

\Subsection{Proof of Lemma~\ref{Leffdimrho}}

Lemma~\ref{LSchur} ensures the following bounds:
\begin{EQA}
	(1 - \rhoIF) \blk\{ \IFT_{\tarpv\tarpv},\IFT_{\nupv\nupv} \}
	\leq 
	\IFT
	& \leq &
	(1 + \rhoIF) \blk\{ \IFT_{\tarpv\tarpv},\IFT_{\nupv\nupv} \}
\label{f6eh3ewfd6ehev65r43tres}
	\\
	(1 - \rhoIF^{2}) \, \IFT_{\tarpv\tarpv}
	\leq 
	\IFTb_{\tarpv\tarpv}
	& \leq &
	\IFT_{\tarpv\tarpv} 
	\, .
\label{ihgt64eghdfyf65ewytgdh}
\end{EQA}
By \eqref{hf8g8hy984u3453fhfh}
\begin{EQA}[rcl]
	\Var^{1/2}(\xivr)
	&=&
	\Var^{1/2}\bigl\{ \IFTb_{\tarpv\tarpv}^{-1/2} (\scorem{\tarpv} \zeta 
		- \IFT_{\tarpv\nupv} \, \IFT_{\nupv\nupv}^{-1} \scorem{\nupv} \zeta) 
	\bigr\}
	\\
	& \leq &
	\Var^{1/2}\bigl\{ \IFTb_{\tarpv\tarpv}^{-1/2} \scorem{\tarpv} \zeta \bigr\}
	+ \Var^{1/2}\bigl\{ \IFTb_{\tarpv\tarpv}^{-1/2} \IFT_{\tarpv\nupv} \, \IFT_{\nupv\nupv}^{-1} \scorem{\nupv} \zeta \bigr\}
	\, .
\label{hf8g8hy984u3453fhfhp}
\end{EQA}
Further, 
\begin{EQA}
	&& \nquad 
	\Var\bigl\{ \IFTb_{\tarpv\tarpv}^{-1/2} \IFT_{\tarpv\nupv} \, \IFT_{\nupv\nupv}^{-1} \scorem{\nupv} \zeta \bigr\}
	=
	\IFTb_{\tarpv\tarpv}^{-1/2} \IFT_{\tarpv\nupv} \, \IFT_{\nupv\nupv}^{-1} \Var(\scorem{\nupv} \zeta)
	\IFT_{\nupv\nupv}^{-1} \, \IFT_{\tarpv\nupv} \, \IFTb_{\tarpv\tarpv}^{-1/2}
	\\
	& \leq &
	\frac{1}{1 - \rhoIF^{2}}\IFT_{\tarpv\tarpv}^{-1/2} \IFT_{\tarpv\nupv} \, \IFT_{\nupv\nupv}^{-1/2} \Var(\scorem{\nupv} \zeta)
	\IFT_{\nupv\nupv}^{-1} \, \IFT_{\tarpv\nupv} \, \IFT_{\tarpv\tarpv}^{-1/2}
	\leq 
	\frac{\rhoIF^{2}}{1 - \rhoIF^{2}} \BBH_{\tarpv}
	\, .
\label{t6dyd67wywehv5e4we4fgdu}
\end{EQA}
Similarly
\begin{EQA}
	\Var\bigl( \IFTb_{\tarpv\tarpv}^{-1/2} \scorem{\tarpv} \zeta \bigr)
	& = &
	\IFTb_{\tarpv\tarpv}^{-1/2} \Var(\scorem{\tarpv} \zeta) \, \IFTb_{\tarpv\tarpv}^{-1/2} 
	\leq 
	\frac{1}{1 - \rhoIF^{2}} \IFT_{\tarpv\tarpv}^{-1/2} \Var(\scorem{\tarpv} \zeta) \, \IFT_{\tarpv\tarpv}^{-1/2}
	\leq 
	\frac{1}{1 - \rhoIF^{2}} \BBH_{\tarpv}
	\, .
\label{7dtyv83ejo42et7gh95}
\end{EQA}
This implies
\begin{EQA}
	\Var\bigl( \xivr \bigr)
	& \leq &
	\biggl( \frac{\rhoIF}{\sqrt{1 - \rhoIF^{2}}} + \frac{1}{\sqrt{1 - \rhoIF^{2}}} \biggr)^{2} 
	\BBH_{\tarpv} 
	=
	\frac{1 + \rhoIF}{1 - \rhoIF} \BBH_{\tarpv} \, ,
\label{dyhdf7y6v6fr6tetd55j}
\end{EQA}
and the assertion follows.


\def\neta{s}
\Section{Nonlinear regression. Tools and proofs}
\label{Stoolscalm}

This section collects some technical statements and proofs of the main results on estimation in nonlinear regression model.

\Subsection{Local concavity}
\label{Swarmstartnl}

Remind that \( \Thetad \) is defined by \eqref{cujduyfdt5c544retd5td54e4} 
with \( \DPNLc^{2} = \nabla \Regrfv(\thetav_{0}) \, \nabla \Regrfv(\thetav_{0})^{\T} \); see \eqref{tsfdc6frdw6red6wt7qagdudwg}.
%
Now we state strong concavity of \( \LL(\upsv) \) with an explicit lower bound on 
\( \IFT(\upsv) = - \nabla^{2} \LL(\upsv) \).
The first technical result bounds the value \( \| \Regrfv(\thetav) - \Regrfv(\thetav_{0}) \| \) over \( \thetav \in \Thetad \).

\begin{lemma}
\label{LnormThetad}
Suppose \nameref{LLM2ref}.
Then for any \( \thetav \in \Thetad \), it holds with \( \DPNLc^{2} = \DPNL^{2}(\thetav_{0}) \)
\begin{EQA}
	\| \Regrfv(\thetav) - \Regrfv(\thetav_{0}) \|^{2}
	& \leq &
	(1 + \dltwbD) \| \DPNLc (\thetav - \thetav_{0}) \|^{2} \, .
\label{vbib7gve4f4eere3gwjhwsde}
\end{EQA}
\end{lemma}

\begin{proof}
Given \( \thetav \in \Thetad \), denote \( \uv = \thetav - \thetav_{0} \).
By definition of \( \Thetad \), it holds \( \| \DPNLc \uv \| \leq \rrdd \).
We now use the representation
\begin{EQA}
	\Regrf_{j}(\thetav) - \Regrf_{j}(\thetav_{0})
	&=&
	\int_{0}^{1} \langle \nabla \Regrf_{j}(\thetav + t \uv), \uv \rangle \, dt
\label{juvdywhdwcsrqtqhggfcvbcdwcd}
\end{EQA}
and in view of \eqref{fchghiu87686574e5rtyyu}, it holds
\begin{EQA}
	\bigl| \Regrf_{j}(\thetav) - \Regrf_{j}(\thetav_{0}) \bigr|^{2}
	=
	\biggl( 
		\int_{0}^{1} \langle \nabla \Regrf_{j}(\thetav + t \uv), \uv \rangle \, dt 
	\biggr)^{2}
	& \leq &
	\int_{0}^{1} \langle \nabla \Regrf_{j}(\thetav + t \uv), \uv \rangle^{2} \, dt 
\label{fqe9f8gefyeuyefuygyde}
\end{EQA}
yielding
\begin{EQA}
	\sum_{j=1}^{\dimq} | \Regrf_{j}(\thetav) - \Regrf_{j}(\thetav_{0}) |^{2}
	& \leq &
	\int_{0}^{1} \| \DPNL(\thetav + t \uv) \uv \|^{2} \, dt
	\leq 
	(1 + \dltwbD) \, \| \DPNLc \uv \|^{2}
\label{c8h8e8e838383j4hr4uhceh}
\end{EQA}
as required.
\end{proof}

\begin{lemma}
\label{LD2Mktk}
Suppose
\nameref{LLM2ref},
\nameref{LLMkref},
and \eqref{8difiyfc54wrboeNL}.
Then for any \( \upsv = (\thetav,\etav) \in \Upsd \) from \eqref{cujduyfdt5c544retd5td54e4},
the matrix \( \IF(\upsv) \) given by 
\begin{EQA}[rcl]
	\IF(\upsv)
	& \eqdef &
	\nabla \Regrfv(\thetav) \, \nabla \Regrfv(\thetav)^{\T} 
		+ \sum_{j=1}^{\dimq} \{ \Regrf_{j}(\thetav) - \eta_{j} \} \, \nabla^{2} \Regrf_{j}(\thetav) 
	\qquad
\label{DPGP2HPGP220}
\end{EQA}
satisfies with \( \smlc \eqdef 2\rrde/\deta^{1/2} \leq 1/2 \)
\begin{EQA}
	(2 - \smlc) \, \DPNL^{2}(\thetav)
	& \leq &
	\IF(\upsv) 
	\leq 
	(2 + \smlc) \, \DPNL^{2}(\thetav) .
\label{pcuc7e4emvyfyemsdcgytgf0}
\end{EQA}
\end{lemma}

\begin{proof}
It suffices to show that
\begin{EQA}
	\biggl| \sum_{j=1}^{\dimq} \bigl\{ \Regrf_{j}(\thetav) - \eta_{j} \bigr\} 
	\, \nabla^{2} \Regrf_{j}(\thetav) \biggr|
	& \leq &
	\smlc \, \DPNL^{2}(\thetav) 
\label{poiuytgfbhnmjklopwsedfgbhN}
\end{EQA}
The structural relation \( \etav_{0} = \Regrfv(\thetav_{0}) \) yields for \( \upsv = (\thetav,\etav) \in \Upsd \)
by \( (1 + \dltwbD)/\sqrt{2} \leq 1 \)
\begin{EQA}
	\| \Regrfv(\thetav) - \etav \|
	&=&
	\| \Regrfv(\thetav) - \Regrfv(\thetav_{0}) \| + \| \etav - \etav_{0} \| 
	\leq 
	\frac{(1 + \dltwbD)^{1/2}}{\sqrt{2}} \rrde + \rrde 
	\leq 
	2 \rrde \, , 
	\qquad
\label{dftgejdtwfrwenanosgenN}
\end{EQA}
and by the Cauchy-Schwarz inequality
\begin{EQA}
	&& \nquad
	\biggl| \sum_{j=1}^{\dimq} \bigl\{ \Regrf_{j}(\thetav) - \eta_{j} \bigr\} \,
		\langle \nabla^{2} \Regrf_{j}(\thetav), \uv^{\otimes 2} \rangle \biggr|
	\\
	& \leq &
	\| \Regrfv(\thetav) - \etav \| 
	\biggl( \sum_{j=1}^{\dimq} \langle \nabla^{2} \Regrf_{j}(\thetav), \uv^{\otimes 2} \rangle^{2} \biggr)^{1/2}
	\leq 
	2 \, \rrde \biggl( \sum_{j=1}^{\dimq} \langle \nabla^{2} \Regrf_{j}(\thetav), \uv^{\otimes 2} \rangle^{2} \biggr)^{1/2}  .
	\qquad	
\label{vjuvygvfyr6r6erg34r5etN}
\end{EQA}
Hence, by 
\eqref{wfegy7r5qrw35edfhgdyfys} of \nameref{LLMkref}
\begin{EQA}
	&& \nquad
	\biggl| \sum_{j=1}^{\dimq} \bigl\{ \Regrf_{j}(\thetav) - \eta_{j} \bigr\} \, 
	\langle \nabla^{2} \Regrf_{j}(\thetav), \uv^{\otimes 2} \rangle \biggr|
	\leq 
	2 \, \deta^{-1/2} \, \rrde \, \| \DPNL(\thetav) \uv \|^{2}
	\leq 
	\smlc \, \| \DPNL(\thetav) \uv \|^{2}
\label{ggtdg2321323213689768768N}
\end{EQA}
and \eqref{poiuytgfbhnmjklopwsedfgbhN} follows.
\end{proof}

Bound \eqref{pcuc7e4emvyfyemsdcgytgf0} implies the following result.

\begin{lemma}
\label{LDFblk2}
Suppose the conditions of Lemma~\ref{LD2Mktk}.
Let \( \IFTb_{\GP,\thetav\thetav}(\upsv) \) and \( \IFTb_{\GP,\etav\etav}(\upsv) \) be given by \eqref{ydsydfyfyedtw3e23w4er}.
It holds  with \( \smlc = 2\rrde/\deta^{1/2} \leq 1/2 \)
\begin{EQA}[rcccl]
	(1 - \smlc) \DPNL^{2}
	& \leq &
	\IFTb_{\thetav\thetav}(\upsv) 
	& \leq & 
	(1 + \smlc) \DPNL^{2}
	\, ,
	\\
	\frac{1 - \smlc}{1 - \smlc/2} \, \Id_{\dimq} 
	& \leq &
	\IFTb_{\etav\etav}(\upsv) 
	& \leq & 
	2 \, \Id_{\dimq} \, .
\label{7dkcuf5dethv7edhjwe}
\end{EQA}
\end{lemma}

\begin{proof}
The bound on \( \IFTb_{\thetav\thetav}(\upsv) \) follows from \eqref{pcuc7e4emvyfyemsdcgytgf0}.
Also 
\begin{EQA}
	\nabla \Regrfv(\thetav)^{\T} \, \IF^{-1}(\upsv) \, \nabla \Regrfv(\thetav)
	& \leq &
	\frac{1}{2 - \smlc} \nabla \Regrfv(\thetav)^{\T} \, \DPNL^{-2}(\thetav) \, \nabla \Regrfv(\thetav)
	\leq 
	\frac{1}{1 - \smlc/2} \, \Id_{\dimq} \,
	\\
	\nabla \Regrfv(\thetav)^{\T} \, \IF^{-1}(\upsv) \, \nabla \Regrfv(\thetav)
	& \geq &
	0 \, .
\label{y67d6tc6e3yhergfteh}
\end{EQA}
This implies \eqref{7dkcuf5dethv7edhjwe}.
\end{proof}

\begin{lemma}
\label{Sconccalm}
Suppose the conditions of Lemma~\ref{LD2Mktk} and let \( \smlc \leq 1/6 \).
Then \( \LL(\upsv) \) is concave in \( \Upsd \) and
for any \( \upsv = (\thetav,\etav) \in \Upsd \)
\begin{EQA}
	\IFT(\upsv) 
	& \geq &
	\begin{pmatrix}
		(2 - \smlc) \DPNL^{2}(\thetav) & - \nabla \Regrfv(\thetav) 
		\\
		- \nabla \Regrfv(\thetav)^{\T} & 2 \Id_{\dimq}  
	\end{pmatrix}	
	\geq 
	\frac{1}{2}
	\begin{pmatrix}
		\DPNL^{2}(\thetav) & 0 \\
		0 & \Id_{\dimq}
	\end{pmatrix}
	\, .
	\qquad
\label{ufcee7wtrdfew5eG}
\end{EQA}
\end{lemma}

\begin{proof}
The first inequality in \eqref{ufcee7wtrdfew5eG} follows from \eqref{pcuc7e4emvyfyemsdcgytgf0}.
Further, 
\begin{EQA}
	&& \nquad
	\begin{pmatrix}
		(2 - \smlc) \DPNL^{2}(\thetav) & - \nabla \Regrfv(\thetav) 
		\\
		- \nabla \Regrfv(\thetav)^{\T} & 2 \Id_{\dimq}  
	\end{pmatrix}	
	- 
	\frac{1}{2}
	\begin{pmatrix}
		\DPNL^{2}(\thetav) & 0 \\
		0 & \Id_{\dimq}
	\end{pmatrix}
	\\
	&=&
	\begin{pmatrix}
		(3/2 - \smlc) \DPNL^{2}(\thetav) & - \nabla \Regrfv(\thetav) 
		\\
		- \nabla \Regrfv(\thetav)^{\T} & \frac{3}{2} \Id_{\dimq}  
	\end{pmatrix}	
	\geq 
	0
\label{dyce8ydeuhdihi3nrv7teu}
\end{EQA}
because of \( 3/2 - \smlc \geq 4/3 \) and
\( \nabla \Regrfv(\thetav) \, \nabla \Regrfv(\thetav)^{\T} = 2 \DPNL^{2}(\thetav) \).
\end{proof}

\begin{lemma}
\label{LDFIFuv}
For any vector \( \uv \in \R^{\dimq} \) and \( \wv = (0,\uv) \), 
it holds with \( \DFM^{2} = \blk\{ \DPNL^{2}, \Id_{\dimq} \} \)
\begin{EQA}
	\| \DFM \, \IFT^{-1} \wv \|
	& \leq &
	2 \| \uv \| \, .
\label{njcfjjcf2sd5gvyj3gfke}
\end{EQA}
If  \( \wv = (\uv,0) \) for \( \uv \in \R^{\dimp} \) then
\begin{EQA}
	\| \DFM \, \IFT^{-1} \wv \|
	& \leq &
	2 \,  \| \DPNL^{-1} \uv \| \, . 
\label{uys7f78g7848rtfhd6fhe}
\end{EQA}
\end{lemma}

\begin{proof}
For \( \wv = (0,\uv) \), identity \( \DFM^{-1} \wv = \wv \) and \eqref{ufcee7wtrdfew5eG} imply
\begin{EQA}
	\| \DFM \, \IFT^{-1} \wv \|
	& = &
	\| \DFM \, \IFT^{-1} \DFM \, \wv \|
	\leq 
	2 \left\| \begin{pmatrix}
		\Id_{\dimp} & 0 \\
		0 & \Id_{\dimq}
	\end{pmatrix} \binom{0}{\uv}
	\right\|
	=
	2 \| \uv \|
	\, 
\label{dcyvyrhe46fg3duwyhy}
\end{EQA}
as claimed in \eqref{njcfjjcf2sd5gvyj3gfke}.
Similarly, with \( \wv = (\uv,0) \), it holds \(  \)
\begin{EQA}
	\| \DFM \, \IFT^{-1} \wv \|^{2}
	& \leq &
	4 \uv^{\T} \DPNL^{-2} \, \DPNL^{2} \, \DPNL^{-2} \uv
	=
	4 \| \DPNL^{-1} \uv \|^{2} \, .
\label{7c7eue3ugftge4ugy3hek}
\end{EQA}
This yields \eqref{uys7f78g7848rtfhd6fhe}.
\end{proof}

%

\Subsection{Local smoothness}
Apart from the basic conditions about linearity of the stochastic component of \( \LL(\upsv) \) and about concavity 
of the expectation \( \E \LL(\upsv) \), 
we need some local smoothness properties of the expected penalized log-likelihood 
\( \E \, \LL(\upsv) \).
We make use of the fact that
the only non-quadratic term in \( \E \, \LL(\upsv) \) is 
\( \fs(\upsv) = - \| \Regrfv(\thetav) - \etav \|^{2}/2 \), 
and the value \( \| \Regrfv(\thetav) - \etav \| \)
is uniformly bouned in \( \Upsd \).

Let us fix for the moment some \( \etav \in \Etad \).
For \( \upsv = (\thetav,\etav) \),
consider \( \fNL(\thetav) = \sum_{j=1}^{\dimq} \fNL_{j}(\thetav) \) with
\(  \fNL_{j}(\thetav) = - | \eta_{j} - \Regrf_{j}(\thetav) |^{2}/2 \), 
\begin{EQA}
    \langle \nabla^{3} \fNL_{j}(\thetav), \uv^{\otimes 3} \rangle 
    & = &
    - 3 \langle \nabla \Regrf_{j}(\thetav),\uv \rangle \,
    	\langle \nabla^{2} \Regrf_{j}(\thetav),\uv^{\otimes 2} \rangle
    - \bigl\{ \Regrf_{j}(\thetav) - \eta_{j} \bigr\} \, \langle \nabla^{3} \Regrf_{j}(\thetav),\uv^{\otimes 3} \rangle ,
    \qquad
\label{d3utd4utd4dt4cm}
    \\
    \langle \nabla^{4} \fNL_{j}(\thetav), \uv^{\otimes 4} \rangle 
    & = &
    - 3 \langle \nabla^{2} \Regrf_{j}(\thetav),\uv^{\otimes 2} \rangle^{2}
    - 4 \langle \nabla \Regrf_{j}(\thetav),\uv \rangle \, 
    	\langle \nabla^{3} \Regrf_{j}(\thetav),\uv^{\otimes 3} \rangle
    \\
    &&
    - \, \bigl\{ \Regrf_{j}(\thetav) - \eta_{j} \bigr\}  
    \langle \nabla^{4} \Regrf_{j}(\thetav),\uv^{\otimes 4} \rangle .
\end{EQA}

The next results explain how local smoothness properties of \( \fNL(\upsv) \) can be characterized under 
\nameref{LLM2ref} and \nameref{LLMkref}.
First we bound the \( \thetav \)-derivatives of \( \fNL(\upsv) \).

\begin{lemma}
\label{Lsmoothcalm}
Assume \nameref{LLM2ref} and \nameref{LLMkref} with \( k=2,3 \). 
Let also \eqref{8difiyfc54wrboeNL} hold.
Then for any \( (\thetav,\etav) \in \Upsd \) and \( \uv \in \R^{\dimp} \) with \( \smlc = 2\rrde/\deta^{1/2} \leq 1/2 \)
\begin{EQA}
	| \langle \nabla^{3} \fNL(\thetav), \uv^{\otimes 3} \rangle |
	& \leq &
	(3 + \smlc) \deta \, \| \DPNL(\thetav) \uv \|^{3} \, .
\label{voedf8efyffte4er4wgwjw}
\end{EQA}
If also \nameref{LLMkref} holds for \( k=4 \), then
\begin{EQA}
	| \langle \nabla^{4} \fNL(\thetav), \uv^{\otimes 4} \rangle |
	& \leq &
	(7 + \smlc) \deta^{2} \, \| \DPNL(\thetav) \uv \|^{4} \, .
\label{voedf8efyffte4er4wgwjw4}
\end{EQA}
\end{lemma}

\begin{proof}
It holds
\begin{EQA}
	&& \nquad
	\langle \nabla^{3} \fNL(\thetav), \uv^{\otimes 3} \rangle 
	=
	- \frac{d^{3}}{dt^{3}} \| \etav - \Regrfv(\thetav + t \uv) \|^{2}/2 \biggl|_{t=0}
	\\
	& = &
	- 3 \sum_{j=1}^{\dimq} \langle \nabla \Regrf_{j}(\thetav),\uv \rangle \,
    	\langle \nabla^{2} \Regrf_{j}(\thetav),\uv^{\otimes 2} \rangle 
	- \sum_{j=1}^{\dimq} \bigl\{ \Regrf_{j}(\thetav) - \eta_{j} \bigr\} \, 
		\langle \nabla^{3} \Regrf_{j}(\thetav),\uv^{\otimes 3} \rangle \, .
	\qquad
\label{hs8d8s8sduhudhuhdsu}
\end{EQA}
By the Cauchy-Schwarz inequality and \eqref{wfegy7r5qrw35edfhgdyfys} of \nameref{LLMkref} with \( k = 2 \)
\begin{EQA}
	&& \nquad
	\biggl| \sum_{j=1}^{\dimq} \langle \nabla \Regrf_{j}(\thetav),\uv \rangle \,
    	\langle \nabla^{2} \Regrf_{j}(\thetav),\uv^{\otimes 2} \rangle \biggr|
	\\
	& \leq &
	\biggl\{ \sum_{j=1}^{\dimq} \langle \nabla \Regrf_{j}(\thetav),\uv \rangle^{2} \biggr\}^{1/2}
	\biggl\{ 
		\sum_{j=1}^{\dimq} \langle \nabla^{2} \Regrf_{j}(\thetav),\uv^{\otimes 2} \rangle^{2} 
	\biggr\}^{1/2}
	\leq 
	\deta \, \| \DPNL(\thetav) \uv \|^{3} .
\label{werd5wr5werw6ed37tr37tr73}
\end{EQA}
Similarly to \eqref{ggtdg2321323213689768768N}, using \eqref{dftgejdtwfrwenanosgenN} and \nameref{LLMkref} with \( k = 3 \)
\begin{EQA}
	&& \nquad
	\biggl| \sum_{j=1}^{\dimq} \bigl\{ \Regrf_{j}(\thetav) - \eta_{j} \bigr\} \, 
		\langle \nabla^{3} \Regrf_{j}(\thetav),\uv^{\otimes 3} \rangle \biggr|
	\\
	& \leq &
	\| \Regrfv(\thetav) - \etav \| \, \,  
	\biggl\{ \sum_{j=1}^{\dimq} \langle \nabla^{3} \Regrf_{j}(\thetav),\uv^{\otimes 3} \rangle^{2} \biggr\}^{1/2} 
	\leq 
	2 \rrdd \, \deta^{2} \, \| \DPNL(\thetav) \uv \|^{3} \, .
\label{079u8760u9jhkigfetr6w}
\end{EQA}
Now \eqref{voedf8efyffte4er4wgwjw} follows by \( 2 \rrdd \, \deta = \smlc \). 
The proof of \eqref{voedf8efyffte4er4wgwjw4} is similar.
\end{proof}

Due to the results of Lemma~\ref{Lsmoothcalm}, regularity of each 
\( \Regrf_{j}(\thetav) \) implies smoothness of 
\( \fNL(\thetav) = - \| \Regrfv(\thetav) - \etav \|^{2}/2 \) 
in \( \thetav \) with \( \etav \) fixed.
Now we check the full dimensional smoothness characteristics of \( \fs(\upsv) \).
As in the case of a general SLS model, consider the third and fourth derivatives of 
\( \fs(\upsv) \) for \( \upsv = (\thetav,\etav) \in \Upsd \).

\begin{lemma}
\label{Lupssmooth}
Assume \nameref{LLM2ref} and \nameref{LLMkref} for \( k = 2,3 \). 
Then for any \( \upsv = (\thetav,\etav) \in \Upsd \) 
and any \( \wv = (\uv,\netav) \) with \( \uv \in \R^{\dimp} \) and \( \netav \in \R^{\dimq} \)
\begin{EQA}
	\bigl| \langle \nabla^{3} \fs(\upsv), \wv^{\otimes 3} \rangle \bigr|
	& \leq &
	(3 + \smlc) \deta \| \DPNL(\thetav) \uv \|^{3} 
	+ 3 \deta \, \| \netav \| \, \| \DPNL(\thetav) \uv \|^{2} 
	\leq 
	\hmax_{3} \, \deta \| \DFM \, \wv \|^{3} \, ,
	\qquad
\label{9dd5de5de4e44wr5w3twuyw2wh}
\end{EQA}
where \( \DFM^{2} = \blk\{ \DPNL^{2}, \Id_{\dimq} \} \) and
\begin{EQA}[c]
	\hmax_{3}
	\eqdef
	(3 + \smlc) (1 + \dltwbD)^{3} + 3 (1 + \dltwbD)^{2} \, .
	\qquad
\label{9dd5de5de4e44wr5w3twuyw2wh4}
\end{EQA}
%
Similarly, under \nameref{LLMkref} with \( k=4 \), 
for any \( \upsv = (\thetav,\etav) \in \Upsd \) 
and any \( \wv = (\uv,\netav) \) 
\begin{EQA}
	\bigl| \langle \nabla^{4} \fs(\upsv), \wv^{\otimes 4} \rangle \bigr|
	& \leq &
	(7 + \smlc) \deta^{2} \| \DPNL(\thetav) \uv \|^{4} 
	+ 4 \deta^{2} \, \| \netav \| \, \| \DPNL(\thetav) \uv \|^{3} 
	\leq 
	\hmax_{4} \, \deta^{2} \, \| \DFM \, \wv \|^{4}
	\, ,
	\qquad
\label{idj4r7dg22www2222hhddh}
\end{EQA}
where 
\begin{EQA}[c]
	\hmax_{4}
	\eqdef
	(7 + \smlc) (1 + \dltwbD)^{4} + 4 (1 + \dltwbD)^{3} \, .
	\qquad
\label{9dd5de5de4e44wr5w3twuyw2wh44}
\end{EQA}
This yields full dimensional conditions \nameref{LLsT3ref} and \nameref{LLsT4ref}
with \( \dltwu_{3} = \hmax_{3} \, \deta \) and \( \dltwu_{4} = \hmax_{3} \, \deta^{2} \).
\end{lemma}

\begin{proof}
Fix any \( \upsv = (\thetav,\etav) \in \Upsd \), \( \wv = (\uv,\netav) \)
with \( \| \DFM \wv \| \leq \rr \).
Then 
\begin{EQA}
    \langle \nabla^{3} \fs(\upsv), \wv^{\otimes 3} \rangle 
	&=&
	\langle \nabla^{3} \fNL(\thetav), \uv^{\otimes 3} \rangle 
	- 3  \sum_{j=1}^{\dimq} \neta_{j} \, \langle \nabla^{2} \Regrf_{j}(\thetav),\uv^{\otimes 2} \rangle.
\label{d3utd4utd4dt4cmn}
\end{EQA}
By 
\nameref{LLMkref} for \( k = 2 \)
\begin{EQA}
	&& \nquad
	\biggl| \sum_{j=1}^{\dimq} \neta_{j} \, 
			\langle \nabla^{2} \Regrf_{j}(\thetav),\uv^{\otimes 2} \rangle 
	\biggr|
	\leq 
	\| \netav \| \biggl( \sum_{j=1}^{\dimq} 
		\langle \nabla^{2} \Regrf_{j}(\thetav),\uv^{\otimes 2} \rangle^{2} 
	\biggr)^{1/2}
	\leq 
	\deta \, \| \netav \| \, \| \DPNL(\thetav) \uv \|^{2} .
\label{p12348ujnwsdcvoihg3}
\end{EQA}
Combining this with 
\eqref{voedf8efyffte4er4wgwjw} of Lemma~\ref{Lsmoothcalm} yields
\begin{EQA}
	| \langle \nabla^{3} \fs(\upsv), \wv^{\otimes 3} \rangle |
	& \leq &
	(3 + \smlc) \deta \| \DPNL(\thetav) \uv \|^{3} 
	+ 3 \deta \, \| \netav \| \, \| \DPNL(\thetav) \uv \|^{2} \, 
\label{qqduwyqwyfettewdquww18d}
\end{EQA}
and \eqref{9dd5de5de4e44wr5w3twuyw2wh} follows.
Further, by \nameref{LLM2ref} and \( \| \DPNL \uv \|^{2} + \| \netav \|^{2} \leq \rr^{2} \), 
\begin{EQA}
	&& \nquad
	(3 + \smlc) \deta \| \DPNL(\thetav) \uv \|^{3} 
	+ 3 \deta \, \| \netav \| \, \| \DPNL(\thetav) \uv \|^{2}
	\\
	& \leq &
	(3 + \smlc) (1 + \dltwbD)^{3} \deta \| \DPNL \uv \|^{3} 
	+ 3 \deta (1 + \dltwbD)^{2} \, \| \netav \| \, \| \DPNL \uv \|^{2}
	\\
	& \leq &
	\bigl\{ (3 + \smlc) (1 + \dltwbD)^{3} + 3 (1 + \dltwbD)^{2} \bigr\} \deta \rr^{3}
	\, ,
\label{2gcv6rtjh9sjeychdkovytf}
\end{EQA}
yielding \eqref{9dd5de5de4e44wr5w3twuyw2wh4}.
Under \nameref{LLMkref} for \( k=4 \), we can use 
\begin{EQA}
    && \nquad
    \langle \nabla^{4} \fs(\upsv), \wv^{\otimes 4} \rangle 
    =
    \langle \nabla^{4} \fNL(\thetav), \uv^{\otimes 4} \rangle 
	- 4  \sum_{j=1}^{\dimq} \neta_{j} \, \langle \nabla^{3} \Regrf_{j}(\thetav),\uv^{\otimes 3} \rangle 
\label{0oijfcyfyr7e474t75738}
\end{EQA}
and \eqref{idj4r7dg22www2222hhddh} follows similarly.
\end{proof}

{\renewcommand{\Section}{\section}
\renewcommand{\Subsection}{\subsection}
\def\GPp{\GP_{\priord}}
\def\GPb{\GP_{1}}
\def\DPp{\DP_{\priord}}
\def\projp{\text{\Large$\pi$}}
\def\Projc{\Proj^{c}}
\def\dimmc{\dimm'}
\def\dimms{\dimm^{*}}
\def\CGP{w}
\def\CGPa{\CGP_{0}}
\def\CGPb{\CGP_{1}}
\def\CGPs{\CGP^{*}}
\def\Fclass{\mathcal{F}}
\def\CONSTGPT{\CONSTi_{\GPT}}

\Chapter{Smooth penalties and minimax rates}
\label{Spriorexample}
This section presents some examples of choosing \( \GP^{2} \) for 
achieving the ``bias-variance trade-off''
and obtaining rate optimal results.
Let pMLE \( \tilde{\upsv}_{\GP} \), its population counterpart \( \upsvs_{\GP} \),
and the background true parameter \( \upsvs \) be given by \eqref{tuGauLGususGE}.
Theorem~\ref{TQFiWibias} yields the following bound for the risk 
\( \riskt_{\QP} \) of \( \tilde{\upsv}_{\GP} \) given by \eqref{7djhed8cjfct534etgdhdy}:
\begin{EQA}
	\E \| \QP (\tilde{\upsv}_{\GP} - \upsvs) \|^{2}
	& \approx &
	\riskt_{\QP}
	\approx
	\dimA_{\QP} + \| \QP \IF_{\GP}^{-1} \GP^{2} \upsvs \|^{2} \, .
\label{bnjmb3ed6w2jh21fg7fj}
\end{EQA}
This suggest to select the operator \( \GP^{2} \) by forcing the ``bias-variance trade-off'' 
\( \| \dimA_{\QP} \asymp \QP \IF_{\GP}^{-1} \GP^{2} \upsvs \|^{2} \).
Later in this section, we illustrate this relation through popular examples of regularization by projection or by roughness penalty.
For any considered choice of penalization \( \GP^{2} \), we assume the conditions of Propositions~\ref{PconcMLEgenc} and \ref{Lvarusetb}
to be fulfilled. 
To simplify the analysis, we also assume \( \VP^{2} = \DPN^{2} = \IF \) and consider two specific choices of \( \QP \):
prediction/response loss with \( \QP = \IF^{1/2} \) and
estimation loss with \( \QP = \Id \). 
In the latter case, we focus on a direct problem with a bounded condition number of \( \IF \).

\Section{Projection estimation: bias-variance trade-off and risk bounds}
\label{SprojGP}
Consider the class of projection estimators given by a set of sub-spaces \( \{ \II_{\dimm} \} \) of the parameter space \( \R^{\dimp} \).
For each \( m \), only projection \( \Proj_{\dimm} \upsv \) on the subspace \( \II_{\dimm} \) is considered 
but there is no any additional penalization.
Formally, this corresponds to the diagonal matrix \( \GP_{\dimm}^{2} \) with \( \dimm \) diagonal elements equal to zero,
and the remaining ones equal to infinity.
Later we everywhere use the sub-index \( \dimm \) in place of \( \GP_{\dimm} \).
It appears that \( \IF_{\dimm}^{-1}(\upsv) \, \GP_{\dimm}^{2} = \{ \IF(\upsv) + \GP_{\dimm}^{2} \}^{-1} \GP_{\dimm}^{2} \) 
for any \( \upsv \in \Ups \) is nothing but the orthogonal projector \( \Projc_{\dimm} = \Id - \Proj_{\dimm} \) 
on the subspace \( \II_{\dimm}^{c} \) of the remaining coordinates \( \dimm+1,\dimm+2,\ldots \).
In particular, for \( \upsv = \upsvs_{\GP} \),
\begin{EQA}
	\IF_{\dimm}^{-1} \, \GP_{\dimm}^{2} \upsvs
	&=&
	\Projc_{\dimm} \upsvs
	=
	\upsvs - \Proj_{\dimm} \upsvs \, .
\label{yeufyy3hfyfhwmdygneum}
\end{EQA}
Similarly, \( \IF_{\dimm}^{-1} \IF = \Proj_{\dimm} \) and \( \IF_{\dimm}^{-1} \IF \, \IF_{\dimm}^{-1} = \IF_{\dimm}^{-1} \). 
As \( \DPN^{2} = \VP^{2} = \IF \), this leads to 
\begin{EQA}[c]
	\dimA_{\QP,\dimm} = \tr(\QP \IF_{\dimm}^{-1} \VP^{2} \IF_{\dimm}^{-1} \QP^{\T}) 
	= 
	\tr (\QP \IF_{\dimm}^{-1} \QP^{\T}) 
	\, .
\label{uef8efhjfjfbufb8874uf}
\end{EQA}
In particular,
\begin{EQA}[c]
	\dimA_{\QP,\dimm} 
	= 
	\begin{cases}
		\tr (\Proj_{\dimm}) = \dimm, 	& \QP = \IF^{1/2} , \\
		\tr(\IF_{\dimm}^{-1}), 		& \QP = \Id .
	\end{cases}
\label{uef8efhjfjfbufb8874ufpr}
\end{EQA}
For the corresponding risk \( \riskt_{\QP,\dimm} \) from \eqref{bnjmb3ed6w2jh21fg7fj}, we obtain by Theorem~\ref{TQFiWibias}
\begin{EQA}
	\riskt_{\QP,\dimm}
	& \approx &
	\begin{cases}
		\dimm + \| \IF^{1/2} \Projc_{\dimm} \upsvs \|^{2} & \QP = \IF^{1/2} , \\
		\tr(\IF_{\dimm}^{-1}) + \| \Projc_{\dimm} \upsvs \|^{2}, 		& \QP = \Id .
	\end{cases}
	 \, .
\label{hgfcuygeuhfui3w23}
\end{EQA}
The optimal (or oracle) choice of \( \dimm \) can be given by minimization of the risk \( \riskt_{\QP,\dimm} \):
\begin{EQA}
	\dimms
	& \eqdef &
	\argmin_{\dimm} \riskt_{\QP,\dimm}
	\, .
\label{d7fjfjjvytv435y6ewee}
\end{EQA}
%
A standard way of obtaining the minimax rate of estimation is based on the 
approximation theory for functional spaces.
One assumes that \( \upsvs \) belongs to a special set \( \Fclass \) like a Sobolev or Besov ball, and 
\begin{EQA}
	\| \upsv - \Proj_{\dimm} \upsv \|
	& \leq &
	\rho_{\dimm} \, , 
	\qquad
	\upsv \in \Fclass \, ,
\label{edydjd7u73ef6fhej7hj}
\end{EQA}
where the \( \rho_{\dimm} \)'s are \( \Fclass \)-specific and decrease to zero as \( \II_{\dimm} \) increase.
As an example, consider a ``smooth'' signal \( \upsvs \) from a Sobolev ball \( \BBB(\smpa,\CGPa) \):
\begin{EQA}
	\BBB(\smpa,\CGPa)
	& \eqdef &
	\Bigl\{ \upsv = (\ups_{j}) \in \R^{\dimp} \colon \sum_{j \geq 1} j^{2\smpa} \ups_{j}^{2} \leq \CGPa \Bigr\} 
\label{BBBsttsjj1sj2S}
\end{EQA}
with \( \smpa > 0 \) and \( \CGPa \asymp 1 \).
Then for any \( \upsv \in \BBB(\smpa,\CGPa) \)
\begin{EQA}
	\| \Projc_{\dimm} \upsv \|^{2}
	& \leq &
	\dimm^{-2\smpa} \sum_{j > \dimm} j^{2\smpa} \ups_{j}^{2}
	\leq 
	\CGPa \, \dimm^{-2\smpa} \, .
\label{ysgduy27828e3h2ewsud}
\end{EQA}
We additionally assume that \( \IF \leq \CONSTIF \, n \Id \), where \( n \) is a scaling parameter meaning the sample size, while
\( \CONSTIF \) is an absolute constant.
For \( \QP = \IF^{1/2} \), it holds
\begin{EQA}
	\| \IF^{1/2} \Projc_{\dimm} \upsvs \|^{2}
	& \leq &
	\CONSTIF \, n \| \Projc_{\dimm} \upsvs \|^{2}
	\leq 
	\CONSTIF \, n \rho_{\dimm}^{2} \, .
\label{ysgduy27828e3h2ewsud}
\end{EQA}
Therefore, the \( \upsvs \)-dependent choice \eqref{d7fjfjjvytv435y6ewee} can be replaced by the \( \Fclass \)-specific choice
\begin{EQA}
	\dimms
	&=&
	\argmin_{\dimm} \{ \dimm + \CONSTIF \, n \rho_{\dimm}^{2} \} \, .
\label{8fifivugvy5435tdjbuuyr}
\end{EQA}
Typically the solution to this problem satisfies the balance relation \( \dimms \asymp \CONSTIF \, n \rho_{\dimms}^{2} \)
leading to the risk \( \riskt_{\QP,\dimms} \asymp \dimms \).
For the case of a Sobolev ball, \( \rho_{\dimm}^{2} = \CGPa \dimm^{-2\smpa} \), 
and the trade-off relation reads as \( \dimm \asymp n \, \dimm^{-2\smpa} \).
This leads to the standard rule of thumb 
\( \dimms \asymp n^{1/(2 \smpa+1)} \) and \( \riskt_{\QP,\dimms} \asymp n^{1/(2\smpa+1)} \).

For the estimation loss with \( \QP = \Id \), the situation is similar as long as 
a \emph{direct} problem is considered and
the condition number \( \CONSTIF = \lambda_{\max}(\IF)/\lambda_{\min}(\IF) \) 
of the Fisher information operator \( \IF = \IF(\upsvs) \) is fixed.
Later we assume that \( n \Id \leq \IF \leq \CONSTIF \, n \Id \), where \( n = \lambda_{\min}(\IF) \). 
As \( \IF \geq n \Id \), we obtain for the value \( \dimA_{\QP,\dimm} \) from \eqref{uef8efhjfjfbufb8874ufpr}
\( \dimA_{\QP,\dimm} = \tr \IF_{\dimm}^{-1} \leq n^{-1} \tr \Proj_{\dimm} \leq n^{-1} \dimm \).
Therefore, the optimal choice of \( \dimm \) can be reduced to minimization of 
\( n^{-1} \dimm + \CONST \rho_{\dimm}^{2} \) which coincides with \eqref{8fifivugvy5435tdjbuuyr}.
For the case of a Sobolev ball with \( \rho_{\dimm}^{2} = \CGPa \dimm^{-2\smpa} \),
this yields \( \dimms \asymp n^{1/(2 \smpa+1)} \) and \( \riskt_{\dimms} \asymp n^{-2\smpa/(2\smpa+1)} \).


\Section{Roughness penalty}
This section explores a more general case of a penalizing family \( \GPT = \{ \GP^{2} \} \). 
We show under rather general conditions that the risk of each \( \tilde{\upsv}_{\GP} \) 
can be decomposed and analyzed as in the case projection estimation with a proper choice
of the projection sub-space.
Assume as earlier that \( \VP^{2} = \IF \).
For any \( \QP \) and any \( \GP^{2} \in \GPT \), it holds
\begin{EQA}[c]
	\dimA_{\QP} = \tr(\QP \IF_{\GP}^{-1} \VP^{2} \IF_{\GP}^{-1} \QP^{\T}) 
	= 
	\tr (\QP \IF_{\GP}^{-1} \IF \, \IF_{\GP}^{-1} \QP^{\T}) 
\label{uef8efhjfjfbufb8GP}
\end{EQA}
and
\begin{EQA}[c]
	\dimA_{\QP} 
	= 
	\begin{cases}
		\tr (\IF_{\GP}^{-1} \IF)^{2} , 	& \QP = \IF^{1/2} , \\
		\tr(\IF_{\GP}^{-2} \IF), 		& \QP = \Id .
	\end{cases}
\label{uef8efhjfjfbufbefi8GP}
\end{EQA}
Similarly
\begin{EQA}
	\bias_{\QP}
	&=&
	\| \QP \IF_{\GP}^{-1} \GP^{2} \upsvs \|^{2}
	=
	\begin{cases}
		\| \IF^{1/2} \IF_{\GP}^{-1} \GP^{2} \upsvs \|^{2} , 	& \QP = \IF^{1/2} , \\
		\| \IF_{\GP}^{-1} \GP^{2} \upsvs \|^{2}, 		& \QP = \Id .
	\end{cases}
\label{d7ewjd6jr4g8ejefgghi4}
\end{EQA}
The aim is to describe these quantities and the related risk bounds in terms of the spectral characteristics of the penalizing matrices \( \GP^{2} \).
Later we assume that each \( \GP^{2} \) fulfills the polynomial growth condition on its spectrum.

\begin{description}
\item[\( \bb{(\GPT)} \)\label{GPTref}]
	\emph{Let \( \GP^{2} \in \GPT \) and \( \gp_{1}^{2} \leq \ldots \leq \gp_{\dimp}^{2} \) be its increasing eigenvalues.
	Then for all \( \dimm < \dimp \) 
	}
\begin{EQA}
	\sum_{j = \dimm+1}^{\dimp} \gp_{j}^{-4}
	& \leq &
	\CONSTGPT \, \dimm \, \gp_{\dimm+1}^{-4} ,
	\qquad
	\sum_{j=1}^{\dimm} \gp_{j}^{4}
	\leq 
	\CONSTGPT \, \dimm \, \gp_{\dimm}^{4} \, .
\label{dsujjue67he3jv7nmeuej}
\end{EQA}

\end{description}

Condition \eqref{dsujjue67he3jv7nmeuej} assumes that \( \gp_{j}^{2} \) grow at least as \( j^{2\smpa} \) for
\( \smpa > 1/4 \).
The constant \( \CONSTGPT \) depends on \( \smpa \) only.

\begin{lemma}
\label{LGPgrowth}
Let \( \IF \leq \CONSTIF \, n \Id \). 
Assume \nameref{GPTref}.
For \( \GP^{2} \in \GPT \), let \( \gp_{1}^{2} \leq \ldots \leq \gp_{\dimp}^{2} \) be its
increasing eigenvalues.
Then for any \( \dimm \)
\begin{EQA}
	\tr (\IF_{\GP}^{-1} \IF)^{2}
	& \leq &
	\Bigl( 1 + \frac{\CONSTGPT \, \CONSTIF^{2} \, n^{2}}{\gp_{\dimm+1}^{4}} \Bigr) \, \dimm 
	\, .
\label{r8uerf6yghe5erjkrf7e3j}
\end{EQA}
In particular, if \( \dimm_{\GP} \) is the largest \( m \) such that \( \gp_{\dimm}^{2} \leq \CONSTIF \, n \) then
\begin{EQA}
	\tr (\IF_{\GP}^{-1} \IF)^{2}
	& \leq &
	(1 + \CONSTGPT) \dimm_{\GP} 
	\, .
\label{r8uerf6yghe5erjkrf7e3jGP}
\end{EQA}
It \( n \Id \leq \IF \leq \CONSTIF \, n \Id \) then 
\begin{EQA}
	\tr (\IF_{\GP}^{-2} \IF)
	& \leq &
	n^{-1} \tr (\IF_{\GP}^{-1} \IF)^{2}
	\, .
\label{8d8d8d8h3rf76fyj34ffg}
\end{EQA}
\end{lemma}

\begin{proof}
As \( \IF \leq \CONSTIF \, n \Id \), it holds by \eqref{dsujjue67he3jv7nmeuej} for any \( \dimm \)
\begin{EQA}
	\tr (\IF_{\GP}^{-1} \IF)^{2}
	& \leq &
	\CONSTIF^{2} \, n^{2} \tr(\CONSTIF \, n \Id + \GP^{2})^{-2}
	\leq 
	\sum_{j=1}^{\dimm} \frac{\CONSTIF^{2} \, n^{2}}{(\CONSTIF \, n + \gp_{j}^{2})^{2} }
	+ \sum_{j=\dimm+1}^{\dimp} \frac{\CONSTIF^{2} \, n^{2}}{(\CONSTIF \, n + \gp_{j}^{2})^{2} }
	\\
	& \leq &
	\dimm + \CONSTIF^{2} \, n^{2} \sum_{j=1+\dimm_{\GP}}^{\dimp} \frac{1}{\gp_{j}^{4}}
	\leq 
	\dimm + \CONSTIF^{2} \, n^{2} \, \CONSTGPT \, \dimm \, \gp_{\dimm+1}^{-4}
\label{uf8f73jr7jewlfgu7e4hjfk}
\end{EQA}
and the first bound follows.
Further, by definition of \( \dimm_{\GP} \)
\begin{EQA}[c]
	\gp_{\dimm_{\GP}+1}^{4}
	\geq 
	\CONSTIF^{2} \, n^{2}
\label{t6dfytfyfyyhe457dfhn}
\end{EQA}
which reduces the second bound to the first one.
\end{proof}

Now we evaluate the bias term using similar arguments. 

\begin{lemma}
\label{LbiasGP}
Assume \( n \Id \leq \IF \leq \CONSTIF \, n \Id \).
Let \( \GP^{2} \) satisfy \eqref{dsujjue67he3jv7nmeuej}.
Then for any \( \dimm \geq 1 \) 
\begin{EQA}
	\| \IF_{\GP}^{-1} \GP^{2} \upsvs \|^{2}
	& \leq &
	\| \Projc_{\dimm} \upsvs \|^{2} + \frac{\CONSTIF \, \CONSTGPT \, \dimm}{n} \, \| \GP \Proj_{\dimm} \upsvs \|^{2} \, ,
\label{67duchuychjjmdyctwbnsr}
\end{EQA}
where \( \Proj_{\dimm} \) is the spectral projector for \( \GP^{2} \), that is, \( \Proj_{\dimm} \) projects onto the subspace \( \II_{\dimm} \)
of the first \( \dimm \) principle components of \( \GP^{2} \).
\end{lemma}

\begin{proof}
The use of \( \IF \geq n \Id \) implies for any \( \dimm \)
\begin{EQA}
	\| \IF_{\GP}^{-1} \GP^{2} \upsvs \|^{2}
	& \leq &
	\| (n \Id + \GP^{2})^{-1} \GP^{2} \upsvs \|^{2}
	\\
	&=&
	\| (n \Id + \GP^{2})^{-1} \GP^{2} \Proj_{\dimm} \upsvs \|^{2} + \| (n \Id + \GP^{2})^{-1} \GP^{2} \Projc_{\dimm} \upsvs \|^{2}
	\\
	& \leq &
	n^{-2} \| \GP^{2} \Proj_{\dimm} \upsvs \|^{2} + \| \Projc_{\dimm} \upsvs \|^{2} \, .
\label{6dfyfhhvui8987urne3yt}
\end{EQA}
If \( \gp_{\dimm}^{2} \leq \CONSTIF \, n \), it holds by the Cauchy-Schwarz inequality and \eqref{dsujjue67he3jv7nmeuej}
\begin{EQA}
	\| \GP^{2} \Proj_{\dimm} \upsvs \|^{2}
	& \leq &
	\| \GP \Proj_{\dimm} \upsvs \|^{2} \sum_{j=1}^{\dimm} \gp_{j}^{2} 
	\leq 
	\| \GP \Proj_{\dimm} \upsvs \|^{2} \, \CONSTGPT \, \dimm \, \gp_{\dimm}^{2}
	\leq 
	\| \GP \Proj_{\dimm} \upsvs \|^{2} \, \CONSTIF \, \CONSTGPT \, \dimm \, n \, .
\label{7ehfciui3e8f8gh8y969gff}
\end{EQA}
This implies \eqref{67duchuychjjmdyctwbnsr}.
\end{proof}

Now we summarize in the case with \( \QP = \IF^{1/2} \).
For \( \QP = \Id \) the conclusion is similar.

\begin{proposition}
\label{PriskGP}
Assume \( n \Id \leq \IF \leq \CONSTIF \, n \Id \) and \( \GP^{2} \) satisfy \nameref{GPTref}.
Let also
\begin{EQA}
	\CONSTIF \, \| \GP \Proj_{\dimm_{\GP}} \upsvs \|^{2}
	& \leq &
	\alp \, n \, .
\label{dujhbr7e7t63tybuwtvhg}
\end{EQA}
If \( \dimm_{\GP} \) is the largest \( m \) such that \( \gp_{\dimm}^{2} \leq \CONSTIF \, n \) then
\begin{EQA}
	\riskt_{\QP}
	& \leq &
	\begin{cases}
		( 1 + \CONSTGPT + \alp \CONSTGPT ) \dimm_{\GP} + n \| \Projc_{\dimm_{\GP}} \upsvs \|^{2} \, , &  \QP = \IF^{1/2} \, , \\
		( 1 + \CONSTGPT + \alp \CONSTGPT ) \dimm_{\GP}/n + \| \Projc_{\dimm_{\GP}} \upsvs \|^{2} \, , &  \QP = \Id \, .
	\end{cases}
\label{ydfjgv7ejed7urkh744u}
\end{EQA}
\end{proposition}

Usually the value \( \alp \) in \eqref{dujhbr7e7t63tybuwtvhg} is small and the term \( \alp \CONSTGPT \) can be ignored
even when \( \| \GP \upsvs \| \) is very large.
For illustration, let us consider the most interesting case when \( \upsvs \) is \( \GPa \)-smooth for \( \GPa^{2} \leq \GP^{2} \), 
that is, \( \GPa \)-smoothness is less restrictive then \( \GP \)-smoothness.

\begin{lemma}
\label{LsmoothGPa}
Let \( \upsvs \) be \( \GPa \)-smooth for some \( \GPa \in \GPT \), that is, \( \| \GPa \upsvs \|^{2} \leq 1 \).
Let also \( \GP \) and \( \GPa \) commute and hence, have the same eigenspaces, and 
\( (\gp_{0,j}^{2}) \) be the ordered eigenvalues of \( \GPa^{2} \).
Moreover, let the ratio \( \gp_{j}^{2}/\gp_{0,j}^{2} \) grow with \( j \).
Then
\begin{EQA}
	n^{-1} \| \GP \Proj_{\dimm_{\GP}} \upsvs \|^{2}
	& \leq &
	\CONSTIF / \gp_{0,\dimm_{\GP}}^{2} \, .
\label{ydfjfd8uuy655tehdf}
\end{EQA}
\end{lemma}

\begin{proof}
As \( \GP \) and \( \GPa \) commute, the same holds for \( \GPa \) and \( \Proj_{\dimm} \). 
Hence,
\begin{EQA}
	\| \GP \Proj_{\dimm} \upsvs \|^{2}
	&=&
	\| \GP \GPa^{-1} \Proj_{\dimm} \GPa \upsvs \|^{2}
	\leq 
	\gp_{m}^{2} / \gp_{0,\dimm}^{2} \, .
\label{8djfcfjcvjvcht5gr7yhe}
\end{EQA}
Applying this bound to \( \dimm = \dimm_{\GP} \) and using \( \gp_{\dimm_{\GP}}^{2} \leq \CONSTIF \, n \) yields the result.
\end{proof}

The right-hand side of \eqref{ydfjfd8uuy655tehdf} is small provided that \( \gp_{0,\dimm_{\GP}}^{2} \) is large when \( \gp_{\dimm_{\GP}}^{2} \approx \CONSTIF \, n \).
Therefore, even a minor smoothness of \( \upsvs \) ensures that the value \( n^{-1} \| \GP \Proj_{\dimm} \upsvs \|^{2} \) is relatively small. 

\medskip
We conclude that a roughness penalty \( \GP^{2} \) satisfying \nameref{GPTref} yields nearly the same risk as the projection estimator 
with a special \( \GP^{2} \)-dependent choice \( \II_{\dimm_{\GP}} \) of the corresponding sub-space.
This reduces the problem of risk minimization to the case of projection estimation considered earlier.

\Section{An example}

Consider a particular example when \( \{ \GP^{2} \} \) is a univariate family of penalizing matrices 
\( \GP^{2} \) of the form \( \GP^{2} = \CGP \GPb^{2} \)
for \( \GPb^{2} = \diag\{ \gp_{1}^{2}, \ldots, \gp_{\dimp}^{2} \} \) fixed.
Everywhere in this section, we assume \( n \Id \leq \IF \leq \CONSTIF \, n \Id \).
Each value \( \CGP \) identifies 
the spectral cut-off value \( \dimm_{\CGP} \) which solves \( \CGP \gp_{\dimm}^{2} \approx \CONSTIF \, n \).
If \( \gp_{j}^{2} = h(j) \) for a strictly increasing function \( h(\cdot) \) then
\begin{EQA}
	\dimm_{\CGP}
	& \approx &
	h^{-1}(\CONSTIF \, n/\CGP) .
\label{d7xcyhjkkf7643ddhrusnsdh}
\end{EQA}
Now we study the bias term beginning from the case when \( \| \GPb \upsvs \|^{2} \leq \CGPb \).
Then \( \| \GP \Proj_{\dimm} \upsvs \|^{2} \leq \| \GP \upsvs \|^{2} \leq \CGP \| \GPb \upsvs \|^{2} \).
This yields the upper bound for the risk:
\begin{EQA}
	\riskt_{\CGP}
	& \approx &
	\dimm_{\CGP} + \| \GP \Proj_{\dimm_{\CGP}} \upsvs \|^{2}
	\leq 
	h^{-1}(\CONSTIF \, n/\CGP) + \CGP \, \CGPb \, .
\label{hdycyuccuje47ctdge6}
\end{EQA}
The optimal/oracle choice \( \CGPs \) of \( \CGP \) is obtained by minimization of this expression 
w.r.t. \( \CGP \) leading to
\begin{EQA}
	\CGPs
	&=&
	\argmin_{\CGP} \bigl\{ h^{-1}(\CONSTIF \, n/\CGP) + \CGP \, \CGPb \bigr\} \, .
\label{usd8w6fjebcje6cjdebdj}
\end{EQA}
Another way of defining the optimal choice is based on \eqref{8fifivugvy5435tdjbuuyr} and 
\eqref{d7xcyhjkkf7643ddhrusnsdh}.
Namely, we define the optimal spectral cut-off value \( \dimms \) by \eqref{8fifivugvy5435tdjbuuyr}
and then identify the corresponding Lagrange multiplier \( \CGP \) by 
\( \CGPs = \gp_{\dimms}^{-2} \CONSTIF \, n \).

For instance, if \( \gp_{j}^{2} = h(j) = j^{2\smpa} \) then \( h^{-1}(j) = j^{1/(2\smpa)} \),
\( \dimm_{\CGP} \approx (\CONSTIF \, n/\CGP)^{1/(2\smpa)} \),
and
\begin{EQA}
	\riskt_{\CGP}
	& \leq &
	\dimm_{\CGP} + \| \GP \Proj_{\dimm_{\CGP}} \upsvs \|^{2}
	\leq 
	(\CONSTIF \, n/\CGP)^{1/(2\smpa)} + \CGP \, \CGPb \, .
\label{hdycyuccuje47ctdge6}
\end{EQA}
This yields
\begin{EQA}
	\CGPs
	& \asymp &
	(\CONSTIF \, n)^{1/(2\smpa+1)} \CGPb^{-2\smpa/(1 + 2\smpa)} ,
	\qquad
	\riskt^{*}
	\asymp
	\CGPs \, \CGPb
	\asymp
	(\CONSTIF \, n \, \CGPb)^{1/(2\smpa+1)} \, .
\label{ucg4e6fhcvje48cjvgetdxbd}
\end{EQA}
The case when \( \upsvs \) is not \( \GPb \)-smooth is a bit more involved
because there is no minimax solution over a class of signals \( \upsvs \).
The \( \upsvs \)-dependent choice of \( \dimms \) follows \eqref{d7fjfjjvytv435y6ewee}
and \( \CGPs = \CONSTIF \, n / (\dimms)^{2\smpa} \).

}
\bibliography{exp_ts,listpubm-with-url}

\end{document}